\documentclass[12pt,a4paper]{amsart}
\usepackage{leftidx}
\usepackage{tikz}
\usetikzlibrary{matrix,positioning,decorations.pathreplacing}
\usepackage{color,subfigure,graphicx}
\usepackage{xspace}
\usepackage{soul}
\usepackage{mathrsfs}

\usepackage{float}
\usepackage[]{hyperref}
\hypersetup{
   colorlinks=true,       % false: boxed links; true: colored links
    linkcolor=blue,          % color of internal links
    citecolor=blue,        % color of links to bibliography
    filecolor=magenta,      % color of file links
    urlcolor=cyan           % color of external links
}

\usepackage{amssymb,amsmath,amsfonts,amsthm,amscd}

\makeatletter
\renewcommand\normalsize{%
    \@setfontsize\normalsize{11.7}{14pt plus .3pt minus .3pt}%
    \abovedisplayskip 10\p@ \@plus4\p@ \@minus4\p@
    \abovedisplayshortskip 6\p@ \@plus2\p@
    \belowdisplayshortskip 6\p@ \@plus2\p@
    \belowdisplayskip \abovedisplayskip}
\renewcommand\small{%
    \@setfontsize\small{9.5}{12\p@ plus .2\p@ minus .2\p@}%
    \abovedisplayskip 8.5\p@ \@plus4\p@ \@minus1\p@
    \belowdisplayskip \abovedisplayskip
    \abovedisplayshortskip \abovedisplayskip
    \belowdisplayshortskip \abovedisplayskip}
\renewcommand\footnotesize{%
    \@setfontsize\footnotesize{8.5}{9.25\p@ plus .1pt minus .1pt}%%
    \abovedisplayskip 6\p@ \@plus4\p@ \@minus1\p@
    \belowdisplayskip \abovedisplayskip
    \abovedisplayshortskip \abovedisplayskip
    \belowdisplayshortskip \abovedisplayskip}
\ifdefined\pdfpagewidth
\else
\fi
\calclayout
\makeatother

\newcommand{\mb}[1]{\ensuremath{\mathbb{#1}}}
\newcommand{\N}{{\mb{N}}}

\newcommand{\R}{{\mb{R}}}

\newcommand{\eps}{\varepsilon}

\newcommand{\D}{\ensuremath{\mathscr D}}

\newcommand{\G}{\ensuremath{\mathcal G}}
\renewcommand{\H}{\ensuremath{\mathcal H}}
\newcommand{\I}{\ensuremath{\mathcal I}}
\renewcommand{\L}{\ensuremath{\mathcal L}}
\newcommand{\M}{\ensuremath{\mathcal M}}
\renewcommand{\O}{\ensuremath{\mathcal O}}

\renewcommand{\k}{\ensuremath{\kappa}}
\newcommand{\y}{\ensuremath{\varrho}}
\newcommand{\Con}{\ensuremath{\mathscr C}}
\newcommand{\Conc}{\ensuremath{\mathscr C}_{c}}
\newcommand{\Cinf}{\ensuremath{\mathscr C^\infty}}
\newcommand{\Cinfc}{\ensuremath{\mathscr C^\infty_{c}}}
\renewcommand{\d}{\ensuremath{\partial}}

\newcommand{\transp}[1]{\ensuremath{\leftidx{^t}{\!#1}{}}}

\newcommand{\rhs}{r.h.s.\@\xspace}

\newcommand{\etc}{etc.\@\xspace}
\newcommand{\eg}{e.g.\@\xspace}
\newcommand{\resp}{resp.\@\xspace}

\newcommand{\pp}{a.e.\@\xspace}

\newcommand{\nhd}{neighborhood\xspace}

\newcommand{\suff}{sufficiently\xspace}

\newcommand{\cst}{{\mathrm{Cst}}}

\newcommand{\Slim}{S}

\DeclareMathOperator{\Char}{Char}

\newcommand{\chart}{\mathcal C}

\newcommand{\hO}{O}
\newcommand{\chdiff}{\phi}
\newcommand{\cdiff}{\phi}

\newcommand{\cdiffL}{\phi_\L}

\newcommand{\tM}{\tilde{\M}}
\newcommand{\hM}{\hat{\M}}
\newcommand{\ThM}{T^* \hM}

\newcommand{\hL}{\hat{\L}}
\newcommand{\ThL}{T^* \hL}

\newcommand{\hatt}{\hat{t}}
\newcommand{\hx}{\hat{x}}

\newcommand{\hxi}{\hat{\xi}}
\newcommand{\hgammaG}{\ensuremath{\leftidx{^{\mathsf G}}{\hat\gamma}{}}}
\newcommand{\chgammaG}{\ensuremath{\leftidx{^{c,\mathsf
        G}}{\hat\gamma}{}}}

\newcommand{\bichar}{bicharacteristic\xspace}
\newcommand{\bichars}{bicharacteristics\xspace}
\newcommand{\bbichar}{broken bicharacteristic\xspace}
\newcommand{\bbichars}{broken bicharacteristics\xspace}
\newcommand{\gbichar}{generalized bicharacteristic\xspace}
\newcommand{\gbichars}{generalized bicharacteristics\xspace}

\newcommand{\gammaB}{\ensuremath{\leftidx{^{\mathsf B}}{\gamma}{}}}

\newcommand{\gammaG}{\ensuremath{\leftidx{^{\mathsf G}}{\gamma}{}}}
\newcommand{\GammaG}{\ensuremath{\leftidx{^{\mathsf G}}{\bar{\gamma}}{}}}
\newcommand{\XG}{\ensuremath{\leftidx{^{\mathsf G}}{X}{}}}

\newcommand{\tgamma}{\tilde{\gamma}}
\newcommand{\tnu}{\tilde{\nu}}
\newcommand{\bld}[1]{\mbox{\boldmath $#1$}}

\newcommand{\Equiv}{\Leftrightarrow}
\newcommand{\imp}{\Rightarrow}

\DeclareMathOperator{\Span}{span}

\newcommand{\Norm}[2]{{\| #1 \|}_{#2}}
\newcommand{\bigNorm}[2]{{\big\| #1 \big\|}_{#2}}

\newcommand{\norm}[2]{{|#1|}_{#2}}
\newcommand{\bignorm}[2]{{\big|#1\big|}_{#2}}

\DeclareMathOperator{\dist}{dist}

\DeclareMathOperator{\supp}{supp}

\newcommand{\scp}[2]{#1 \cdot  #2}

\newcommand{\Bigdup}[2]{\Big\langle #1 , #2 \Big\rangle}
\newcommand{\bigdup}[2]{\big\langle #1 , #2 \big\rangle}
\newcommand{\dup}[2]{\langle #1 , #2 \rangle}

\newcommand{\ovl}[1]{\overline{#1}}

\newcommand{\unitfunction}[1]{\bld{1}_{#1}}

\newcommand{\Rdp}{\R^d_+}

\let \div \relax
\DeclareMathOperator{\div}{div}
\newcommand{\nablag}{\nabla_{\!\! g}}
\newcommand{\divg}{\div_{\! g}}

\DeclareMathOperator{\length}{length}

\newcommand{\para}[1]{\leftidx{^\parallel}{#1}{}}
\newcommand{\compressed}[1]{\leftidx{^c}{#1}{}}

\newcommand{\cT}{\compressed{T}}
\newcommand{\cTL}{\cT^*\L}

\newcommand{\cHb}{\compressed{\Hb}}
\newcommand{\cphi}{\compressed{\phi}}
\newcommand{\cy}{\compressed{\y}}
\newcommand{\cgamma}{\compressed{\gamma}}
\newcommand{\cdist}{\compressed{\!\dist}}

\newcommand{\TM}{T^* \M}
\newcommand{\TL}{T^* \L}

\newcommand{\dTM}{\d(\TM)}
\newcommand{\dTL}{\d(\TL)}
\newcommand{\pdTL}{\para{\dTL}}
\newcommand{\pdTM}{\para{\dTM}}

\newcommand{\pTL}{\para{\TL}}
\newcommand{\pTM}{\para{\TM}}

\newcommand{\pT}{\para{T}}
\newcommand{\pE}{\para{\mathcal E}}
\newcommand{\pG}{\para{\mathcal G}}
\newcommand{\pH}{\para{\mathcal H}}
\newcommand{\pEb}{\para{\mathcal E}_\d}
\newcommand{\pGb}{\para{\mathcal G}_\d}
\newcommand{\pHb}{\para{\mathcal H}_\d}

\newcommand{\ppi}{\pi_{\parallel}}
\newcommand{\pxi}{\para{\xi}}
\newcommand{\py}{\para{\y}}
\newcommand{\pgamma}{\para{\gamma}}
\newcommand{\pe}{\para{e}}

\newcommand{\sd}{{\mathrm{d}}}
\newcommand{\gl}{{\mathrm{3}}}
\newcommand{\sg}{{\mathrm{g}}}

\newcommand{\sdGb}{\G^{\sd}_\d}
\newcommand{\glGb}{\G^{\gl}_\d}
\newcommand{\sgGb}{\G^{\sg}_\d}

\newcommand{\Hb}{\H_\d}
\newcommand{\Gb}{\G_\d}

\DeclareMathOperator{\Hamiltonian}{H}
\newcommand{\Hp}{\Hamiltonian_p}
\newcommand{\HpG}{\Hamiltonian_p^{\G}}
\newcommand{\Hpz}{\Hp z}
\newcommand{\Hppz}{\Hp^2 z}
\newcommand{\Hz}{\Hamiltonian_z}
\newcommand{\Hzp}{\Hz p}
\newcommand{\Hzzp}{\Hz^2 p}

\newcommand{\n}{\mathsf n}
\newcommand{\nx}{\n_x}
\newcommand{\nxs}{\n^*_x}
\newcommand{\nxsC}{\n^{*, \chart}_x}

\newcommand{\jp}[1]{\langle #1 \rangle}

\newcommand{\U}{\mathscr U}

\newcommand{\tchi}{\tilde{\chi}}

\theoremstyle{plain}
\newtheorem{theorem}{Theorem}[section]
\newtheorem{proposition}[theorem]{Proposition}
\newtheorem{lemma}[theorem]{Lemma}
\newtheorem{corollary}[theorem]{Corollary}
\newtheorem{assumption}[theorem]{Assumption}

\newtheorem{theorembis}{Theorem}
\newtheorem{corollarybis}[theorembis]{Corollary}
\newcounter{theorembiss}

\theoremstyle{definition}
\newtheorem{definition}[theorem]{Definition}

\newtheorem{remark}[theorem]{Remark}

\newcommand{\et}{\ensuremath{\text{and}}}
\newcommand{\avec}{\ensuremath{\text{with}}}
\newcommand{\si}{\ensuremath{\text{if}}}

\newcommand{\dans}{\ensuremath{\text{in}}}

\newcommand{\where}{\ensuremath{\text{where}}}

\newcommand*\quot[2]{{^{\textstyle #1}\big/_{\textstyle #2}}}
\numberwithin{equation}{section}


\begin{document}
\title[Measure and continuous vector fields] 
{Measure and continuous vector field at a boundary II:
  geodesics and support propagation}

%% authors
\author{Nicolas Burq}
\address{Nicolas Burq. Laboratoire de Math\'ematiques d'Orsay, Universit\'e
  Paris-Sud, Universit\'e Paris-Saclay, B\^atiment~307, 91405
  Orsay Cedex \& CNRS UMR 8628 \& Institut Universitaire de France}
 \email{nicolas.burq@math.u-psud.fr}
\author{Belhassen Dehman}
  \address{Belhassen Dehman. Universit\'e de Tunis El Manar, Facult\'e
  des Sciences de Tunis, 2092 El  Manar \& Ecole Nationale d'Ing\'enieurs de Tunis, ENIT-LAMSIN, B.P. 37, 1002 Tunis, Tunisia. }
\email{Belhassen.Dehman@fst.utm.tn}
 \author{J\'er\^ome Le Rousseau}
 \address{J\'er\^ome Le Rousseau.
   Universit\'e Sorbonne Paris Nord, Laboratoire Analyse, G\'eom\'etrie et Applications, LAGA, CNRS, UMR 7539, F-93430, Villetaneuse, France.}
\email{jlr@math.univ-paris13.fr}
\date{\today}

\begin{abstract}
Nonnegative measures that are solutions to a transport equation with
{\em continuous} coefficients have been widely studied.  Because of
the low regularity of the associated vector field, there is no natural
flow since nonuniqueness of integral curves is the general rule. It
has been known since the works by L.~Ambrosio~\cite{Ambrosio:08} and L.~Ambrosio and G.~Crippa~\cite{Ambrosio-Crippa:08,AmCr14}
 that such measures can be described
as a {\em superposition} of $\delta$-measures supported on integral
curves. In this article, motivated by some {\em observability}
questions for the {\em wave equation}, we are interested in such
transport equations in the case of {\em domains with
  boundary}. Associated with a wave equation with
$\Con^1$-coefficients are \bichars that are integral curves of a
continuous Hamiltonian vector field.  We first study in
details their behaviour in the presence of a boundary and define their
natural generalisation that follows the laws of geometric optics. Then,
we introduce a natural class of transport equations with a source term on
the boundary, and we prove that any {\em nonnegative} measure
satisfying such an equation has a union of maximal
\gbichars for support. This result is a weak form of the {superposition principle}
in the presence of a boundary. With its companion article
\cite{BDLR1}, this study completes the proof of wave observability
generalizing the celebrated result of Bardos, Lebeau, and
Rauch~\cite{BLR:92} in a low regularity framework where coefficients
of the wave equation (and associated metric) are $\Con^1$ and the
boundary and the manifold are $\Con^2$.
\end{abstract} 
\maketitle
\tableofcontents
%%%%%%%%%%%%%%%%%
% section
%%%%%%%%%%%%%%%%%
\section{Introduction}

We are interested in studying the properties of
nonnegative measures solutions of ordinary differential equations.
If $X$ is a {\em smooth} vector field, in say $\R^d$, and if a measure $\mu$
solves the free transport equation 
\begin{align}
  \label{eq: free transport equation}
  \transp{X} \mu =0,
\end{align}
then 
\begin{align}
  \label{statement: measure transport}
  \text{\em the measure} \  \mu \ \text{\em is invariant along the integral
  curves of}\  X.
\end{align}
Such a transport equation   arises naturally in many
contexts. Our initial motivation lies in the observation of waves. Understanding propagation
properties of measures fufilling an equation like \eqref{eq: free transport equation} has become a key point to obtain an observability
inequality, that is, an estimate of the energy by a ``recording'' of the
solution in a restricted domain for some time $T>0$, and, as a corollary, an exact controllability property for the wave
equation. We refer  to \cite{BLR:92,BG:1997} for this topic where microlocal techniques allow one to find the proper geometrical
condition for wave observability to hold in connection with measure transport
along the geodesic flow. This condition known as the Geometric
Control Condition (GCC) roughly states: 
\begin{center}
  {\em every geodesic  should reach the observation region in the
    prescribed time $T>0$.}
\end{center}
In this framework, the vector field that plays the role of $X$ in \eqref{statement:
  measure transport} is the Hamiltonian vector field
$\Hp$ associated with the principal symbol $p$ of the wave operator. It is
precisely the vector field that generates {\em \bichars} in phase-space and  geodesics are
their projection on the base space. 

In the case of a domain with a boundary, generalized geodesics have
to be  considered. They follow the
reflections laws of optics at boundaries and more complex behavior when
reaching the boundary tangentially. In a smooth setting, their mathematical study in
connection with the propagation of singularities for waves was carried
out in the work by R.~Melrose and J.~Sj\"{o}strand
\cite{MS:78}; see also \cite[Chapter 24]{Hoermander:V3}. Analyzing a transport equation for measures becomes more subtle if the boundary in
encountered. In the transport equation we consider there is  a source term associated with a
second nonnegative measure supported in the boundary. It takes the form
\begin{align}
    \label{eq: Gerard-Leichtnam equation-intro}
    \transp{\Hp} \mu  = -\int_{\y \in \pHb \cup \pGb} 
    \frac{\delta_{\y^+} - \delta_{\y^-}}
    {\dup{\xi^+- \xi^-}{\nx}_{T_x^*\M, T_x\M} 
    } \ d \nu (\y), 
  \end{align}
where $\nu$ is a nonnegative measure on (the cotangent bundle of) the
boundary. All terms in \eqref{eq: Gerard-Leichtnam equation-intro} are
explained in Sections~\ref{sec: Geometrical setting-intro} and
\ref{sec: main result} below.  Such an equation was first derived in
the work of P.~G\'erard and \'E.~Leichtnam~\cite{GL:1993}. In the case
of smooth coefficients, the transport of the measure $\mu$ can still
be well understood along the generalized geodesic flow; see
\cite{Lebeau:96} where a slightly different measure is introduced. In
our main results below, we consider a generalized version of this
equation by adding a term $f \mu$ with $f$ a continuous real function.

\medskip In the framework of wave observability, the issue of limited
smoothness of coefficients has received little attention in works
relying on microlocal techniques as those cited above. The main reason
is that using such techniques is known to consume many
derivatives. In \cite{BDLR0}, the derivation of
an observability estimate is carried out for waves in the case of
rough ($\Con^1$) coefficients on a manifold {\em without} boundary.  A
transport equation for a nonnegative measure is also at the heart of
the argument, yet with a vector field $X= \Hp$ which is only {\em
  continuous}.  This transport equation takes the simple form of
\eqref{eq: free transport equation} since there is no boundary in the
setting of \cite{BDLR0}. For a regularity as low as
Lipschitz, or $\log$-Lipschitz, for the vector field $X$, the
Cauchy-Lipschitz formula and its extensions, yield a unique flow for
$X$ and one obtains the same invariance result as in \eqref{statement: measure
  transport}.  If now $X$ is a $\Con^0$-vector field, while the Peano
theorem yields the existence of integral curves for $X$, one cannot
guarantee the uniqueness of those integral curves and the existence of
a flow. This compromises the extension of the transport
result in \eqref{statement: measure transport}. Still, a remarkable quantitative result is achieved in the work
of L.~Ambrosio and G.~Crippa \cite{Ambrosio:08,Ambrosio-Crippa:08,AmCr14}: the measure $\mu$ solution
to \eqref{eq: free transport equation} can be written as a sum of
positive measures, each defined as a constant times $\delta_\gamma$,
that is, the measure supported by an integral curve $\gamma$ of $X$
and invariant along $\gamma$. The sum is defined by means of a
nonnegative measure on the space of all continuous curves. The result of
\cite{AmCr14} is thus a {\em superposition principle}.

The result of \cite{AmCr14} relies on a smoothing procedure
 of the vector field $X$ as in \cite{Evans-Gangbo:99} that we have not been
able to extend in the presence of a boundary.
%at a
%boundary in the case of the Hamiltonian vector field\footnote{Or
%  rather the 
%extension of $\Hp$ that differs from $\Hp$ at the boundary in the so-called glancing
%region. This extension is properly recalled in what follows.}  that interests us. 
To prepare the
work exposed in the present article, in \cite{BDLR0}  we proved a
weaker version of the superposition principle of \cite{AmCr14}, that
is,
\begin{align*}
  %\label{eq: weak version Ambrosio}
 \supp \mu \ \text{\em is a union of maximal integral
curves of} \ X,
 \end{align*}
if $\mu$ is a solution to the transport equation~\eqref{eq: free transport equation}. Our proof scheme is much different from that of
\cite{AmCr14} and is inspired by the microlocal techniques used in
\cite{MS:78}.  In the present article we  extend this latter proof scheme to regions at the
boundary. We prove that if a nonnegative mesure $\mu$ 
fulfills a transport equation of the form given in \eqref{eq: Gerard-Leichtnam equation-intro}
then 
\begin{align}
  \label{statement: measure support transport}
  \supp \mu \ 
  \text{\em  is a union of maximal \gbichars}.
\end{align}
From this result, an observability estimate for waves is deduced in
the companion article \cite{BDLR1} in the case of
$\Con^1$-coefficients for the wave operator under a geometric control
condition, using a transport equation of the form of \eqref{eq:
  Gerard-Leichtnam equation-intro} also proven in \cite{BDLR1}. The
main result of the present article, Theorem~\ref{theorem: measure
  support propagation} below that states \eqref{statement: measure
  support transport} in details, is thus the missing link to complete
the proof of wave observability in the case of $\Con^1$-coefficients
and the presence of a ($\Con^2$) boundary thus generalizing the
celebrated result of Bardos, Lebeau, and Rauch~\cite{BLR:92}.

\medskip
Note that the regularity we consider, that is $\Con^1$-coefficients
for the wave operator, stands as a limit case, as far as the
underlying geometry is concerned, since the Hamiltonian vector field
$\Hp$ that generates \bichars is then only continuous. If considering
coefficients with lower regularity, one obtains a Hamiltonian vector
field $\Hp$ that may not be continuous. Then, the mere existence of
\bichars and geodesics is not guaranteed (and the meaning
of~\eqref{eq: free transport equation} neither)

\subsection{Towards a superposition principle}
With the result on the measure support stated in \eqref{statement:
  measure support transport} and proven here, a quantitative version
of this result in the spirit of the superposition principle of
L.~Ambrosio and G.~Crippa \cite{Ambrosio:08,Ambrosio-Crippa:08,AmCr14}
appears now to be a very interesting open question.  A first step in
this direction is given here: we prove that the Dirac-measure
$\delta_{\GammaG}$ supported on a maximal \gbichar $\GammaG$ obeys a
transport equation of the form of \eqref{eq: Gerard-Leichtnam
  equation-intro}. Not having this property would ruin any hope to
obtain a superposition principle, since any sum of such Dirac-measures
automatically fulfills an equation of this type. A precise statement
and a proof is given in Appendix~\ref{sec: proof prop mesure single
  bichar}.

\subsection{Perspectives in the presence of a flow}

The regularity level we consider here yields a continuous Hamiltonian
vector field. As mentionned above, this suffices for \gbichars to
exist but lack of uniqueness prevents the existence of a flow {\em in
  general}. One cannot however exclude that uniqueness holds true and
a flow exists in some particular cases. Does the result presented here
improve in such cases? We foresee that it does in the following strong
form: a measure solution to the transport equation~\eqref{eq:
  Gerard-Leichtnam equation-intro} is actually transported along the
\gbichar flow. This is the subject of an ongoing work \cite{BDLR3}.

A second perspective concerns the case where the Hamiltonian vector
field $\Hp$ is continuous and moreover lies in the Sobolev class
$W^{1,1}$. Then, as $\Hp$ is naturally divergence free, this fits the
setting of the celebrated article of R.J.~DiPerna and
P.L.~Lions~\cite{DiPerna-Lions}. Then, a flow exists in a weak sense,
with flow properties fulfilled almost everywhere. Can one obtain a
weak form of transport along such flow for (a large class of) measures
solutions to \eqref{eq: Gerard-Leichtnam equation-intro}? This stands
as a very interesting and natural open question. Extension to more
general vector fields as in \cite{Ambrosio:04} is of great interest also. 

%%%%%%%%%%%%%%%%%
% sub-section
%%%%%%%%%%%%%%%%%
\subsection{Outline and notation}

The present article is organized as follows. In Section~\ref{sec:
  Geometrical setting-intro}, the geometrical notions necessary to
the statement of our results are presented and our main result, that is, a
precise formulation of the propagation of the measure support given in
\eqref{statement: measure support transport}, is stated in  Section~\ref{sec: main
  result}.  Section~\ref{sec:
  measure support propagation1} covers the case of a transport equation in two cases: 
(1) away from any boundary and (2) a boundary transverse to the
vector field. A full treatment of the boundary requires the
introduction of additional geometrical notions; this is done in
Section~\ref{sec: setting at the boundary}. The proof of the main
result on the propagation of the measure support is carried out in
Section~\ref{sec: Propagation of the measure support}. A mass property
of the boundary measure $\nu$ that can be deduced from the transport
equation \eqref{eq: Gerard-Leichtnam equation-intro} is proven in
Section~\ref{sec: proof prop: no mass on Gd G3}.

In Appendix~\ref{sec: proof prop mesure single bichar} we prove that
the linear measure supported by a single \gbichar fulfills an equation
of the form of \eqref{eq: Gerard-Leichtnam equation-intro}.
Appendix~\ref{sec: Existence and continuity properties of generalized
  bicharacteristics} provides a proof of the existence of
\gbichars. Despite nonuniqueness a continuity property is also
given. In Appendix~\ref{sec proof prop: quasi-normal coordinates} we
present the quasi-normal coordinates that are often used near the
boundary.

\medskip
Often for the sake of concision, the Einstein summation
convention for repeated indices are used. 
In local coordinates $B(x,r)$ denotes the Euclidean open ball centered at $x$ with radius~$r$.

%%%%%%%%%%%%%%%%%
% subsection
\subsection{Acknowledgements}
This research was partially supported by Agence Nationale de la
Recherche through project ISDEEC ANR-16-CE40-0013 (NB), by the
European research Council (ERC) under the European Union’s Horizon
2020 research and innovation programme, Grant agreement 101097172 -
GEOEDP (NB), and by the Tunisian Ministry for Higher Education and
Scientific Research within the LR-99-ES20 program (BD). The authors
acknowledge GE2MI (Groupement Euro-Maghr\'ebin de Math\'ematiques et
de leurs Interactions) for its support.
They also thank F.~Santambroggio for drawing their attention to
the work by L. Ambrosio and
G. Crippa~\cite{Ambrosio:08,Ambrosio-Crippa:08,AmCr14} at some early
stage of their study.

%%%%%%%%%%%%%%%%%
% sub-section
%%%%%%%%%%%%%%%%%
\section{Geometrical setting I}
\label{sec: Geometrical setting-intro}

Consider a $\Con^{2}$ compact connected $d$-dimensional manifold
$\M$ with boundary equipped with a $\Con^1$-Riemannian  metric $g$. 
In this section we introduce some geometrical notions that are necessary for the
statement of our main results. Many more details are available in
Section~\ref{sec: setting at the boundary}.

A simple example that fits our present setting is that of 
a bounded open subset $\Omega$ of $\R^d$ with a $\Con^{2}$-boundary,
that is, with the boundary given locally by $\varphi(x)=0$ with
$\varphi \in \Con^{2}(\R^d)$ and $d \varphi \neq 0$. Then $\M = \Omega \cup \d\Omega$ and one can simply consider the Euclidean
metric.  In the spirit of this simple example, we consider an open
$d$-dimensional manifold\footnote{The manifold $\tM$ can be
  constructed by embedding $\M$ in $\R^{2d}$ thanks to the
  Whitney theorem \cite{Whitney:36}.} $\tM$ such that $\M \subset \tM$ and 
extend the metric $g$ to a \nhd of $\M$ is a $\Con^1$-manner.

\subsection{Local coordinates}
\label{sec: Local coordinates}
Equip a compact \nhd $\hM$ of $\M$ in $\tM$ with a finite
$\Con^2$-atlas. A local chart is denoted $(\hO, \chdiff)$ with $\hO$
an open subset of $\hM$ and $\chdiff$ a one-to-one map from $\hO$ onto an open
subset of $\R^d$. Charts can be chosen so that 
\begin{align}
  \label{eq: local chart boundary}
  &\chdiff(\hO \cap \M) = \chdiff(\hO) \cap \{x_d\geq 0\}
    \ \text{is an open subset of}\  \ovl{\Rdp}, \\
  &\chdiff(\hO \cap \d\M) = \chdiff(\hO) \cap \{x_d=0\},
    \  \et \  
    \chdiff(\hO \setminus \M) = \chdiff(\hO) \cap \{x_d<0\},
    \notag
\end{align}
if $\hO \cap \d\M \neq \emptyset$. Denote the local coordinates by $x=(x', x_d)$ with $x'
\in \R^{d-1}$. Note that $\M$ being compact it
contains its boundary $\d\M$.

In a  local chart, the metric $g$ is
given by $g_x = g_{ij} (x) dx^i \otimes dx^j$, where  $g_{ij}\in
\Con^1(\chdiff(\hO))$. We use the classical notation $(g^{ij}(x) )_{i,j}$ for  the
inverse of  $(g_{ij}(x) )_{i,j}$. 
The metric $g_x = (g_{ij}(x) )_{i,j}$ provides an inner  product
on $T_x \M$. The metric $g_x^* = g^{ij}(x) d \xi_i \otimes d \xi_j$ provides an
inner product
on $T^*_x \M$, denoted $g^*_x(\xi,\tilde{\xi})$, for $\xi, \tilde{\xi} \in T^*_x
\M$. Define the associated norm
\begin{align*}
  \norm{\xi}{x} = g^*_x(\xi,\xi)^{1/2}.
\end{align*}
In this introductory section, near a boundary point, local coordinates
are chosen according to the following proposition as they simplify the
exposition of some geometrical notions.
%%%%%%%%%%%%%%%%%%%%%%%%
% proposition          %
%%%%%%%%%%%%%%%%%%%%%%%%
\begin{proposition}[quasi-normal geodesic coordinates]
  \label{prop: quasi-normal coordinates}
  Suppose $m^0 \in \d\M$. There exists a
  $\Con^2$-local chart $(\hO, \chdiff)$ such that $m^0 \in \hO$,  $\chdiff(m) = (x',z)$, with $x' \in \R^{d-1}$ and $z \in \R$, and 
  \begin{enumerate}
   \item $\chdiff(\hO\cap \M) =  \{ z\geq 0\} \cap \chdiff(\hO)$, $\chdiff(\hO\cap
      \d\M) =  \{ z= 0\} \cap \chdiff(\hO)$,  and $\chdiff(\hO\setminus \M) =  \{ z<0\} \cap \chdiff(\hO)$ ;
    \item at the boundary, the representative of the metric has the form 
      \begin{align*}
        g(x',z=0) = \sum_{1\leq i, j \leq d-1} g_{i j}(x',z=0)  dx^i \otimes dx^j
        + |dz|^2.
      \end{align*}
  \end{enumerate}
\end{proposition}
In other words the matrix of $g=(g_{ij})$ has the block-diagonal form
   {\em at the boundary}
   \begin{equation}
     \label{eq: structure metric quasi-geodesic coordinates}
      g (x',z=0) = \begin{tikzpicture}[
  every left delimiter/.style={xshift=.75em,yshift=0.4em},
  every right delimiter/.style={xshift=-.75em,yshift=0.4em},
  style1/.style={
  matrix of math nodes,
  every node/.append style={text width=#1,align=center,minimum height=1.8ex},
  nodes in empty cells,
  left delimiter=(,
  right delimiter=)
  },
baseline=(current  bounding  box.center)
]
\matrix[style1=0.27cm] (1mat)
{
  & &  & \\
  &  & &\\
  &  & &\\
  & &  &  \\
  & &  &  \\
};
\node[font=\normalsize] 
  at (1mat-1-4.north) {$0$};
  \node[font=\normalsize] 
  at (1mat-2-4.north) {$\vdots$};
\node[font=\normalsize] 
  at (1mat-3-4) {$0$};
\node[font=\normalsize] 
  at (1mat-4-4.south) {$1$};
 \node[font=\normalsize] 
  at (1mat-4-1.south) {$0$};
\node[font=\normalsize] 
  at (1mat-4-2.south) {$\ \cdots$};
  \node[font=\normalsize] 
  at (1mat-4-3.south) {$\ 0$};
  \node[font=\Huge] 
  at (1mat-2-2) {$*$};
\end{tikzpicture}.
    \end{equation}
  Naturally, the same form holds for $g_x^* =
  (g^{ij}(x))$ at the boundary. 
One deduces that 
  \begin{align*}
    g_{j d}(x',z)  = z h_{j d}(x',z)  \ \ \et \ \ 
    g_{d d}(x',z)  = 1 +  z h_{d d}(x',z),
  \end{align*}
  for some continuous functions $h_{j d}$, $j=1, \dots, d$.

  Proposition~\ref{prop: quasi-normal coordinates} can be found in
  \cite{BM:21} with a different regularity level. We provide a proof
  in Appendix~\ref{sec proof prop: quasi-normal coordinates} based on
  that of \cite{BM:21} with a generalization to other levels of
  regularity.
  %%%%%%%%%%%%%%%%%%%%%%%%
% remark               %
%%%%%%%%%%%%%%%%%%%%%%%%
\begin{remark}
  \label{remark: no normal geodesic coodinates}
  Because of the low regularity of $g$ and $\M$ one {\em
    cannot choose} normal geodesic coordinates, that is, local
  coordinates for which $g_{j d} = g_{d j} = 0$ for $j \neq d$ and
  $g_{d d}=1$ near a point $m^0$ of the boundary. The
  coordinates that Proposition~\ref{prop: quasi-normal coordinates}
  provides only have this property in a \nhd of $m^0$ {\em within} the
  boundary $\d\M$.
\end{remark}

\medskip One sets $\L = \R \times \M$ and
$\hL = \R \times \hM$ .  From a local chart $(\hO, \chdiff)$
in the atlas for $\hM$ one defines a map
$\cdiffL: (t,m) \mapsto (t, \chdiff(m))$ from $\O = \R\times \hO$ onto
$\R\times \chdiff(\hO)$, yielding a local chart $(\O, \cdiffL)$  for $\hL$ and thus a
finite atlas.

For $x = \chdiff(m)$, $m \in \hO \cap \M$, denote by $v = (v',v^d)$
and $\xi =(\xi',\xi_d)$ the associated coordinates in $T_m \M$ and
$T_m^* \M$, with $v', \xi'\in \R^{d-1}$ and $v^d, \xi_d \in \R$.  We
write $T_x\M$ and $T_x^*\M$ by abuse of notation.  In what follows, it
will be convenient to write $z$ in place of $x_d$, in particular for
the local coordinates given by Proposition~\ref{prop: quasi-normal
  coordinates}. Accordingly we shall denote the associated cotangent
variable $\xi_d$ by the letter $\zeta$, that is, $\xi= (\xi', \zeta)$.
We however do not change the notation for the associated tangent
variable $v^d$.  With local charts at the boundary given by
Proposition~\ref{prop: quasi-normal coordinates}, if $x \in \d\M$ and
$v \in T_x\d\M$ then $v = (v',0)$ and we use the bijective map
$(\xi',0) \mapsto \xi'$ to parameterize $T_x^* \d\M$.

Also classically set
\begin{align*}
  T \M = \bigcup_{x \in \M} \{x\} \times T_x \M, 
  \quad \TM = \bigcup_{x \in \M} \{x\} \times T_x^* \M\\
  \big(\text{\resp}\ T \hM = \bigcup_{x \in \M} \{x\} \times T_x \hM, 
  \quad \ThM = \bigcup_{x \in \M} \{x\} \times T_x^* \hM\big).
\end{align*}
With $\M$ containing its boundary $\d\M$, one sees that $ T \M$ (\resp
$\TM$) contains  $\{x\} \times T_x \M$ (\resp $\{x\} \times T_x^* \M$)
for $x \in \d\M$. We denote by $\dTM$ the boundary of $\TM$ that is
the set of $(x,\xi)$ with $x \in \d\M$. 
In the local coordinates, $\dTM$ is
given by $\{ z=0\}$ and $\TM$ by $\{ z\geq 0\}$.

In the associated local chart on $\L$, the representative of $(t,m)
\in \L$ is $(t,x) = (t,x',z)$.  We shall use the letter $\y$ to denote
an element of $\TL$, that is, $\y = (t,x; \tau, \xi)$ with $(t,x) \in
\L$, $\tau \in \R$ and $\xi \in T^*_x\M$.  Classically, we write $\TL
\setminus 0$ for the set of points $\y = (t,x; \tau, \xi)$ with
$(\tau, \xi)\neq0$.  The boundary $\dTL$ is the set of points $\y=
(t,x;\tau,\xi)$ such that $x \in \d\M$. Note that $\dTL$ is locally
given by $\{ z=0\}$ and $\TL$ is locally given by $\{ z\geq 0\}$.

\subsection{Wave operator and bicharacteristics}
\label{sec: wave operator}

On the manifold $\M$ consider the elliptic operator $A = A_{\k,g}= \k^{-1} \divg (\k
\nablag)$, that is, in local coordinates
\begin{align*}
  A f = \k^{-1} (\det g)^{-1/2} 
  \sum_{1\leq i, j \leq d}  \d_{x_i}\big( 
  \k  (\det g)^{1/2} g^{i j}(x)\d_{x_j} f
  \big).
\end{align*}
Its principal symbol is simply
$a(x,\xi) = - g^*_x(\xi,\xi)  = - g_x^{i j} \xi_i \xi_j = - \norm{\xi}{x}^2$. Note
that for $\k=1$, one has $A = \Delta_g$, the Laplace-Beltrami
operator associated with $g$ on $\M$. Together with $A$ consider
the wave operator $ P_{\k,g} = \d_t^2 - A_{\k,g}$. Its
principal symbol in a local chart is given by 
\begin{align*}
  p(\y) = -\tau^2 + \norm{\xi}{x}^2.
\end{align*}
Note that $p(\y)$ is smooth in the variables $(\tau,\xi)$ and $\Con^1$
in $x$. 

For a function $f$ of the variable $\y$, the Hamiltonian vector field
$\Hamiltonian_f$ is defined by $\Hamiltonian_f(h)  = \{ f, h\}$, where
$\{.,.\}$ is the Poisson bracket. 
In local coordinates one has 
\begin{align}
  \label{eq: Hp}
  \Hp (\y) &= \d_\tau p (\y) \d_t 
             + \nabla_\xi p(\y) \cdot \nabla_x
             -\nabla_x p (\y) \cdot \nabla_\xi \\
           &=  - 2 \tau \d_t + 2 g^{ij}(x) \xi_i \d_{x_j} 
             - \d_{x_k} g^{ij}(x) \xi_i \xi_j  \d_{\xi_k}.\notag
\end{align}
Recall the following definition. 
%%%%%%%%%%%%%%%%%%%%%%%%
% definition           %
%%%%%%%%%%%%%%%%%%%%%%%%
\begin{definition}
  \label{def: bichar}
  Suppose $V$ is an open subset of $\TL \setminus \dTL$ and $J\subset \R$ is an interval. A $\Con^1$-map $\gamma: J \to V \cap \Char p$ is called a
  \bichar in $V$ if 
  \begin{align}
    \label{eq: def: bichar-intro}
    \frac{d}{d s} \gamma(s) = \Hp \big( \gamma(s) \big), 
    \qquad  s \in J. 
  \end{align}
  It is called {\em maximal} in $V$ if it cannot be extended by
  another \bichar also valued in~$V$.
\end{definition}
\begin{remark}
  \label{rem:def-bichar}
   \leavevmode  
\begin{enumerate}
\item
  \label{rem: constant tau}
  If $\gamma(s) = \big(t(s), x(s), \tau(s), \xi(s)\big)$ is a \bichar observe
  that $\tau(s)$ is constant because of the form of $\Hp$.
  Since $\gamma(s) \in \Char p$ one has 
  \begin{align*}
    \norm{\xi(s)}{x(s)} = |\tau(s)| 
  \end{align*}
  also constant along a \bichar.
\item With $\tau(s) = \cst \neq 0$ then 
  $t(s)$ is an affine function of $s$ since $d t / d s = - 2 \tau(s)$.  
  If $V = \TL \setminus \dTL$ and $\gamma$ is maximal with $S^+ = \sup J < +\infty$
  then 
  \begin{align*}
    \lim_{{s \to S^+} \atop {s \in J}} \gamma(s) \in \dTL.
  \end{align*}
  This is part of Lemma~\ref{lemma: limit of a maximal bichar from the
    interior} below.
\item
  \label{rem: bichar in hL}
  If needed we also call \bichar a map $\gamma: J \to \ThL \cap
    \Char p$, such that \eqref{eq: def: bichar-intro} holds, using
    the extension of $\Hp$ in $\ThL$. This is then a \bichar above
    $\hL$ and we mention it explicitly.
\end{enumerate}
\end{remark}
Note that $\transp{\Hp} f (\y) = 2 \tau \d_t f (\y)  - 2 \d_{x_j}
  \big(g^{ij}(x) \xi_if (\y)  \big) 
  + \d_{\xi_k} \big(\d_{x_k} g^{ij}(x) \xi_i \xi_j f (\y) \big)$ and 
deduce 
\begin{align*}
  %\label{eq: transp Hp}
  \transp{\Hp}  = - \Hp.
\end{align*} 
Recall also that 
\begin{align}
  \label{eq: derivative along bichar}
  \Hp f (\gamma(s)) = \frac{d}{d s} f (\gamma(s)),
  \ \ \si \ \gamma \ \text{is a \bichar}. 
\end{align}

\subsection{A partition of the cotangent bundle at the boundary}
\label{sec: A partition of the cotangent bundle at the boundary-intro}

Denote  by $\pdTL \subset \dTL$ the bundle of points $\y = (\y',0)= 
(t,x',z=0,\tau,\xi',0) \in \TL$ for $\y' =
(t,x',z=0,\tau,\xi') \in T^* \d\L$. 
 Identifying $\y'$ and $(\y',0)$ as presented
above  thanks to the chosen  local coordinates  allows one to
indentify $\pdTL$ and $T^* \d\L$.

Denote by $\ppi$ the map from $\dTL$ into $\pdTL$ given by
\begin{align*}
  \ppi (t,x',z=0,\tau,\xi',\zeta) = (t,x',z=0,\tau,\xi',0).
\end{align*}
%%%%%%%%%%%%%%%%%%%%%%%%
% definition           %
%%%%%%%%%%%%%%%%%%%%%%%%
\begin{definition}[elliptic, glancing, and hyperbolic regions]
  \label{def: E', H', G'-intro}
  One partitions $\pdTL$  into three homogeneous regions.
  \begin{enumerate}
  \item The elliptic region $\pEb = \pdTL \cap \{ p >0\}$; 
    if
    $\y \in \pEb$ it is called an elliptic point.          	
  \item The glancing region $\pGb = \pdTL \cap \{ p =0\}$;   if $\y
    \in \pGb$ it is called a glancing point.
  \item The hyperbolic region $\pHb =\pdTL \cap \{ p < 0\}$; if
    $\y \in \pHb$ it is called a   hyperbolic point.
  \end{enumerate}
\end{definition}
Since $p (\y) = -\tau^2 + \zeta^2 + g_x(\xi',\xi')_x$ by \eqref{eq: structure metric quasi-geodesic coordinates} if $\y \in
\dTL$, one has the following properties:
\begin{enumerate}
  \item If $\y \in \pEb$ then $\ppi^{-1}\big( \{ \y\}\big) \cap \Char p
    = \emptyset$.
  \item  If $\y \in \pGb$ then $\ppi^{-1}\big( \{ \y\}\big) \cap \Char p
    = \{ \y\}$.
  \item If $\y = (t,x',z=0,\tau,\xi',0) \in \pHb$ then $\ppi^{-1}\big( \{ \y\}\big) \cap \Char p
    = \{ \y^-, \y^+\}$, where
    \begin{align}
      \label{eq: relevement H-G-intro}
      \y^\pm = (t,x',z=0,\tau,\xi^\pm), \ \ \where \ \xi^\pm=(\xi',\zeta^\pm)
       \   \avec \  \zeta^\pm = \pm \sqrt{- p(\y)}.
    \end{align}
\end{enumerate}
Associated with the previous partition of $\pdTL$  is a partition of
$\Char p \cap \dTL$. Indeed, if $\y \in \Char p \cap \dTL$ then
$\ppi(\y) \in \pdTL$ and 
$p \big( \ppi(\y)\big ) \leq 0$. 
Note that having $\y \in \Char p \cap \dTL$ and $p \big( \ppi(\y)\big
) = 0$ is equivalent to having $\y \in \pGb$.
%%%%%%%%%%%%%%%%%%%%%%%%
% definition           %
%%%%%%%%%%%%%%%%%%%%%%%%
\begin{definition}[partition of $\Char p$ at the boundary]
  \label{def: H, G-intro}
  One partitions $\Char p \cap \dTL$  into two homogeneous regions
  $\Gb$ and $\Hb$:
  \begin{enumerate}
  \item $\Gb = \pGb$; $\y \in \Gb \ \Equiv \ \y \in \Char p$ and $\ppi(\y) = \y$.
  \item $\y \in \Hb$ if $\y \in \Char p$ and $\ppi(\y) \in \pHb$. It
    is also called a hyperbolic point. If $\y =
    (t,x',z=0,\tau,\xi',\zeta)$ one says that $\y  \in \Hb^+$ if $\zeta>0$ and
    $\y  \in \Hb^-$ if $\zeta<0$.
    \end{enumerate}
\end{definition}
Thus, if $\y \in \pHb$ then $\ppi^{-1}\big( \{ \y\}\big) \cap \Char p
= \{ \y^-, \y^+\}$ with $\y^+ \in \Hb^+$ and $\y^- \in \Hb^-$, with
$\y^\pm$ as given in \eqref{eq: relevement H-G-intro}.

Introducing the following involution on $\dTL$
\begin{align}
  \label{eq: involution-intro}
  \Sigma(t,x',z=0,\tau,\xi',\zeta) = (t,x',z=0,\tau,\xi',-\zeta),
\end{align}
one finds that $\Sigma (\y^-) = \y^+$ if $\y \in \pHb$. Thus, $\Sigma$
is a one-to-on map from  $\Hb^-$ onto $\Hb^+$.

\subsection{Glancing region, gliding vector field, and generalized bicharacteristics}
\label{sec: glancing region, generalized bichar-intro}

Recall that we denote by $z$ the variable $x_d$.  One computes
\begin{align*}
  \Hpz (\y) = \Hpz (x,\xi)= 2 g^{dj}(x) \xi_j.
\end{align*}
Observe that $\Hpz$ is a $\Con^1$-function. Note that 
$\Hpz_{|z=0} = 2 \zeta$ in the present local coordinates. Hence, $\pGb= \Gb=\{
z=\Hpz=p=0 \}$ and $\Hb^\pm =\{ z=p=0, \ \Hpz \gtrless 0\}$ locally. With \eqref{eq: derivative along bichar} this means that a
\bichar going through a point $\y\in \Hb$ has a contact of order
exactly one with the boundary: it is transverse to $\dTL$. A \bichar
going through a point $\y\in \Gb$ has a contact of order greater than or
equal to two:  it is tangent to $\dTL$.

One can further compute $\Hppz$.  It is continuous and gives the
following partition of $\Gb$. 
%%%%%%%%%%%%%%%%%%%%%%%%
% definition           %
%%%%%%%%%%%%%%%%%%%%%%%%
\begin{definition}[partition of $\Gb$]
  \label{def: diffrative}
  Introduce
  \begin{align*}
    &\sdGb = \{ \y \in \Gb; \Hppz (\y) >0\},\\
    &\glGb =\{ \y \in \Gb; \Hppz (\y) =0\}, \\ 
    &\sgGb =\{ \y \in \Gb; \Hppz (\y) <0\}.
  \end{align*}
    One calls $\sdGb$ the diffractive set, $\sgGb$ the gliding set.
    One calls $\glGb$ the glancing set of order three: if $\y^0 \in
    \glGb$ a \bichar that goes through $\y^0$ has a contact with the
    boundary of order greater than or equal to three.
\end{definition}
\ref{def: diffrative}
On $\pdTL$ one defines 
\begin{align} 
  \label{sec: def gliding vector field into}
    \HpG (\y)  = \Big(\Hp + \frac{\Hppz}{\Hzzp}\Hz\Big)  (\y) ,
  \end{align}
referred to as the gliding vector field. Note that in the present
coordinates one has $\Hzzp_{|z=0}=2$. More explainations on $\HpG$ are
given in Section~\ref{sec: description glancing region}.  In turn, one
defines the following vector field on $\TL$
\begin{align}
  \label{eq: def XG}
  \XG(\y)=
  \begin{cases}
    \Hp (\y) & \text{if}\ \y \in \TL\setminus \sgGb,\\
    \HpG (\y) & \text{if}\ \y \in \sgGb,
  \end{cases}
\end{align} 
that is, $\XG = \Hp + \unitfunction{\sgGb} (\HpG - \Hp) $.
%%%%%%%%%%%%%%%%%%%%%%%%
% definition           %
%%%%%%%%%%%%%%%%%%%%%%%%
\begin{definition}[\gbichar]
	\label{def: generalized bichar-intro}
  Suppose $J \subset \R$ is an interval, $B$ a discrete subset of $J$, and
  \begin{align*}
    \gammaG: J\setminus B 
    \to \Char p \cap \TL.
  \end{align*}
  One says that $\gammaG$ is a \gbichar if  the
  following properties hold:
  \begin{enumerate}
  \item For $s \in J \setminus B$, $\gammaG(s) \notin \Hb$ and the map $\gammaG$ is differentiable
    at $s$ with
	\begin{equation*}
		\frac{d}{d s} \gammaG(s) = \XG\big( \gammaG(s)\big).
	\end{equation*}
  \item If $S\in B$, then $\gammaG (s) \in \TL \setminus \dTL$ for $s
    \in J\setminus B$ \suff
    close to $S$ and moreover
    \begin{enumerate}
      \item  if $[S-\eps,S] \subset J$ for some $\eps>0$, then $\gammaG(S^-) = \lim_{s \to S^-} \gammaG(s) \in
      \Hb^-$;
       \item if $[S,S+\eps] \subset J$ for some $\eps>0$, then $\gammaG(S^+) = \lim_{s \to S^+} \gammaG(s) \in
      \Hb^+$;
    \item and if $[S-\eps,S+\eps]\subset J$ for some $\eps>0$, then
    $\gammaG(S^+) = \Sigma \big(\gammaG(S^-)\big)$.
    \end{enumerate}   
  \end{enumerate}
\end{definition}
Recall that $\TL$ contains its boundary $\dTL$; as a result a
\gbichar $\gammaG (s)$ may lie in the boundary for $s$ in
some interval. Details on \gbichars are given in Section~\ref{sec: Broken and generalized bicharacteristics}.

\medskip
When one refers to a (generalized) \bichar one often means the points visited in
$\TL$ by
$s\mapsto \gammaG(s)$ as $s$ varies, that is, 
\begin{align*}
	\{ \gammaG(s); \ s \in J \setminus B\}.
\end{align*}
Observe however that this set may not be a closed set if $B\neq
\emptyset$ as its intersection with $\Hb$ is empty.  
Consequently we rather use its closure to describe the set of reached points.
%%%%%%%%%%%%%%%%%%%%%%%%
% definition           %
%%%%%%%%%%%%%%%%%%%%%%%%
\begin{definition}[\gbichar]
  \label{def: generalized bichar 2-intro}
	By \gbichar one also refers to 
	\begin{align*}
	\GammaG = \ovl{\{ \gammaG(s); \ s \in J \setminus B\}} 
	= \{ \gammaG(s); \ s \in J \setminus B\} 
	\cup \bigcup_{s\in B} \{ \gammaG(s^-), \gammaG(s^+)\}.
\end{align*}
\end{definition}

\bigskip
The following theorem states that for every point of  $\TL$ one can find a maximal \gbichar that goes through this point.
%%%%%%%%%%%%%%%%%%%%%%%%
% theorem              %
%%%%%%%%%%%%%%%%%%%%%%%%
\begin{theorem}
  \label{theorem: existence bichar}
  Suppose $J\setminus B \ni s \mapsto \gammaG(s) = (t(s), x(s),
  \tau(s), \xi(s))$ is a \gbichar. If $\gammaG$ is maximal
  then $J = \R$. Moreover, $t(\R) = \R$ if $\tau(s) = \cst \neq 0$.
  
  If $\y^0 \in \Char p \cap \TL$ there exists a maximal \gbichar $s \mapsto \gammaG(s)$ with $s \in \R \setminus B$ 
  such that $\gammaG(0) = \y^0$ if $\y^0 \notin \Hb$ and 
$\gammaG(0^\pm) = \y^0$ if $\y^0 \in \Hb^\pm$.
\end{theorem}
If it does not create confusion, by abuse of notation, we 
sometimes write $\gammaG(0) = \y^0$ even in the
case $\y^0 \in \Hb$ with the understanding that  $\gammaG(0^\pm) = \y^0$ if $\y^0 \in \Hb^\pm$.

Note that there is no uniqueness of such a  maximal \gbichar because of the limited smoothness of $\XG$.
The result of Theorem~\ref{theorem: existence bichar} is classical in the case of smooth coefficients; see
\cite{MS:78} or \cite[Section 24.3]{Hoermander:V3}. Here, in
the case of the present limited smoothness a possible proof of
Theorem~\ref{theorem: existence bichar} follows quite
closely the arguments developed in what
follows. Instead of duplicating a quite long proof, we chose an
argument that allows one to consider Theorem~\ref{theorem: existence
  bichar} as a consequence of our main result, Theorem~\ref{theorem:
  measure support propagation} below.
We refer to Appendix~\ref{sec: proof: theorem: existence
  bichar} for this proof.

Despite the lack of uniqueness, some form of continuity related to all
\bichars passing through one point holds. For $\y^0 = (t^0, x^0,
\tau^0, \xi^0) \in \Char p
\cap \TL$ and $T>0$
introduce
\begin{align*}
  \Gamma^T(\y^0)
  = \{|t-t^0| \leq T \} \cap \bigcup_{\y^0 \in \GammaG} \GammaG,
\end{align*}
that is, the union of all \gbichar that pass through $\y^0$, restricted to the  time interval $[t^0-T, t^0+T]$.  
%\begin{multline*}
%  \Gamma^T(y^0)  = \{ \y =(t,x,\tau,\xi) \in \Char p \cap \TL; \ 
%  |t-t^0| \leq T, \ \y =\gammaG(s), \ \et \ \y^0 =\gammaG(0) \\
%  \pour \ \gammaG \ \text{a maximal \gbichar}\}.
%\end{multline*} 
%%%%%%%%%%%%%%%%%%%%%%%%
% proposition          %
%%%%%%%%%%%%%%%%%%%%%%%%
\begin{proposition}
  \label{prop:continuity-flow}  
  Suppose $\y^0 \in \Char p \cap \TL \setminus
  0$ and $T>0$. 
  \begin{multline*}
    \forall \eps>0, \ \exists \delta>0, \ \forall \y^1,\y \in \Char p
    \cap \TL, \\
    \dist(\y^1, \y^0) \leq \delta \ \et \ 
    \y \in \Gamma^T(\y^1)
    \ \ \imp \ \
    \dist(\y, \Gamma^T(\y^0) \leq \eps.   
  \end{multline*}
\end{proposition}
A proof is given in Appendix~\ref{sec: proof: continuity properties
  bichar}.

\section{Main result and open questions}
\label{sec: main result}
Our main result concerns the description of the support of a measure
density\footnote{The word `density' is omitted in what follows, yet
  $\mu$ has the density property and thus acts on continuous functions
on $\hL$.} $\mu$ defined on $\U$ an open subset of  $T^* \hL$. Recall
that $\hL$ is a local extension of $\L$; see Section~\ref{sec:
  Geometrical setting-intro}.

A first assumption made on $\mu$ is the following.
%%%%%%%%%%%%%%%%%%%%%%%%
% Assumption         %
%%%%%%%%%%%%%%%%%%%%%%%%
\begin{assumption}
  \label{assumption: first properties of the semiclassical measure}
  The measure $\mu$ is nonnegative and supported in $\U \cap \Char p 
  \cap \TL
   \setminus 0$.
\end{assumption}
In particular, $\mu$ vanishes in a \nhd of $(\tau,\xi)=0$ and in 
 $\U \setminus \TL$. 
% Near the boundary, with the chart $\mathcal C = (\hO, \chdiff)$ used above, this means
% that 
% \begin{equation*}
%   \chdiff\big(\supp(\mu) \cap \hO\big) \subset
% \chdiff\big(\Char p \cap \hO\big) \cap  \{ z\geq 0\}.
% \end{equation*} 

A second assumed property is the following one.
%%%%%%%%%%%%%%%%%%%%%%%%
% Assumption          %
%%%%%%%%%%%%%%%%%%%%%%%%
\begin{assumption}
  \label{assumption: Gerard-Leichtnam equation}
  One has, in the sense of distributions,  
  \begin{align}
    \label{eq: Gerard-Leichtnam equation}
    \transp \Hp \mu =   f \mu  - \int_{\y \in \pHb \cup \pGb} 
    \frac{\delta_{\y^+} - \delta_{\y^-}}
    {\dup{\xi^+- \xi^-}{\nx}_{T_x^*\M, T_x\M} 
    } \ d \nu (\y) \qquad \dans \ \U,
  \end{align}
  with $f$ a continuous real function on $T^* \hL$ and $\nu$ a nonnegative measure on $\pdTL$ and 
   where $\y^\pm$ and $\xi^\pm$ are as given in 
   \eqref{eq: relevement H-G-intro}. 

   Here, $\nx$ stands for  the unitary inward
    pointing normal vector in the sense of the metric; it is recalled at the end of Section~\ref{sec: Musical isomorphisms, normals}. 
\end{assumption}
%%%%%%%%%%%%%%%%%%%%%%%%
% remark               %
%%%%%%%%%%%%%%%%%%%%%%%%
\begin{remark}
  \label{remark: integrand GL equation on G-intro}
  If $\y \in \pGb$ then $\y^-$ and $\y^+$ coincide with $\y$ and $\xi^+= \xi^-$.  The value of
  the integrand in \eqref{eq: Gerard-Leichtnam equation} thus requires
  some explanation in this case.  In fact, first consider
  $\y^0 = (\y^{0\prime}, 0) \in \pHb$ with $\y^{0\prime}=(t^0,x^{0\prime},z=0, \tau^0, \xi^{0\prime})$. Then $\y^{0,\pm} \neq \y^0$ and \eqref{eq: relevement
    H-G-intro} give
  $\xi^{0,+} - \xi^{0,-} = 2 \zeta^+ dz$, yielding
  $\dup{\xi^{0,+} - \xi^{0,-}}{\n_{x^0}}_{T_x^*\M, T_x\M} = 2 \zeta^+$ since
  $\nx= \d_z$ in the considered coordinates. With a
  $\Con^1$-test function $q(\y)$ one has 
  \begin{align*}
    \dup{\delta_{\y^{0,+} } - \delta_{\y^{0,-} }}{q}
    &= q \big(\y^{0\prime},\zeta^+\big)  
      -  q \big(\y^{0\prime},-\zeta^+\big) .
  \end{align*}
  The integrand is thus 
\begin{align*}
    \frac{q \big(\y^{0\prime},\zeta^+\big)  
      -  q \big(\y^{0\prime},-\zeta^+\big) }{2 \zeta^+}.
\end{align*}
  If now a sequence $(\y^{(n)} )_n\subset \pHb$ converges to $\y \in \pGb$
then 
\begin{align}
  \label{eq: understanding GL-glancing}
  \frac{\dup{\delta_{\y^{(n),+}} - \delta_{\y^{(n),-}}}{q}}
  {\dup{\xi^{(n),+}-\xi^{(n),-}}{\nx}_{T_x^*\M, T_x\M}}
  \to \d_{\zeta}  q (\y).
\end{align}
The integrand in \eqref{eq: Gerard-Leichtnam equation} for $\y\in
 \pGb$ is thus to be understood as the derivative with respect to the
 variable $\zeta$ at $\zeta=0$.
Note that this interpretation is very coordinate dependent. We
give a more geometrical interpretation using more intrinsic
coordinates in Section~\ref{sec: On the measure equation}.
\end{remark}

Our main result states that if a point lies in the support of $\mu$
solution to \eqref{eq: Gerard-Leichtnam equation}, then there exists a maximal
\gbichar initiated at this point contained in $\supp \mu$.
%%%%%%%%%%%%%%%%%%%%%%%%
% theorem              %
%%%%%%%%%%%%%%%%%%%%%%%%
\begin{theorem}
  \label{theorem: measure support propagation}
 Suppose $\mu$ is a measure density that fulfills
 Assumptions~\ref{assumption: first properties of the semiclassical
   measure} and \ref{assumption: Gerard-Leichtnam equation}.  Suppose
 $\y^0 \in \supp \mu$. There exists a maximal \gbichar $s \mapsto
 \gammaG(s)$ with $s \in J \setminus B$, with $B$ a discrete subset
 of $\R$, such that
  \begin{align*}
    \y^0 \in \GammaG \subset \supp \mu, 
  \end{align*}
	with $\gammaG$ and $\GammaG$ as in Definitions~\ref{def: generalized bichar-intro} and \ref{def: generalized bichar 2-intro}.
  In other words, the support of $\mu$ is a union of maximal
  \gbichars in $\U$.

  If $\U =T^*\hL$ then $J=\R$ for  each maximal
  \gbichar. 
\end{theorem}

As mentionned in the introductory section, this is an
extension of the superposition principle of L.~Ambrosio and G.~Crippa
\cite{Ambrosio:08,Ambrosio-Crippa:08,AmCr14}, yet in a nonquantitative form: here, we only describe
the geometry of the support of the measure and not the measure
itself. A very natural open question is the following:
\begin{center}
  {\em 
  Is there an extension of the quantitative superposition principle of
  \cite{AmCr14}\\ for a nonnegative measure that fulfills both Assumptions~\ref{assumption: first properties of the semiclassical measure}
  and \ref{assumption: Gerard-Leichtnam equation}?}
\end{center}
In other words,  can one
write such a mesure $\mu$ as a ``sum'' of
positive measures, each defined as a constant times $\delta_{\tiny \GammaG}$,
where $\GammaG$ is a \gbichar as in Definition~\ref{def:
  generalized bichar-intro}?

\medskip
A first natural question is the following: for a \gbichar
$\GammaG$, is the measure $\delta_{\tiny \GammaG}$ well defined and does it
fulfill a transport equation of the form of \eqref{eq:
  Gerard-Leichtnam equation}? Despite the fact that $\gammaG(s)$ can
have an infinite number of points of discontinutity,  for $s\in B$, that can accumulate,
one can answer positively these question. This is
done in the beginning of Appendix~\ref{sec: proof prop mesure single
  bichar}. 

\medskip
Note that the nonquantitative superposition principle of
Theorem~\ref{theorem: measure support propagation} suffices for the purpose
of the companion article \cite{BDLR1} towards the derivation
of observability estimate for the wave equation in the case of
$\Con^1$-coefficients and a $\Con^2$-boundary.

\bigskip
\medskip
An important consequence of Assumptions~\ref{assumption: first properties of the semiclassical measure} and \ref{assumption: Gerard-Leichtnam equation} is also the following property of the measure $\nu$.
%%%%%%%%%%%%%%%%%%%%%%%%
% proposition          %
%%%%%%%%%%%%%%%%%%%%%%%%
\begin{proposition}
  \label{prop: no mass on Gd G3}
  There exists $C>0$ such that 
   $|\tau| \geq C>0$ in $\supp \nu \cap (\pHb \cup \pGb)$.
  One  has $\dup{\nu}{\unitfunction{\sdGb \cup \glGb}}=0$ and
  $\dup{\mu}{ \unitfunction{\sdGb}}=0$, that is, the measure $\nu$
  has no mass on $\sdGb \cup \glGb$ and the measure $\mu$ has no mass on $\sdGb$.
\end{proposition}
This proposition is due to N.~Burq and P.~G\'erard~\cite{BG:1997} in
the case of smooth coefficients. The proof requires refinements in the
present low regularity setting. It is given in Section~\ref{sec: proof
  prop: no mass on Gd G3}.
%%%%%%%%%%%%%%%%%
% section
%%%%%%%%%%%%%%%%%
\section{Transport equation, measure support propagation away from
  or across boundaries}
\label{sec: measure support propagation1}
%%%%%%%%%%%%%%%%%
% sub-section
%%%%%%%%%%%%%%%%%
\subsection{Support propagation away from boundaries}

Suppose $\Omega$ is an open subset of a $\Con^2$ $d$-dimensional manifold. 
Denote by $^1\D'(\Omega)$ and $^1\D^{\prime,0}(\Omega)$ the spaces
of density  distributions and density Radon measures on
$\Omega$. 

Consider a continuous vector field $X$ and a continuous real function $f$
on $\Omega$ and suppose $\mu$ is a nonnegative measure density
on~$\Omega$. Assume that $\mu$ is such that $\transp{X}\mu = f\mu$ in the
sense of distributions, that is,
\begin{equation}
  \label{EDO}  
  \bigdup{\transp{X}\mu}{a}_{^1\D'(\Omega), \Cinf_c(\Omega)}
  = \bigdup{\mu}{X a}_{^1\D^{\prime,0}(\Omega), \Con^0_c(\Omega)}=
  \bigdup{\mu}{f a}_{^1\D^{\prime,0}(\Omega), \Con^0_c(\Omega)}, 
\end{equation}
for $a \in \Cinfc(\Omega)$.
If $f$ vanishes and $X$ is moreover Lipschitz, one concludes that $\mu$ is invariant
along the flow that $X$ generates. However, if $X$ is not Lipschitz,
there is no such flow in general. Yet, integral curves do exist by the
Cauchy-Peano theorem.

Away from any boundary a precise statement associated with
\eqref{statement: measure support transport} is given in the following theorem.
%%%%%%%%%%%%%%%%%%%%%%%%
% theorem              %
%%%%%%%%%%%%%%%%%%%%%%%%
\begin{theorem}
\label{th: ODE}
On $\Omega$, suppose $X$ is a continuous vector field, $f$ is a continuous
real function, and $\mu$ is a nonnegative density measure  that is solution to $\transp{X}\mu =  f\mu$ in the sense of
distributions.  Then, the support of $\mu$ is a union of maximally
extended integral curves of the vector field $X$.

In other words, if $m^0\in \Omega$ is in  $\supp \mu$, then there exist an
 interval $I$ in $\R$ with  $0\in I$  and a $\Con^1$ curve $\gamma: I
 \to  \Omega$ that cannot be extended such that $\gamma(0) = m^0$ and 
\begin{equation*}
  \frac{d }{ds} \gamma(s) = X( \gamma (s)), \qquad s \in I,
\end{equation*} 
and 
$\gamma(I) \subset \supp \mu$.
\end{theorem}
This theorem and its proof can be found in \cite{BDLR0}. We decided to 
reproduce the argument here as it increases the readability of the present article for the following two reasons
:
\begin{enumerate}
  \item the proof of Theorem~\ref{th: ODE}  is much simpler than the
argument we develop below  to understand
the structure of the support of $\mu$ at a boundary if fulfilling the
more general equation \eqref{eq: Gerard-Leichtnam equation};
\item  the
techniques used at the boundary, despite their high level of
technicality, are in the same spirit as those in the proof of
Theorem~\ref{th: ODE}. In particular, some of the cutoff functions
introduced in the proof of Theorem~\ref{th: ODE} are used further
in the article. 
\end{enumerate}

The strategy of the proof of Theorem~\ref{th: ODE} is very much inspired by the Melrose and
Sj\"ostrand approach to the propagation of singularities~\cite{MS:78}
and relies on careful choices of test functions allowing one to
construct sequences of points in the support of the measure relying on
nonnegativity\footnote{of the measure in our case and of some
  operators for Melrose and Sj\"ostrand, via the G\r{a}rding
  inequality.}. Then, a limiting procedure leads to the conclusion, in
the spirit of the classical proof of the Cauchy-Peano theorem.

Theorem~\ref{th: ODE} is stated on an open subset of a smooth
manifold. Yet, its result is of local nature. Using a
local chart one may assume that $\Omega$  is an open subset of $\R^d$
instead without any loss of generality.

The proof of Theorem~\ref{th: ODE} is
made of two steps that are stated in the following propositions.
%%%%%%%%%%%%%%%%%%%%%%%%
% proposition               %
%%%%%%%%%%%%%%%%%%%%%%%%
\begin{proposition}
\label{prop: invariant}
Suppose $X$ is a $\Con^0$-vector field on $\Omega$ an open subset of
$\R^d$. For a closed set $F$ of $\Omega$, the following two
properties are equivalent.
\begin{enumerate}
\item \label{one} The set $F$ is a union of maximally extended integral curves of the vector field $X$.
\item \label{two} For any compact  $K\subset \Omega$ where the vector field $X$ does not vanish, 
\begin{align*}
  \forall \eps>0, \ \exists \delta_0>0, \  \forall x\in K\cap F,  
  \ \forall \delta \in [- \delta_0, \delta_0], 
  \ \ B\big(x+ \delta X(x), \delta \eps\big) \cap F \neq \emptyset.
\end{align*}
\end{enumerate}
\end{proposition}
%%%%%%%%%%%%%%%%%%%%%%%%
% proposition               %
%%%%%%%%%%%%%%%%%%%%%%%%
\begin{proposition}
  \label{prop: symboles}
  On $\Omega$ an open subset of $\R^d$, suppose $X$ is a
  $\Con^0$-vector field and $f$ is a continuous real function.  Consider a
  nonnegative measure $\mu$ on $\Omega$ solution to $\transp{X}\mu = f \mu$
  in the sense of distributions.
Then, the closed set $F= \supp \mu$ satisfies the second property in
Proposition~\ref{prop: invariant}.
\end{proposition}

\medskip
\begin{proof}[Proof of Proposition~\ref{prop: invariant}]
 First, we prove that Property~\eqref{one} implies
  Property~\eqref{two} and consider a compact set $K$ of $\R^d$ such
  that $K \subset \Omega$ and $K \cap F \neq \emptyset$. 

  There exists $\eta>0$ such that $K \subset K_\eta \subset
  \Omega$ with $K_\eta = \{ x \in \Omega;\ \dist(x, K) \leq
  \eta\}$. One has $\Norm{X}{} \leq C_0$ on $K_\eta$ for some
  $C_0>0$. Suppose $x \in K$ and $\gamma(s)$ is a maximal integral
  curve defined on an interval $]a,b[$, $a, b \in \ovl{\R}$ and such
  that $0 \in ]a,b[$ and $\gamma(0) =x$.  If $b <
\infty$ then there exists $s^1 \! \! \in \, ]0,b[$ such that $\gamma(s^1)
\notin K_\eta$. Since $\gamma(s) \in K_\eta$ if $s <  \eta /C_0$, one finds
  that   $b \geq \eta /C_0$. Similarly, one has $|a| \geq \eta /C_0$.
  Consequently, there exists $\Slim>0$ such that any
  maximal integral curve $\gamma(s)$ of the vector field $X$ with
  $\gamma (0) \in K$ is defined for $s \in I= (- \Slim, \Slim)$. 
  
  Pick $x \in K\cap F$. According to
  the Property~\eqref{one}, there exists
  \begin{align*}
    \gamma: I \rightarrow F \ \text{such that} \ 
    \dot{\gamma}(s) = X(\gamma(s)) 
     \ \text{and} \  \gamma(0) =x.
    \end{align*}
    By uniform continuity of the vector field $X $ in a compact
    neighborhood of $K$ one has 
\begin{align*}
  \gamma(s) = \gamma(0) + \int_0^s \dot{\gamma}(s) ds
  = \gamma(0) + \int_0^s X(\gamma(s)) ds
  = x+ s X(x) + r(s),
\end{align*}
for $s\in (-S, S)$, 
where $\lim_{s\rightarrow 0} \Norm{r(s)}{} / s =0$,
{\em uniformly} with respect to $x$. One deduces that for any  $\eps
>0$ there exists $0< \delta_0< S$ such that $\Norm{r(s)}{}< s \eps $  for any $s\in (- \delta_0, \delta_0)$,
which implies 
\begin{align*} 
  F \ni \gamma(s) \in B\big(x+ sX(x), s\eps\big).
\end{align*}

\medskip Second, we prove that Property~\eqref{two} implies
Property~\eqref{one}. It suffices to prove that for any $x\in F$ there
exist an open interval $I \ni 0$ and an integral curve 
\begin{align*}  
  \gamma: I \rightarrow F 
  \ \text{such that} \ 
  \dot{\gamma}(s) =X( \gamma(s))
   \ \text{and} \ \gamma(0) =x. 
\end{align*}
Then, the standard continuation argument shows that this local
integral curve included in $F$ can be extended to a maximal integral
curve also included in $F$.

If $X(x)=0$, then the trivial integral curve
$\gamma(s) = x$, $s\in \R$, is included in $F$. As a consequence,
one assumes $X(x) \neq 0$ and one picks a compact
neighborhood $K$ of $x$ containing  $B(x,\eta)$ with $\eta >0$ and where, for some $0 < c_K < C_K$, 
\begin{align*} 
 c_K\leq \Norm{ X(y)}{} \leq C_K, \quad y \in K.
\end{align*}

Let $n \in \N^*$. Set $x_{n,0} = x$ and  $\eps = 1 /n$ and apply
Property~\eqref{two}. One deduces that there exist $0< \delta_n \leq
1 / n$ and a point 
\begin{align*}  
x_{n,1}\in F \cap B \big(x_{n,0}+\delta_n X(x_{n,0}), {\delta_n} / n\big). 
\end{align*}
If $x_{n,1}\in K$ one can perform this construction again, yet starting
from $x_{n,1}$ instead of $x_{n,0}$. If a sequence of points
$x_{n,0}, x_{n,1}, \dots, x_{n,L^+}$ is obtained in this manner one has
\begin{align}  
  \label{eq: construction points - away from boundary}
x_{n,\ell+1} \in F \cap B\big( x_{n,\ell}+\delta_n X(x_{n, \ell}),
  {\delta_n} /n\big), \qquad \ell = 0, \dots, L^+-1.
\end{align}
%%%%%%%%%%%%%%%%%%%%%%%
% figure
%%%%%%%%%%%%%%%%%%%%%%%
\begin{figure}
  \begin{center}
    \subfigure[Iterative construction of the curve $\gamma_n$. \label{fig: construction gamma n-a}]
    {\resizebox{7cm}{!}{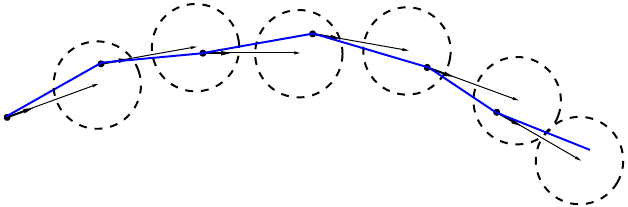}}
    \quad 
    \subfigure[Convergence of $\gamma_n$ as $n$ increases. \label{fig: construction gamma n-b}]
    {\resizebox{7cm}{!}{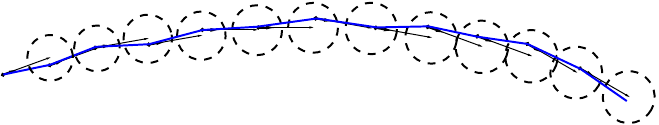}}
    \caption{Construction and convergence of the sequence $(\gamma_n)_n$.}
  \label{fig: construction gamma n}
  \end{center}
\end{figure}
One can
carry on the construction as long as $ x_{n,L^+} \in K$.
The same construction for $\ell\leq 0$ can be performed, with the
property 
\begin{align}  
  \label{eq: construction points - away from boundary bis}
x_{n,\ell-1} \in F \cap B\big(x_{n,\ell} - \delta_n X(x_{n, \ell}),
  {\delta_n} / n\big), \qquad |\ell| = 0, \dots, L^--1.
\end{align}
Having $\Norm{X}{}\leq C_K$ on $K$ and $B(x,\eta) \subset K$ ensures that one can construct the sequence at least for 
\begin{align*}  
L^+ = L^- = L_n 
   = \Big\lfloor \frac{ \eta}{\delta_n (C_K + 1)}\Big\rfloor +1
  \leq \Big\lfloor \frac{ \eta}{\delta_n (C_K + 1/ n)} \Big\rfloor +1,
\end{align*}
where $\lfloor . \rfloor$ denotes the floor function. 
With the constructed  points $x_{n,\ell}$, $|\ell| \leq L_n$,
define the following continuous curve $\gamma_n(s)$ for
$|s|\leq L_n\delta_n$:
\begin{align*}  
\gamma_n (s) = x_{n,\ell}  + (s-\ell\delta_n) \frac{x_{n,\ell+1} - x_{n,\ell}} { \delta_n}
 \  \text{for}\  s\in [\ell \delta_n, (\ell + 1) \delta _n) \  \text{and}\  |\ell|
  \leq L_n-1.
\end{align*}
This curve and its construction is illustrated in Figure~\ref{fig:
  construction gamma n-a}. 
Note that $\gamma_n(s)$ remains in a compact set, uniformly with
respect to $n$. In this compact set $X$ is uniformly continuous. 

Set $\Slim = \eta /(C_K+1)$. Since $\Slim \leq L_n\delta_n$, in fact,  we
only consider the function $\gamma_n(s)$ for $|s| \leq
\Slim$ in what follows.  
Note that since $x_{n,\ell} \in F$ for $|\ell| \leq L_n$ then one has 
\begin{align}
  \label{eq: discrete curve almost in F}
  \dist\big(\gamma_n (s), F\big)\leq  \delta_n ( C_K+ 1 / n), 
  \qquad |s|\leq \Slim.
\end{align}
From \eqref{eq: construction points - away from boundary}, for $\ell
\geq 0$ and  $s\in (\ell \delta_n, (\ell + 1) \delta _n)$, one has 
\begin{align*}  
  \dot{\gamma}_n(s)= \frac{x_{n,\ell+1} - x_{n,\ell}}{\delta_n} =
X(x_{n,\ell}) + \mathcal{O}( 1 / n).
\end{align*}
Similarly, from \eqref{eq: construction points -
  away from boundary bis}, for $\ell
\leq 0$ and  $s\in ((\ell-1) \delta_n, \ell\delta _n)$, one has
\begin{align*}  
  \dot{\gamma}_n(s)= \frac{x_{n,\ell} - x_{n,\ell-1}}{\delta_n} =
X(x_{n,\ell}) + \mathcal{O}( 1 / n).
\end{align*}
In any case, using the uniform continuity of the
vector field $X$, one finds
\begin{align*}  
  \dot{\gamma}_n(s) 
  = X( \gamma_n(s)) + e_n(s),
\end{align*}
where the error $|e_n|$ goes to zero {\em uniformly} with respect to
$|s|\leq\Slim$ as $n\rightarrow + \infty$.

\medskip
Since the curve $\gamma_n$ is continuous, one finds
\begin{equation}\label{eq.appro}
  \gamma_n(s) 
  = \gamma_n(0) + \int_0^s \dot{\gamma}_n(\sigma)d\sigma 
  = x+  \int_0^s X({\gamma}_n(\sigma)) d\sigma 
  + \int_0^s e_n(\sigma) \, d \sigma.
\end{equation}
Let now $n$ grow to infinity. With \eqref{eq.appro}, 
the family of curves
$(s \mapsto \gamma_n(s), |s| \leq\Slim)_{n\in \mathbb{N}^*}$ is
equicontinuous and pointwise bounded; by the Arzel\`a-Ascoli theorem one can extract a
subsequence $(s \mapsto \gamma_{n_p})_{p \in \mathbb{N}}$ that
converges uniformly to a curve $\gamma(s), |s| \leq \Slim$. Convergence is illustrated in Figure~\ref{fig:
  construction gamma n-b}.  Passing
to the limit $n_p\rightarrow+ \infty$ in~\eqref{eq.appro} one finds that
$\gamma(s)$ is solution to
\begin{align*}  
\gamma(s) = x+  \int_0^s X({\gamma}(\sigma))d\sigma.
\end{align*}
From \eqref{eq: discrete curve almost in F}, for any $|s|
\leq \Slim$, there exists $(y_p)_p \subset F$ such that 
$\lim_{p\rightarrow + \infty} y_p = \gamma(s)$.
Since $F$ is closed one concludes that $\gamma(s) \in F$. 
\end{proof}
  
\medskip
\begin{proof}[Positivity  argument and proof of Proposition~\ref{prop:
    symboles}]
 Consider a compact set $K$ where the vector field $X$ does not
  vanish.  By continuity of the vector field there exist $0< c_K\leq  C_K$ such that  $0< c_K \leq\Norm{ X(x)}{}
  \leq C_K$, for all $x \in K$. 

Consider $x^0\in K\cap \text{supp}(\mu)$. By performing a
rotation and a dilation of coefficient $\Norm{ X(x^0)}{}\in [c_K,
C_K]$, one
can assume that $X(x^0) = (1,0, \dots, 0)\in \R^d$. One writes $x
= (x_1, x')$ with $x' \in \R^{d-1}$. 

%%%%%%%%%%%%%%%%%%%%%%%
% figure
%%%%%%%%%%%%%%%%%%%%%%%
\begin{figure}
  \begin{center}
    \subfigure[\label{fig: function chi}]
    {\resizebox{5.5cm}{!}{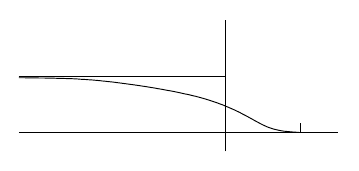}}
    \quad 
    \subfigure[\label{fig: function beta}]
    {\resizebox{5.5cm}{!}{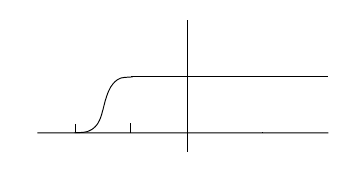}}
    \caption{The two localization functions $\chi$ (a) and $\beta$ (b)
      used to
      build the test function $q$.}
  \label{fig: localization functions}
  \end{center}
\end{figure}
\medskip Let $\chi \in \Cinf(\R)$ be given by 
\begin{equation}
  \label{eq: def chi}
  \chi (s) = \bld{1}_{s<1}\, \exp (1/(s-1)),
\end{equation} 
and 
$\beta \in \Cinf(\R)$ be such that
\begin{equation}
  \label{eq: def beta}
  \beta \equiv 0 \ {\text on}\   ]-\infty, -1],\quad 
 \beta'>0 \ {\text{on}}\  ]-1, -1/2[,\quad
 \beta \equiv 1\ {\text{on}}\  [-1/2, +\infty[.
\end{equation} 
These two functions are represented in Figure~\ref{fig: localization functions}.
Then set 
\begin{align} 
  \label{eq: def q eps delta x - test function}
  q = e^{A x_1} (\chi \circ v) (\beta \circ w), 
  \quad g =e^{A x_1} (\chi' \circ v) (\beta \circ w)  X v, 
  \quad h = e^{A x_1} (\chi \circ v) (\beta' \circ w) X w, 
  \end{align}
with $A>0$ meant to be chosen \suff large below and 
\begin{align*}
  &v(x)= 1/2-\delta^{-1} (x_{1}- x^0_{1})
  +8 (\eps \delta)^{-2} \Norm{x'- x^{0\prime}}{}^{2}\\
  &\et \quad 
  w(x)=2\eps^{-1} \big(1  - \delta^{-1} (x_{1}-x^0_{1}) \big),  
\end{align*}
for $\eps>0$ and $\delta>0$ both meant to be chosen small in what follows.
One has
\begin{align*}
  X q  =  g + h + A (X x_1)q.
\end{align*}

\medskip
The function $q$ is compactly supported. Indeed, in
the support of $\beta\circ w$ one has $w\geq -1$ implying
\begin{align*}
  x_1- x^0_1 \leq \delta (1+ \eps/2),
\end{align*}
while on the support of $\chi \circ v$ one  has $v \leq 1$ which gives
\begin{equation*}
  - 1/2+8  (\eps \delta)^{-2} \Norm{x'- x^{0\prime}}{}^{2}
  \leq \delta^{-1} (x_{1}- x^0_{1}).
\end{equation*}
On the supports of $q$ and $(\chi'
\circ v) (\beta \circ w)$ one thus finds
\begin{equation}
  \label{proche}
  - \delta /2  \leq  x_1- x^0_1 \leq \delta ( 1+ \eps/2)
  \ \  \text{and} \ \ 
  8  (\eps \delta)^{-2} \Norm{x'- x^{0\prime}}{}^{2} \leq 3/2 + \eps/2.
\end{equation}
Similarly, on the support of $\beta' \circ w$ one has $-1\leq w \leq
-1/2$ yielding
\begin{align*}
  \delta ( 1+ \eps /4) \leq x_1- x^0_1 \leq \delta (1+ \eps/2),
\end{align*} 
which implies that on the support of
$h$ one has
 \begin{equation}\label{proche2}
   \delta ( 1+ \eps /4)  \leq x_1- x^0_1 
   \leq \delta ( 1+ \eps/2)
   \ \  \text{and} \ \ 
   8  (\eps \delta)^{-2} \Norm{x'- x^{0\prime}}{}^{2} \leq 3/2 + \eps/2.
\end{equation}
In particular, in the case $\eps \leq 1$, one finds
 \begin{equation}\label{proche3}
   \supp h 
   \subset B\big(x^0 + \delta X(x^0), \eps \delta\big).
\end{equation}
These estimations of the supports of $q$ and
$h$ are illustrated in Figure~\ref{fig: test function supports}.
%%%%%%%%%%%%%%%%%%%%%%%
% figure
%%%%%%%%%%%%%%%%%%%%%%%
\begin{figure}
  \begin{center}
    \subfigure[Support of $q$. \label{fig: q eps delta support}]
    {\resizebox{7.5cm}{!}{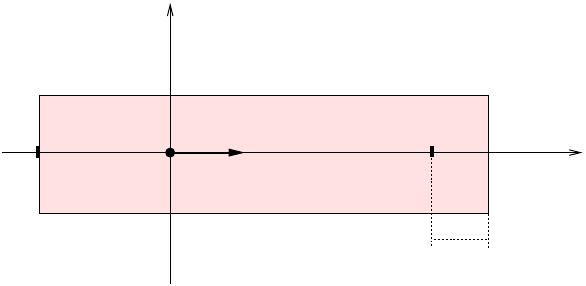}}
    \quad 
    \subfigure[Support of $h$. \label{fig: g eps delta support}]
    {\resizebox{7.5cm}{!}{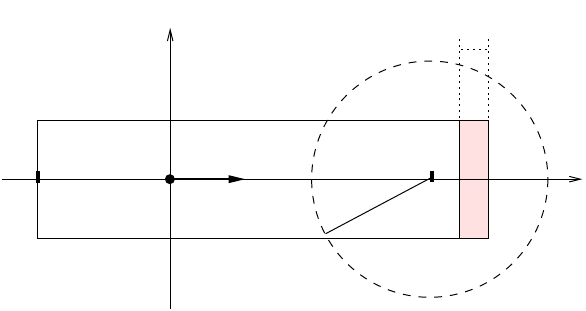}}
    \caption{Estimation of the test function supports in the case $\eps \leq 1$.}
  \label{fig: test function supports}
  \end{center}
\end{figure}
%%%%%%%%%%%%%%%%%%%%%%%%
% lemma                %
%%%%%%%%%%%%%%%%%%%%%%%%
\begin{lemma}
  \label{lem: positivite g}
  For any $0 < \eps \leq 1$ there exists $\delta_0>0$ such that for
  any $x^0\in K$ and $0 < \delta \leq \delta_0$
  \begin{enumerate}
    \item the function $g$ is nonnegative and is positive in a \nhd of
      $x^0$.
    \item $X x_1 \geq 1/2$ in $\supp q$.
  \end{enumerate}
\end{lemma}
\begin{proof}
  Consider $0 < \eps \leq 1$. One has
  $g = (\chi' \circ v) (\beta \circ w) X v$. Since
  $\beta \geq 0$ and $\chi'<0$ it suffices to prove that $X v(x) \leq 0$
  for $x$ in the support of $(\chi' \circ v) (\beta \circ w)$ for $\delta>0$
  chosen sufficiently small, uniformly with respect to $x^0 \in K$.

  Write
  \begin{align}
    \label{eq: decomposition X pres x0}
    X(x) - X(x^0) = \alpha^1(x, x^0) \d_{x_1} + \alpha'(x, x^0) \cdot \nabla_{x'},
  \end{align}
  with $\alpha^1(x, x^0)\in \R$ and $\alpha'(x, x^0)\in \R^{d-1}$. 
  By \eqref{proche},  for $x \in \supp (\chi' \circ v) (\beta \circ w)$ one
  has $\Norm{ x - x^0}{} \lesssim \delta$. From the uniform continuity of
  $X$ in any compact set one concludes that 
  \begin{align}
    \label{eq: smallness Xx - Xx0}
    |\alpha^1(x, x^0)| + \Norm{\alpha'(x, x^0)}{} = o(1) \ \
    \text{as} \  \delta \to 0^+,
  \end{align}
  uniformly\footnote{Observe that the change of variables made above
    for $X(x^0) = (1,0,\dots,0)$ does not affect uniformity since the
    dilation is made by a factor in $[c_K,C_K]$.} with respect to $x^0\in K$ and $x \in \supp (\chi' \circ v) (\beta \circ w)$.
  Using that $X(x^0) = \d_{x_1}$ and the form of $v$ given above, one writes 
  \begin{align*}
    X v(x) &= \big(X(x) v\big)(x) 
    = \big(\d_{x_1} v + \big(X(x) - X(x^0)\big) v\big)(x) \\
    &= - \delta^{-1}\Big( 1 + \alpha^1(x, x^0)  - 16 \eps^{-1}(\eps \delta)^{-1} \alpha'(x,x^0) \cdot (x' -
      x^{0\prime})\Big).
  \end{align*}
  Using again \eqref{proche}, one thus finds for $x \in \supp (\chi'
  \circ v) (\beta \circ w)$  
  \begin{align*}
    \big| \alpha^1(x, x^0) 
    -  16 \eps^{-1}(\eps \delta)^{-1} \alpha'(x,x^0) \cdot (x' - x^{0\prime})
    \big| 
    \lesssim \big| \alpha^1(x, x^0) \big| + \eps^{-1} \Norm{\alpha'(x,x^0)}{}.
  \end{align*}
  With $\eps$ fixed above and with \eqref{eq: smallness Xx - Xx0} one
  finds that $X v(x) \sim - \delta^{-1}$ as $\delta \to 0^+$
  uniformly with respect to $x^0\in K$ and
  $x \in \supp (\chi' \circ v) (\beta \circ w)$.

  One also has
  $g(x^0) = - \delta^{-1}\chi'(1/2) \beta(2
  \eps^{-1})>0$ and thus $g$ is positive in a \nhd
  of $x^0$, concluding the first part. 

  With \eqref{eq: decomposition X pres x0} one has $X x_1 = 1 +
  \alpha^1(x, x^0)$ and thus $X x_1\geq 1/2$ if $\Norm{x-x^0}{}\leq
  \eta$, with $\eta>0$ \suff small, uniformly in $x^0 \in K$. The
  estimate of $\supp q$ in \eqref{proche} gives $\Norm{x-x^0}{}\leq
  \eta$ for $\delta_0$ chosen \suff small. This gives the
  second part.
\end{proof}

We are now in a position to conclude the proof of
Proposition~\ref{prop: symboles}. Note that it suffices to prove the result
for $0 < \eps \leq 1$.  Choose $\delta_0>0$ as given by
Lemma~\ref{lem: positivite g}.
With ~\eqref{EDO}, for $0 < \delta \leq \delta_0$, one has
\begin{equation*}%\label{somme}
  0= \bigdup{\mu}{(X-f) q} = \dup{\mu}{g} + \dup{\mu}{h} +  \bigdup{\mu}{(A (X x_1) - f)q}.
\end{equation*}
By Lemma~\ref{lem: positivite g}, $(A (X x_1) - f)q\geq 0$ for $A \geq
2\sup_K |f|$, implying $\bigdup{\mu}{(A (X x_1) - f)q}\geq 0$.  By
Lemma~\ref{lem: positivite g}, $g\geq 0$ and $g$ is positive in a \nhd
of $x^0$. As $x^0 \in \supp \mu$ one finds $\dup{\mu}{g}>0$.
Consequently, 
$\dup{\mu}{h}\neq 0$. 
By the support estimate for $h$ given in \eqref{proche3} the conclusion follows: $\supp \mu \cap
B\big(x^0 + \delta X(x^0), \eps \delta\big)\neq \emptyset$. 
\end{proof}

%%%%%%%%%%%%%%%%%
% sub-section
%%%%%%%%%%%%%%%%%
\subsection{Support propagation at a boundary in the transverse case}
As above consider $X$ a $\Con^0$-vector field in an open subset $\Omega$
of a $\Con^2$ $d$-dimensional manifold. 
Suppose $\I$ is a $\Con^1$-hypersurface and $x^1 \in \I \cap \Omega$. In a local
chart $(\O, \cdiff)$ with $\O\subset \Omega$ \nhd of
$x^1$, $\I$ is given by $\varphi(x) =0$ with $d \varphi(x^1) \neq
0$ for some $\Con^1$-function $\varphi$. Assume that $X$ is transverse to $\I$ at $x^1$, meaning that $X
\varphi (x^1) = d \varphi(x^1) \big(X(x^1)\big) \neq 0$. 
This property remains true in a bounded \nhd $V \Subset \O$ of $x^1$.
Set 
\begin{align*}
  V^+ = V \cap \{ \varphi >0\}, \quad V^- = V \cap \{ \varphi <0\}.
\end{align*}
Consider a nonnegative measure $\mu$ that is solution in $V$ to the following
transport equation with a single-layer potential
\begin{align}
  \label{eq: equation t X mu everywhere}
  \transp{X}\mu =  f \mu 
  + \tilde{\mu} \otimes \delta_\I,
\end{align}
where $\tilde{\mu}$ is a measure on $\I$. Recall that
$\delta_\I = |\nabla \varphi|\, \varphi^*\delta$ (see \eg
\cite[Theorem 6.1.5]{Hoermander:V1}). 

The main result of this section is the following proposition. 
%%%%%%%%%%%%%%%%%%%%%%%%
% theorem                  %
%%%%%%%%%%%%%%%%%%%%%%%%
\begin{theorem}
  \label{th: propagation from boundary - transport equation}
  Suppose  $X$ is transverse to the hypersurface $\I$ in $V$ and $\mu$ is a
  nonnegative measure that vanishes in $V^-$ and is solution to
  \eqref{eq: equation t X mu everywhere}.  If
  $x^* \in \I \cap \supp \mu$, then there exists an integral curve of $X$,
  $s\mapsto \gamma(s)$, such that $\gamma(0) = x^*$ and
  \begin{enumerate}
    \item if $X \varphi>0$, then there exists $\Slim>0$ such that $\{
      \gamma(s)\}_{s\in [0,\Slim[} \subset \supp \mu$;  
    \item  if $X \varphi<0$, then there exists $\Slim>0$ such that $\{
      \gamma(s)\}_{s\in ]-\Slim,0]} \subset \supp \mu$.
  \end{enumerate}
\end{theorem}
In other words, a half integral curve of $X$ initiated at $x$ is {\em
  locally} contained in $\supp \mu$. This half integral curve is
naturally located in $\ovl{V^+}$.  The theorem is based on the
following proposition whose proof is given below. 
%%%%%%%%%%%%%%%%%%%%%%%%
% proposition          %
%%%%%%%%%%%%%%%%%%%%%%%%
\begin{proposition}
  \label{prop: discrete propagation of support - boundary}
  Suppose $\alpha = \pm 1$ and $K \subset V$  is a compact set such that $ \alpha
  X \varphi > 0$ on $K \cap \I$. Then, 
  \begin{multline*}
    \forall \eps>0, \ \exists \delta_0>0, \  \forall x^0 \in K
    \cap \ovl{V^+} \cap \supp \mu,   \ \forall \delta \in [0,
    \delta_0], \\
    B\big(x^0+  \delta \alpha
    X(x^0), \delta \eps\big) \cap \supp \mu \neq \emptyset.
    \end{multline*}
\end{proposition}

%%%%%%%%%%%%%%%%%%%%%%%%
% remark               %
%%%%%%%%%%%%%%%%%%%%%%%%
\begin{remark}
  \label{remark: smooth tansverse case}
  In the case $X$ is $\Con^1$ one can use its flow to find coordinates $x=(x',x_d)$
  such that the hypersuface $\I$ is locally given by $\{ x_d=0\}$ and
  $X = \d_{x_d}$. Then the result of Proposition~\ref{prop: discrete
    propagation of support - boundary} becomes obvious. Here,  we give
  a more involved proof that applies to the present case of a continuous
  vector field. 
\end{remark}

\begin{proof}[Proof of Theorem~\ref{th: propagation from boundary -
    transport equation}]
  On the one hand, Proposition~\ref{prop: discrete propagation of
    support - boundary} is the counterpart of Proposition~\ref{prop:
    symboles}, and it can be used to adapt the proof of
  Theorem~\ref{th: ODE} and obtain a proof of Theorem~\ref{th:
    propagation from boundary - transport equation}. On the other
  hand, one can use the result of Proposition~\ref{prop: discrete
    propagation of support - boundary} in conjunction with the result
  of Theorem~\ref{th: ODE} to obtain a little shorter proof of
  Theorem~\ref{th: propagation from boundary - transport equation}. We
  choose this second strategy here.

  We treat the case $X \varphi>0$ here; the other case can be treated
  similarly.  Suppose $x^0\in \I \cap \supp \mu$. Consider a bounded
  \nhd $W^0$ of $x^0$ where $X \varphi \geq C_0>0$. Consider $n \in \N$ and
  set $\delta =\eps =1/(n+1)$. By
  Proposition~\ref{prop: discrete propagation of support - boundary}
  there exists $x_n^1 \in B(x^0 +  X(x^0)/n, 1/n^2)$ such that $x_n^{1} \in
  \supp \mu$. For $n$ chosen \suff large, one has $x_n^1 \in V^+$. In
  $V^+$ one has $\transp X \mu =  f \mu$. One then applies
  Theorem~\ref{th: ODE}: there exists a {\em maximal}
  integral curve $\tgamma_n: ]S_1^n, S_2^n[ \to V^+$ of the vector field
  $X$, with $S_1^n < 0 <S_2^n$ with $\tgamma_n(0) = x_n^1$, that lies in $\supp \mu$.
  Since $X$ is bounded in $V^+$ (recall that $\ovl{V}$ is chosen compact), one finds that there exists $\Slim>0$
  such that $\Slim < S_2^n$ for $n$ chosen \suff large.

  Set $s_n^{1} = \Norm{x_n^1- x^0}{} / \Norm{X(x^0)}{}$. One has
  $s_n^{1} \to 0$ as $n \to \infty$. Then define $\gamma_n:
  [0,\Slim]\to V^+$ by 
  \begin{align*}
    \gamma_n (s) =
    \begin{cases}
      x^0 + s \frac{ x_n^1- x^0}{s_n^{1}}
      & \text{if} \  0\leq s \leq s_n^{1}, \\
      \tgamma_n(s - s_n^{1})
      & \text{if} \  s_n^{1} \leq s \leq \Slim.
    \end{cases}
  \end{align*}
  We then follow the arguments in the proof of Proposition~\ref{prop:
    invariant}. One has
  \begin{align*}  
    \dot{\gamma}_n(s) 
  = X( \gamma_n(s)) + e_n(s),
  \end{align*}
  for $s \in [0,\Slim]$ where the error $|e_n|$ goes to zero {\em uniformly} with respect to
$|s|\leq\Slim$ as $n\rightarrow + \infty$.
  Naturally, one has $e_n(s) =0$ for $s \in
  [s_n^{1},\Slim]$. Equation~\ref{eq.appro} is also valid and the
  Arzel\`a-Ascoli theorem applies leading to a limit curve that
  fulfills the sought requirements. 
  \end{proof}

%%%% proof of proposition
\begin{proof}[Proof of Proposition~\ref{prop: discrete propagation of
    support - boundary}]
  We prove the result in the case $X \varphi>0$, that is $\alpha =1$, on $K \cap \I$. The proof in
  the case $X \varphi<0$ can be written {\em mutatis mutandis}. 

  In $K \cap \I$, having $X \varphi>0$ simply means having $X$ pointing
  towards $\{ \varphi >0\}$. In fact, in $K \cap \I$ one has $X \varphi \geq C_0>0$ for some
  $C_0>0$. Hence, in a bounded open  \nhd $W_1$ in $\R^d$ of $K
  \cap \I$ one has $X \varphi \geq C_0/2$. Introduce also $W_2$ a bounded
  open  \nhd  of $\ovl{W_1}$ where  $X \varphi \geq C_0/4$.
  
  Since $K \setminus W_1$ is compact, with the result of
  Proposition~\ref{prop: symboles}, it suffices to consider the case
  $x^0 \in K \cap W_1 \cap \supp \mu$.  Suppose $x^0$ is such a point.
  In the compact set $\ovl{W_1} \cap K$ one has $c_1 \leq \| X \| \leq C_1$.
  As in the proof of Proposition~\ref{prop: symboles}, performing a rotation
and a dilation of coefficient $\| X(x)\|\in [c_1, C_1]$, one can assume
that $X(x^0) = (1,0, \dots, 0)\in \R^d$.  By abuse of notation, we still
use the letter $\varphi$ for the function used to define the
hypersurface $\I$.

  From the equation \eqref{eq: equation t X mu everywhere} satisfied
  by $\mu$ one has
  \begin{align}
    \label{eq: transp X mu vanishing on varphi q}
    \dup{(\transp{X}-f)\mu}{\varphi q}_{\D^{\prime 1}(\R^d),\Con^1_c(\R^d)}=0,
    \qquad q \in \Con^1_c(\R^d), 
  \end{align}
  since $\varphi \delta_\I=0$. 
  
  Consider $0< \eps\leq 1$ (observe that this case is sufficient for the
  conclusion to hold).  For $\delta>0$, one uses the function
  $\varphi q$ as a test function, with
  $q$ as defined in \eqref{eq: def q eps delta x -
    test function}.  With \eqref{eq: transp X mu vanishing on varphi
    q} one finds
\begin{align}
  \label{eq: action on test function - boundary 1}
  0 &= \dup{(\transp{X}-f)\mu}{\varphi q}_{\D^{\prime 1}(\R^d),
  \Con^1_c(\R^d)}
      = \dup{\mu}{X (\varphi q)- f q}_{\D^{\prime 0}(\R^d), \Con^0_c(\R^d)} \\
  &= 
    \dup{\mu}{\varphi  g}_{\D^{\prime 0}(\R^d),\Con^0_c(\R^d)}  
    + \dup{\mu}{\varphi h}_{\D^{\prime
    0}(\R^d),\Con^0_c(\R^d)}
    +\bigdup{\mu}{(X \varphi) q}_{\D^{\prime 0}(\R^d), \Con^0_c(\R^d)}
    \notag \\
    &\quad 
    +\bigdup{\mu}{(A (X x_1)-f) \varphi q}_{\D^{\prime 0}(\R^d), \Con^0_c(\R^d)}, 
    \notag
\end{align}
recalling that
$X q = g + h + A (X x_1)q$ with
$g$ and $h$ given in \eqref{eq: def q eps delta x - test function}. 

In $W_2$ one has $X \varphi \geq C_0/4$.  If
$\delta_0>0$ is chosen \suff
small,  one has $\supp q \subset W_2$
uniformly with respect to $x^0 \in W_1$ for $0 < \delta\leq
\delta_0$.  As $q\geq 0$ one finds $(X \varphi)
q\geq 0$. Moreover $\big( (X \varphi) q\big)(x^0)
>0$ and since $x^0 \in \supp \mu$ this yields 
\begin{align*}
  %\label{eq: action on test function - boundary 2}
  \bigdup{\mu}{(X \varphi) q}_{\D^{\prime 0}(\R^d), \Con^0_c(\R^d)}  
    >0. 
\end{align*}
By Lemma~\ref{lem: positivite g},  for $\delta_0$
chosen \suff small and $A$ \suff large, one has
$g\geq 0$ if $0< \delta\leq \delta_0$ and $Ax_1-f\geq 0$.  As $\varphi
\geq 0$ in $\supp \mu$ one thus
finds
\begin{align}
  \label{eq: action on test function - boundary 3}
  \dup{\mu}{\varphi  g}
  +\bigdup{\mu}{(A (X x_1)-f) \varphi q}_{\D^{\prime 0}(\R^d), \Con^0_c(\R^d)}
  \geq 0, \quad \text{for} \ 0 < \delta \leq \delta_0.
\end{align}
From \eqref{eq: action on test function - boundary 1}--\eqref{eq:
  action on test function - boundary 3} one obtains
$\dup{\mu}{\varphi h}_{\D^{\prime 0}(\R^d),\Con^0_c(\R^d)}  <0$. 
By the support estimate for $h$ given in
\eqref{proche3} the conclusion follows.
\end{proof}

\medskip
The following lemma is not of direct use  in this section but it is used in
another section below. 
%%%%%%%%%%%%%%%%%%%%%%%%
% lemma                %
%%%%%%%%%%%%%%%%%%%%%%%%
\begin{lemma}
  \label{lem: support goes to the trace - transverse X}
  Suppose $\mu$ is a nonnegative measure solution to \eqref{eq: equation t
    X mu everywhere} that vanishes in $V^-$.  Suppose $x \in \I$. Then
  $x \in \supp \mu$ if and only if $x\in \supp(\tilde{\mu}\otimes \delta_\I)$
\end{lemma}
\begin{proof}
  If $x \notin \supp \mu$ then $\mu$ vanishes in a \nhd of
  $x$ and from \eqref{eq: equation t
    X mu everywhere} one has $\tilde{\mu} \otimes \delta_\I$ vanishing
  in that \nhd. This gives $x \notin \supp(\tilde{\mu}\otimes \delta_\I)$.

  Suppose
  now that $x\notin \supp(\tilde{\mu}\otimes \delta_\I)$. Then, from \eqref{eq:
    equation t X mu everywhere}, the equation fulfilled by $\mu$ is
  $\transp{X}\mu =f \mu$ locally near $x$. By Theorem~\ref{th: ODE}, if
  $x \in \supp \mu$ then there exist $\Slim>0$ and a integral curve
  $\gamma(s)$ of $X$ such that $\gamma(0)= x$ and
  $\{ \gamma(s)\}_{s\in ]-\Slim,\Slim[} \subset \supp \mu$. Yet, as $X$ is
  transverse to $\I$ here, half of the integral curve lies in $V^-$
  where $\mu$ vanishes; this gives a contradiction. Hence, if
  $x \notin \supp(\tilde{\mu}\otimes \delta_\I)$ then $x\notin \supp \mu$.
\end{proof}

%%%%%%%%%%%%%%%%%
% section
%%%%%%%%%%%%%%%%%
\section{Geometrical setting II}
\label{sec: setting at the boundary}

Here, we carry on with the introduction of the geometrical notions to be used
in the subsequent sections. In Section~\ref{sec: Geometrical setting-intro} we used the
quasi-normal geodesic coordinates of Proposition~\ref{prop:
  quasi-normal coordinates} to obtain a `straight path' towards the
necessary notions for the statement of the main result in
Theorem~\ref{theorem: measure support propagation}: glancing and
hyperbolic regions, hamiltonian vector field and gliding vector field,
\bichars and \gbichars. In the present section we provide
additional results and notions. Yet, we do not rely on quasi-normal
geodesic coordinates for the following two reasons: (1) the
simplifications provided by such coordinates at the boundary hide some
of the geometrical properties, and more important, (2) we wish to 'push' the definition of
the glancing and hyperbolic regions and gliding vector field away from the boundary to ease
arguments in the proof of Theorem~\ref{theorem: measure support
  propagation}: extending $\pdTL$ away from $\dTL$ we obtain a
foliation of $\TL$. Since the advantageous structure of quasi-normal
geodesic coordinates is lost away from the boundary it is better to
work in arbitary coordinates from the beginning. Note that due to the
considered low regularity of the coefficients, the foliation we
introduce is not a geometrical object in the sense that it depends on
the chosen coordinates and on an extension of the conormal vector
field. Yet, this foliation is only used in a single local chart
in what follows. 

In  Section~\ref{sec: Geometrical setting-intro}, $T_x^*\d\L$ was
identified with the set of conormal vectors $\xi = (\xi',0)$. This is
not natural in general.  In fact, if $m^0 \in \d\M$ set $x^0 =
\chdiff(m^0)$, for a local chart $(\hO, \chdiff)$. One has  $x^0 =(x^{0\prime}, 0)$. 
The injection $\d\M \to \M$ yields a  natural
injection of $T_{m^0} \d\M$ into $T_{m^0} \M$ and  the {\em surjection} of
$T_{m^0}^* \M$ into $T_{m^0}^* \d\M$ by duality that take the form 
$v' \mapsto (v',0)$ and $(\xi',\zeta)  \mapsto \xi'$ respectively in the considered local coordinates
for $v' \in T_x \d\M$ and  $(\xi',\zeta) \in T^*_x \M$.
As most often done, $T_x \d\M$ is naturally viewed as a linear 
subspace of $T_x \M$ and in the chosen coordinates,  $v' \in  T_x
\d\M$ identifies with  $v= (v',0) \in  T_x \d\M$. However, $T^*_x \d\M$
is identified with the set of covectors orthogonal to the unit
vector field $\nx$ at the boundary. In general such covectors do not
take the form $(\xi',0)$; yet, they do in quasi-normal geodesic coordinates.

Note that if a notation appearing in  Section~\ref{sec: Geometrical
  setting-intro} is used  in what follows, say $\pdTL$, it denotes
the same object.

\begin{remark}
  In what follows, since the metric $g$ is also defined in $\hM
  \setminus \M$ we also
consider \bichars that leave or enter $\TL$. 
To avoid possible confusion we write $\Char p \cap \TL$ or
$\Char p \cap \{ z \geq 0\}$ if only considering characteristic points in
the cotangent bundle $\TL$ and not the extension made outside $\M$ and
$\L$.
\end{remark}

%%%%%%%%%%%%%%%%%
% sub-section
%%%%%%%%%%%%%%%%%
\subsection{Musical isomorphisms, normal and conormal vectors}
\label{sec: Musical isomorphisms, normals}

Consider a point $(x,v) \in  T \M$. As is done classicaly, denote by $v^\flat$,
and $(x,v)^\flat$ by  extension, the unique element  of $T_x^*\M$ such
that
\begin{align*}
  \dup{v^\flat}{u}_{T_x^* \M, T_x \M} = g_x(v, u), \qquad u \in T_x \M.
\end{align*}
In local coordinates, this reads $(v^\flat)_i = g_{ij} (x) v^j$, $1\leq i \leq d$. 
One thus obtains a map $\flat: v \mapsto v^\flat$ from $T_x \M$ into
$T_x^*\M$, and by extension from $T \M$ into
$\TM$. With the invertibility of $g_x = (g_{ij}(x) )_{i,j}$ one
readily sees that $\flat$ is an isomorphism. Moreover, one has
$g_x^*(v^\flat, v^\flat) = g_x (v,v)$, meaning that $\flat$ is  an isometry. 

The inverse isometry is denoted by $\sharp$. One has
$\sharp: \xi \mapsto \xi^\sharp$ from $T_x^* \M$ onto
$T_x\M$, and by extension from $\TM$ onto
$T\M$. One has 
\begin{align*}
  \dup{\xi^\sharp}{\omega}_{T_x \M, T_x^* \M} = g^*_x(\xi, \omega), \qquad \omega \in T^*_x \M,
\end{align*}
and in local coordinates 
\begin{align*}
  (\xi^\sharp)^i =  g^{ij} (x) \xi_j, \quad 1\leq i \leq d.
\end{align*}
One can also write, for $\xi \in T_x^* \M$ and $v \in T_x\M$
\begin{align*}
  \dup{\xi}{v}_{T_x^* \M, T_x \M} 
  = g_x(\xi^\sharp, v) 
  = g_x^*(\xi, v^\flat)
  = \dup{\xi^\sharp}{v^\flat}_{T_x \M, T_x^* \M} .
\end{align*} 
From~\eqref{eq: Hp} one finds
\begin{align*}
  %\label{eq: Hp 2}
  \Hp (\y) = - 2 \tau \d_t + 2 (\xi^\sharp)^j \d_{x_j} 
  - \d_{x_k} g^{ij}(x) \xi_i \xi_j  \d_{\xi_k}.
\end{align*}

For $x = (x',0) \in \d\M$, denote by $\nx \in T_x \M$ the unitary
{\em inward}~\footnote{Here, we choose the inward direction for $\nx$
  to be consistent with having $z>0$ if $(x',z,\xi) \in \TM$ as far
  as sign are concerned.}
pointing normal vector to $\d\M$,  meaning that $g_x(\nx, \nx)=1$ and
$g_x (\nx, v) =0$ for all $v \in T_x \d\M$ and $\nx^d >0$. 
Set $\nxs = \nx^\flat \in T_x^* \M$. One has  $g_x^*(\nxs, \nxs)
=1$ and 
\begin{align*}
 \dup{\nxs}{v}_{T^*_x \M, T_x\M} =0, \quad v \in T_x \d\M. 
\end{align*}
In a local chart at the boundary as in \eqref{eq: local chart boundary} this gives $\nxs =
\big(0,\dots,0,  (g^{dd}(x))^{-1/2}\big)$, that is, $\nxs =
(g^{dd}(x))^{-1/2}  d x_d$. One  deduces
\begin{align*}
  \nx^i
  = \big( (\nxs)^\sharp\big)^i
  = g^{ij}(x) (\nxs\big)_j
  = (g^{dd})^{-1/2}  g^{id}(x).
\end{align*}
Note that $\nx$ and $\nxs$ have $\Con^1$ regularity. 

%%%%%%%%%%%%%%%%%%%%%%%%
% remark               %
%%%%%%%%%%%%%%%%%%%%%%%%
\begin{remark}
  \label{remark: extension unitary normal vectors}
  We insist on the fact that $\nx$ and $\nxs$ are here defined on the
  boundary only. A natural extension away from the boundary would use normal geodesic coordinates that are not
  available here (see Remark~\ref{remark: no normal geodesic
    coodinates}) or a local geodesic flow but the latter may not
  exist due to the potential lack of uniqueness of geodesics here.
 % We do not know of any {\em geometric} way to extend
 %  $\nx$ and $\nxs$ with the same $\Con^1$ regularity. 
\end{remark}

%%%%%%%%%%%%%%%%%
% sub-section
%%%%%%%%%%%%%%%%%
\subsection{Partition of the cotangent bundle revisited}
\label{sec: A partition of the cotangent bundle at the boundary}

Consider a local chart at the boundary as in \eqref{eq: local chart
  boundary}.
As above, denote $z=x_d$ with $x=(x',z)$ and, accordingly, the
associated cotangent variables read $\xi = (\xi', \zeta)$. 
For $x \in \d\M$, denote by $\pT^*_x \M$ the orthogonal of $\nxs$ in the
sense of $g_x^*$, that is, 
\begin{align*}
  \pT^*_x \M = \{ \xi \in T^*_x \M; \ g_x^*(\xi, \nxs)=0\}.
\end{align*}
One has $\pT^*_x \M  =\flat (T_x \d\M)$ and 
$T^*_x \M  =  \pT^*_x \M \oplus \Span(\nxs)$.
Denote by $\ppi$ the orthogonal projection onto $\pT_x^* \M$.
For $\xi \in T_x^* \M$ set $\pxi = \ppi (\xi)$
that reads 
\begin{align}
  \label{eq: def pxi}
  \pxi = \xi - g_x^*\big( \xi , \nxs) \nxs.
\end{align}
If $\y =(t,x,\tau,\xi)\in \TL$  the following computations in the
considered local coordinates are useful
  \begin{align}
    \label{eq: form Hp z}
    &\Hzzp (\y) = \Hzzp (x) = 2 g^{d d}(x) \neq 0,\\
    &\Hpz (\y) = \Hpz (x,\xi)= 2 g^{dj}(x) \xi_j
      = \big( 2  \Hzzp (x))^{1/2} g^*_x(\xi, \nxs),\nonumber
  \end{align}
  yielding, for $\y \in \dTL$, that is, $x \in \d \M$, 
  \begin{align*}
    \pxi =  \xi -  \alpha  \Hpz (x,\xi) \nxs, \quad
    \ \ \text{with} \ \alpha (x)= \big( 2  \Hzzp (x))^{-1/2},
  \end{align*}
  as  $\nxs = \big(0,\dots,0,  2\alpha(x)\big)$.
  In local coordinates this gives 
\begin{align}
  \label{eq: ppi local coordinates}
  \pxi &= \Big(\xi_1, \dots, \xi_{d-1}, \zeta -
         \frac{\Hpz}{\Hzzp}(x,\xi)\Big)
         =\Big(\xi_1, \dots, \xi_{d-1},
         -  \frac{2}{\Hzzp} (x) \sum_{j=1}^{d-1} g^{j d}(x)\xi_j
        \Big).
\end{align}
Above was mentionned the surjection $T_x^* \M$ into
$T_x^*\d\M$. Consider the map 
\begin{align*}
 \pT^*_x \M &\to T_x^*\d\M\\
 (\xi',\zeta) &\mapsto \xi',
\end{align*}
One finds that it is an isomorphism, giving a geometrical
identification of $ T_x^*\d\M$ as a subspace of $T_x \M$. However, we
keep the notation $\pT^*_x \M$ to avoid any possible confusion.

One also denotes by $\Sigma$  the orthogonal symmetry 
with respect to $\pT_x^*\M$,
that is, 
\begin{align*}
  %\label{eq: def sxi}
  \Sigma(\xi) &= \xi - 2 g_x^*\big( \xi , \nxs) \nxs\\
       &=  \xi - 2  \alpha \Hpz(x,\xi) \nxs
         = \pxi -  \alpha \Hpz(x,\xi) \nxs. \notag
\end{align*}
In local coordinates this gives 
  $\Sigma(\xi) = \Big(\xi_1, \dots, \xi_{d-1}, \zeta - 2\frac{\Hpz}{\Hzzp}(x,\xi)\Big)$.
One has $\Sigma(\xi) + \xi = 2\, \pxi$. 

Accordingly, for $x \in \d\M$ set
\begin{align*}
  \pT^*_{t,x}\L = \{ (\tau,\xi)\in T^*_{t,x}\L; \ \xi \in \pT_x^* \M\},
\end{align*}
and 
\begin{align*}
   \pdTM = \bigcup_{x\in \d\M} \{ x\} \times \pT_x^* \M, 
  \quad 
  \pdTL = \bigcup_{(t,x)\in \d\L} \{ (t,x)\} \times   \pT^*_{t,x}\L,
\end{align*}
and for $\y = (t,x,\tau, \xi) \in \d\L$ one writes $\ppi(\y) = \py = (t,x,\tau, \pxi)$ and
$\Sigma(\y) = (t,x,\tau, \Sigma(\xi))$.

Naturally, for $\y \in \dTL$ one has $\ppi(\y) = \ppi(\Sigma(\y)) = \py$ and  
\begin{align*}
  \y \in \pdTL 
  \ \ \Equiv \ \ \py = \y
  \ \ \Equiv \ \ \Sigma(\y) = \y.
\end{align*}

From what is written above for $x \in \d\M$ one has 
\begin{align}
  \label{eq: form xi sxi}
    &\xi =  \pxi +  (\alpha \Hpz) (\y)\, \nxs 
      = \pxi +  \alpha(x)  \Hpz (x,\xi) \, \nxs , \\
    &\Sigma(\xi) =  \pxi +  (\alpha \Hpz) (\Sigma(\y)) \, \nxs
      = \pxi -   \alpha (x)  \Hpz (x,\xi) \, \nxs \notag
\end{align}
yielding
  \begin{align}
    \label{eq: change sign Hp z Sigma}
    \Hpz (\Sigma(\y)) = - \Hpz (\y), \qquad \y \in \dTL.
  \end{align}
From~\eqref{eq: form Hp z} one also has the following
characterization  of $\pdTL$ 
  \begin{align}
    \label{eq: characterization pT*L}
    \pdTL = \{ z=\Hpz=0\}.
  \end{align}
  Note that for any $\Con^2$-function $\phi$ with $d \phi \neq 0$
  such that $\d\M = \{ \phi=0\}$ locally, one finds 
  \begin{align*}
     z=\Hpz=0 \ \ \Equiv \ \  \phi = \Hp \phi=0,
  \end{align*}
  meaning that the use of \eqref{eq: characterization pT*L} 
  made below is not coordinate dependent. 

  One can observe  that $\pdTL$ is a symplectic submanifold of
  $\TL$.

  \medskip
  Note that one has
  \begin{align}
  \label{eq: d xid Hpz}
  \d_{\zeta}\Hpz = - \Hz \Hpz = \Hz \Hzp = \Hzzp,
\end{align}
and 
\begin{align}
 \label{eq: pxi indep xid}
  \Hz \pxi = - \d_{\zeta} \pxi =0.
\end{align}

  \medskip Suppose $\y \in \dTL$. Since $\ppi$ maps $\dTL$ into $\pdTL$,
  then $d \ppi(\y)$, its differential at $\y$, maps $T_\y \dTL$ into
  $T_{\py} \pdTL$. 
 %%%%%%%%%%%%%%%%%%%%%%%%
% lemma                %
%%%%%%%%%%%%%%%%%%%%%%%%
\begin{lemma}
  \label{lemma: kernel diff ppi}
  Suppose $\y \in \dTL$. 
  One has $\ker (d \ppi(\y)) = \Span(\Hz)(\y).$
\end{lemma}
In local coordinates one has $\Hz = - \d_{\zeta} = (0,\dots, 0, -1) \in
\R^{2 d+2}$ at $\y \in \dTL$. 
%%%% proof of lemma
\begin{proof}
  Consider $v \in T_\y \dTL$,  that is, $v \in
  \Span\{\d_t, \d_{x_i}, \d_\tau,\d_{\xi_j}\}$, $i=1,\dots,d-1$,
  $j=1,\dots, d$. From the form of
  $\ppi$ given in \eqref{eq: ppi local coordinates} one has 
  \begin{align*}
    d \ppi (\y) (v) 
    = v -  d (\Hpz/\Hzzp) (\y) (v)\d_{\zeta}
  \end{align*}
  If one has $d \ppi (\y) (v)=0$ then one sees that 
  $v \in \Span(\Hz)$. 

  Conversely, if $v = \d_{\zeta}$ one has 
  $d \ppi (\y) (v) = \big(1 - d (\Hpz/\Hzzp) (\y) (\d_{\zeta})\big) \d_{\zeta}$
  and, using that $\Hzzp(\y)$ is independent of $\xi$, 
  \begin{align*}
   d (\Hpz/\Hzzp) (\y) (\d_{\zeta}) = \d_{\zeta} (\Hpz/\Hzzp)(\y)
    = \frac{\d_{\zeta} \Hpz(\y)}{\Hzzp(x)} = \frac{\Hzzp}{\Hzzp}(x) = 1, 
  \end{align*}
by \eqref{eq: d xid Hpz},
  implying  $d \ppi (\y) (v) =0$. 
\end{proof}
If $x \in \d\M$ observe with what precedes that the maps
\begin{align}
  \label{eq: parametrization cotangent bundles}
  (x,\xi) \mapsto (x,\pxi, \Hpz (x,\xi)),
            \ \ \text{and} \ \ 
  \y \mapsto (\py, \Hpz (\y)),
\end{align}
yield natural parametrizations of $\dTM$ and $\dTL$ that is
used in what follows. Observe however that these coordinates are only
$\Con^1$.

%%%%%%%%%%%%%%%%%%%%%%%%
% definition           %
%%%%%%%%%%%%%%%%%%%%%%%%
\begin{definition}[outward, inward, and tangentially pointing vectors]
  \label{def: pointing outward or inward}
  Consider $\y = (t,x,\tau,\xi) \in \TL$ with $x \in \d\M$. 
  One says that 
  \begin{enumerate}
    \item $\y$ points strictly outward if $\Hpz (\y)<0$;
    \item $\y$ points tangentially if $\Hpz (\y)=0$, that is, $\y \in
      \pTL$;
    \item $\y$ points strictly inward if $\Hpz (\y)>0$.
  \end{enumerate}
\end{definition}
Set $v = \xi^\sharp$.  One has $v^d = g^{dj} \xi_j = \Hpz(\y)/2 =
\alpha^{-1} g_x^*(\xi, \nxs)/2$. The terminology of  Definition~\ref{def: pointing outward or
  inward} is thus related to the sign of $v^d$ (and not that of
$\zeta$) and, as we shall see
below, to the behavior of \bichar that goes through $\y$ if
moreover $\y \in \Char p$; see Lemma~\ref{lem: direction bichar hyperbolic pt}.
With
\eqref{eq: change sign Hp z Sigma} one has the
following properties. 
%%%%%%%%%%%%%%%%%%%%%%%%
% lemma                %
%%%%%%%%%%%%%%%%%%%%%%%%
\begin{lemma}
  \label{lem: Property Sigma}
  Consider $\y \in \dTL$. One has
  \begin{enumerate}
    \item $\y$ strictly points inward if and only if 
  $\Sigma(\y)$ strictly points outward;
\item $\y$ points tangentially if and only if $\Sigma(\y) = \y$. 
    \end{enumerate}
\end{lemma}

\bigskip
As in Section~\ref{sec: A partition of the cotangent bundle at the boundary-intro}
the vector bundle $\pdTL$ is written as the union of the three
bundles $\pEb$, $\pGb$, $\pHb$. Their definition is identical to that
given in Definition~\ref{def: E', H', G'-intro}, which we reproduce
here to ease reading. 
%%%%%%%%%%%%%%%%%%%%%%%%
% definition           %
%%%%%%%%%%%%%%%%%%%%%%%%
\begin{definition}[elliptic, glancing, and hyperbolic regions]
  \label{def: E', H', G'}
  One partitions $\pTL$ into three homogeneous regions.
  \begin{enumerate}
  \item The elliptic region $\pEb= \pdTL \cap \{ p > 0\}$; if
    $\y\in \pEb$ it is called an elliptic point.          	
  \item The glancing region $\pGb = \pdTL \cap \{ p = 0\}$; if $\y
    \in \pGb$ it is called a glancing point.
  \item The hyperbolic region $\pHb = \pdTL \cap \{ p < 0\}$; if
    $\y\in \pHb$ it is called a   hyperbolic point.
  \end{enumerate}
\end{definition}
%%%%%%%%%%%%%%%%%%%%%%%%
% lemma                %
%%%%%%%%%%%%%%%%%%%%%%%%
\begin{lemma}
  \label{lemma: relevement pG pH}
  The sets $\pEb$, $\pGb$ and $\pHb$ are also characterized by 
  \begin{align*}
  &\y \in  \pEb \ \ \Equiv \ \
        \y \in \pdTL \ \  \text{and} \ \ 
    \ppi^{-1}(\{ \y\}) \cap \Char p=\emptyset,\\
  &\y \in  \pGb \ \ \Equiv \ \
     \y \in \pdTL \ \  \text{and} \ \ 
    \ppi^{-1}(\{ \y\}) \cap \Char p=\{ \y\}, \\
  &\y \in  \pHb \ \ \Equiv \ \
     \y \in \pdTL \ \  \text{and} \ \ 
      \ppi^{-1}(\{ \y\}) \cap \Char p=\{ \y^-, \y^+\},
\end{align*}
where, in the last case, $\y^\pm = (t,x,\tau,\xi^\pm)$ if $\y=
(t,x,\tau,\xi)$ with 
\begin{align}
  \label{eq: relevement pG pH}
  \xi^+ = \xi + \lambda  \nxs
 \ \ \text{and} \ \ 
  \xi^- = \xi - \lambda \nxs,
\end{align} 
with $\lambda = \big(- p (\y)\big)^{1/2} = \big(\tau^2 - |\xi|_x^2 \big)^{1/2}$.
\end{lemma}
For $\y \in \pHb$ the notation $\y^\pm$ is used in what follows
with the definition given in this lemma. By \eqref{eq: form xi sxi} one has 
\begin{align*}
  %\label{eq: Hpz lambda}
 \alpha(x) \Hpz (\y^+) =\lambda >0 
  \ \ \text{and} \ \ 
  \alpha(x)  \Hpz (\y^-) = -  \lambda < 0,  
\end{align*}
with $\alpha(x) =\big( 2 \Hzzp (x)\big)^{-1/2}$,
that is, $\y^+$ points {\em inward} and $\y^-$ points {\em outward} in
the sense given in Definition~\ref{def: pointing outward or inward}.
%%%% proof of lemma
\begin{proof}[Proof of Lemma~\ref{lemma: relevement pG pH}]
If $\y = (t,x; \tau, \xi)\in \dTL$, then $\tilde{\y} \in \ppi^{-1}(\{ \y\})$
reads $\tilde{\y} = (t,x; \tau, \tilde{\xi})$ with $\tilde{\xi}  = \xi
+ \lambda \nxs$ for some $\lambda \in \R$.  And one has $\norm{\tilde{\xi}}{x}^2
= \norm{\xi}{x}^2 + \lambda^2$ and $p(\tilde{\y}) = p (\y)+ \lambda^2$.
If $\y \in \pE$ one has $p (\y) >0$ and thus no choice
of $\lambda \in \R$ can yield $p(\tilde{\y}) =0$.
If $\y \in \pG$ one has $p (\y) = 0$ and  thus
the only choice of $\lambda \in \R$ to have $p(\tilde{\y}) = 0$ is $\lambda=0$.
If $\y \in \pH$ one has $p (\y) < 0$ and  thus
one has $p(\tilde{\y}) = 0$ if and only if $\lambda= \pm \big(- p (\y)\big)^{1/2}$. 
\end{proof}
The following definition is counterpart to Definition~\ref{def: H, G-intro} yet not
coordinate dependent. 
%%%%%%%%%%%%%%%%%%%%%%%%
% definition           %
%%%%%%%%%%%%%%%%%%%%%%%%
\begin{definition}[partition of $\Char p \cap \dTL$]
  \label{def: G H}
  Set
  \begin{align*}
    \Gb = \{ \y \in \dTL; \ p(\y)=0 \ \text{and}\  \Hpz(\y)=0\},
  \end{align*}
  and $\Hb = \Hb^+ \cup \Hb^-$ with 
  \begin{align*}
    \Hb^\pm = \{ \y \in \dTL; \ p(\y)=0 \ \text{and}\  \Hpz(\y) \gtrless 0\},
  \end{align*}
  that is, the set of characteristic points at the boundary
  that point tangentially ($\Gb$), strictly inward ($\Hb^+$), and
  outward ($\Hb^-$). Recall that $\dTL$ is locally $\{ z=0\}$. 

  Together $\Gb$ and $\Hb$ (\resp $\Gb$, $\Hb^+$, and $\Hb^-$) form a
  partition of $\Char p \cap \dTL$.
\end{definition}
The index $\d$ in Definition~\ref{def: G H} expresses that only
boundary points are considered, that is $z=0$. At
places we use extensions of these sets  away from the boundary. Yet,
this is done in local charts only and not in a geometrically
invariant way; see Section~\ref{sec: local extension away from boundary}.

The sets $\Hb^\pm$ are connected and open in $\Char p \cap \dTL$.
The set $\Gb$ is a connected and closed subset of
$\Char p \cap \dTL$. 

We now make the connexion
between $\pGb$ and $\Gb$ on the one hand, and $\pHb$ and $\Hb$ on the
other hand.  One has $\pGb = \Char p \cap \pdTL$. Since $\pdTL$ is
charaterized by $\Hpz=0$ and $z=0$, see \eqref{eq: characterization pT*L}, one finds that $\pGb = \Gb$ in
fact. 

Consider $\y \in \pHb$ and $\y^\pm$ as given by Lemma~\ref{lemma:
  relevement pG pH}. One has $\y^\pm$ in $\Char p$ and 
\begin{align}
  \label{eq: connection Hpz lambda}
  \alpha(x)\Hpz (\y^\pm) =   g^*_x(\xi^\pm, \nxs) = \pm  \lambda
\end{align}
with $\lambda>0$ as in Lemma~\ref{lemma:
  relevement pG pH}. 
Thus $\y^+ \in \Hb^+$ and $\y^- \in \Hb^-$. One also has $\ppi(\y^\pm) =
\y$. 

One has thus obtained the following proposition.
%%%%%%%%%%%%%%%%%%%%%%%%
% proposition               %
%%%%%%%%%%%%%%%%%%%%%%%%
\begin{proposition}
  \label{prop: relevement H G}
  One has $\ppi(\Gb) = \Gb= \pGb$ and $\ppi(\Hb)  = \pHb$ and, conversely, 
  \begin{align*}
    \ppi^{-1} (\pGb) \cap \Char p = \Gb = \pGb, \qquad \ppi^{-1} (\pHb) \cap \Char p = \Hb.
  \end{align*}
  
  One has $\Sigma(\Gb) = \Gb$ and $\Sigma(\Hb^+) = \Hb^-$.
  \begin{enumerate}
  \item  If $\y \in \Gb = \pGb$ then $\pi^{-1} (\{ \y\}) \cap \Char p =\{ \y\}$
    and $\Sigma(\y) = \y$. 
  \item 
  If $\y \in \pHb$ then
  $\pi^{-1} (\{ \y\}) \cap \Char p =\{ \y^+, \y^-\}$ with
  $\y^+ \in \Hb^+$ and $\y^- \in \Hb^-$, as given in Lemma~\ref{lemma:
    relevement pG pH}, and $\Sigma (\y^+) = \y^-$.
  Conversely, if $\y\in \Hb^\pm$ then $\Sigma(\y)
  \in \Hb^\mp$ and $\ppi(\y) = \ppi \big(\Sigma(\y)\big) \in \pHb$.  
  \end{enumerate}
\end{proposition}
By extension, if $\y\in \Hb$ one also says that $\y$ is a
hyperbolic point. The set $\Hb$ is also
called the hyperbolic region. 

%%%%%%%%%%%%%%%%%%%%%%%%
% remark               %
%%%%%%%%%%%%%%%%%%%%%%%%
\begin{remark}
  \label{remark: regularity Hpz}
The function $\Hpz$ is key in the description of the regions $\Gb$ and
$\Hb^\pm$. Having $p$ only $\Con^1$ one may think that $\Hpz$ is only
$\Con^0$. However, above we computed
\begin{align}
  \label{eq: formula Hp z}
  \Hpz (\y) = 2 (\xi^\sharp)^d = 2 g^{dj}(x) \xi_j.
\end{align}
One thus sees that $\Hpz$ is in fact a $\Con^1$-function of $\y$.
\end{remark}
If $\y = (t,x,\tau,\xi) \in\Char p \cap \dTL= \Hb^\pm \cup \Gb$ with
what we wrote above one has
\begin{align}
  \label{eq: decomposition point in char p}
  \xi = \pxi + \alpha \Hpz(\y) \nxs,
  \quad \text{with} \ \
   \alpha \Hpz (\y) = \begin{cases}
    \big( - p (\py)\big)^{1/2}>0 & \text{if}\ \y \in \Hb^+,\\
    0 & \text{if}\ \y \in \Gb,\\
    - \big( - p (\py)\big)^{1/2}<0 & \text{if}\  \y \in \Hb^-.
  \end{cases}
\end{align}

%%%%%%%%%%%%%%%%%%%%%%%%
% remark               %
%%%%%%%%%%%%%%%%%%%%%%%%
\begin{remark}
  \label{remark: limit r pm -> r}
  Note that $\ovl{\Hb^\pm} = \Hb^\pm \cup \Gb$, implying $\ovl{\Hb} = \Hb
  \cup \Gb = \Char p \cap \dTL$. If 
  $(\y^{n})_{n \in \N} \subset \pHb$ converges to
  $\y \in \Gb$, then $(\y^{n})^\pm \to \y$ as $n \to +\infty$. 
\end{remark}

\begin{figure}
  \begin{center}
    \resizebox{3.5cm}{!}{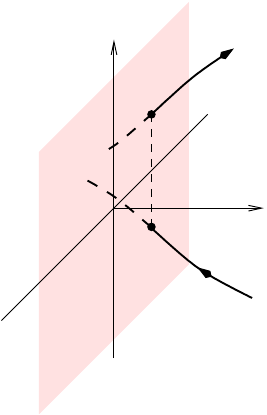}
    \caption{Two \bichars; one going through $\y^+ \in \Hb^+$ and one
      through $\y^- \in \Hb^-$. Here $\y^+ = \Sigma(\y^-)$.}
  \label{fig: bichar H}
  \end{center}
\end{figure}

\medskip The notions of Definition~\ref{def: pointing outward or
  inward} and the description given in Proposition~\ref{prop:
  relevement H G} are of importance because of the following result,
using that for a bicharacteristic $\gamma(s)$, the value of $\Hp f
(\gamma(s))$, is equal to the derivative of $s \mapsto f(\gamma(s))$.
%%%%%%%%%%%%%%%%%%%%%%%%
% lemma                %
%%%%%%%%%%%%%%%%%%%%%%%%
\begin{lemma}
  \label{lem: direction bichar hyperbolic pt}
  Consider $\y  \in \Char p$. Denote by
  $\gamma(s) = (t(s), x(s), \tau(s), \xi(s))$
  a \bichar with $x(s) = (x'(s), z(s))$ such that $\gamma(0) = \y$.  
  \begin{enumerate}
  \item If $\y \in \Gb$ then $\frac{d}{d s} z_{|s=0}=(\Hpz)(\gamma(0))=0$;
  \item If $\y\in \Hb^\pm$ then $\frac{d}{d s} z_{|s=0}=(\Hpz)(\gamma(0))\gtrless 0$.
  \end{enumerate}
\end{lemma}
Thus, for $\y \in \Gb$, a glancing point at the boundary, any \bichar that goes through $\y$ is
tangent to $\dTL$.  For $\y \in \Hb$, any \bichar that
goes through $\y$ is transverse to $\dTL$, either entering $\TL$
if $\y \in \Hb^+$, or exiting $\TL$ if $\y\in \Hb^-$. This is
illustrated in Figure~\ref{fig: bichar H}. However, the
geometry of a \bichar that goes through a glancing point needs to be
further described. This is the subject of Section~\ref{sec: description glancing region}.
%%%%%%%%%%%%%%%%%%%%%%%
% figure
%%%%%%%%%%%%%%%%%%%%%%%

\bigskip
We conclude this section by noting that the glancing set $\pGb=\Gb$ is a
submanifold away from $0$ (in the cotangent variable). The
singularity at $0$ comes from its natural conic structure. 
%%%%%%%%%%%%%%%%%%%%%%%%
% proposition               %
%%%%%%%%%%%%%%%%%%%%%%%%
\begin{proposition}[submanifold property of $\Gb$]
  \label{prop: G submanifolds}
  The set $\pdTL$ is a $\Con^1$-submanifolds of $\TL$ of codimension
  two respectively.
Away from $(\tau,\xi) =(0,0)$, the set $\Gb$ is a
  $\Con^1$-submanifold of $\TL$ of codimension three defined by
  $p=z=\Hpz=0$ (and thus a $\Con^1$-submanifold of $\pdTL$
  of codimension one).
\end{proposition}
%%%% proof of proposition
\begin{proof}
  The first result is clear since $\pdTM$ is the orthogonal of
  $\nxs$ in $\{z=0\}$, the metric is $\Con^1$, and $\nxs$ and $d z$ are linearly independent.

  \medskip The second result amounts to proving that $p$, $z$, and
  $\Hpz$ have differentials of rank three at a point $\y = (t,x, \tau, \xi) \in \Gb$
  with $(\tau,\xi)\neq (0,0)$.  Let us consider such a point.
  
  Observe that if $p(\y) =0$ then $(\tau,\xi) =(0,0) \ \Equiv \ \tau
  =0 \ \Equiv \xi=0$. One thus has $\tau \neq 0$. Assume that $\alpha d z(\y)  + \beta d p(\y)  + \gamma d(\Hpz)(\y)
  =0$ for some $\alpha, \beta,\gamma \in \R$.
  As $\d_\tau z = \d_\tau \Hpz=0$, and  $\d_\tau p = 2 \tau \neq 0$  one obtains $\beta=0$. 
  Since $\d_{\zeta} z =0$ and $\d_{\zeta} \Hpz = \Hzzp\neq 0$ one
  also has $\gamma=0$. Finally, as $d z\neq 0$ one finds $\alpha=0$.  
\end{proof}

%%%%%%%%%%%%%%%%%
% sub-section
%%%%%%%%%%%%%%%%%
\subsection{Some local extensions away from the boundary}
\label{sec: local extension away from boundary}

Above, the glancing set $\Gb$ and the hyperbolic sets $\Hb$ are defined
at the boundary, in some geometrical fashion, that is, independently
of the chosen local coordinates. In
Section~\ref{sec: Propagation of the measure support} it is
convenient to ``push'' the notions of glancing set and the hyperbolic
sets away from $\dTL$, that is $\{z=0\}$. Yet, as mentionned above there is no
geometrical way to extend  $\nxs$ away from $\d\M$ in a $\Con^1$
geometrical fashion. Still, the
construction of  \gbichars performed in
Section~\ref{sec: Propagation of the measure support} only relies
on local arguments. Here, we thus extend the formentionned
notions away from the boundary, yet only in a fixed local chart.

Consider a local chart $\chart =(\hO,\chdiff)$ at the boundary 
as in \eqref{eq: local chart boundary} where
the boundary is given by $\{ z=0\}$. 

Extend $\nxs$ to be equal to
$\nxsC = \big(0,\dots,0, (g^{dd}(x))^{-1/2}\big)$, that is,
$\nxsC = (g^{dd}(x))^{-1/2} d x_d \in T_x^* \M$ above all points $x$ of
the chart. (The use of the notation  $T_x^* \M$ is quite abusive but
we have now been \suff clear that the extension is not geometrical by
any means.)

For any $x$ in the chart one can set 
\begin{align*}
  \pT^*_x \M = \big(\nxsC\big)^\perp
\end{align*}
and for $\xi \in T_x^* \M$ one can define $\ppi$ and 
$\pxi$ as in \eqref{eq: def pxi}.  For $\y =(t,x,\tau,\xi) \in \TL$ one defines $\py =(t,x,\tau,\pxi) $
and the definition of $\pT^*_{(t,x)} \L$ follows 
similarly as above. 
Then Set 
\begin{align*}
  \pTM = \bigcup_{x \in \chdiff(\hO)} \{x\} \times \pT^*_x \M, 
   \ \ 
  \pTL = \bigcup_{(t,x) \in I \times \chdiff(\hO)} \{(t,x)\} \times \pT^*_{(t,x)} \L.
\end{align*}
One has $\pTM\cap \{z=0\} = \pdTM$ and $\pTL\cap \{z=0\} = \pdTL$. 
In the local chart, one can then define
\begin{align*}
 \pE = \pTL \cap \{ p > 0\}, \ \ 
\pG = \pTL \cap \{ p = 0\}, \ \ \pH = \pTL \cap \{ p < 0\},
\end{align*}
thus extending the elliptic, glancing and hyperbolic regions away from
the boundary. The characterization result of Lemma~\ref{lemma:
  relevement pG pH} extends  {\em mutatis mutandis}. 
If one sets
\begin{align*}
    \G = \{ \y \in \TL; \ p(\y)=0, \ \text{and}\  \Hpz(\y)=0\},
  \end{align*}
  and $\H = \H^+ \cup \H^-$ with 
  \begin{align*}
    \H^\pm = \{ \y \in \TL; \ p(\y)=0, \ \text{and}\  \Hpz(\y) \gtrless 0\}.
  \end{align*}
One has
\begin{align*}
  \Char p \cap \TL = \G \cup \H = \G \cup \H^+ \cup \H^-, 
\end{align*}
and 
\begin{align*}
    \Gb = \G \cap \dTL,  \ \ \Hb =  \H \cap \dTL,
   \ \  \text{and} \ \ \Hb^\pm = \H^\pm \cap \dTL.
  \end{align*}
If $\y \in  \pHb$ then $\ppi^{-1}(\{ \y\}) \cap \Char p=\{ \y^-,
\y^+\}$ and 
formulae~\eqref{eq: relevement pG pH} extend: with $\y=
(t,x,\tau,\xi)$ one has $\y^\pm = (t,x,\tau,\xi^\pm)$ with 
\begin{align}
  \label{eq: relevement pG pH-extended}
  \xi^+ = \xi + \lambda  \nxs
 \ \ \text{and} \ \ 
  \xi^- = \xi - \lambda \nxs,
\end{align} 
with $\lambda = \big(- p (\y)\big)^{1/2} = \big(\tau^2 - |\xi|_x^2 \big)^{1/2}$.

 The result of Proposition~\ref{prop: relevement H G} also extends:
\begin{align*}
    \ppi^{-1} (\pG) \cap \Char p = \G = \pG, \qquad \ppi^{-1} (\pH) \cap \Char p = \H,
  \end{align*}
and so does \eqref{eq: decomposition point in char p}:
if $\y = (t,x,\tau,\xi) \in\Char p \cap \TL$ one has
\begin{align}
  \label{eq: decomposition point in char -ext}
  \xi = \pxi + \alpha \Hpz(\y) \nxs,
  \quad \text{with} \ \
   \alpha \Hpz (\y) = \begin{cases}
    \big( - p (\py)\big)^{1/2}>0 & \text{if}\ \y \in \H^+,\\
    0 & \text{if}\ \y \in \G,\\
    - \big( - p (\py)\big)^{1/2}<0 & \text{if}\  \y \in \H^-.
  \end{cases}
\end{align}

In \eqref{eq: parametrization cotangent bundles},
a natural parametrizations of $\dTM$ and $\dTL$ is mentionned. It extends to $\TM$
and $\TL$ with what is introduced above:
\begin{align}
  \label{eq: parametrization cotangent bundles extension}
  \begin{array}{ll}
    \TM &\!\!\to \pTM \times \R\\
  (x,\xi) &\!\!\mapsto (x,\pxi, \Hpz (x,\xi))
\end{array}
            \ \ \text{and} \ \ 
  \begin{array}{ll}
    \TL &\!\!\to \pTL \times \R\\
  \y &\!\!\mapsto (\py, \Hpz (\y)),
    \end{array}
\end{align}
With $\pTM$ given by $\{\Hpz=0\} \cap \TM$ and $\pTL$ given by
$\{\Hpz=0\} \cap \TL$,  these two maps are in fact $\Con^1$ local diffeomorphisms by the
first part of the following proposition obtained by adapting 
the proof of Proposition~\ref{prop: G submanifolds}.
%%%%%%%%%%%%%%%%%%%%%%%%
% proposition               %
%%%%%%%%%%%%%%%%%%%%%%%%
\begin{proposition}
  \label{prop: G submanifolds extension}
  The set $\pTL$ is a  $\Con^1$-submanifold of $\TL$ of codimension
  one, in a local chart.
  Away from $(\tau,\xi) =(0,0)$, the set $\G$ is a $\Con^1$-submanifold of $\TL$ of codimension two
  defined  by $\Hpz=p=0$ (and thus a $\Con^1$-submanifold of $\pTL$
  of codimension one) in a local chart.
\end{proposition}

A use we make of the parametrizations given in \eqref{eq:
  parametrization cotangent bundles extension} is through the
following lemma.
%%%%%%%%%%%%%%%%%%%%%%%%
% lemma                %
%%%%%%%%%%%%%%%%%%%%%%%%
\begin{lemma}
  \label{lemma: derivability zeta}
  Suppose $s \mapsto \y(s)  \in \TL$.
Assume that $s\mapsto \py(s)$ and $s \mapsto  \Hpz\big(\y(s)\big)$
are both differentiable at $s=s_0$. Then $s \mapsto \y(s)$ is also
  differentiable at $s=s_0$.

Assume moreover that $\y(s_0) \in \pTL$ and $\frac{d}{d s}
\Hpz\big(\y(s)\big)_{|s=s_0}=0$. Then, $\frac{d}{d s}\y(s)_{|s=s_0} = \frac{d}{d s}\py(s)_{|s=s_0}$.
\end{lemma}
\begin{proof}
  Write $\y(s) = \big( t(s),x(s), \tau(s), \xi(s)\big)$.  The first
  part is a consequence of Proposition~\ref{prop: G submanifolds
    extension}.  With \eqref{eq: decomposition point in char -ext} one
  has
  \begin{align*}
    \frac{d}{d s}\xi (s)_{|s=s_0}
    &=  \frac{d}{d s}\pxi(s)_{|s=s_0}
    + \frac{d}{d s} \Hpz \big(\y(s)\big)_{|s=s_0} \alpha(x(s_0)) n^*_{x(s_0)} \\
    &\quad + \Hpz \big(\y(s_0)\big)  \frac{d}{d s} \big(\alpha(x(s)) n^*_{x(s)} \big)_{|s=s_0}.
  \end{align*}
  If $\y(s_0) \in \pTL$ then $\Hpz \big(\y(s_0)\big) =0$, which
  yields the second result. 
\end{proof}

In the local extension, $\ppi$ maps $\TL$ into $\pTL$. If $\y \in \TL$
then, $d \ppi(\y)$, its differential at $\y$ maps linearly $T_\y \TL$
into $T_{\py} \pTL$. With a proof similar to that of Lemma~\ref{lemma:
  kernel diff ppi} one has the following result. 
 %%%%%%%%%%%%%%%%%%%%%%%%
% lemma                %
%%%%%%%%%%%%%%%%%%%%%%%%
\begin{lemma}
  \label{lemma: kernel diff ppi-bis}
  Consider $\y \in \TL$. 
  One has $\ker (d \ppi(\y)) = \Span(\Hz)(\y).$
\end{lemma}
\begin{proof}
 The only difference concerns the choice of $v \in T_\y \TL$ whereas
 one chose  $v \in T_\y \dTL$ in the proof Lemma~\ref{lemma:
  kernel diff ppi}. Thus one has $v \in
  \Span\{\d_t, \d_{x_i}, \d_\tau,\d_{\xi_j}\}$, $i=1,\dots,d$,
  $j=1,\dots, d$. 
  The remainder of the proof is unchanged. 
\end{proof}

In Section~\ref{sec: Propagation of the measure  support} below, we
denote by $v \in \R^{2d+1}$ local coordinates for $\para{\TL}$ and use
the variable $\vartheta$ to denote the value of $\Hpz$. Then, according to  
\eqref{eq: parametrization cotangent bundles extension}, $(v,\vartheta)$ gives
$\Con^1$-coordinates on $\TL$. One has 
\begin{align*}
  \d_{\zeta} = \sum_{j=1}^{2 d+1} (\d_{\zeta} v_j) \d_{v_j} 
  + (\d_{\zeta} \vartheta) \d_\vartheta.
\end{align*}
Note that \eqref{eq: d xid Hpz} and \eqref{eq: pxi indep xid} also
hold here.
One thus has $\d_{\zeta}
v_j=0$, yielding 
\begin{align}
  \label{eq: change of variables xid zeta}
  \d_{\zeta} = - \Hz = \Hzzp\, \d_\vartheta.
\end{align}

%%%%%%%%%%%%%%%%%
% sub-section
%%%%%%%%%%%%%%%%%
\subsection{Partition of  the glancing region $\Gb$ and gliding vector field}
\label{sec: description glancing region}
By Definition~\ref{def: G H} and Proposition~\ref{prop: G submanifolds}, the glancing region $\Gb$ is the
$\Con^1$-submanifold of $\TL$ locally given by $z=\Hpz=p=0$.
By \eqref{eq: formula Hp z} $\Hpz$ is a $\Con^1$-function. As in
Section~\ref{sec: glancing region, generalized bichar-intro} it makes
sense to compute $\Hppz$ yielding a $\Con^0$-function and we recall
the partition of $\Gb = \sdGb \cup \glGb  \cup \sgGb$ made in Definition~\ref{def: diffrative}:
 \begin{itemize}
   \item $\sdGb = \{ \y \in \Gb; \Hppz (\y) >0\}$,  the diffractive
    set;
    \item 
      $\glGb =\{ \y \in \Gb; \Hppz (\y) =0\}$, the glancing sets of
      order three, meaning that a bicharacteristics that goes through
      a point of  $\glGb$ has
    a contact with the boundary of order at least equal to three;
   \item    $\sgGb =\{ \y \in \Gb; \Hppz (\y) <0\}$, the gliding set.
      \end{itemize}

As $\Hp^2z$ is a $\Con^0$-function, the set $\glGb$
defined by  $z=\Hpz=\Hp^2z=p=0$ is a  $\Con^0$-submanifolds of $\TL$.

\begin{figure}
  \begin{center}
    \resizebox{3.5cm}{!}{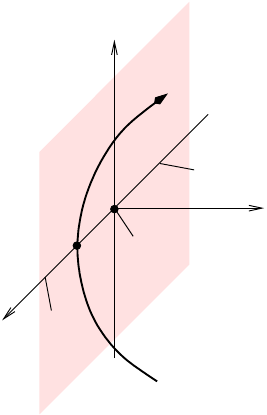}
    \caption{A \bichar going through a diffractive point ($\sdGb$).}
  \label{fig: bichar G+}
  \end{center}
\end{figure}

%%%%%%%%%%%%%%%%%%%%%%%%
% lemma                %
%%%%%%%%%%%%%%%%%%%%%%%%
\begin{lemma}
  \label{lem: bichar diffractive points}
  Consider  $\y \in \sdGb$. Denote by
  $\gamma(s) = (t(s), x(s), \tau(s), \xi(s))$
  a \bichar above $\hL$ with $x(s) = (x'(s), z(s))$ such that $\gamma(0) = \y$.   Then, 
  $\frac{d}{d s}z_{|s=0}=0$ and $\frac{d^2}{d s^2}z_{|s=0} >0$,
  meaning that for some $S>0$,  $\gamma(s)\in \TL \setminus \dTL$ for $s\in
]-S, S[\setminus \{0\}$.  
\end{lemma}
Here, $\gamma(s)$ is not a \bichar in the sense of Definition~\ref{def:
  bichar} as it encounters a point in $\dTL$ at $s=0$. It is understood as in Remark~\ref{rem:def-bichar}-\eqref{rem: bichar in hL}.
The behavior of a \bichar above $\hL$ going through a point of $\sdGb$ is
illustrated in Figure~\ref{fig: bichar G+}.
%%%% proof of lemma
\begin{proof}
  The result follows as $\frac{d}{d s }z_{|s=0} = \Hp
  z(\y) =0$ and 
  $\frac{d^2}{d s^2}z_{|s=0} = \Hppz(\y)>0$.
\end{proof}

%%%%%%%%%%%%%%%%%%%%%%%
% figure
%%%%%%%%%%%%%%%%%%%%%%%

\bigskip The map $\ppi$ maps $\dTL$ onto $\pdTL$ and acts as a
projection onto $\pdTM$ in the cotangent $\xi$ variable; see the
beginning of Section~\ref{sec: A partition of the cotangent bundle at
  the boundary}. For $\y\in \dTL$ the differential $d \ppi (\y)$ maps
$T_\y \dTL$ into $T_{\py} \pdTL$.  For $\y \in \pdTL$ one has $\y =
\ppi(\y)$ and $\Hp (\y) \in T_\y \dTL$ since $\Hpz(\y)=0$. 
One may thus set
\begin{align}
  \label{eq: definition HpG}
  \HpG (\y) = d \ppi (\y) \big( \Hp(\y)\big) \in T_\y \pdTL.
\end{align}
%%%%%%%%%%%%%%%%%%%%%%%%
% lemma                %
%%%%%%%%%%%%%%%%%%%%%%%%
\begin{lemma}
  \label{lemma: properties HpG}
  Consider $\y \in \pdTL$. 
  One has 
  \begin{align}
    \label{eq: definition HpG-bis}
    \HpG (\y)  = \Big(\Hp + \frac{\Hppz}{\Hzzp}\Hz\Big)  (\y)  .
  \end{align}
\end{lemma}
As already seen above $\Hzzp$ does not vanish; this makes
formula~\eqref{eq: definition HpG-bis} sensible.
One calls $\HpG$  the {\em gliding} vector field, here defined above
$\pdTL$. Below we extend its definition at every point of $\TL$; yet
this definition is only local and depends on the considered
coordinates. 

This definition of $\HpG$ is more satisfactory than that given in
Section~\ref{sec: glancing region, generalized bichar-intro}. There no detail was given; we only
intended to be able to state our main result, that is,
Theorem~\ref{theorem: measure support propagation}.
%%%% proof of lemma
\begin{proof}
  Since $\ker( d \ppi(\y)) = \Span(\Hz)$ by Lemma~\ref{lemma: kernel
    diff ppi} one has $\Hp^\G  (\y) = \Hp  (\y)  +
  \lambda \Hz  (\y) $ for some $\lambda\in \R$. Since $\HpG  (\y)
  \in T_\y \pdTL$  and $\pdTL$ is defined by $z=\Hpz=0$ one has
  \begin{align*}
    0 = \HpG (\Hpz) (\y) = \big(\Hppz + \lambda \Hz \Hpz) (\y)
    =\big(\Hppz - \lambda \Hzzp \big) (\y).
  \end{align*}
  Since $\Hzzp \neq 0$ one finds that $\lambda = (\Hppz / \Hzzp)(\y)$, hence the given formula for $\Hp^\G(\y)$.
\end{proof}
%%%%%%%%%%%%%%%%%%%%%%%%
% lemma                %
%%%%%%%%%%%%%%%%%%%%%%%%
\begin{lemma}
  \label{lemma: properties HpG-bis}
  If $\y \in \Gb$ then $\HpG(\y) \in T_\y \Gb$.
\end{lemma}
Note that the tangent space $T_\y \Gb$ makes sense since $\Gb$ is
$\Con^1$-manifold by Proposition~\ref{prop: G submanifolds}.

For $\y \in \Gb$, observe that $\Hp(\y)$ and $\HpG(\y)$ coincide if
and only if $\y \in \glGb$, that is $\Hppz=0$.

 %%%% proof of lemma
\begin{proof}
  Since $\y \in \pdTL$ one has $\HpG(\y) \in T_\y \pdTL$ by the definition
  of $\HpG$ in \eqref{eq: definition HpG}. With $\pdTL$ given by
  $z=\Hpz=0$ on has $\HpG z (\y) = \HpG \Hpz (\y) =0$.

  On $\pdTL$, with \eqref{eq: definition HpG-bis} one computes
  \begin{align*}
    \HpG p
    = \frac{\Hppz}{\Hzzp} \Hzp 
    = -  \frac{\Hppz}{\Hzzp} \Hpz=0.
  \end{align*}
  Hence, the result since $\Gb$ is given by $z=\Hpz=p=0$.
\end{proof}
A fairly important remark is the following one. 
\begin{remark}
  \label{rem: non vanishing hamiltonian vector field}
  Observe that $\Hp(\y) \neq0$ if $\tau \neq 0$, where as above
  $\y =(t,x,\tau,\xi) \in \TL$. In
  fact $\Hp(\y) t = - 2 \tau$. 
  If  $\y \in \Char p$, one has
  \begin{align*}
    \tau \neq 0 \ \  \Equiv \ \  
    \xi\neq 0 \ \ \Equiv \ \  
    (\tau,\xi)  \neq (0,0). 
  \end{align*}
  Hence, if $\y \in \Char p$ is such that $ (\tau,\xi) \neq  (0,0)$ then 
  $\Hp(\y) \neq 0$. Considering the form of $\HpG$ given above one also
  has  $\HpG(\y) t = - 2 \tau$. One thus finds  $\HpG(\y) \neq0$
  if $\y \in \Char p$ is such that $(\tau,\xi) \neq  (0,0)$.
\end{remark}

\bigskip
In the same framework as in Section~\ref{sec: local extension away
  from boundary},  in a local chart, we extend the definition of
$\HpG$ away from $\pdTL$ by setting
\begin{align}
  \label{eq: definition HpG-extended}
   \HpG = \Hp + \Big( \frac{\Hppz}{\Hzzp} - \frac{\Hp
  \Hzzp}{(\Hzzp)^2} \Hpz \Big)\Hz.
\end{align}
On $\pdTL$ where 
$\Hpz=0$  formula \eqref{eq: definition HpG-extended}  coincides
with \eqref{eq: definition HpG-bis}.  
Observe that 
\begin{align}
  \label{eq: HpG z = Hp z}
   \HpG z = \Hpz.
\end{align}

The reason for formula
\eqref{eq: definition HpG-extended} is as follows. 
Consider a \bichar $\gamma(s) \in \Char p \cap \TL$ and set
$\pgamma(s) = \ppi (\gamma(s)) \in \pTL$.  
If $\gamma(s) = (t(s),x(s),\tau, \xi(s))$ then
$\pgamma(s) = (t(s),x(s),\tau, \pxi(s))$
and one has
\begin{align*}
  \frac{d}{d s}\pgamma(s)
  = d \ppi(\gamma(s)) \big(\frac{d}{d s} \gamma(s)\big)
  = d \ppi(\gamma(s)) \big(\Hp( \gamma(s))\big)
 \in T_{\pgamma(s)} \para{\TL}.
\end{align*}
If $\gamma(s) \in \Gb \subset \pdTL$ then $\pgamma(s)= \gamma(s)$
and $d \ppi(\gamma(s)) \big(\Hp( \gamma(s))\big) =
\HpG\big(\gamma(s)\big) \in T_{\gamma(s)} \pdTL$ by the definition
of $\HpG$ in $\pdTL$ introduced in \eqref{eq: definition HpG}. 

More generally, for  $\gamma(s)\in \pTL$ one has $\gamma(s) =
\pgamma(s)$ and thus one has  $d \ppi(\gamma(s)) \big(\Hp( \gamma(s))\big) \in
T_{\gamma(s)} \pTL$. The proof of Lemma~\ref{lemma: properties HpG}
applies, with Lemma~\ref{lemma: kernel diff
  ppi} replaced by Lemma~\ref{lemma: kernel diff
  ppi-bis}, and yields 
\begin{align*}
  d \ppi(\gamma(s)) \big(\Hp( \gamma(s))\big) \big( \gamma(s)\big) 
  = \Big(\Hp + \frac{\Hppz}{\Hzzp}\Hz\Big)  \big( \gamma(s)\big) , 
  \end{align*}
which coincides with \eqref{eq: definition HpG-extended} since
$\Hpz=0$ on $\pTL$. 

Yet, if $\gamma(s)\in \H$, with $\H$ locally defined as
in Section~\ref{sec: local extension away from boundary}, then $\gamma(s) \neq
\pgamma(s)$ implying that $d \ppi(\gamma(s)) \big(\Hp(
\gamma(s))\big) \notin T_{\gamma(s)} \TL$. 
Local coordinates can be of some help however. One has $\pxi(s) =
\xi(s) - (\alpha \Hpz) \big(\gamma(s)\big) \nxs$ yielding 
\begin{align*}
  \pgamma(s) = \gamma(s) -  (\alpha \Hpz) \big(\gamma(s)\big) \nxs,
\end{align*}
here identifying $\nxs$ with $(0,0,0,\nxs)$ in the $(t,x,\tau,\xi)$
variables. Recalling that the $\zeta$-componant of $(\alpha \Hpz)
\nxs$ is $\Hpz/\Hzzp$ while other componants are zero one obtains
\begin{align*}
  \frac{d}{d s}\pgamma(s) 
  &= \frac{d}{d s}\gamma(s) 
    - \frac{d}{d s}\big((\Hpz/\Hzzp) \big(\gamma(s)\big)\big) \d_{\zeta}\\
  &= \Hp\big(\gamma(s)\big)
    + \frac{ \Hp\Hpz}{\Hzzp}  \big(\gamma(s)\big) \Hz
    - \frac{ \Hp \Hzzp}{(\Hzzp)^2}  \Hpz\big(\gamma(s)\big) \Hz\\
  &= \Hp\big(\gamma(s)\big)
    + \Big( \frac{ \Hppz}{\Hzzp} 
    - \frac{ \Hp \Hzzp}{(\Hzzp)^2}  \Hpz\Big)\big(\gamma(s)\big) \Hz.
\end{align*}
In the local coordinates considered here it is thus not difficult to identify a tangent vector at
$\gamma(s)$ with a tangent vector at $\pgamma(s)$. With this
identification one may
write 
\begin{align}
 \label{eq: dynamique pgamma}
  \frac{d}{d s}\pgamma(s) 
  &=\HpG\big(\gamma(s)\big),
\end{align}
with the understanding that $\HpG\big(\gamma(s)\big)$ is computed at
$\gamma(s)$ 
according to \eqref{eq: definition HpG-extended} and viewed as a tangent vector at $\pgamma(s)$. We 
use this notation in what follows, with the necessary care since
$\HpG\big(\gamma(s)\big)$ may not be equal to
$\HpG\big(\pgamma(s)\big)$. In fact the following result holds. 
%%%%%%%%%%%%%%%%%%%%%%%%
% lemma                %
%%%%%%%%%%%%%%%%%%%%%%%%
\begin{lemma}
  \label{lemma: HpG rho neq HPG prho}
  One has $\HpG(\y) =\HpG\big(\py\big)$ if and only if $\py = \y$,
  that is, $\y\in\pTL$. 
\end{lemma}
%%%% proof of lemma
\begin{proof} Suppose $\HpG(\y) =\HpG\big(\py\big)$. With
  \eqref{eq: HpG z = Hp z} one has $\Hpz (\y) = \HpG z(\y)
  =\HpG z\big(\py\big) = \Hpz\big(\py\big) =0$. Hence the conclusion.
\end{proof}

 If one is
not inclined to perform this abuse of notation one should stick with 
\begin{align}
  \label{eq: dynamique pgamma-rigorous}
  \frac{d}{d s}\pgamma(s) 
  &=d \ppi\big(\gamma(s)\big) \Big(\!\Hp\big(\gamma(s)\big)\!\Big).
\end{align}
Note that neither \eqref{eq: dynamique pgamma} nor \eqref{eq: dynamique pgamma-rigorous} is to be viewed as a differential equation for
$\pgamma(s)$, since the right-hand-side is a function of $\gamma(s)$
and not $\pgamma(s)$.

%%%%%%%%%%%%%%%%%%%%%%%%
% remark               %
%%%%%%%%%%%%%%%%%%%%%%%%
\begin{remark}
  \label{remark: non vanishing HpG}
  Observe that Remark~\ref{rem: non vanishing hamiltonian vector
    field} applies to the local extension of $\HpG$ away from the
  boundary: if $\y \in \Char p \cap \TL$ is such that
  $ (\tau,\xi) \neq (0,0)$ then $\HpG(\y) \neq 0$.
\end{remark}

%%%%%%%%%%%%%%%%%
% sub-section
%%%%%%%%%%%%%%%%%
\subsection{The compressed cotangent bundle}
\label{sec: The compressed cotangent bundle}

The symmetry $\Sigma$ introduced in Section~\ref{sec: A partition of
  the cotangent bundle at the boundary} acts as an involution on
$\dTL$ that leaves $\pdTL$ invariant.
For $(t,x) \in \d\L$, one sets
\begin{align*}
  \cT_{(t,x)}^* \L = \quot{T_{(t,x)}^* \L}{\Sigma}, 
\end{align*}
called the {\em compressed} cotangent space at $(t,x)$, and 
\begin{align*}
 \cTL
  = \bigcup_{(t,x) \in \d\L} \{(t,x)\} \times \cT_{(t,x)}^*\L 
  \ \ \cup \ 
  \bigcup_{(t,x) \in \L\setminus \d\L} \{(t,x)\} \times T_{(t,x)}^*\L ,
\end{align*}
called the  {\em compressed} cotangent bundle.
The quotient with respect to $\Sigma$ is thus only performed above the
boudary $\d\L$, that is $\{ z=0\}$ in local coordinates.
This allows one to identify  a point $\y \in \Hb^\pm$ with
$\Sigma(\y) \in \Hb^\mp$. This turns out to be usefull in the
construction of \gbichars in what follows, allowing such
\bichars to be continuous across hyperbolic points.

Denote by $\cphi$ the associated quotient map $\TL \to
\cTL$. One has $\cphi(\Hb^+) = \cphi(\Hb^-)$.  It acts as the
identity on $\TL \setminus \dTL$ and on $\pdTL$. For $\y
\in \TL \setminus \Hb$ one thus writes $\cphi(\y) = \y$ by
abuse of notation.
Set 
\begin{align*}
  \cHb = \cphi(\Hb) = \cphi(\Hb^+) = \cphi(\Hb^-). 
\end{align*}

One can endow $\cTL$ with a natural metric inherited from that of
$\TL$.  In fact, given two points
$\cy^0, \cy^1\in \cTL$ 
consider a path $\gamma(s) = [0,1] \setminus B\to \TL$,
with $B = \{ s_1, \dots, s_k\} \subset [0,1]$, with $s_1 < s_2 < \cdots
< s_k$, such that 
the limits 
\begin{align*}
  %\label{eq: path distance compress cotangent bundle 1}
  \gamma(s_n^-) = \lim_{s \to s_n^-} \gamma (s) \ \ \text{and} \ \ 
  \gamma(s_n^+) = \lim_{s \to s_n^+} \gamma (s),
  \qquad n =1, \dots, k,
\end{align*}
exist and $\gamma (s_n^-),  \gamma (s_n^+) \in \dTL \setminus \pdTL$ with moreover 
\begin{align*}
 % \label{eq: path distance compress cotangent bundle 2}
  \Sigma(\gamma (s_n^-)) = \gamma (s_n^+),
\end{align*}
and 
\begin{align*}
  %\label{eq: path distance compress cotangent bundle 3}
  \cphi(\gamma (0)) = \cy^0, \ \  \cphi(\gamma (1)) = \cy^1.
\end{align*}
One sets 
\begin{align*}
  \length(\gamma ) = \ell_{[0,s_1]} \gamma
  + \sum_{n=1}^{k-1} \ell_{[s_n,s_{n+1}]} \gamma 
  +\ell_{[s_k,1]} \gamma,
\end{align*}
where $\ell_{[a,b]} \gamma$ stands for the length of $\gamma$ for
$s\in [a,b]$. One then defines 
\begin{align}
  \label{eq: distance compressed cotangent bundle}
  \cdist ( \cy^0, \cy^1) = \inf \length(\gamma ),
\end{align}
where the infimum is computed over all paths fulfilling the above conditions.

The quotient map is continuous: for some $C>0$ one has 
\begin{align}
  \label{eq: cont quotient map}
  \cdist  \big( \cphi(\y^0), \cphi(\y^1) \big) \leq C  \Norm{\y^0 -
  \y^1}{}, \quad \y^0, \y^1\in
  \TL.
\end{align}
% This can be seen in a local chart using the variables $(u,z,\vartheta)$ with
% $(u,z) = \py$ and $\vartheta= \Hpz$. 

With the metric structure now given on $\cTL$ note that the path
$\cgamma: [0,1] \to \cTL$ given by $\cgamma (s)= \cphi(\gamma(s))$
is continuous if $\gamma$ is as described previously.

The following lemma follows from what is above. 
%%%%%%%%%%%%%%%%%%%%%%%%
% lemma                %
%%%%%%%%%%%%%%%%%%%%%%%%
\begin{lemma}
  \label{lemma: relevement chemin fibre compresse}
Suppose $J$ is an interval, $\cgamma: J \to \cTL$ is continuous, and 
$B \subset J$ is a discrete set, that is, made of isolated points, such that $\cgamma(s) \in\cHb$ if and only if $s\in B$. Then, there exists a unique map $\gamma: J\setminus B \to
\TL$ such that $\cphi(\gamma(s)) = \cgamma(s)$ if $s \in
J\setminus B$ and $\gamma$ is continuous away from points in
$B$. Moreover,  for $S \in B$ the limits 
\begin{align*}
  %\label{eq: path distance compress cotangent bundle 1}
  \gamma(S^-) = \lim_{s \to S^-} \gamma (s) \in \Hb^- \ \ \text{and} \ \ 
  \gamma(S^+) = \lim_{s \to S^+} \gamma (s) \in \Hb^+
\end{align*}
exist and $\Sigma (\gamma(S^-) ) = \gamma(S^+)$.

In addition, if $\tilde{J} \subset J$ is an interval such that
$\gamma(s)$ lies in one local chart for $s \in \tilde{J}$, then $s
\mapsto \pgamma(s) = \ppi(\gamma)(s)$ can be defined for $s \in \tilde{J} \setminus
B$ (see Section~\ref{sec: local extension away from boundary}) and
moreover extended to the whole
interval $\tilde{J}$ in a continuous manner.
\end{lemma}

% For a measure $\mu$ defined on $\TL$ 
% denote by $\csupp(\mu)$ the image by $\cphi$ of $\supp \mu$. 
% Denote also by $\cChar(p)$ the image by $\cphi$ of $\Char p$. 
% Note that both are closed subsets of $\cTL$. 

%%%%%%%%%%%%%%%%%
% sub-section
%%%%%%%%%%%%%%%%%
\subsection{Broken and generalized bicharacteristics}
\label{sec: Broken and generalized bicharacteristics}
Consider $\gamma: J \to \TL \setminus \dTL$ a maximal \bichar,
with $J= ]S_0,S_1[$ (see Definition~\ref{def: bichar}).
Set $\gamma(s) = (t(s), x(s), \tau(s),
\xi(s))$. Recall that $\tau(s) = |\xi(s)| = \cst$ and consider the only
interesting case of a nonzero value for $\tau(s)$.  Concerning the potential limit at $s = S_0^+$ and
$s = S_1^-$ one has the following result.
%%%%%%%%%%%%%%%%%%%%%%%%
% lemma                %
%%%%%%%%%%%%%%%%%%%%%%%%
\begin{lemma}
  \label{lemma: limit of a maximal bichar from the interior}
  Assume $S_0^+ < +\infty$ (\resp $S_0^- > - \infty$).
  The limit of $\gamma(s)$ as $s \to S_0^+$ (\resp $s \to S_1^-$) exists and
  \begin{align*}
    \lim_{s \to S_0^+} \gamma(s)  \in  \Hb^+ \cup \glGb \cup \sdGb \qquad
    \text{(\resp}\ \lim_{s \to S_1^-} \gamma(s)  
    \in \Hb^- \cup \glGb \cup \sdGb \text{)}.
  \end{align*} 
\end{lemma}
In other words, the limit lies in $\Char p \cap \dTL$ yet away
from $\sgGb$.  
%%%% proof of lemma
\begin{proof}
  Along a \bichar the value of  $\tau$ is constant since $\Hp \tau =0$. As the
  \bichar lies in $\Char p$ where $\norm{\xi}{x}^2 = \tau^2$, one finds
  that $\xi$ remains bounded. Thus $\gamma(s)$ lies in a compact set.
  Since $\Hp$ is continuous and thus bounded there, one finds that
  $\gamma(s)$ has a single accumulation point as $s \to S_1^-$.
  Hence, $\y^1 = \lim_{s \to S_1^-} \gamma(s)$ exists and belongs to
  $\dTL \cap \Char p$. One sets $\gamma(S_1)= \y^1$. Assume that
  $\y^1 \in \sgGb$. In local coordinates as in \eqref{eq: local chart
    boundary} one has
  $\gamma(s) = \big(t(s),x'(s), z(s), \tau(s),\xi'(s),\zeta(s)\big)$ for
  $s \in [S_1-\eps,S_1]$ for some $\eps>0$. Naturally, one has
  $z(s) \geq 0$. Moreover, $z(S_1) = \Hpz(\y^1)=0$ and $\Hppz(\y^1)
  <0$. Since $\Hpz (\gamma(s)) = \frac{d}{d s} z(s)$ and
  $\Hppz\big(\gamma(s)\big) = \frac{d^2}{d s^2} z(s)$ along a
  \bichar one concludes that $z(s)<0$ for $s \in [S_1-\eps',S_1]$ for
  some $\eps'>0$, a contradiction.
\end{proof}

The following property is elementary yet important in the generalization of
the notion of \bichars. It states the existence  of some interval of uniform minimal size for all \bichars originating from a \nhd of a hyperbolic point. 
%%%%%%%%%%%%%%%%%%%%%%%%
% lemma                %
%%%%%%%%%%%%%%%%%%%%%%%%
\begin{lemma}
  \label{lemma: limit of a bichar in H}
  Consider $\y^0 \in \Hb^-$ (\resp $\Hb^+$). There exists $S_0>0$ and an
  open \nhd $V^0$ of $\y^0$ in $\Hb^-$ (\resp $\Hb^+$) such that for
  any maximal \bichar $\gamma(s)$ in $\TL \setminus \dTL$
  defined on $]s_1,s_2[$, with $s_1< s_2$ and
  $\lim_{s \to s_2^-} \gamma(s) \in V^0$ (\resp
  $\lim_{s \to s_1^+} \gamma(s) \in V^0$), one has $s_2 - s_1 \geq S_0$.
%Moreover, the value of $S_0$ depends continuously on $\y^0$.
\end{lemma}
%This property is quite meaningful in the present context of nonuniqueness
%of maximal \bichars.
%%%% proof of lemma
\begin{proof}
  In the proof we only treat the case $\y^0 \in \Hb^-$. The case $\y^0
  \in \Hb^+$ can be treated similarly. Near $\y^0$ use local
  coordinates as in \eqref{eq: local chart boundary}.
  
  One has $\Hpz(\y^0) <0$. Thus, there exists a $2(d+1)$-dimensional
  ball $\mathcal B_2$ of radius $2 R$ centered at $\y^0$ such that
  $\Hpz <- C_0$  and $\Norm{\Hp}{}\leq C_1$ in
  $V=\mathcal B_2 \cap \{ z\geq 0\}$ for some $C_0>0$ and $C_1>0$.
  Set $V^0= \mathcal B_1 \cap \Hb \subset \Hb^-$ with $\mathcal B_1$
  the ball of radius $R$ centered at $\y^0$.

Suppose $\gamma$ is a maximal \bichar defined on $]s_1,s_2[$ as in the
  statement, with $\y^1 = \lim_{s \to s_2^-} \gamma(s) \in V^0$.  Set
  $\gamma(s_2) = \y^1$ by continuous extension.

  By continuity $\gamma(s)\in V$ for $s \in [s_2-\delta,s_2]$ for some
  $\delta>0$. 
  If $\gamma(s) \in V$ one has
  \begin{align*}
    %\label{eq: broken bichar near H-}
    \frac{d}{d s} z(s)= \Hpz(\gamma(s))
    \leq -C_0.
  \end{align*}
  Thus, $z(s)>0$ for $s<s_2$ and $\gamma([s,s_2]) \subset  V$ implying, first, that the
  {\em maximal} \bichar cannot cease to exist if $\gamma(s)$ remains in
  $V \setminus \d\L$, and second, that $\gamma(s)$ can only exit $V$
  through $\d \mathcal B_2 \cap \{ z>0\}$.  Set
  $S^- = \inf \{ s_1\leq S <s_2; \ \gamma([S,s_2])\subset V\}$.  One has 
  $s_2-s_1 \geq s_2-S^- \geq R /C_1$.
\end{proof}

\bigskip
A first generalization of \bichars is that of \bbichars. It is
composed of a sequence (finite or not) of \bichars that are connected
at hyperbolic points at the boundary and follow the optics
law of reflection\footnote{The optics
law of reflection, often called Descartes' law or
Snell's law, goes back in fact to the works of Euclid (ca. 300 BC) and
Hero of Alexandria (ca. 60 AD).} at such points.
%%%%%%%%%%%%%%%%%%%%%%%%
% definition           %
%%%%%%%%%%%%%%%%%%%%%%%%
\begin{definition}[broken \bichar]
	\label{def: broken bichar}
  Suppose $J \subset \R$ is an interval, $B$ a discrete subset of $J$, and
  \begin{align*}
    \gammaB: J\setminus B 
    \to \Char p  \cap \TL.
  \end{align*}
  One says that $\gammaB$ is a  \bbichar if  the
  following properties hold:
  \begin{enumerate}
  \item for $s \in J \setminus B$, $\gammaB(s) \notin \Hb \cup \glGb
    \cup \sgGb$ and   the map $\gammaB$ is differentiable
    at $s$ with
	\begin{equation*}
		\frac{d}{d s} \gammaB(s) = \Hp\big( \gammaB(s)\big).
	\end{equation*}
  \item  If $S\in B$, then $\gammaB (s) \in \TL \setminus \dTL$ for $s
    \in J\setminus B$ \suff
    close to $S$ and moreover
    \begin{enumerate}
      \item  if $[S-\eps,S] \subset J$ for some $\eps>0$, then $\gammaB(S^-) = \lim_{s \to S^-} \gammaB(s) \in
      \Hb^-$;
       \item if $[S,S+\eps] \subset J$ for some $\eps>0$, then $\gammaB(S^+) = \lim_{s \to S^+} \gammaB(s) \in
      \Hb^+$;
    \item and if $[S-\eps,S+\eps]\subset J$ for some $\eps>0$, then
    $\gammaB(S^+)=\Sigma \big(\gammaB(S^-)\big)$.
    \end{enumerate}   
\end{enumerate}   
\end{definition}
%%%%%%%%%%%%%%%%%%%%%%%
% figure
%%%%%%%%%%%%%%%%%%%%%%%
\begin{figure}
  \begin{center}
    \resizebox{10cm}{!}{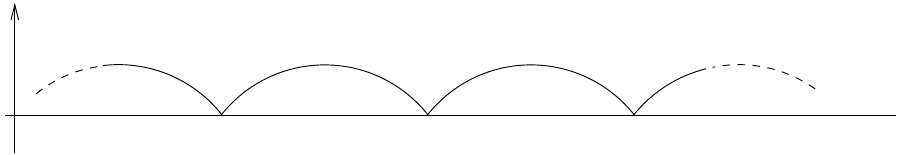}
    \caption{A \bbichar in a \nhd of part of the boundary.}
  \label{fig: broken bichar} 
  \end{center}
\end{figure}
Figure~\ref{fig: broken bichar} sketches what a \bbichar may look
like.

\begin{remark}
  \label{rem: def broken bichar}

  \leavevmode  
  \begin{enumerate}
  \item 
    \label{rem: broken bichar pgammaB C0}
  Near a point $S \in J \setminus B$, $s \mapsto \gammaB(s)$ is
  $\Con^1$.  At a point $S \in B$, $s \mapsto \gammaB(s)$ is
  discontinuous. However, if one sets $\py(s) = \ppi
  \big(\gammaB(s)\big)$ in local coordinates  (see Section~\ref{sec: local extension away from boundary}) one
  sees that $\py(s)$ is $\Con^0$. In particular the $z$-component is
  continuous and vanishes at $s=S^\pm$ for $S \in B$. 

  If one considers the map $\cphi$ introduced in Section~\ref{sec: The
    compressed cotangent bundle} one sees that $s\mapsto
  \cphi(\gammaB(s))$ takes values in the compressed cotangent bundle
  and can be extended to the whole interval $J$ as a continuous
  function. At this stage, the map $s\mapsto
  \cphi(\gammaB(s))$ is however not needed for the understanding of the
  behavior of \bbichars.

  \item Similarly to what is observed in Remark~\ref{rem:def-bichar}-\eqref{rem: constant tau} one has, for $s \in J \setminus B$,
    \begin{align*}
    \norm{\xi(s)}{x(s)} = |\tau(s)|
  \end{align*}
   constant along a \bbichar.
  \item Observe that one allows a \bbichar to reach points in $\sdGb$;
there, the tangent vector $\frac{d}{d s} \gammaB(s)$ is also given by
$\Hp\big(\gammaB(s)\big)$.  

\item \label{rem: isolated points Broken Bichar}Points of
$B$ are naturally isolated because of Lemma~\ref{lemma: limit of a
  bichar in H}. In fact, points of $B$ can
only accumulate at the boundary of $J$ as stated in the following
lemma.
\end{enumerate}
\end{remark}

%%%%%%%%%%%%%%%%%%%%%%%%
% lemma                %
%%%%%%%%%%%%%%%%%%%%%%%%
\begin{lemma}
  \label{lemma: broken bichar - closure B}
  Suppose $\gammaB: J\setminus B \to \TL$
  is a \bbichar. One has 
  \begin{enumerate}
  \item $\ovl{B} \setminus B \subset \d J \setminus J $;
  \item if $S \in \ovl{B} \setminus B$ then the limit
    \begin{align*}
      \lim_{{s \to S}\atop s \in J \setminus B}\gammaB(s) \ \ \text{exists}
    \end{align*}
    and belongs to $\glGb \cup \sgGb$.
    \item If $S \in \d J \setminus \ovl{B}$ then the limit 
    \begin{align*}
      \lim_{{s \to S}\atop s \in J \setminus B}\gammaB(s) \ \ \text{exists}
    \end{align*}
  and does not belong to $\sgGb$.
  \end{enumerate}
\end{lemma}
%%%% proof of lemma
\begin{proof}
  Consider $S \in \ovl{B} \setminus B$. Since $d x / d s$ is bounded, $x(s)$ has a limit as
  $s \to S$. One may thus use a single local chart for $s$ \suff close
  to $S$. This allows one to use the local extension of $\HpG$  away
  from $\dTL$ as well as the extention of $\ppi$ and $\pTL$ given in 
  Section~\ref{sec: local extension away from boundary}.
  As in the beginning of the proof
  Lemma~\ref{lemma: limit of a maximal bichar from the interior}, one
  sees that $\gammaB(s)$ lies in a compact set for $s$ near
  $S$. Thus
  $s \mapsto \gammaB(s)$ has at least one accumulation point as
  $J\setminus B  \ni s \to S$. We prove that this accumulation point is unique.

  Parameterize $\gammaB(s)$ by $\py(s) = \ppi(\gammaB(s)) = (t(s), x'(s), z(s), \tau(s), \pxi(s))$ and $\Hpz(\gammaB(s))$; see
  \eqref{eq: parametrization cotangent bundles}. By
  Remark~\ref{rem: def broken bichar}-\eqref{rem: broken bichar
    pgammaB C0} $\py(s)$ is a continuous
  function for $s \in J$ and with \eqref{eq:
    dynamique pgamma}--\eqref{eq:
    dynamique pgamma-rigorous} one has
  \begin{align*}
    \frac{d}{d s} \py(s)  = \HpG(\gammaB(s)), \quad \text{for} \ s \in
    J \setminus B.
  \end{align*}
As $ \HpG(\gammaB(s))$ is bounded as $J \setminus B  \ni s \to S$, then
  $s\mapsto \py(s)$ has a single
  accumulation point as $J \setminus B  \ni  s \to S$, denoted by $\py^{0} = (t^0,
  x^{\prime 0}, z^0, \tau^0,\pxi^0)$. By the continuity of $\py(s)$
  the same holds if one allows $J \ni s \to S$.  

  Consider $(s_n)_n \subset B$ such that $s_n \to S$. Since
  $\big(\y(s_n^\pm) \big)_n \subset \Hb$, then $\big(\py(s_n^\pm)
  \big)_n \subset \pHb$, and one finds
  $\py^0 \in \ovl{\pHb} = \pHb \cup \pGb$. In particular $z^0=0$. 
  
  As $\gammaB(s) \in \Char p$, with Lemma~\ref{lemma: relevement pG
    pH} and \eqref{eq: connection Hpz lambda} (and its extension~\eqref{eq: relevement pG pH-extended} in
  Section~\ref{sec: local extension away from boundary}), the function  $\Hp(\gammaB(s))$ is
  such that 
  \begin{align*}
    \begin{cases}
    \alpha(x(s)) \Hpz \big(\gammaB(s) \big) =  \lambda(s) \ \text{or} -  \lambda(s) 
    &\text{if} \ \ s \in J \setminus B,\\
     \alpha(x(s)) \Hpz \big (\gammaB(s^\pm) \big) = \pm   \lambda(s) 
    &\text{if} \ \  s \in B.
\end{cases}
\end{align*}
with $\lambda(s)  = \big(-p(\py(s))\big)^{1/2}$.
If $J \setminus B \ni s \to S$, since $\py(s) \to \py^0$, one finds
that $\lambda(s)$ has a limit, and thus $\Hpz
\big(\gammaB(s) \big)$ has {\em at most} two accumulation points. The
resulting potential accumulation points for $\gammaB(s)$, as $J\setminus B \ni s \to S$,  lie in $\ppi^{-1}(\{\py^0\}) \cap \Char p$
and
\begin{itemize}
  \item either there are two distinct accumulation points in $\Hb$,  image of
    one another by the symmetry map $\Sigma$,
    \item or there is a single accumulation point in $\Gb$.
\end{itemize}
(See Lemma~\ref{lemma: relevement pG pH} and Proposition~\ref{prop:
  relevement H G}.)  Our claim is that $\py^0 \in \Gb$ and thus
$\gammaB(s)$ has a single accumulation point $\y^0 = \py^0$, meaning
that $\gammaB(s)$ has a limit as $J\setminus B \ni s \to S$. In
particular $\Hp z(\y^0)=0$.

Proceed by contradiction.  Suppose there are two distinct accumulation
points $\y^{1,-} \in \Hb^-$ and $\y^{1,+}\in \Hb^+$ with
$\Sigma(\y^{1,-}) = \y^{1,+}$. Consider the sequence
$(s_n)_n$ introduced above. The sequence $\gammaB(s_n^-)$ has an accumulation
point that lies in $\ovl{\Hb^-}$; since it is then
accumulation point of  $\gammaB(s)$ this point is necessarily
$\y^{1,-}$. Hence, $\gammaB(s_n^-)$ has a unique accumulation point,
meaning it converges to $\y^{1,-}$.  
By Lemma~\ref{lemma: limit of a
    bichar in H} there exists a \nhd $V^1$ of $\y^{1,-}$ in $\Hb^-$, such
  that any maximal \bichar in $\TL$ initiated at a point of $V^1$
  exists for $s$ in an interval of minimal length $\ell^0>0$. This is
  in contradiction with having $\gammaB(s_n^-) \to \y^1$ and
  $|s_{n+1} - s_n| \to 0$.  This proves our claim.

%  is one in $\Hb$ for $\gammaB(s)$ as $J\setminus B \ni s \to S$, and thus none in $\Gb$. Consider the sequence
% $(s_n)_n$ introduced above. Then $\gammaB(s_n^-)$ has an accumulation
% point $\y^{1} \in \ovl{\Hb^-}$ and
% $\gammaB(s_n^+)=\Sigma(\gammaB(s_n^-)\big)$ has
% $\Sigma(\y^{1}) \in \ovl{\Hb^+}$ for accumulation points. Then, $\y^1$
% and $\Sigma(\y^{1})$ are accumulation points of $\gammaB(s)$ as
% $J \setminus B \ni s \to S$. One has $\y^1 \neq \Sigma(\y^{1})$ since
% otherwise $\y^1 \in \Gb$, which is exluded here. Thus, $\y^1\in \Hb^-$
% and $\Sigma(\y^{1})\in \Hb^+$ are the only two (distinct) accumulation
% points of $\gammaB(s)$  as $J\setminus B \ni s \to S$. Since $\y^1 \notin \ovl{\Hb^+}$ and $\Sigma(\y^{1}) \notin \ovl{\Hb^-}$ 
% one moreover has $\gammaB(s_n^-) \to \y^1$ and
%   $\gammaB(s_n^+) \to \Sigma(\y^1)$ (uniqueness of the accumulation
%   point for each sequence).  By Lemma~\ref{lemma: limit of a
%     bichar in H} there exists a \nhd $V^1$ of $\y^1$ in $\Hb^-$, such
%   that any maximal \bichar in $\TL$ initiated at a point of $V^1$
%   exists for $s$ in an interval of minimal length $\ell^0>0$. This is
%   in contradiction with having $\gammaB(s_n^-) \to \y^1$ and
%   $|s_{n+1} - s_n| \to 0$.  This proves our claim. 

\medskip Assume now that $\y^0 \in \sdGb$. Then, in addition to having
$\Hpz(\y^0)=0$, one has $\Hppz(\y^0) >0$. One considers an open
\nhd $W$ of $\y^0$ where $\Hppz> 0$. There exists $\delta>0$ such
that $\gammaB(s) \in W$ for
$s \in J \cap ]S-\delta, S+ \delta[ \setminus B$ and
$\gammaB(s^\pm) \in W$ for $s \in B \cap ]S-\delta, S+ \delta[$.  
Consider $r \in B \cap ]S-\delta, S[$. The following argument is
similar if one starts from $r \in B \cap ]S, S+\delta[$.  
Set $R = \inf B \cap ]r,S]$. Since points in $B$ are isolated (see
Remark~\ref{rem: def broken bichar}-\eqref{rem: isolated points Broken
  Bichar}), one has $r < R$. From the continuity of $z(s)$ one has $z(R)=0$.
Yet, let us consider the \bbichar
$\gammaB(s)$ for $s \in ]r,R[\subset J$. One has $z(r^+)=0$ and $z'(r^+)=\Hpz\big(
\gammaB(r^+)\big) >0$ since $\gammaB(r^+) \in \Hb^+$. Since $]r,R[
\subset J$ and $\gammaB(s) \in W$ for $s \in ]r,R[$ one has $z''(s) = \Hppz
\big(\gammaB(s) \big) >0$ for $s \in ]r,R[$. Thus $z(s)$ increases on $]r,R[$ and
$z(R) >0$. A contradiction with $z(R)=0$. Thus $\y^0$ cannot be in $\sdGb$.

    In conclusion, for $S \in \ovl{B}\setminus B$,  we have found
    \begin{align*}
      \y^0 =\lim_{J \setminus B \ni s\to S} \gammaB(s) \in \glGb\cup \sgGb.
    \end{align*}
    With the definition of a \bbichar that takes values in $\TL \setminus (\Hb \cup \glGb \cup  \sgGb)$ one finds that $S \notin
    J$. Thus $S \in \d J \setminus  J$.

    \bigskip Finally, consider the case of a point $S \in \d J
    \setminus \ovl{B}$. Thus there exists $\eps>0$ such that $B \cap
    ]S- \eps, S+ \eps[ = \emptyset$.  With the same argument as in the begining of
    Lemma~\ref{lemma: limit of a maximal bichar from the interior} one
    finds that 
    \begin{align*}
      \y^0 = \lim_{{s \to S}\atop s \in J \setminus B}\gammaB(s) \ \ \text{exists}.
    \end{align*}
    If $\y^0 \in \dTL$ then arguing as in the proof of  Lemma~\ref{lemma: limit
      of a maximal bichar from the interior} one finds that $\y^0 \notin
    \sgGb$. 
\end{proof}

\begin{figure}
  \begin{center}
    \resizebox{12cm}{!}{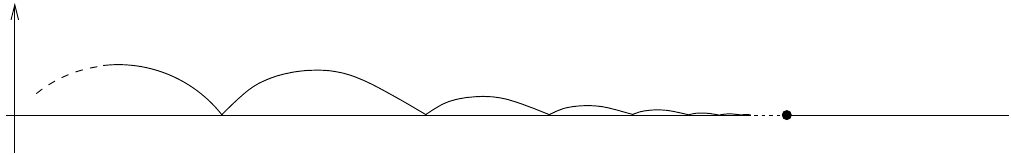}
    \caption{A \bbichar in the case $B$ has an
      accumulation point.}
    \label{fig: maximal broken bichar} 
  \end{center}
\end{figure}

Figure~\ref{fig: maximal broken bichar}  illustrates the behavior of a
\bbichar $\gammaB(s)$ as $s \to S$ with $S$ an accumulation of
$B$.   
%%%%%%%%%%%%%%%%%%%%%%%
% figure
%%%%%%%%%%%%%%%%%%%%%%%
For a \bbichar defined on $J \setminus B$ with $J = ]S_0, S[$ if
$\lim_{s \to S^-} \gammaB \in \dTL$ and $S \notin \ovl{B}$ this
means that the \bbichar is in fact a simple \bichar in a \nhd of
$S^-$ and Lemma~\ref{lemma: limit of a maximal bichar from the
  interior} applies (this is also covered by the third item of Lemma~\ref{lemma: broken bichar - closure B}). If $S \in \ovl{B}$ then by Lemma~\ref{lemma:
  broken bichar - closure
  B} 
one has $\lim_{s \to S^-} \gammaB(s)\in \glGb \cup \sgGb$. The following
lemma yields a finer understanding of the behavior of the \bbichar as it reaches such a  limit point. 
%%%%%%%%%%%%%%%%%%%%%%%%
% lemma                %
%%%%%%%%%%%%%%%%%%%%%%%%
\begin{lemma}
  \label{lemma: limit broken bichar in G-G0}
  Consider a local chart at the boundary $(\hO,\chdiff)$ as in
  \eqref{eq: local chart boundary} and $W$ a bounded open subset of
  $\TL$ that lies above $\hO$.  Suppose 
  $\y^0 \in \big(\glGb \cup \sgGb\big) \cap W$. Let also $S>0$,
  $J=]0,S[$, $B$ a discrete subset of $J$, and
  $\gammaB: J\setminus B \to \TL$ be a \bbichar. Assume that
  $lim_{s \to S^-} \gammaB(s)$ exists in $\TL$ and that this limit is
  $\y^0$.  If one sets $\gammaB(S) = \y^0$, then, $\gammaB$ is
  differentiable at $s=S^-$ and
	\begin{align*}
		\frac{d}{d s} \gammaB (S^-) = \HpG(\y^0).
   \end{align*}	
   In fact, there exists some $C>0$ uniform with respect to
   $\y^0$, such that 
	\begin{align}
	\label{eq: precise null-derivability z reach G-G0}
          0\leq z \big(\gammaB(s)\big) \leq \int_s^S
          \big|\Hpz\big( \gammaB(\sigma)\big) \big|\, d\sigma,
	\end{align}	
        and 
	\begin{align}
	\label{eq: precise null-derivability Hp z reach G-G0}
          \big|\Hpz \big(\gammaB(s)\big)\big| \leq C \int_s^S
         \big| 
          \big(\Hppz\big) \big( \gammaB(\sigma)\big)
	- \big(\Hppz\big) \big( \gammaB(S)\big)
          \big|\, d\sigma,
	\end{align}
        for $s \in J \setminus B$ such that $\gammaB(s) \in W$.
\end{lemma}

%%%% proof of lemma
\begin{proof}
Consider $S_0\in ]0,S[$ such that $\gammaB(s)\in W$ for $s \in
[S_0, S]\setminus B$. Set $\py(s) = (t(s),x'(s),z(s),
\tau(s),\pxi(s)) = \ppi\big(\gammaB(s)\big)$ and $\py^0 =
\ppi\big(\y^0\big)$ (with the local extension introduced in Section~\ref{sec: local extension away from boundary}). 
Along the \bbichar, $s\mapsto \py(s)$ is continuous and piecewise
$\Con^1$. One has {\em piecewisely} on $[S_0, S[$
\begin{align*}
  %\label{eq: tangent Hp' broken bichar}
  \frac{d}{d s} \py(s)
  = \HpG \big(\gammaB(s)\big),
\end{align*}
by \eqref{eq: dynamique pgamma}--\eqref{eq:
    dynamique pgamma-rigorous} since a \bbichar is a regular \bichar away from $s
\in B$.
For $s, s' \in [S_0, S[$, one thus has
$\py(s') - \py(s) = \int_s^{s'} \HpG \big(\gammaB(\sigma)\big) d \sigma$.
One has $\py(s') \to  \py^{0} = \py(S)$ as $s' \to S^-$. Dominated
convergence yields
\begin{align}
  \label{eq: py limit point broken bichar}
  \py(S) - \py(s) = \int_s^{S} \HpG \big(\gammaB(\sigma)\big) d \sigma.
\end{align}
As $\gammaB(s) \to \y^0$ one has  $\HpG \big(\gammaB(s)\big) =
\HpG(\y^0) + o(1)$ as $s \to S^-$, yielding
\begin{align*}
  \py(S) - \py(s) = (S-s) \HpG(\y^{0}) + o(S-s).
\end{align*}
This implies that $\py(s)$ is differentiable at $s = S^-$ and
$\frac{d}{d s} \py(S^-) = \HpG(\y^{0})$. 

Since the $\d_z$ component of $\HpG$ is $(\Hpz) \d_z$ and $z(S)=0$, with
\eqref{eq: py limit point broken bichar} one has
\begin{align*}
  - z(s) = \int_s^{S} \Hpz \big(\gammaB(\sigma)\big) d \sigma, 
\end{align*}
yielding \eqref{eq: precise null-derivability z reach G-G0}.

\medskip
If one proves
\begin{align}
  \label{eq: null-derivability Hpz reach G-G0}
  \Hpz \big(\gammaB(s)\big)= o(1) (s-S), 
  \quad \text{for} \ s \to S^- \ \text{in} \  J \setminus B,
\end{align}	
then one concludes that $\gammaB$ is differentiable at $s = S^-$ by
Lemma~\ref{lemma: derivability zeta} with $\frac{d}{d s} \gammaB (S^-) = \HpG(\y^0)$.

\bigskip
We now prove 
\eqref{eq: null-derivability Hpz reach G-G0} by proving the estimate
for  $\Hpz\big(\gammaB(s)\big)$ in \eqref{eq: precise null-derivability Hp z reach G-G0}. Consider the function 
\begin{align}
  \label{eq: function g - limit in glancing}
	g(s) = \frac{1}{2} \big(\Hpz\big)^2 
	\big( \gammaB(s)\big)
	- \big(\Hppz\big) \big( \gammaB(S)\big)  
	z(s).
\end{align}
Since $ \big(\Hppz\big) \big( \gammaB(S)\big)\leq 0$ and $z(s) \geq
0$ one has $g(s) \geq 0$. Despite having $\big( \Hpz\big) \big(
\gammaB(s)\big) $ discontinuous across any point of $B$, observe that
$g(s)$ can be  extended to $[S_0, S]$ as a continuous function. One
has $g(S)=0$.
Moreover, between two points of $B$ one has
\begin{align*}
  \frac{d}{d s} g(s) = \big(\Hpz\big)\big( \gammaB(s)\big) \eps(s,S), 
  \ \ \text{with} \ \ 
  \eps(s,S) = \big(\Hppz\big) \big( \gammaB(s)\big)
	- \big(\Hppz\big) \big( \gammaB(S)\big).
\end{align*}
One finds 
\begin{align*}
  \big| \frac{d}{d s} g(s)\big|
  \leq  \big| \big(\Hpz\big)\big( \gammaB(s)\big)\big| \, |\eps(s, S)|
  \lesssim g^{1/2}(s)|\eps(s, S)|,
  \qquad s \in [S_0, S[ \setminus B.
\end{align*}
Classically, with $a>0$, one replaces $g$ by $g^a= g+a>0$ for which
one has $|\frac{d}{d s} g^a(s)| \lesssim g^a(s)^{1/2} |\eps(s,S)|$ leading to
\begin{align*}
  \big|\frac{d}{d s} \big(g^a(s)^{1/2}\big)\big| \lesssim |\eps(s,S)|, 
  \qquad s \in [S_0, S[ \setminus B.
\end{align*}
Consider  $S_0 \leq s \leq s' < S$. Note that $B \cap [s,s']$ is finite since $B$ does not have any accumulation
point in $]0,S[$ by Lemma~\ref{lemma: broken bichar - closure B}. Note
also that $\eps(s,S)$ is bounded  on $[S_0,S]$. As a result
$(g^a)^{1/2}$ is Lipschitz and thus absolutely continuous on
$[s,s']$ and one finds 
\begin{align*}
  \big| g^a(s)^{1/2}- g^a(s')^{1/2}\big|
  \leq \int_s^{s'}\big|\frac{d}{d s} \big(g^a(s)^{1/2}\big)\big| \, d \sigma
  \lesssim \int_s^{s'} |\eps(\sigma,S)|\, d \sigma.
\end{align*}
Letting  $a \to 0^+$ 
one obtains 
\begin{align*}
  \big| g(s)^{1/2}- g(s')^{1/2}\big|
  \lesssim \int_s^{s'} |\eps(\sigma,S)|\, d \sigma.
\end{align*}
Since $s\mapsto\eps(s,S)$ is integrable on $[S_0,S]$ and $g(s') \to 0$ as $s' \to
S^-$, one finally obtains 
\begin{align*}
  0\leq |\Hpz \big( \gammaB(s)\big)\big| \lesssim g^{1/2}  (s)\lesssim \int_s^S
  |\eps(\sigma,S)|\, d \sigma.
\end{align*}
This is estimate \eqref{eq: precise null-derivability Hp z reach G-G0}.
Observe that $\eps(s,S) =o(1)$ as $s \to S^-$, because of the 
continuity of the function $\Hppz$. Hence, the estimate for $\Hpz \big( \gammaB(s)\big)$ in \eqref{eq: null-derivability Hpz reach G-G0} follows. 
\end{proof}
\begin{remark}
  If $\M$ is a $\Con^3$-manifold and $g$ a $\Con^2$-metric then $\Hpz$
  is $\Con^2$ and one can prove that a \bbichar
  $\gammaB(s)$  as above cannot
  reach a point in $\sgGb$ as a limit 
  as $s \to S^-$ even if $S \in
  \ovl{B}$. The proof  goes by contradiction as follows using the setting of the
  proof of Lemma~\ref{lemma: limit broken bichar in
    G-G0}. Consider $\y^0 \in \sgGb$ and a \bbichar that
  converges to $\y^0$ as $s \to S^-$. Consider $V$ a bounded \nhd of $\y^0$
  where $\Hppz \leq -C_0<0$ and where $\gammaB(s)$ lies for $s \in
  [S- \eps, S[$ for some $\eps>0$. With $p \in \Con^2$ one finds that $\Hppz \in \Con^1$. Set
  $g(s) = \Hpz\big(\gammaB(s) \big)^2/2 - z(s) \Hppz
  \big(\gammaB(s) \big)$. One has $g$ continuous and nonnegative, and
  moreover $g(s)=0$ if and only if $z(s) =  \Hpz\big(\gammaB(s)
  \big)=0$. Since moreover $\Hppz \big(\gammaB(s)
  \big)<0$ one finds that $g(s)=0$ if and only if $\gammaB(s) \in \sgGb$. Away from the points of $B$ one finds 
  $\frac{d}{d s} g(s) = - z(s) \Hp^3 z \big(\gammaB(s) \big)$. Since $\Hp^3 z $ is
  continuous and thus bounded in $V$ one obtains
  \begin{align*}
    \big| \frac{d}{d s} g(s)\big| \lesssim z (s) \lesssim - z(s) \Hppz
    \big(\gammaB(s) \big) \lesssim g(s). 
  \end{align*}
  Since $g$ is absolutely continuous on $[s,s']$ for $s < s' <S$ it follows that 
  \begin{align*}
    |g(s') - g(s)|  \leq  \int_{s}^{s'} \big| \frac{d}{d \sigma}
    g(\sigma) \big| \, d \sigma 
    \lesssim \int_s^{s'} g(\sigma) d \sigma.
    \end{align*}
    Letting $s'
  \to S^-$ one has  $g(s') \to 0$ and thus $0 \leq g(s) \lesssim
  \int_s^S  g(\sigma) d \sigma$.  The Gr\"onwall inequality yields
$g\equiv 0$, meaning that $\gammaB(s) \in \sgGb$ for $s\in [S-\eps,S]$;
a contradiction since we have considered a \bbichar on $ [S-\eps,S[$.

  Here, we have not been able to exclude that a \bbichar reaches
  $\sgGb$ in a limiting process. The proof that we have just recalled
  does not apply due to the lack of smoothness of $\Hppz$. This question remains open. 
  Observe however that in Lemma~\ref{lemma: limit broken bichar in
    G-G0} we obtain that $\frac{d}{d s} \gammaB(S) = \HpG\big(
  \gammaB(S) \big)$ in both cases $\y^0 \in \glGb$ and $\y^0\in \sgGb$.
  This suffices for the analysis we carry out in what follows.

  Note also that no thorough study of nonuniqueness issues at boundary
  has been carried out for $\Con^k$ coefficients, $k\geq 2$, up to our
  knowledge.
\end{remark}

The following proposition is a consequence of Lemma~\ref{lemma: limit
  broken bichar in G-G0} and generalizes part of it.
%%%%%%%%%%%%%%%%%%%%%%%%
% proposition              %
%%%%%%%%%%%%%%%%%%%%%%%%
\begin{proposition}
  \label{prop: limit generalized bichar in G-G0}
  Consider $\y^0 \in \glGb \cup \sgGb$. Let also $S>0$,
  $J=[0,S[$, $B$ a discrete subset of $J$, and
  $\gamma: J\setminus B  \to \TL \cap \Char p$
  be a map such that 
  \begin{align*}
    \lim_{{s \to S^-} \atop s \notin B} \gamma(s) = \y^0.
  \end{align*}
  One sets $\gamma(S) = \y^0$. Assume moreover that 
  \begin{enumerate} 
  \item if $s_0 \in B$ then  $\gamma(s)$ is a
    \bbichar for $s \in [s_0-\delta, s_0[\cup ]s_0, s_0 +
    \delta]$ for some $\delta>0$;
  \item if $s_0 \in [0,S[\setminus B$ and $\gamma(s_0) \notin \glGb \cup \sgGb$ then $\gamma(s)$ is a
    \bbichar for $s \in [s_0 -\delta, s_0 + \delta]\setminus B$
    for some $\delta >0$;
  \item $\pgamma(s) = \ppi\big(\gamma(s)\big)$ is differentiable at $s= S^-$ and 
  \begin{align}
    \label{eq: tangential differentiability approaching G}
    \frac{d}{d s} \pgamma(S^-) = \HpG (\y^{0}).
  \end{align}
  
  \end{enumerate}
Then, $\gamma(s)$ is differentiable at $s=S^-$ and 
\begin{align}
  \label{eq: differentiability approaching G}
    \frac{d}{d s} \gamma(S^-) = \HpG(\y^0).
\end{align}
\end{proposition}
%%%% proof of proposition
\begin{proof}
  For $s$ near $S^-$, use a single chart as in the proof of Lemma~\ref{lemma: limit broken bichar in G-G0}.
 Write $\gamma(s) = \big (t(s), x'(s), z(s), \tau(s),
  \xi(s) \big)$. One has $\pgamma(s) = \big (t(s), x'(s), z(s), \tau(s),
  \pxi(s) \big)$.
  
  With \eqref{eq: tangential differentiability approaching G} and
  Lemma~\ref{lemma: derivability zeta}, if one proves that $s \mapsto
  \Hpz(s)$ is differentiable at
  $s=S^-$ and  $\frac{d}{d s} \big(\Hpz(\gamma(s)\big)_{|s=S^-}=0$
  then $s \mapsto \gamma(s)$ is differentiable at $s=S^-$ one obtains \eqref{eq:
    differentiability approaching G}.

%   By Proposition~\ref{prop: G submanifolds}, $\Gb$ is a manifold
%   defined by $p= z= \Hpz=0$. Here, one has
%   $p\big(\gamma(s)\big)=0$. 

%   From the form of $\HpG$ given in \eqref{eq: definition HpG-extended}
%   one has $\HpG z = \Hpz$ and thus by \eqref{eq: tangential
%     differentiability approaching G}
% one has
%   \begin{align*}
%     \frac{d}{d s} z (S^-) =  \HpG z (\y^0) = \Hpz (\y^0) =0.
%   \end{align*} 
%   If one proves that
%   $\frac{d}{d s} \big(\Hpz(\gamma(s)\big)_{|s=S^-}=0$ then one has
%   $\frac{d}{d s} \gamma(S^-) \in T_{\y^0} \Gb \subset T_{\y^0} \pdTL$.  This implies that
%   $\frac{d}{d s} \gamma(S^-) = d \ppi(\y^0) \big( \frac{d}{d s}
%   \gamma(S^-)\big)$.  However, one has
%   $d \ppi(\y^0) \big( \frac{d}{d s} \gamma(S^-)\big) = \frac{d}{d s}
%   \pgamma(S^-)$. Hence with \eqref{eq: tangential differentiability approaching G} one concludes that \eqref{eq: differentiability
%     approaching G} holds.

    Below we prove that 
    \begin{align}
      \label{eq: differentiability approaching G-2}
      \Hpz\big(\gamma(s)\big) = o (S-s) \quad
      \text{as} \  s \to S^-  \ \text{in} \ [0,S[ \setminus B.
    \end{align}
    Consider $s \in [0,S[ \setminus B$. One says that $s \in
    B_{0}$ if $\gamma(s) \in \glGb\cup \sgGb$; then $z(s) =0$ and $\Hpz \big(
    \gamma(s)\big)=0$. One thus only needs to prove \eqref{eq:
      differentiability approaching G-2} for $s\to S^-$  in  $[0,S[ \setminus (B
    \cup B_{0})$. 
    
    Consider $\eps>0$. By continuity of $\Hppz$, there exists
    $0 < S_\eps < S$ such that
    \begin{align*}
      \big| 
          \big(\Hppz\big) \big( \gamma(s)\big)
	- \big(\Hppz\big) \big( \gamma(S)\big)
      \big|\leq \eps,
    \end{align*}  
      if $s \in [S_\eps, S] \setminus B$.  Consider
    $s \in [S_\eps, S[ \setminus (B \cup B_{0})$;
    by assumption $\gamma$ is locally a \bbichar.  Define $s_1$
    as the supremum of the connected component of $s$ in
    $[S_\eps,S[ \setminus B_{0}$. Note that $s_1$ can be equal to
    $S$. By continuity one has $\gamma(s_1)
    \in \glGb \cup \sgGb$
    and on the interval $[s,s_1]$ one faces the situation described in
    Lemma~\ref{lemma: limit broken bichar in G-G0}.  One thus has
  \begin{align*}
    0\leq \big|\Hpz \big(\gamma(s)\big)\big| 
    &\leq C
    \int_s^{s_1} 
          \big|\big(\Hppz\big) \big( \gamma(\sigma)\big)
	- \big(\Hppz\big) \big( \gamma(s_1)\big)
          \big|
    d\sigma\\
   &\leq  2 C \eps  (s_1 -s) \lesssim (S-s) \eps,
  \end{align*}
 meaning \eqref{eq: differentiability approaching G-2} holds.
\end{proof}

\bigskip
The notion of \bbichar is not sufficient to understand the
propagation of the support of measures as points in $\glGb$ and $\sgGb$
are not considered. Moreover, a sequence of \bbichars may
converge to a curve that is not a \bbichar. This leads to the
introduction  of \gbichars as in Definition~\ref{def: generalized bichar-intro}.
With the notion of  \bbichars introduced above one may also
write the definition of  \gbichars as follows. 
%%%%%%%%%%%%%%%%%%%%%%%%
% definition           %
%%%%%%%%%%%%%%%%%%%%%%%%
\begin{definition}[\gbichar]
	\label{def: generalized bichar}
  Suppose $J \subset \R$ is an interval, $B$ a discrete subset of $J$, and
  \begin{align*}
    \gammaG: J\setminus B 
    \to \Char p \cap \TL.
  \end{align*}
  One says that $\gammaG$ is a \gbichar if  the
  following properties hold:
  \begin{enumerate}
  \item for $s \in J \setminus B$, $\gammaG(s) \notin \Hb$ and the map $\gammaG$ is differentiable
    at $s$ with
	\begin{equation*}
		\frac{d}{d s} \gammaG(s) = \XG\big( \gammaG(s)\big).
	\end{equation*}
  \item if $s^0\in B$, then $\gammaG$ is a \bbichar
    on an interval $[s^0-\eps, s^0+\eps] \setminus \{s^0\}$ for some $\eps>0$.
  \end{enumerate}
\end{definition}
Recall the definition of the vector field $\XG$ in \eqref{eq: def XG}.
%%%%%%%%%%%%%%%%%%%%%%%%
% remark               %
%%%%%%%%%%%%%%%%%%%%%%%%
\begin{remark}
  \label{rem: def generalized bichar}

  \leavevmode  
  \begin{enumerate}
  \item
 From Lemma~\ref{lemma: properties HpG} one sees that $\XG$ is
  continuous at points of $\TL \setminus \sgGb$, in particular at
  points in $\glGb$. It is however
  discontinuous at points of $\sgGb$. In fact if $\y \in \sgGb$ then
  $\XG(\y) = \HpG(\y) \neq \Hp(\y)$ since $\Hppz<0$ and in any \nhd of $\y$ there are
  points in $\Hb$ where $\XG = \Hp$.  
  Yet, note that restricted to $\Gb$ the vector field $\XG$ is continuous.

\item Similarly to what is observed in Remarks~\ref{rem:def-bichar}-\eqref{rem: constant tau}  and \ref{rem: def broken bichar} one has, for $s \in J \setminus B$,
    \begin{align*}
    \norm{\xi(s)}{x(s)} = |\tau(s)|
  \end{align*}
   constant along a \gbichar.
\end{enumerate}
\end{remark}

%%%%%%%%%%%%%%%%%%%%%%%%
% lemma                %
%%%%%%%%%%%%%%%%%%%%%%%%
\begin{lemma}
  \label{lemma: generalized bichar locally gliding}
  Suppose $s^0 \in J \setminus \ovl{B}$ and $\gammaG(s^0)\in \sgGb$. Then,
  $\gammaG(s)\in \sgGb$ for $s$ in a \nhd of $s^0$.
\end{lemma}
%%%% proof of lemma
\begin{proof}
  Proceed by contradiction and 
  assume that there
  exists a sequence $s_n \in J\setminus \ovl{B}$ that converges to
  $s^0$ such that $\gammaG(s_n) \notin \sgGb$.  If a
  subsequence $s_{n_k}$ is such that $\gammaG(s_{n_k})\in \glGb$, then
  \begin{align*}
    0= \Hppz\big( \gammaG(s_{n_k})\big)
    \ \mathop{\longrightarrow}_{k \to +\infty} \
    \Hppz\big( \gammaG(s^0)\big)< 0,
  \end{align*}
  a contradiction.  Thus, there exists a subsequence
  $s_{n_k}$ with $\gammaG(s_{n_k})\notin \glGb\cup \sgGb$. Then, the
  \gbichar is a \bbichar in a \nhd of $s_{n_k}$. Consider the part of
  the \gbichar that is a maximal \bbichar. It ceases to exist at a
  point between $s^0$ and $s_{n_k}$. There, by Lemma~\ref{lemma: limit
    broken bichar in G-G0} it reaches a point in $\glGb$ (here one is
  away from hyperbolic points). One is thus back to assuming that
  there is a subsequence $r_n$ that converges to $s^0$ with
  $\gammaG(r_n) \in \glGb$ leading to a contradiction.
\end{proof}

If $s \mapsto \gammaG(s)$  is a \gbichar it is obviously discontinuous 
across points  $s\in B$ and continuous otherwise. The following lemma states $\Con^1$-regularity away from points in $\ovl{B}$.
%%%%%%%%%%%%%%%%%%%%%%%%
% lemma                %
%%%%%%%%%%%%%%%%%%%%%%%%
\begin{lemma}
  \label{lemma: C1 generalized bichar}
  If $s^0 \in J \setminus \ovl{B}$ then $s \mapsto \gammaG(s)$ is
  $\Con^1$ in a \nhd of $s^0$.
\end{lemma}
%%%% proof of lemma
\begin{proof}
  There exists a \nhd $W$ of $s^0$ with $W \cap \ovl{B} = \emptyset$. If one proves that $\gammaG(s)$ is $\Con^1$ at $s^0$ the same holds for any point in $W$, hence the result. 

  First, assume that $\gammaG(s^0)\notin \sgGb$. At such point the
  vector fields $\XG$ is continuous; see the first part of
  Remark~\ref{rem: def generalized bichar}. It follows that
  $\gammaG(s)$ is $\Con^1$ at $s^0$. Second, assume that
  $\gammaG(s^0)\in \sgGb$. Then, $\gammaG(s)\in \sgGb$ for $s$ in a
  \nhd of $s^0$ by Lemma~\ref{lemma: generalized bichar locally
    gliding}.  There, one has $\XG\big(\gammaG(s)\big) =
  \HpG\big(\gammaG(s)\big)$. Since $\HpG$ is continuous this yields
  the result.
\end{proof}

%%%%%%%%%%%%%%%%%%%%%%%%
% remark               %
  %%%%%%%%%%%%%%%%%%%%%%%%
  \begin{remark}
    \label{remark: only diff at B bar setminus B}

   \leavevmode  
    \begin{enumerate}
    \item
 For $s^0 \in \ovl{B}\setminus B$, the definition only states
the differentiability of $\gammaG(s)$ at $s=s^0$. In particular, the
definition does not imply that the derivative of $\gammaG$  is
continuous near such a point.
\item
  In a local chart with the notation of Section~\ref{sec: local extension away from boundary}, if one sets $\py(s) = \ppi \big(\gammaG(s)\big)$ one
  sees that $\py(s)$ is $\Con^0$. In particular the $z$-component is
  continuous and vanishes at $s=S^\pm$ for $S \in B$. 
If one considers the map $\cphi$ introduced in Section~\ref{sec: The
    compressed cotangent bundle} one sees that $s\mapsto
  \cphi(\gammaG(s))$ takes values in the compressed cotangent bundle
  and can be extended to the whole interval $J$ as a continuous
  function. This aspect is used in the construction of a
  \gbichar in what follows.
  \end{enumerate}
\end{remark}

%%%%%%%%%%%%%%%%%%%%%%%%
% lemma                %
%%%%%%%%%%%%%%%%%%%%%%%%
\begin{lemma}
  \label{lemma: B bar structure}
  One has the following equivalence
  \begin{align}
  \label{eq: caracterisation B bar}
    s \in \ovl{B}\setminus B \  \ \Equiv  \  \
    \exists (s_n)_n \subset B, \ \text{with} \ s_{n+1} \notin \{s_{k}; k
    \leq n\} \
    \text{such that} \  s_n \to s.
  \end{align}
  Moreover, $\ovl{B}\setminus B$ is a closed set. 
\end{lemma}
% proof of lemma
\begin{proof}
  The ``$\Rightarrow$'' part of \eqref{eq: caracterisation B bar} is
  straightforward. The ``$\Leftarrow$'' part is a consequence of
  $B$ being a discrete set. Assume now that $s
\in \ovl{\ovl{B} \setminus B}$. Then there exists $(s_n)_n \subset
\ovl{B} \setminus B$ such that $s_n \to s$. One can construct a
sequence $(s'_n)_n \subset B$ with $s'_{n+1} \notin \{s'_{k}; k
  \leq n\}$
such that $s'_n \to  s$. One has $s \in
\ovl{B}\setminus B$ by \eqref{eq: caracterisation B bar}.
\end{proof}
%%%%%%%%%%%%%%%%%%%%%%%%
% lemma                %
%%%%%%%%%%%%%%%%%%%%%%%%
\begin{lemma}
  \label{lemma: generalized bichar - derivative endpoint}
  Suppose $\gammaG: J \setminus B \to \TL$
 is a \gbichar. 
  \begin{enumerate}
    \item If $s \in J \cap \ovl{B} \setminus B$ then $\gammaG(s) \in
      \glGb\cup \sgGb$.\\[-8pt]
    \item Set $S = \sup J$
  and suppose that $S< +\infty$. 
  Then, the limit
  \begin{align}
    \label{eq: derivative generalized bichar endpoint1}
    \y^0 = \lim_{{s\to S} \atop {s \in J \setminus B}} \gammaG(s)
  \end{align}
  exists. If one sets $\gammaG(S) = \y^0$ then $\gammaG(s)$ is
differentiable at $s= S^-$ and 
\begin{align}
  \label{eq: derivative generalized bichar endpoint2}
    \frac{d}{d s} \gammaG(S^-) = \XG(\y^0).
\end{align}
Then, Part (1) applies if $S \in \ovl{B} \setminus B$. 
A similar result holds for $S= \inf J$.
  \end{enumerate}
\end{lemma}
%%%% proof of lemma
\begin{proof}
  {\bfseries Part (1):} $s \in J \cap \ovl{B} \setminus B$ .  The
  argument of the proof of Lemma~\ref{lemma: broken bichar - closure
    B} can be applied.

  \medskip
  \noindent
{\bfseries Part (2):} $S = \sup J < +\infty$. The result is clear if $S \in J
\setminus B$ using the first part. Suppose $S \notin  J
\setminus B$.
First, consider the case $S \notin \ovl{B}$. Then one has $\gammaG(s)
\in \Char p \cap \TL \setminus \Hb$ on an interval of the form $[S-
\eps, S[$ for some $\eps>0$, and 
\begin{align*}
  \frac{d}{d s} \gammaG(s) = \XG\big(\gammaG(s)\big), 
\quad s \in [S-\eps, S[.
\end{align*}
Since $\gammaG(s)$ remains bounded, then $\XG\big(\gammaG(s)\big)$
remains bounded yielding the existence of the limit
\begin{align*}
    \y^0 = \lim_{{s\to S} \atop {s \in [S-\eps, S[}} \gammaG(s).
  \end{align*}
On the one hand if $\y^0 \in \TL \setminus \Hb$, then $s \mapsto \XG\big(\gammaG(s)\big)$ is
continuous on $[S-\eps, S[$ and has a limit at $s=S^-$; the proof follows the arguments of Lemma~\ref{lemma: C1 generalized bichar}. This implies that $\gammaG(s)$ is
differentiable at $s= S^-$ and \eqref{eq: derivative generalized
  bichar endpoint2} holds.
On the other hand if $\y^0 \in \Hb$ then $\y^0 \in \Hb^-$, meaning
that locally the \gbichar is a regular \bichar reaching a
hyperbolic point. This again gives    that $\gammaG(s)$ is
differentiable at $s= S^-$ and \eqref{eq: derivative generalized
  bichar endpoint2} holds.

\medskip
Second, consider the case $S \in \ovl{B}$. The argument of the proof
of Lemma~\ref{lemma: broken bichar - closure B} can be applied  {\em
  mutatis mutandis} yielding the existence  of the limit $\y^{0}$ in \eqref{eq:
  derivative generalized bichar endpoint1}.
Then, the proof of Lemma~\ref{lemma: limit broken bichar in G-G0}
can be applied with
some slight modifications. One finds with the same arguments that
\begin{align*}
  \frac{d}{d s} \py(S^-) = \HpG(\y^{0}).
\end{align*}
One introduces the same function $g$ as in \eqref{eq: function g -
  limit in glancing}.
At points in $J \setminus B$ the function $g$ is differentiable, even at
point in $\ovl{B} \setminus B$. In fact, if $ s^1\in J \setminus\ovl{B}$ on
finds 
\begin{multline}
  \label{eq: derivee g - gammaG}
  \frac{d}{d s} g(s^1) = \big(\Hpz\big)\big( \gammaG(s^1)\big) \eps(s^1,S), \\
 \text{with} \ \ 
  \eps(s^1,S) = \big(\Hppz\big) \big( \gammaG(s^1)\big)
	- \big(\Hppz\big) \big( \gammaG(S)\big).
\end{multline}
as in the proof of Lemma~\ref{lemma: limit broken bichar in G-G0}. 
At a point $s^1 \in J \cap \ovl{B} \setminus B$ one has $\gammaG(s^1) \in
\glGb\cup \sgGb$ by the first part of the lemma. Hence, $g(s^1) =
\frac{d}{d s} g(s^1) =0$ since $z(s^1) = 0$,  $\frac{d}{d s} z(s^1)=
\HpG z (s^1)= \Hpz (s^1)=0$, and $\frac{d}{d s} \Hpz (s^1) = \HpG \Hpz
(s^1)=0$ (recall that $\HpG \Hpz=0$ on $\pdTL$ since $\HpG$ is tangent
to $\pdTL$).
Hence \eqref{eq: derivee g - gammaG} holds in $J \setminus B$. Setting
$g^a = g+a$ for $a>0$ one finds
\begin{align*}
  \big|\frac{d}{d s} \big(g^a(s)^{1/2}\big)\big| \lesssim |\eps(s,S)|, 
  \qquad s \in J \setminus B.
\end{align*}
Since $B$ is at most countable, one finds that $(g^a)^{1/2}$ is
Lipschitz thus absolutely continuous on $J$. The remainder of the
proof of Lemma~\ref{lemma: limit broken bichar in G-G0} then applies. 
One finds that $\gammaG(s)$ is
differentiable at $s= S^-$ and \eqref{eq: derivative generalized
  bichar endpoint2} holds.
\end{proof}

\medskip
For a \gbichar the notation $\GammaG$ was introduced in
Definition~\eqref{def: generalized bichar 2-intro}. If one considers
the map $\cphi$ associated with the compressed cotangent bundle
introduced in Section~\ref{sec: The compressed cotangent bundle}
observe that
  \begin{align*}
    \cphi(\GammaG) 
    = \ovl{\cphi \{ \gammaG(s); \ s \in J \setminus B\}}
    \ \ \text{and} \ \ 
    \GammaG = \cphi^{-1} \big( \cphi(\GammaG) \big).
  \end{align*}

\subsection{On the measure equation}
\label{sec: On the measure equation}

The measure equation of  Assumption~\ref{assumption: Gerard-Leichtnam
  equation}, is one of the two hypothesis in our main result,
Theorem~\ref{theorem: measure support propagation}:
\begin{align*}
    \transp \Hp \mu  = f \mu - \int_{\y \in \pHb \cup \pGb} 
    \frac{\delta_{\y^+} - \delta_{\y^-}}
    {\dup{\xi^+- \xi^-}{\nx}_{T_x^*\M, T_x\M} 
    } \ d \nu (\y) \qquad \dans \ \U.
  \end{align*} 
With what precedes one sees that this equation is  of geometrical
nature. 

In Section~\ref{sec: main result} we gave an interpretation of the
integrand for $\y \in \pGb=\Gb$ using that $\Gb \cup \Hb = \ovl{\Hb}$. 
Here a similar interpretation can be done yet using the geometrical
$\Con^1$-parametrization~\ref{eq: parametrization cotangent bundles}
of the cotangent bundle. 

We thus  use the
  parametrization $(\py, \Hpz(\y))$ of $\TL$ and denote by $\vartheta$
  the variable $\Hpz (\y)$. If $\y = (t,x,\tau,\xi)\in \pH$ then with
  Lemma~\ref{lemma: relevement pG pH} one has $\y^\pm =
  (t,x,\tau,\xi^\pm) \in \H^\pm$ with
  \begin{align*}
    \xi^+ = \xi + \alpha(x) \vartheta^+ \nxs
    \ \ \text{and} \ \ 
    \xi^- = \xi + \alpha(x)\vartheta^- \nxs
  \end{align*}
  and $\vartheta^+ >0$ and 
  $\vartheta^- = - \vartheta^+$.
  Then
  \begin{align}
    \label{eq: computation denominator GL}
    \dup{\xi^+- \xi^-}{\nx}_{T_x^*\M, T_x\M}  
    = \alpha(x) (\vartheta^+- \vartheta^-) = 2 \alpha(x) \vartheta^+.
    \end{align}

If $q = q(\py, \vartheta)$ is a $\Con^1$-test function then, for $\y\in \pHb$
one has 
\begin{align*}
  \frac{\dup{\delta_{\y^+} - \delta_{\y^-}}{q}}
  {\dup{\xi^+-\xi^-}{\nx}_{T_x^*\M, T_x\M}}
  = \frac{ q(\y, \vartheta^+) - q(\y, - \vartheta^+)}
  {2 \alpha(x) \vartheta^+}.
\end{align*}
If now a sequence $(\y^{(n)} )_n\subset \pHb$ converges to $\y = \py \in \pGb$
then 
\begin{align}
  \label{eq: remark: integrand GL equation on G}
  \frac{\dup{\delta_{\y^{(n)+}} - \delta_{\y^{(n)-}}}{q}}
  {\dup{\xi^{(n)+}-\xi^{(n)-}}{\nx}_{T_x^*\M, T_x\M}}
  \ \to \ \frac{1}{\alpha(x)}\d_\vartheta q (\y,0).
\end{align}
 Up to the factor $1/\alpha$, the integrand in \eqref{eq: Gerard-Leichtnam equation} for $\y\in
 \pGb$ is thus to be understood as the derivative with respect to the
 variable $\vartheta$ at $\vartheta=0$.
If one considers the quasi-normal geometric coordinates given by
Proposition~\ref{prop: quasi-normal coordinates} that are used in
Section~\ref{sec: Geometrical setting-intro} then one has $\vartheta = 2
\zeta$ and $\alpha=1/2$ since $\alpha = (2 \Hzzp)^{-1/2}$ and $\Hzzp=2$
in those coordinates. One thus recovers the observation made in Remark~\ref{remark: integrand GL equation on G-intro}.

%%%%%%%%%%%%%%%%%%%%%%%%
% remark               %
%%%%%%%%%%%%%%%%%%%%%%%%
\begin{remark}
Observe that since $\Sigma (\y^+)  =  \y^-$, under
Assumption~\ref{assumption: Gerard-Leichtnam equation} one finds
\begin{align}
  \label{eq: anti-invariance Hp mu by involution}
  \Sigma_* (\transp \Hp \mu - f \mu)  = - (\transp \Hp \mu - f \mu).
\end{align}
\end{remark}

%%%%%%%%%%%%%%%%%
% section
%%%%%%%%%%%%%%%%%
\section{Propagation of the measure  support}
\label{sec: Propagation of the measure  support}

This section is devoted the proof of our main result, Theorem~\ref{theorem: measure support propagation}.

In what follows, we consider all possible situations for $\y^0$:
$\y^0 \in \TL \setminus \dTL$, that is, a point away from the
boundary (Section~\ref{sec: propagation Away from the boundary}),
$\y^0 \in \Hb$, that is, a hyperbolic point at the boundary
(Section~\ref{sec: propagation hyperbolic points}),
$\y^0 \in \glGb\cup \sgGb$, that is, an order-3-glancing  or a gliding
point (Section~\ref{sec: General glancing points}), and
$\y^0 \in \sdGb$, that is a diffractive point
(Section~\ref{sec: diffractive points}).  In all situations,
we prove a local result. These results put together yield the
proof of Theorem~\ref{theorem: measure support propagation}; see
Section~\ref{sec: proof theorem: measure support propagation}. Working
locally allows one to use a single chart. This permits to use the
extension away from the boundary of notions only geometrically
meaningfull at the boundary $\pT_\y^* \L$, $\G$, $\H$, $\HpG$, \etc;
see Section~\ref{sec: local extension away from boundary}. We 
make use of these extension without mentionning it in what follows. 

%%%%%%%%%%%%%%%%%
% sub-section
%%%%%%%%%%%%%%%%%
\subsection{Away from the boundary}
\label{sec: propagation Away from the boundary}

With Assumption~\ref{assumption: Gerard-Leichtnam equation},
away from the boundary, one has $\transp{\Hp} = f \mu$ locally.  Then,
Theorem~\ref{th: ODE} gives the following result.
%%%%%%%%%%%%%%%%%%%%%%%%
% proposition          %
%%%%%%%%%%%%%%%%%%%%%%%%
\begin{proposition}
  \label{prop: propagation Away from the boundary}
  Suppose $\y^0 \in \supp \mu \cap \big(\TL \setminus
  \dTL\big)$. If $V$ is an open \nhd of $\y^0$ in $\TL$ such
  that $V \cap \dTL = \emptyset$ then there exists a maximally
  extended \bichar $\gamma(s)$ in $V$, for $s \in J$, 
  an open interval of $\R$, with $0 \in J$, such that 
  \begin{align*}
    \gamma(0) = \y^0 \ \ \text{and} \ \ 
    \{ \gamma(s); \ s \in J\} \subset \supp \mu.
  \end{align*}
\end{proposition}

\medskip
%%%%%%%%%%%%%%%%%
% sub-section
%%%%%%%%%%%%%%%%%
\subsection{Hyperbolic points}
\label{sec: propagation hyperbolic points}

Consider $\y^0 \in \Hb^\pm$. Consider a local chart as in \eqref{eq: local chart boundary}. As $\Hb^\pm$ is an open subset of $\dTL \cap \Char p$ one can
choose $V^0$ an open subset of $T^*\R^{1+d}$ such that $\y^0 \in V^0$
and
\begin{align*}
  \big( V^0 \cup \Sigma (V^0) \big)
  \cap \dTL \cap \Char p \subset \Hb^\pm.
\end{align*}

If $q\in \Con^1_c(\R^{2d+2})$ is such that $\supp q \subset V^0$,  then, from \eqref{eq: Gerard-Leichtnam equation} one deduces
\begin{align*}
  \bigdup{\transp{\Hp} \mu }{q} = \bigdup{\mu }{fq}
  \mp \int_{\y \in \pHb \cup \pGb} 
    \frac{q(\y^\pm)} 
    {\dup{\xi^+- \xi^-}{\nx}_{T_x^*\M, T_x\M} 
    } d \nu (\y), 
\end{align*}
 and $\xi^+- \xi^-\neq 0$ on the support of $\y \mapsto q(\y^\pm)$ that does
 not meet $\pGb$. One thus
 finds that the measure $\mu$ is solution in $V^0$ to an equation of
 the form 
 \begin{align*}
   \transp{\Hp} \mu = f \mu + \tilde{\mu} \otimes \delta_{z=0},
 \end{align*}
 where $\tilde{\mu}$ is a measure on $\{z=0\} = \dTL$ in $V^0$.
The following lemma is the result of Lemma~\ref{lem: support goes to
  the trace - transverse X}
translated into the present setting. 
%%%%%%%%%%%%%%%%%%%%%%%%
% lemma                %
%%%%%%%%%%%%%%%%%%%%%%%%
\begin{lemma}
  \label{lem: support goes to the trace - hyperbolic point}
Suppose $\y \in \Hb$.  Then,  $\y \in \supp \mu$ if and only if
  $\y \in \supp(\tilde{\mu} \otimes \delta_{z=0})$.
\end{lemma}
Similarly, one has $\transp{\Hp} \mu = f \mu + \hat{\mu} \otimes
\delta_{z=0}$ in $\Sigma(V^0)$, with  $\hat{\mu}$ a measure on $\{z=0\}\cap \Sigma(V^0)$.
 From \eqref{eq: anti-invariance Hp mu by involution} one has
 \begin{align}
  \label{eq: supp tilde mu involution}
  \Sigma_* \big(\tilde{\mu} \otimes \delta_{z=0})
  = - \hat{\mu}\otimes \delta_{z=0}
  \ \ \dans \ \ \Sigma (V^0).
\end{align}
Lemma~\ref{lem: support goes to the trace - hyperbolic point} 
and \eqref{eq: supp tilde mu involution} give the following
proposition.
%%%%%%%%%%%%%%%%%%%%%%%%
% proposition                %
%%%%%%%%%%%%%%%%%%%%%%%%
\begin{proposition}
  \label{prop: propagation rho- to rho+}
  Suppose $\y \in \Hb$. 
  Then, $\y \in \supp \mu$ if and only if $\Sigma (\y) \in \supp \mu$.
\end{proposition}

%The same analysis can be performed in $\Sigma(V^0)$ that is an open \nhd of $\Sigma (\y) \in
%\H^{\mp}$. 
%As $\Sigma_*  \delta_{z=0} = \delta_{z=0}$,   one then obtains
%$\Sigma_*   \tilde{\mu} =  - \tilde{\mu}$, which gives
%\begin{align}
%  \label{eq: supp tilde mu involution}
%  \Sigma \big(\supp(\tilde{\mu})  \cap V^0\big) = \supp(\tilde{\mu})  \cap \Sigma (V^0).
%\end{align}

\medskip
By Lemma~\ref{lem: direction bichar hyperbolic pt}, 
for $\y \in \Hb^+$ if $\gamma(s)$ is a \bichar that goes
through $\y$ at $s=0$, then $\{ \gamma(s)\}_{s\in ]0,S[}$ lies in 
$\TL \setminus \dTL$ while $\{ \gamma(s)\}_{s\in ]-S,0[}$ lies
in the complement of $\TL$,
 for some $S>0$. For $\y\in \Hb^-$, this is the opposite.
%%%%%%%%%%%%%%%%%%%%%%%%
% definition           %
%%%%%%%%%%%%%%%%%%%%%%%%
\begin{definition}
  \label{def: half bichar}
For $\y \in \Hb^+$ (\resp $\Hb^-$) and $\gamma(s)$ a \bichar that goes
through $\y$ at $s=0$  call $\gamma(s)$,
for $s\geq 0$ (\resp $s\leq 0$), a {\em half} \bichar initiated at $\y$.
\end{definition}
For $\y\in \Hb$ a half \bichar is locally contained in
$\TL$, that is, in $\{z \geq 0\}$.
%%%%%%%%%%%%%%%%%%%%%%%%
% definition           %
%%%%%%%%%%%%%%%%%%%%%%%%
\begin{definition}
  \label{def: incoming and outcoming bichar at hyperbolic point}
  Suppose $F$ is a closed set of $\TL$ and 
  $\y  \in F \cap \Hb^+$ 
  (\resp $F \cap \Hb^-$).
  One says that a half \bichar $\gamma(s)$ initiated at $\y$ is
locally contained in $F$ if for some $S>0$ one has $\{
\gamma(s)\}_{s\in [0,S[} \subset F$ (\resp $\{ \gamma(s)\}_{s\in ]-S,0]} \subset F$).
\end{definition}

\medskip
With the notion introduced in Definitions~\ref{def: half bichar} and
\ref{def: incoming and outcoming bichar at hyperbolic point} we now
state the following result.
%%%%%%%%%%%%%%%%%%%%%%%%
% proposition                        %
%%%%%%%%%%%%%%%%%%%%%%%%
\begin{proposition}
  \label{prop: propagation from rho+, rho-}
  Suppose $\y^0 \in \Hb^\pm$ and $V^0$ is open subset of $T^* \R^{1+d}$
  such that $V^0 \cap \dTL \cap \Char p \subset \Hb^\pm$.  If
  $\y^0 \in \supp \mu$, then there exists a half \bichar initiated at
  $\y^0$ that is locally contained in $\supp \mu$. Moreover, this
  half \bichar can be chosen maximally extended in $V^0$.
\end{proposition}
%%%% proof of proposition
\begin{proof}
  Applying Theorem~\ref{th: propagation from boundary - transport
    equation} to $\y^0$ gives the existence of a half \bichar
  $\{ \gamma(s)\}_{s\in [0,S[}$, with $S>0$, initiated at $\y^0$
  contained in $\supp \mu$.  Since
  $\gamma(S/2) \in \TL \setminus \dTL$, combining with the
  result of Proposition~\ref{prop: propagation Away from the boundary} 
  one can obtain a maximally extended such half \bichar in $V^0$.
\end{proof}
\bigskip
From Propositions~\ref{prop: propagation rho- to rho+} and \ref{prop:
  propagation from rho+, rho-} one deduces the following result. 
%%%%%%%%%%%%%%%%%%%%%%%%
% theorem              %
%%%%%%%%%%%%%%%%%%%%%%%%
\begin{corollary}
  \label{th: hyperbolic point}
  Suppose $\y^0 \in \Hb$. There is no  half \bichar initiated at $\y^0$
  locally contained in $\supp \mu$ if and only if there is no half
  \bichars initiated at  $\Sigma(\y^0)$ locally contained in $\supp \mu$.
\end{corollary}
\setcounter{theorembiss}{\arabic{theorem}}
Equivalently, this result reads as follows.
%%%%%%%%%%%%%%%%%%%%%%%%
% theorem              %
%%%%%%%%%%%%%%%%%%%%%%%%
\begin{corollarybis}
%\begin{corollary}
  \label{th: hyperbolic point -bis}
 Suppose $\y^0 \in \H$. If there exists a
  half \bichar initiated at $\y^0$ locally contained in $\supp \mu$ 
  then $\{ \y^0, \Sigma (\y^0)\} \subset \supp \mu$ and there exists  a
  half \bichar initiated at  $\Sigma (\y^0)$ locally contained in $\supp \mu$.
%  \end{corollary}
\end{corollarybis}
In other words, if a \bichar locally contained in $\supp \mu$
hits the boundary at a hyperbolic point, then 
the support of $\mu$ is transported along at least one
\gbichar through that hyperbolic point and its image by $\Sigma$.

\medskip
\subsection{Glancing points}

In Section~\ref{sec: General glancing points}, we provide
a propagation tool for general glancing points to 
be used for points in $\sgGb$ and $\glGb$. For points in $\sdGb$, that
is, diffractive points, we provide a more useful propagation tool in
Section~\ref{sec: diffractive points}.

%%%%%%%%%%%%%%%%%%%%%%%%
% lemma                %
%%%%%%%%%%%%%%%%%%%%%%%%
\begin{lemma}
  \label{lemma: zeta localizaton - char p -G}
 Suppose $\y^0 =(t^0,x^0,\tau^0,\xi^0) \in \Gb$, that is,
 $\Hpz(\y^0)=0$, and $V^0$ is a
  bounded \nhd of $\y^0$ in $\TL$ that lies in a local chart. There exists
  $C_0>0$ such that
  \begin{align*}
    %\label{eq: local dependence of zeta upon other variables}
    |\Hpz(\y)| \leq C_0 \Norm{\py - \py^0}{}^{1/2},
    \qquad \y \in \Char p \cap V^0.
  \end{align*}
\end{lemma}
%%%% proof of lemma
\begin{proof}
  With \eqref{eq: decomposition point in char -ext}, if $\y =(t,x,\tau,\xi) \in \Char p \cap V^0$   one has
  \begin{align*}
    \alpha(x)^2 \Hpz (\y)^2 = \lambda^2 = \tau^2  - |\pxi|_x^2 = -
    p(\py).
  \end{align*}
  As $\y^0 \in \Gb$, one has $p (\y^0)=p(\py^0)  = 0$, yielding
  \begin{align*}
    0 \leq \alpha(x)^2 \Hpz (\y)^2 
    =p(\py^0) - p(\py)
    \lesssim \Norm{\py - \py^0}{},
  \end{align*}
  since $V^0$ is bounded and $p$ is $\Con^1$. Since $\alpha(x)^{-1} = (2 \Hzzp)^{1/2}$ is
 bounded on $V^0$ the result follows.
\end{proof}
\subsubsection{General glancing points}
\label{sec: General glancing points}
We prove the following Proposition. 
%%%%%%%%%%%%%%%%%%%%%%%%
% proposition          %
%%%%%%%%%%%%%%%%%%%%%%%%
\begin{proposition}
  \label{prop: propagation glancing region}
  Suppose $K$ is a compact set of $\TL \setminus 0$ that lies in a
  local chart. There exists $C_0>0$
  such that 
  \begin{multline}
    \label{eq: py propagation glancing region}
    \forall \eps>0, \ 
    \exists \delta_0>0, \
      \forall \y^0   \in K \cap \Gb \cap \supp \mu, \
      \forall \delta\in ]0,\delta_0], \
    \exists \y \in \supp \mu 
    \\
    \text{such that}\ \ 
    \py \in 
    B\big(\y^{0}+ \delta \HpG (\y^{0}), C_0 \delta \eps\big),
  \end{multline}
and $|\Hpz(\y)| \leq C_0  \delta^{1/2}$.

\end{proposition}
The notation $\y^{0}+ \delta \HpG (\y^{0})$ is here to be understood
in the local coordinates where such computation makes sense.
%%%% proof of proposition
\begin{proof}
We make some preliminary remarks:
\begin{enumerate}
\item Since $\Hpz(\y^0)=0$, the estimation of $|\Hpz(\y)|$ follows from Lemma~\ref{lemma: zeta localizaton - char p -G}
  as one has
$\Norm{\py - \py^0}{} \lesssim \delta$. We thus prove the existence of $\y$
such that  \eqref{eq: py propagation glancing region} holds.
\item It suffices to consider $0 < \eps \leq 1$. 
\end{enumerate}

 Consider $\y^0 \in K \cap \Gb$. Then, one has $\py^0 = \y^0$. 
On $K \cap \supp \mu$ one has 
$0< c_K\leq \Norm{\HpG}{} \leq C_K$ by Remark~\ref{rem: non vanishing
  hamiltonian vector field} since $\supp \mu \subset \Char p$ by
Assumption~\ref{assumption: first properties of the semiclassical measure}. One has $\HpG(\y^0) \in T_{\y^0} \Gb \subset T_{\y^0}  \pTL$. On $\pTL$, by performing a rotation and a dilation of
scale factor $\Norm{\HpG(\y^0)}{} \in [c_K, C_K]$, one can assume that
$\HpG(\y^0)= (1,0, \dots,0)$. Since $\HpG z(\y^0) = \Hpz(\y^0) =0$,
the above transformation can be chosen  not affecting the $z$ variable.
One may thus use  $u \in \R^{2d}$ such that $(u,z)$ are 
new coordinates on $\pTL$ and $\HpG (\y^{0} ) = \d_{u_1}$, for $u=(u_1, u')$
with $u_1\in \R$ and $u'\in \R^{2d-1}$. 
With $\vartheta = \Hpz(\y)$, one can use  
$(u,z,\vartheta)$ as local  coordinates on $\TL$. Since $\y^0 \in \Gb$
the coordinates of 
$\y^0$ are of the form $(u_1^0,u^{0\prime},z=0,\vartheta=0)$.

Consider the functions $\chi$ and $\beta$ as in \eqref{eq: def
  chi}--\eqref{eq: def beta} and also a convex function $j \in \Cinf(\R)$
such  that
\begin{align*}
  j\equiv 0 \ \text{on}\  ]-\infty,1/2], 
  \quad j'>0 \ \text{on}\  ]1/2, +\infty[, 
  \quad j(s) > s \ \text{for}\  s\geq 1. 
\end{align*}
Observe that these properties imply 
\begin{align}
  \label{eq: support propery f}
  \big( a \geq 1 \ \text{and}\ j(s) \leq a\big) \ \ 
  \imp s \leq a.
\end{align}
A possible choice is simply $j(s) = \alpha s
\unitfunction{2s-1>0}e^{1/(1-2s)}$ with $\alpha>e$. 
Pick also $\psi \in \Cinfc(\R)$ nonnegative such that 
\begin{align*} 
  \psi \equiv 1 \ \text{in}\  [-R,R], \quad \psi \equiv 0  \
  \text{in}\ \R \setminus [-R-1,R+1],
\end{align*}
for some $R>0$ to be set below.
Then set 
\begin{align*} 
  q (u,z,\vartheta) = e^{A u_1}(\chi \circ v) (u,z)\, (\beta \circ w)
  (u) \, \psi(\vartheta) ,
\end{align*}
with 
\begin{align*}
  &v(u,z)= 1/2-\delta^{-1} (u_{1}- u^0_{1})
  +8 (\eps \delta)^{-2} \Norm{u'- u^{0\prime}}{}^{2} 
  + j \big( 2 (\delta\eps)^{-1} z\big),
 \intertext{and} 
  &w(u)=2\eps^{-1} \big(1  - \delta^{-1} (u_{1}-u^0_{1}) \big).  
\end{align*}
Arguing as in the proof of Proposition~\ref{prop: symboles} one finds
\begin{equation}
  \label{eq: supp q}
  - \delta /2  \leq  u_1- u^0_1 \leq \delta ( 1+ \eps/2),
\end{equation}
and 
\begin{equation*}
  8  (\eps \delta)^{-2} \Norm{u'- u^{0\prime}}{}^{2} 
  + j \big(2 (\delta\eps)^{-1} z\big) 
  \leq 3/2 + \eps/2,
\end{equation*}
in $\supp q$. The second inequality implies
\begin{equation}
  \label{eq: supp q 2}
  8 \Norm{u'- u^{0\prime}}{}^{2} \leq (\eps \delta)^{2}\big(3/2 +
  \eps/2\big)\leq 2 (\eps \delta)^{2},
\end{equation}
and 
\begin{equation}
  \label{eq: supp q 3}
  z \leq \frac{\delta\eps}{2} \big(3/2 + \eps/2\big) \leq \delta\eps ,
\end{equation}
using \eqref{eq: support propery f}.  One thus finds that $\supp q
\cap \{ z \geq 0\}$ is compact. Since $\supp \mu \subset \{ z \geq
0\}$ by Assumption~\ref{assumption: first properties of the
  semiclassical measure}, then the action of $\transp \Hp \mu - f\mu$
on the test function $q$ makes sense.

In $\Char p$, with $|u-u^0|$ and $z$ bounded as above, $|\vartheta|$
is also bounded by means of Lemma~\ref{lemma: zeta localizaton - char
  p -G}. Thus, for $R>0$ chosen \suff large one finds that
\begin{align}
  \label{eq: flatness psi in supp mu}
  \psi \equiv 1\  \text{in}\
  \text{a \nhd of the $\vartheta$-projection of}\ 
  \supp q \cap \Char p.
\end{align}

Applying the measure equation~\eqref{eq: Gerard-Leichtnam equation} of Assumption~\ref{assumption: Gerard-Leichtnam equation} to $q$, one first considers the contribution of the \rhs of
\eqref{eq: Gerard-Leichtnam equation} 
associated with hyperbolic points. 
Consider $\y \in \pHb$. One has $\{ \y^+, \y^-\} = \ppi^{-1}\{\y\}\cap \Char p$, and the two
points only differ by their $\vartheta$-component: $\Hpz(\y^+) = -\Hpz(\y^-)$; see Lemma~\ref{lemma:
  relevement pG pH}, Proposition~\ref{prop: relevement H G} and
\eqref{eq: decomposition point in char p}. With the form of $\psi$ and
the value of $R$ chosen above, one has 
\begin{align*}
  q\big( \y^+\big) 
  = q\big(\y^-\big) 
\end{align*}
 by \eqref{eq: flatness psi in supp mu}.
Hence, the contribution of $\pHb$ to the integral in
\eqref{eq: Gerard-Leichtnam equation} vanishes with the present choice
of test function.  

\medskip Second,  one considers the contribution of the \rhs of
\eqref{eq: Gerard-Leichtnam equation} 
associated with glancing points. 
Consider $\y \in \pGb$.  The contribution to
the integrand is described in \eqref{eq: remark: integrand GL equation on G}, that is, a differentiation of the test function with
respect to $\vartheta$. Since 
\begin{align*}
  \d_{\vartheta} q (u,z,\vartheta) = e^{A u_1}(\chi \circ v) (u,z)\, (\beta \circ w)
  (u) \, \psi'(\vartheta) 
\end{align*}
with \eqref{eq: flatness psi in supp
  mu} one finds that that $\d_{\vartheta} q$
vanishes in $\pGb = \Gb \subset \Char p$. 
One concludes that
\begin{align}
  \label{eq: t Hp mu=0}
  \dup{\transp{\Hp} \mu}{q} = \dup{\mu}{\Hp q} =\dup{\mu}{f q}.
\end{align}

\medskip
One has $\Hp q = g + h+ j + A (\Hp u_1)q$, with 
\begin{align*} 
  &g (u,z,\vartheta) = e^{A u_1}(\chi' \circ v) (u,z)\, (\beta \circ w)
  (u) \, \psi(\vartheta) \Hp v (u,z,\vartheta),\\
  &h (u,z,\vartheta) = e^{A u_1}(\chi \circ v) (u,z)\, (\beta' \circ w)
  (u) \, \psi(\vartheta) \Hp w (u,z,\vartheta),\\
  &j (u,z,\vartheta) = e^{A u_1}(\chi \circ v) (u,z)\, (\beta \circ w)
  (u) \, \Hp \psi (u,z,\vartheta).
\end{align*}
By \eqref{eq: flatness psi in supp mu}, as $\supp \mu \subset
\Char p$ by Assumption~\ref{assumption: first properties of the
  semiclassical measure},  one has
\begin{align}
  \label{eq: mu vanish on j}
  \dup{\mu}{j }=0.
\end{align} 
The support properties \eqref{eq: supp q}--\eqref{eq: supp q 3}
naturally hold also for $\supp\big( g\big)$ and
$\supp h$. Moreover, arguing as in the 
proof of Proposition~\ref{prop: symboles}, for $0< \eps\leq 1$, one has
\begin{equation}
  \label{eq: supp h - G}
   \y = (u, z,\vartheta) \in \supp h \ \ \imp \ \ 
   u \in  B\big(u^0 + \delta \HpG(\y^0), \eps \delta\big)
   \ \ \text{and} \ \ z \leq \eps \delta.
\end{equation}
Recall that $\HpG(\y^0) = \d_{u_1}$.

%%%%%%%%%%%%%%%%%%%%%%%%
% lemma                %
%%%%%%%%%%%%%%%%%%%%%%%%
\begin{lemma}
  \label{lem: positivite g - G}
  For any $0 < \eps \leq 1$ there exists $\delta_0>0$ such that for
  any $\y^0\in K\cap \Gb \cap \supp \mu$ and $0 < \delta \leq \delta_0$
  \begin{enumerate}
    \item the function
  $g$ is nonnegative in $\supp \mu$
  and is positive in a \nhd
  of $\y^0$.
  \item $\Hp u_1 \geq 1/2$ in $\supp q \cap \supp \mu$.
  \end{enumerate}
\end{lemma}
Consider  $\delta_0>0$ as given by Lemma~\ref{lem: positivite g -
  G} and $0 < \delta \leq \delta_0$. By \eqref{eq: t Hp mu=0} and
\eqref{eq: mu vanish on j} one finds
\begin{equation*}%\label{eq: somme-G}
  0= \dup{ (\transp{\Hp} -f) \mu}{q}
  = \dup{\mu}{g}
  +  \dup{\mu}{h}
  + \bigdup{\mu}{\big(A (\Hp u_1)-f\big)q}.
\end{equation*}
One has $(A (\Hp u_1) - f)q \geq 0$ in $\supp(\mu)$ for $A \geq
2\sup_K |f|$, implying
\begin{align*}
  \bigdup{\mu}{\big(A (\Hp u_1)-f\big)q}\geq 0.
\end{align*}
With $g\geq 0$, $\y^0 \in \supp \mu$ and $g(\y^0)>0$ one obtains  $\dup{\mu}{h}\neq 0$, meaning that $\supp \mu \cap \supp h \neq \emptyset$. With \eqref{eq: supp h - G} one 
  concludes the proof of Proposition~\ref{prop: propagation glancing region}.
\end{proof}

The proof of Lemma~\ref{lem: positivite g - G} is very close to that
of Lemma~\ref{lem: positivite g}. However, some details need to be handled
carefully. We thus provide a complete proof 
\begin{proof}[Proof of Lemma~\ref{lem: positivite g - G}]
  Consider $0 < \eps \leq 1$.  One has
  \begin{align*}
    g (u,z,\vartheta)
    = e^{A u_1}(\chi' \circ v) (u,z)\
    (\beta \circ w) (u)\,  \psi(\vartheta) \Hp v (u,z,\vartheta).
  \end{align*}
  Since
  $\beta \geq 0$, $\chi'<0$, and $\psi \geq 0$,  it suffices to prove
  that $\Hp v (u,z,\vartheta) \leq 0$
  for $(u,z,\vartheta)$ in $\supp q \cap \supp \mu$ for $\delta>0$
  chosen sufficiently small, uniformly with respect to $\y^0 \in K$.
  Since $v$ is independent of $\vartheta$ and $\Hz = - \d_{\zeta} = -
  \Hzzp\, \d_\vartheta$  one has $\Hp v = \HpG v$ by \eqref{eq: change of variables xid zeta}; see
  \eqref{eq: definition HpG-extended}.

  If $\y \in \supp g \cap \supp \mu$
  then $\Norm{\py - \py^{0}}{}\lesssim \delta$ by \eqref{eq: supp
    q}--\eqref{eq: supp q 3} with in particular
  $0 \leq z \lesssim \delta$ since $z\geq 0$ in $\supp \mu$ by
  Assumption~\ref{assumption: first properties of the
  semiclassical measure}.  With Lemma~\ref{lemma: zeta localizaton - char p -G}
  one has $|\vartheta| \lesssim \Norm{\py - \py^{0}}{}^{1/2}$. Thus $\vartheta= o(1)$ as $\delta \to 0^+$ implying $\Norm{\y  -
  \y^0}{}= o(1)$.

  Write 
  \begin{align*}
    \HpG(\y)  - \HpG(\y^0) \in  \alpha^1(\y, \y^0) \d_{u_1} + \alpha'(\y,
    \y^0) \cdot \nabla_{u'} + \gamma(\y, \y^0)\d_z + \Span \{ \d_\vartheta\},
  \end{align*}
  with $\alpha^1(\y, \y^0) \in \R$, $\alpha'(\y, \y^0) \in \R^{2d-1}$,
  and  $\gamma(\y, \y^0)\in \R$. From the uniform continuity of
  $\HpG$ in any compact set one concludes that 
  \begin{align}
    \label{eq: smallness HpG - HpG0}
    |\alpha^1(\y, \y^0) | + \|  \alpha'(\y, \y^0) \| + |\gamma(\y, \y^0)|= o(1) \ \
    \text{as} \  \delta \to 0^+,
  \end{align}
  uniformly with respect to $\y^0\in K\cap \Gb \cap \supp \mu$ and $\y \in \supp q  \cap
  \supp \mu$. Using that $\HpG(\y^{0}) =
  \d_{u_1}$ and the form of $v$, one writes
  \begin{align*}
    \HpG v(\y) 
    &= \big(\d_{u_1} v + (\HpG(\y)  - \HpG(\y^0) ) v\big) (\y) = -
      \delta^{-1} + r
  \end{align*}
  with 
   $r =  - \delta^{-1}\alpha^1(\y, \y^0) 
      +16 (\eps \delta)^{-2} \alpha'(\y, \y^0)\cdot (u' -
      u^{0\prime})
      + 2 (\delta\eps)^{-1}\gamma(\y, \y^0) j' \big( 2 (\delta\eps)^{-1} z\big)
      $.
  With \eqref{eq: supp q 2} and \eqref{eq: supp q 3} one finds
  \begin{align*}
    |r|  \lesssim \delta^{-1} \big| \alpha^1(\y, \y^0) \big|
    +  (\delta\eps)^{-1}\| \alpha'(\y, \y^0)  \| +  (\delta\eps)^{-1} |\gamma(\y, \y^0)|,
  \end{align*}
  using that $|j'(s)|\lesssim 1$ if $s\leq 2$.
  With $\eps$ fixed above and with \eqref{eq: smallness HpG - HpG0} one
  has $r = \delta^{-1}o(1)$. One thus  finds
$\HpG v (\y) \sim - \delta^{-1}$ as $\delta \to 0^+$
  uniformly with respect to $\y^0\in K\cap \Gb \cap \supp \mu$ and
  $\y \in \supp q \cap \supp \mu$.

\medskip
  One has
  $g(\y^0) = - \delta^{-1}\chi'(1/2) \beta(2
  \eps^{-1})>0$ and thus $g$ is positive in a \nhd
  of $\y^0$. 

  \medskip
  Note that $\Hp u_1 = \HpG u_1 = 1 + \alpha'(\y,\y^0)$. With
  \eqref{eq: smallness HpG - HpG0} one finds tht $\Hp u_1 \geq 1/2$
  for $\delta$ \suff small, uniformly in with respect to $\y^0\in
  K\cap \Gb \cap \supp \mu$ and $\y \in \supp q \cap \supp \mu$.
\end{proof}

\subsubsection{Diffractive points}
\label{sec: diffractive points}
For points in $\sdGb$ we rely on the following result.
%%%%%%%%%%%%%%%%%%%%%%%%
% proposition          %
%%%%%%%%%%%%%%%%%%%%%%%%
\begin{proposition}
  \label{prop: propagation diffractive points}
  Suppose $\y^0  \in \sdGb \cap \supp \mu$. There exist $C>0$ and
  $\delta_0>0$ such that for all $0 < \delta \leq \delta_0$ there
  exist $\y_p \in \supp \mu \cap \H^-$ and $\y_f \in \supp \mu \cap \H^+$ 
  such that 
  \begin{multline*}
    \Norm{\py_p - \py^{0}}{} 
    +   \Norm{\py_f - \py^{0}}{}
    \leq C \delta,
    \quad- C \delta^{1/2}  \leq \Hpz (\y_p) \leq - \delta^2,\\
    \text{and} \ \ \delta^2  \leq \Hpz (\y_f) \leq  C\delta^{1/2}.
  \end{multline*}
\end{proposition}
In the statement, the subscript $p$ stands for {\em past} and $f$ for {\em
  future}, in view of the use we make of these two points in the
construction of a \bichar contained in $\supp \mu$ that goes through
a point $\y^0 \in \sdGb \cap \supp \mu$ (see Section~\ref{sec: proof theorem: measure support propagation}).

%%%% proof of proposition
\begin{proof}[Proof of Proposition~\ref{prop: propagation diffractive points}]
  If $\y^0 = (t^0,x^{0\prime},z^0,\tau^0, \xi^0)$ one 
  has $z^0=0$. 
  Introduce the following two functions
  \begin{align*}
    \phi_0(\y) = \Norm{\py-\py^{0}}{}^2 
    \quad \phi(\y) = - \Hpz(\y) + \phi_0(\y).
  \end{align*}
  Note that $\phi_0$ is independent of the variable $\vartheta = \Hpz(\y)$. 
  
  Here, we write the proof of the existence of the point $\y_f$ as in
  the statement of the proposition. For the existence of the point $\y_p$ one simply
  changes $\phi(\y)$ into $\Hpz(\y) + \phi_0(\y)$ and the proof follows
  {\em mutatis mutandis}. 

  One has $\Hp \phi(\y) = - \Hppz(\y) + \Hp \phi_0(\y)$. Since
  $|\Hp \phi_0(\y)| \lesssim \Norm{\py-\py^{0}}{}$ and
  since $\Hppz(\y^0)>0$ as $\y^0 \in \sdGb$, there exist a \nhd $V^0$
  of $\y^0$ in $\TL$ and $C_0>0$ such that 
  \begin{align}
    \label{eq: small nhd - diffractive}
    \Hp \phi(\y) \leq - C_0<0, \qquad \y \in V^0.
  \end{align}
  
  Consider $\chi  \in \Cinf(\R)$  as given by \eqref{eq: def chi}
%\begin{equation}
%  \label{eq: def tilde chi}
%  \tchi (s) = \bld{1}_{s<1}\, \exp (1/(s-1)) \exp(-s),
%\end{equation} 
and suppose 
  $\psi \in \Cinfc(\R)$ nonnegative, with $\supp \psi \subset [-3,3]$
  and  such that  $\psi  \equiv 1$  on  a \nhd  of  $[-2,2]$. Introduce the following family of test functions, $\delta>0$,
  \begin{align*}
    %\label{eq: test function - diffractive}
    q(\y) = \psi^2(\phi_0(\y)/\delta^2) \ \chi(\phi(\y)/\delta^2)
    \exp(- \phi(\y)/\delta^2).
  \end{align*}
  First, consider the support of $q$. 
  In $\supp q$ one has 
  \begin{align}
    \label{eq: support of q - diffractive}
    \Norm{\py-\py^{0}}{}^2 \leq 3 \delta^2, \qquad \phi(\y)\leq \delta^2.
  \end{align}
  In $\supp q \cap \Char p$, by
  Lemma~\ref{lemma: zeta localizaton - char p -G}, from \eqref{eq:
    support of q - diffractive} one has 
  $|\Hpz(\y)|\lesssim \delta^{1/2}$. Consequently, for $\delta$
  chosen \suff small one has
  $\supp q \cap \Char p   \subset V^0$,
   meaning that
   \begin{align}
     \label{eq: support of q - diffractive 2}
     \Hp \phi(\y) \leq - C_0<0
     \quad  \text{if}\ \y \in \supp q \cap
     \Char p,
   \end{align}
   by \eqref{eq: small nhd - diffractive}.

   Second, compute $\Hp q = h_1 + h_2 + h_3$ with 
   \begin{align*}
     %\label{eq: Hp q- diffractive-h1}
     &h_1(\y) = 2 \delta^{-2} \big( \psi' \psi\big)
     (\phi_0(\y)/\delta^2)\
     \chi(\phi(\y)/\delta^2)\ \exp(- \phi(\y)/\delta^2)
     \Hp \phi_0(\y),\\
     %\label{eq: Hp q- diffractive-h2}
     &h_2 (\y) =  \delta^{-2} \psi^2(\phi_0(\y)/\delta^2) \ 
       \chi'(\phi(\y)/\delta^2)\ \exp(- \phi(\y)/\delta^2)
       \Hp \phi(\y),\\
       &h_3 (\y) =  - \delta^{-2} q(\y)
       \Hp \phi(\y).
   \end{align*}
   Consider the support of $h_1$. Since
   $\supp \psi' \cap \R_+ \subset [2,3]$, one has
   $2 \delta^2 \leq \phi_0(\y)\leq 3 \delta^2$ in $\supp h_1$.  As $\phi(\y)\leq \delta^2$
   by \eqref{eq: support of q - diffractive}, one concludes
   that
   \begin{align*}
    \Hpz (\y) = \phi_0(\y) - \phi(\y)\geq \delta^2 \ \
     \text{if}\ \y \in \supp\big(h_1).
  \end{align*}

  \bigskip With the functions $h_1$ and $h_2$ defined above one has
  \begin{align}
    \label{eq: decompositon <mu, Hp q> - diffractive}
    \dup{\transp{\Hp} \mu - f \mu}{q}=
    \dup{\mu}{h_1} +\dup{\mu}{h_2} +\dup{\mu}{h_3 - f q}.
  \end{align}
  Below, we prove the following lemma.
  %%%%%%%%%%%%%%%%%%%%%%%%
  % lemma                %
  %%%%%%%%%%%%%%%%%%%%%%%% 
  \begin{lemma}
    \label{lemma: sign Hp q - diffractive}
    One has $\dup{\transp{\Hp} \mu - f \mu}{q} \leq 0$,  $\dup{\mu}{h_2} > 0$, and $\dup{\mu}{h_3 - f q}\geq 0$ for $\delta$ chosen \suff small. 
  \end{lemma}
  From \eqref{eq: decompositon <mu, Hp q> - diffractive} and
  this lemma one concludes that $\dup{\mu}{h_1} < 0$, meaning 
   that there exists $\y_f\in \supp \mu \cap \supp h_1$.
  The above analysis yields $\y_f$  as in the statement of the proposition.
\end{proof}
%%%% proof of lemma
  \begin{proof}[Proof of Lemma~\ref{lemma: sign Hp q - diffractive}]

     First, consider the action of $\mu$ on $q$.
     With \eqref{eq: Gerard-Leichtnam equation} in Assumption~\ref{assumption: Gerard-Leichtnam equation} one has
  \begin{align*}
    \dup{\transp{\Hp} \mu -f \mu}{q} = - \int_{\y \in \pHb \cup \pGb} 
    \frac{q (\y^+) - q (\y^-)}
    {\dup{\xi^+- \xi^-}{\nx}_{T_x^*\M, T_x\M} 
    }\  d \nu (\y).
  \end{align*}
   Since the measure $\nu$ is nonnegative and since $\ovl{\pHb} = \pHb
   \cup \pGb$ it suffices to prove that the integrand is nonnegative
   on $\pHb$ to conclude that $\dup{\transp{\Hp} \mu -f \mu}{q}\leq 0$.
   
   Consider $\y \in \pHb$. By Lemma~\ref{lemma: relevement pG pH} one has
   \begin{align*}
     \dup{\xi^+- \xi^-}{\nx}_{T_x^*\M, T_x\M} 
     =  2 \lambda>0.
   \end{align*}
   Next, as $\y^\pm\in \Hb^\pm$, one
   has $\Hpz(\y^+) - \Hpz(\y^-) >0$ from Definition~\ref{def: G H}. As $\phi_0 (\y^+)=
   \phi_0(\y^-) = \phi_0(\y)$ 
   one finds
   \begin{align*}
     \phi(\y^+)- \phi(\y^-) 
     = - \Hpz(\y^+) + \Hpz(\y^-) <0.
   \end{align*}
   Using now that $s \mapsto \tchi(s) = \chi(s) \exp(-s)$ is a nonincreasing function, 
  one obtains
   \begin{align*}
     q(\y^+) - q (\y^-) 
     = \psi^2(\phi_0(\y)/\delta^2)
     \Big(
     \tchi(\phi(\y^+)/\delta^2) - \tchi(\phi(\y^-)/\delta^2)     
     \Big) \geq 0,
   \end{align*}
   implying that $\dup{\transp{\Hp}\mu -f \mu}{q}\leq 0$.

   \medskip
    Second, consider the action of $\mu$ on $h_2$.
    Since $\chi'\leq 0$ and $\Hp \phi \leq 0$ in $\supp q\cap
     \Char p$  by \eqref{eq: support of q - diffractive 2}
     one finds that $h_2 \geq 0$. 
     Consider now the value of $h_2$ at $\y^0$:
     \begin{align*}
       h_2 (\y^0) &=  \delta^{-2} \psi^2(\phi_0(\y^0)/\delta^2) \ 
       \chi'(\phi(\y^0)/\delta^2)\ \exp(\phi(\y^0)/\delta^2)\, 
       \Hp \phi(\y^0)\\
       &= \delta^{-2} \psi^2(0) \chi'(0) \Hp \phi(\y^0).
     \end{align*}
     Since $\psi(0)=1$, $\chi'(0) <0$ and $\Hp \phi(\y^0) <0$ by
     \eqref{eq: support of q - diffractive 2} one
     obtains $h_2 (\y^0)>0$. Since $\y^0 \in \supp \mu$ and
     $\mu$ is nonnegative, one
     concludes that $\dup{\mu}{h_2} > 0$.

     One has $h_3 - f q = (- \delta^{-2} \Hp \phi(\y)-f) q$. With
     \eqref{eq: support of q - diffractive 2} one finds that $-
     \delta^{-2} \Hp \phi(\y)-f \geq 0$ in $\supp \mu$ for $\delta$
     chosen \suff small. This gives $\dup{\mu}{h_3 - f q}\geq 0$.
  \end{proof}

  \bigskip With Proposition~\ref{prop: propagation diffractive points}
  one has the following result.
  %%%%%%%%%%%%%%%%%%%%%%%%
  % propositoin              %
  %%%%%%%%%%%%%%%%%%%%%%%% 
  \begin{proposition}
    \label{prop: bichar diffractive point}
    Suppose $\y^0 \in \sdGb \cap \supp \mu$. There exist $S>0$ and a
    local \bichar above $\hL$, $\gamma: [-S,S] \to \ThL$, such that $\gamma(0) =
    \y^0$, $\gamma(s) \in \TL \setminus \dTL$ for $s \neq 0$,
    and 
    \begin{align*}
      \Gamma = \{ \gamma(s); \ s \in [-S,S] \} \subset \supp \mu.
    \end{align*}
  \end{proposition}
  %%%% proof of proposition
  \begin{proof}
    In a local chart consider a ball $\mathcal B$ of radius $R$ centered at $\y^0$ where
    $\Hppz \geq C>0$. There exist $c_0, C_0$ both positive such that 
    \begin{align*}
      c_0 \leq \Norm{\Hp(\y)}{}\leq C_0, \quad \y \in \mathcal B.
    \end{align*}
    For $\delta>0$ chosen small, by
    Proposition~\ref{prop: propagation diffractive points}
    there exists
    $\y_\delta \in \mathcal B
    \cap\supp \mu \cap \H^+$ such that
    \begin{align*}
    \Norm{\py_\delta- \py^{0\prime}}{}
    \lesssim \delta, \qquad
    \delta^2 \lesssim \Hpz(\y)  \lesssim \delta^{1/2}.
    \end{align*}
   (This point is
    denoted $\y_f$ in Proposition~\ref{prop: propagation diffractive points}.)
    Consequently, for some $C_1>0$ one has 
    \begin{align*}
    \Norm{\y_\delta- \y^0}{} \leq C_1 \delta^{1/2}.
    \end{align*}
    One either has $\y_\delta \in \Hb^+$ or
    $\y_\delta \in \H^+\setminus \dTL$.  In either case, by
    Propositions~\ref{prop: propagation Away from the boundary} and
    \ref{prop: propagation from rho+, rho-} there exists a maximally
    extended \bichar $\gamma_\delta$ in $\mathcal B \cap \TL$, defined on a
    interval of the form $]S_\delta^-, S_\delta^+[$ that is moreover
    contained in $\supp \mu$ and such that
    $\gamma_\delta(0) = \y_\delta$.  We only consider this
    \bichar on the interval $[0,S_\delta^+[$. Since $\Hppz >0$ in
    $\mathcal B$, then $\Hpz(\gamma_\delta(s))$ increases  as $s\in
    [0,S_\delta^+[$ increases. Since $\Hpz (\y_\delta)>0$ one finds
    that $\Hpz$ remains positive and hence $z$ increases  along
    the bicaracteristic. Consequently, if $S_\delta^+ < \infty$ then
    the maximaly extended \bichar leaves $\mathcal B \cap \TL$ for $s = S_\delta^+$, yet
    not through the boundary $\{z=0\}$, that is, $\dTL$. 

    Choose $ \delta_0>0$ such that $C_1 \delta_0^{1/2} < R/2$.  Then,
    for $0 < \delta \leq \delta_0$ one has
    $\Norm{\y_\delta - \y^0}{}< R/2$ and $\dist(\y_\delta, \d \mathcal
    B) > R/2$ implying that $S_\delta^+ > R / (2 C_0)$. Set $S = R / (2 C_0)$.
    
    One has 
    \begin{align*}
      \gamma_\delta(s) = \y_\delta 
      +  \int_0^s \Hp({\gamma_\delta}(\sigma)) d\sigma, 
      \qquad s \in [0, S].
    \end{align*}
    If one lets $\delta$ vary in  $]0,\delta_0]$, 
    the set of \bichars $\gamma_\delta(s)$ is
    equicontinuous on $[0, S]$.  For $\delta\to 0$, by the Arzel\`a-Ascoli theorem one
    can extract a subsequence that converges uniformly to a curve
    $\gamma(s)$ with $[0, S]$. From the continuity of $\Hp$, one has
    \begin{align*}
      \gamma(s) = \y^0 +  \int_0^s \Hp({\gamma}(\sigma)) d\sigma,
      \qquad s \in [0, S],
    \end{align*}
   that is, a \bichar that goes through $\y^0$ at $s=0$. Moreover,
   $\gamma(s) \in \supp \mu$ for $s \in [0, S]$ since $\supp \mu$ is
   a closed set. Arguing as
   above one has $z>0$ along this \bichar if $s>0$. 

   The argument can be used {\em mutatis mutandis} to construct the sought
   \bichar for $s \in [-S,0[$.
  \end{proof}

\subsection{Final construction of \gbichar
  in the measure support}
\label{sec: proof theorem: measure support propagation}

All possible cases listed in the beginning of Section~\ref{sec: Propagation of the measure  support} need to be considered. Cases are then used sequentially, in various orders,
depending on the \gbichar that is constructed along the proof. 

\bigskip
\noindent
{\bfseries Case 1: $\bld{\y^0}$ is away from the boundary.} if $\y^0
\in \supp \mu \cap \big(\TL \setminus \dTL\big)$ and if one sets $V=
\U \cap \TL \setminus \dTL$, Proposition~\ref{prop: propagation Away
  from the boundary} yields a maximal \bichar $\gamma(s)$, with $s \in
J= ]S_1, S_2[$ that lies in $\supp \mu\cap V$ and such that
       $\gamma(0) = \y^0$.

       One says that $\gamma$ leaves $\U$ at $s=S_2^-$ if $S_2
       <\infty$ and the limit point at $s=S_2^-$ is in $\d\U$ (see Lemma~\ref{lemma: limit of a
         maximal bichar from the interior}), with
       the same notation at $s=S_1^+$.      
       If, on the one hand, either $S_1 = -\infty$ or $\gamma$ leaves
       $\U$ at $s=S_1^+$, and, on the other hand, $S_2 = +\infty$ or
       $\gamma$ leaves $\U$ at $s=S_2^-$, then one has obtained a maximal
       \gbichar contained in $\supp \mu$ that goes through $\y^0$. If
       however, for instance $S_2 < +\infty$ and $\gamma$ does not leave $\U$ at $s=S_2^-$, then $\y^1 = \lim_{s\to
         S_2^-} \gamma(s)$ exists by Lemma~\ref{lemma: limit of a
         maximal bichar from the interior} and moreover $\y^1 \in
       \Hb^- \cup \glGb \cup \sdGb \setminus \d\U$.  Since $\supp \mu$
       is a closed subset of $\U$ one has $\y^1 \in \supp \mu$.  Note
       that if $\y^1 \in \sdGb$ one can set $\gamma(S_2) = \y^1$ and
       one has $\frac{d}{d s} \gamma (S_2) = \Hp (\y^1)$. Now if $\y^1
       \in \glGb$ one can also set $\gamma(S_2) = \y^1$ and one has
       $\frac{d}{d s} \gamma (S_2) = \Hp (\y^1) = \HpG(\y^1)$ by
       Lemma~\ref{lemma: properties HpG} as $\Hppz(\y^1)=0$.

   One may then
  consider  Cases 2, 3 and 4 below to carry on the construction of a
  \gbichar contained in $\supp \mu$, with $\y^1$ being now
  the point where one initiates this \gbichar. 

  Naturally, the same
  reasoning is applied in the case $-\infty < S_1$ and $\gamma$ does not leave $\U$ at $s=S_1^+$. Then
  $\y^1= \lim_{s\to S_1^+}  \gamma(s) \in \supp \mu \cap (\Hb^+ \cup \glGb \cup \sdGb) \setminus \d\U$. 

\bigskip  
\noindent
{\bfseries Case 2:  $\bld{\y^0}$ is a hyperbolic point.}
Consider $\y^0 \in \supp \mu\cap \Hb$. Then
$\Sigma(\y^0) \in \supp \mu \cap \Hb$ by Proposition~\ref{prop:
  propagation rho- to rho+}.  Set $\y^{0\pm} \in \H^\pm$ so as to have
$\{ \y^{0-}, \y^{0+}\} = \{\y^0, \Sigma(\y^0)\}$. By
Proposition~\ref{prop: propagation from rho+, rho-}, there exists $S>0$
and a locally defined \bbichar
$\gammaB: [-S,0[\cup ]0,S] \to \TL \setminus \dTL$ contained in
$\supp \mu$ such that
\begin{align*}
	\lim_{s\to 0^+}\gammaB(s) = \y^{0+} 
	\ \  \text{and} \ \ 
	\lim_{s\to 0^-}\gammaB(s)  = \y^{0-}.
\end{align*} 
Moreover $\gammaB(-S), \gammaB(S) \in \supp \mu\cap \TL \setminus \dTL$,
that is, the constructed local \bbichar yields endpoints away from the boundary. To
carry on with the construction of the \gbichar, one needs then to consider Case~1.

\bigskip 
\noindent
{\bfseries Case 3:  $\bld{\y^0}$ is a diffractive point.} Consider 
  $\y^0 \in \supp \mu\cap \sdGb$.  One may then apply
  Proposition~\ref{prop: bichar diffractive point}: there
  exist $S>0$ and a local \bichar above $\hL$, $\gamma: [-S,S] \to \TL$, contained
  in $\supp \mu$ such that $\gamma(0) = \y^0$,
  $\gamma(s) \in \TL \setminus \dTL$ for $s \neq 0$. In
  particular, $\gamma(-S)$ and $\gamma(S)$ are in $\supp \mu \cap \TL\setminus \dTL$.  To
  carry on with the construction of the \gbichar, one needs
  then to consider Case~1.

\bigskip 
\noindent
{\bfseries Intermezzo:  Construction of a maximal \bbichar.}
  Observe that the construction proposed up to this point may imply going back and
forth between Cases~1, 2 and 3 if no point in $\supp \mu \cap \glGb$
is reached. One then obtains a \bbichar
that lies in $\supp \mu$. One sees that the existence of such a
\bbichar yields the existence of a maximal \bbichar contained in $\supp \mu$ 
by means of classical arguments; see for example \cite{Demailly:96}. 
%%%%%%%%%%%%%%%%%%%%%%%%
% proposition          %
%%%%%%%%%%%%%%%%%%%%%%%%
\begin{proposition}
  \label{prop: maximal broken bichar}
  Suppose  $\y^0 \in \supp \mu \cap \TL \setminus (\glGb \cup
  \sgGb)$. There exists a maximal \bbichar  $s \mapsto
  \gammaB(s)$ defined for $s \in J \setminus B$ with $0 \in J = ]S_1, S_2[$
  and $B$ as in Definition~\ref{def: broken bichar} and such that 
  \begin{enumerate}
  \item if $\y^0 \notin \Hb$ then $0\notin B$ and $\gammaB(0) = \y^0$;
    \item if $\y^0 \in \Hb^+$ (\resp $\Hb^-$) then $0\in B$ and $\gammaB(0^+) =
      \lim_{s\to 0^+} \gammaB(s) = \y^0$ (\resp $\gammaB(0^-) =
      \lim_{s\to 0^-} \gammaB(s) = \y^0$);
    \item for all $s \in J \setminus B$, $\gammaB(s) \in \supp \mu$;
      \item for all $S \in B$, $\gammaB(S^\pm) =
      \lim_{s\to S^\pm} \gammaB(s) \in \supp \mu\cap \Hb^\pm$.
\end{enumerate}
\end{proposition}
If, on the one hand, either $S_1 = -\infty$ or $\gamma$ leaves
       $\U$ at $s=S_1^+$, and, on the other hand, $S_2 = +\infty$ or
       $\gamma$ leaves $\U$ at $s=S_2^-$, then one has obtained
  a  maximal \gbichar contained in $\supp \mu$
  that goes through $\y^0$. If $\U = T^*\hL$ then $J=\R$.

  \medskip
  Consider now for instance the case $S_2 < +\infty$ and $\gamma$ does not leave $\U$ at $s=S_2^-$. If
  $S_2 \notin \ovl{B}$ then $\gammaB$ is a maximal \bichar near $S_2$ implying
  that $\y^1= \lim_{s \to S_2^-} \gammaB(s)$ exists and
  belongs to $\supp \mu \cap (\Hb^- \cup \glGb \cup \sdGb)$ by Lemma~\ref{lemma: limit
    of a maximal bichar from the interior}. One can discard having $\y^1 \in
  \Hb^- \cup \sdGb$ since  Cases~2 and 3 above allow one to further
  extend the \bbichar contained in $\supp \mu$ contradicting its maximality. Thus, if
  $S_2 \notin \ovl{B}$ one has $\y^1 \in \supp \mu \cap \glGb$.

  If now $S_2 \in \ovl{B}$ then 
  $\y^1 = \lim_{s\to S_2^-} \gammaB(s)$ exists by Lemma~\ref{lemma:
    broken bichar - closure B} with moreover
  $\y^1 \in \supp \mu \cap (\glGb \cup \sgGb)$. 

  Starting from $\y^1 \in \supp \mu \cap (\glGb \cup \sgGb)$ obtained in either cases, one
  now considers Case~4 below to carry on
  the construction of a \gbichar contained in $\supp \mu$,
  with $\y^1$ being now the point where one initiates the \gbichar.

  Naturally, the same
  reasoning is applied in the case $-\infty < S_1$ and $\gamma$ does not leave $\U$ at $s=S_1^+$. Then,
  $\y^1= \lim_{s\to S_1^+}  \gammaB(s) \in \supp \mu \cap (\glGb \cup \sgGb)$. 

\bigskip
\noindent
{\bfseries Case~4: $\bld{\y^0}$ is an order-3-glancing point or a gliding point.}

The following lemma gives the existence of local \gbichar
that goes through a point in $\supp \mu \cap (\glGb \cup \sgGb)$.
  %%%%%%%%%%%%%%%%%%%%%%%%       
  %       lemma                %
  %%%%%%%%%%%%%%%%%%%%%%%%       
        \begin{lemma}
          \label{lemma: construction generalized bichar near glancing
            point}
          There exists $S_0>0$ such that 
         for any $\y^0 = (t^0,x^0, \tau^0, \xi^0) \in \supp \mu \cap
          (\glGb \cup \sgGb)$ there is a \gbichar 
          $\gammaG: J\setminus B   \to \TL$ where $J = [-S, S]$, with
          $S = S_0 / |\tau^0|$, and $B$ a discrete subset of $J$, and such that 
          $\gammaG(0) = \y^0$ and $\GammaG \subset \supp \mu$.
        \end{lemma}
        We recall that the notation $\GammaG$ is intoduced in Definition~\ref{def: generalized bichar 2-intro}.
        The proof of Lemma~\ref{lemma: construction generalized bichar
          near glancing point}  is quite lengthy. We thus rather first conclude the
        proof of Theorem~\ref{theorem: measure support propagation}
        and below, in Section~\ref{sec: Local construction of a
          generalized bichar}, we proceed with the proof of
        Lemma~\ref{lemma: construction generalized bichar near
          glancing point}.

\bigskip
\noindent {\bfseries Conclusion of the construction of a maximal 
  \gbichar.}
 
With the four cases treated above, given
$\y^0 \in \supp \mu \subset \U \cap \TL $ there exists a local \gbichar that goes through this point (with the understanding of a limit
if the point $\y^0 \in \Hb$) and is contained in $\supp \mu$. This
yields the existence of a maximal \gbichar with the same
properties by means of classical arguments; see for example
\cite{Demailly:96}.

Suppose now that $\U = T^*\hL$ and 
that  $\gammaG(s)$ is such a maximal \gbichar defined for
$s \in J \setminus B$ with $0 \in J = ]S_1, S_2[$ and $B$ as in
Definition~\ref{def: generalized bichar}. Suppose that $S_2 <
+\infty$. Then Lemma~\ref{lemma: generalized bichar - derivative
  endpoint} implies that the limit
  \begin{align*}
    \y^1 = \lim_{{s\to S_2} \atop {s \in J \setminus B}} \gammaG(s)
  \end{align*}
exists and is in  $\supp \mu$ as it is a closed set  and moreover 
\begin{align*}
  \frac{d}{d s} \gammaG(S_2^-) = \XG(\y^1).
  \end{align*}
Yet, with the above argument there exists a local \gbichar that goes through $\y^1$ allowing one to extend $\gammaG(s)$
for $s > S_2$ contradicting its maximality: thus $S_2 = +\infty$. The
same reasonning gives $S_1 = - \infty$. 
This concludes the proof of
  Theorem~\ref{theorem: measure support propagation}.
  \hfill \qedsymbol \endproof

\subsection{Local construction for a
  point in $\glGb$ or $\sgGb$.}
\label{sec: Local construction of a generalized bichar}

In this section we prove Lemma~\ref{lemma: construction generalized
    bichar near glancing point}.

First, we construct locally a continuous curve $\gammaG$ that goes through
$\y^0$ and is contained in $\supp \mu$ and, second, we prove that it is a \gbichar.

\medskip
{\bfseries Local setting}

For any $x \in \d\M$ and $R>0$ consider the closed Riemannian ball $B_g(x,R)
= \{ \tilde{x}\in \M; \
\dist_g(x,\tilde{x}) \leq R\}$. In fact, since $\M$ is compact, one can choose $R_0>0$ \suff small so 
that, for any $x\in \d\M$ there exists a local chart $(\hO,\chdiff)$ of $\M$
such that $B_g(x,R_0) \subset \hO$. 

Suppose now $\y^0 = (t^0, x^0, \tau^0, \xi^0) \in \glGb \cup \sgGb$ 
and $(\hO,\chdiff)$ is a local chart chosen as above, that is, such that $B_g(x^0, R_0) \subset \hO$. 
One has $|\xi^0|_{x^0} = |\tau^0|$ since $\y^0 \in \Char p$. 
If a \gbichar $\gammaG(s) = (t(s), x(s), \tau(s), \xi(s))$
going through $\y^0$ at $s=0$ is constructed then $\tau(s) = \tau^0$
and thus $|\xi(s)|_{x(s)} = |\tau^0|$ as pointed
out in Remark~\ref{rem: def generalized bichar}. 
Thus consider the following compact set
\begin{align*}
  K^0 = \big\{ \y =(t,x,\tau,\xi) \in \TL ; x \in B_g(x^0, R_0), \
  |\tau | + |\xi|_x \leq 4 |\tau^0|\big\}.
\end{align*}
The constructed \gbichar $\gammaG(s)$ will be in $K^0$ for
$|s|$ sufficiently small.  Note that ``room'' is made in $K^0$ in the
cotangent directions for an
interative process to be carried out while remaining in $K^0$. 
 
In the local coordinates associated with $(\hO,\chdiff)$ the hamilitonian
vector field $\Hp$ is as
given in \eqref{eq: Hp}. The componant acting in the $x$-directions is
given by $v(\y)= 2 g^{ij} \xi_i \d_{x_j}$. Observe that the same holds
for the gliding vector field $\HpG$ by \eqref{eq: definition
  HpG-extended}. 
Because of the form of $K^0$  there exists $C_0>0$ such that 
\begin{align*}
 \Norm{v(\y)}{} \leq C_0 |\tau^0| \ \ \text{if} \ \y \in K^0.
\end{align*}
Thus, if for some $s \in \R$ one has $\dist_g(x^0,x(s)) \geq R_0$ then $|s| \geq R_0
/ (C_0 |\tau_0|)$. 
Thus, set  $S_{\max} = S_0 /|\tau_0|$ with $S_0 = R_0 /
(2C_0)$. The \gbichar is construced in an iterative
process for  $s \in  [-S_{\max}, S_{\max}]$ and the choice of $S_{\max}$
ensures that one remains in $K^0$ within that process.

Also, since one remains in $K^0$, the same local chart $(\hO,\chdiff)$ can be used. 
We also use local coordinates as $(u,z,\vartheta)$ as introduced in
the proof of Proposition~\ref{prop: propagation glancing region}. We
assume that it can be used in the whole local chart. This can be
assumed from the beginning by refining the atlas. Recall that $(u,z)$
provides coordinates for $\pTL$ with $u \in \R^{2d}$ and $\vartheta = \Hpz$.
Below, we will go back and forth between the $(t,x,\tau,\xi)$ and
$(u,z,\vartheta)$ coordinates. By abuse of notation we write $\y = (u,z,\vartheta)$ or
$\y=(t,x,\tau,\xi)$.
Here, our starting point is of the forms $\y^0 = (t^0, x^1, \tau^1,
\xi^1)$ and $\y^0 = (u^0,z=0,\vartheta=0)$.

For a point $\y^\ell$, $\ell \in \N$,  constructed below $(t^\ell, x^\ell, \tau^\ell,
\xi^\ell)$ and $(u^\ell, z^\ell, \vartheta^\ell)$ refer to its
coordinates in the two variable systems.

\medskip
{\bfseries Construction of $\bld{\gammaG}$}

Consider $n \in \N^\ast$ and $\eps
= 1/n$. Use $\delta_0>0$ as given by Proposition~\ref{prop:
  propagation glancing region} using $K=K^0$ therein. Set $\delta_n =
\min (\delta_0, 1/n)$.
In the local coordinates $(u,z,\vartheta)$ we construct a {\em piecewise}
continuous curve $\gamma_n$ initiated at $\y^0$.  

Consider $\y^1 \in \supp \mu$ as given by Proposition~\ref{prop:
  propagation glancing region} with $\delta = \delta_n$ therein. 
Since $\HpG \tau =0$ note that in the   $(t,x,\tau,\xi)$ coordinates
one has 
\begin{align}
  \label{eq: increase tau}
|\tau^1 -\tau^0| \leq \eps \delta_n \leq 1/n^2.
\end{align}
Now, in the $(u,z,\vartheta)$
coordinates,  one defines the following
affine curve 
\begin{align*}
  \gamma_n(s) = \y^0 + \frac{s}{\delta_n} (\y^1 - \y^0) 
  \ \ \text{for} \ s \in [0,\delta_n] \ \ \text{and} \ \ S_1 = \delta_n.
\end{align*}

One then faces two options to further construct $\gamma_n$. 
\begin{enumerate}
\item If $\y^1 \in \glGb \cup \sgGb$, like $\y^0$, one picks a second point $\y^2$
also according Proposition~\ref{prop: propagation glancing region},
yet starting from $\y^1$,  and one
further constructs $\gamma_n$ on the interval $[S_1, S_1+ \delta_n]$ in
some affine manner as above.
\begin{align*}
  \gamma_n(s) = \y^1 + \frac{s-S_1}{\delta_n} (\y^2 - \y^1) 
  \ \ \text{for} \ s \in [S_1,S_1+\delta_n] \ \ \text{and} \ \ S_2 =
  S_1 + \delta_n.
\end{align*}
One carries on with this iteration yielding points $\y^1, \dots, \y^k
\in \supp \mu \cap (\glGb \cup \sgGb)$ and $S_\ell = \ell \delta_n$, $\ell =1, \dots,
k$, until either $S_k > S_{\max}$, meaning one is done with the
construction of $\gamma_n$ for $s \in [0, S_{\max}]$, or $\y^k
\in \supp \mu \cap \TL \setminus (  \glGb \cup \sgGb)$, in which case one turns to the second
construction option just below.

Observe that $(k-1) \delta_n \leq S_{\max}$ and iterating
estimate~\eqref{eq: increase tau} one has 
\begin{align}
  \label{eq: increase tau-bis}
  | \tau^\ell -\tau^0| \leq \eps \ell \delta_n
  \leq \eps S_{\max}+ \eps \delta_n
  =S_0 / (n |\tau^0|) + 1/n^2, 
\quad \ell = 1, \dots, k.
\end{align}
This means that for $n$ chosen \suff large one has $\y^1, \dots, \y^k
\in K^0$. 

\item 
If after one or serveral steps using the above piecewise affine
construction $\y^\ell \in \supp \mu\cap \TL \setminus (\glGb \cup \sgGb)$ one constructs a
maximal \bbichar in $K^0$ initiated at $\y^\ell \in K^0$
according to Proposition~\ref{prop: maximal broken bichar}. Maximality
is only understood in the future here, that is, for $s \geq S_\ell$:
this maximal
\bbichar is defined on $[S_\ell, S_{\ell+1}[ \setminus B_{[S_\ell, S_{\ell+1}[}$,
with $B_{[S_\ell, S_{\ell+1}[}$ a discrete subset  of $[S_\ell,
S_{\ell+1}[$. It may happen that $S_\ell \in  B_{[S_\ell,
  S_{\ell+1}[}$ if $\y^\ell \in \Hb$: then the maximal \bbichar
enters $\TL \setminus \dTL$ through $\y^\ell$ if $\y^\ell \in \Hb^+$
or $\Sigma( \y^\ell)$ if  $\y^\ell \in \Hb^-$ (see Case~2 in Section~\ref{sec: proof theorem: measure support propagation}).

Two instances may occur: (a) the maximal \bbichar is such that
$S_{\ell+1} >S_{\max}$ and one is done with the construction of
$\gamma_n$, or (b) it reaches a point
$\y^{\ell+1} \in \glGb \cup \sgGb$ at $s =S_{\ell+1}$. Starting, from
this point $\y^{\ell+1}$ one reinitiates the construction with the
first option above until one reaches $s = S_{\max}^-$. Note that
$\tau$ remains constant along the \bbichar just
constructed. Hence, with \eqref{eq: increase tau-bis} one has the
estimate
\begin{align*}
 |\tau^{\ell+1} -\tau^0| = |\tau^{\ell} -\tau^0| \leq \eps \ell \delta_n.
\end{align*}
\end{enumerate}

Note that $\gamma_n$ contains at most
$\lfloor S_{\max} / \delta_n\rfloor +1$ affine pieces given by
Proposition~\ref{prop: propagation glancing region} as in item (1)
above. In particular,  simlarly to \eqref{eq: increase tau-bis} one
finds that 
\begin{align*}
  %\label{eq: increase tau-ter}
|\tau^{k} -\tau^0| \leq \eps k \delta_n 
\leq \eps S_{\max} + \eps \delta_n\leq  S_0 / (n |\tau^0|) + 1/n^2,
\end{align*}
for any $\y^k$ that is an endpoint of an
affine piece constructed above. Thus,  for $n$ chosen \suff large the constructed curve
$\gamma_n$ remains in $K^0$ as announced above.

With the above construction of the curve
$\gamma_n$ one obtains an alternating sequence of affine pieces and  maximal \bbichars.
Following a part of $\gamma_n$ made by a \bbichar, one finds an affine
part (unless that \bbichar ends the construction of $\gamma_n$). Hence, the number of \bbichars that constitutes $\gamma_n$
is also finite. Denote by $m_n$ this number. One has
$m_n \leq \lfloor S_{\max} /\delta_n\rfloor +1$.  If \bbichars
compose $\gamma_n$, that is, if $m_n \geq 1$, set
$M_n = \{1, \dots, m_n\}$ and for each \bbichar set
$B_{n,j}$, $j \in M_n$, to be the discrete set of points $s$ where the $j$th
\bbichar is discontinuous, that is, at hyperbolic points.  Each
point of $B_{n,j}$ is isolated. Yet, if $\# B_{n,j} = \infty$, points
of $B_{n,j}$ accumulate to some $s \notin B_{n,j}$. In such case, 
recall that $\gamma_n(s^-) \in \glGb \cup \sgGb$ (see Lemma~\ref{lemma: broken bichar - closure B}).  Define
$B_n = \cup_{j \in M_n} B_{n,j}$. It contains all the points where
$\gamma_n$ is discontinuous (if any) corresponding to hyperbolic
points at the boundary.

The $j$th \bbichar is defined for $s \in [S, S'[ \setminus B_{n,j}$. Set $\sigma_{n,j}^0=S$.  It may
happen that $\sigma_{n,j}^0 \in B_{n,j}$.
Index the (at most countable) other ordered
elements of $B_{n,j}$ as follows:
\begin{align*}
  \sigma_{n,j}^0< \sigma_{n,j}^1 < \sigma_{n,j}^2 <  \cdots < \sigma_{n,j}^\ell
  < \cdots,
\end{align*}
with $1\leq \ell \leq L_{n,j} = \# B_{n,j}$. 
If $L_{n,j}= +\infty$ set 
\begin{align*}
  \sigma_{n,j}^\infty = \sup_{0 \leq \ell \leq L_{n,j}} \sigma_{n,j}^\ell.
\end{align*}
If $L_{n,j}< \infty$ set also $\sigma_{n,j}^{L_{n,j}+1} = \sigma_{n,j}^\infty$ to be the
value $s$ of the endpoint of the $j$th \bbichar. Using the index value
$L_{n,j}+1$ can be useful in summations in what follows.
With the
maximal \bbichar defined on $[S,S'[\setminus B_{n,j}$ one has
$\sigma_{n,j}^\infty=S'$.
If $j < m_n$ then
$\sigma_{n,j}^\infty < S_{\max}$. Note however that one may have $\sigma_{n,m_n}^\infty > S_{\max}$ since the
$m_n$th maximal \bbichar may carry one beyond $s = S_{\max}$. 
One has
\begin{multline*}
  \sigma_{n,1}^0 <\sigma_{n,1}^1 < \cdots < \sigma_{n,1}^\infty 
  < \sigma_{n,2}^0 < \cdots < \sigma_{n,2}^{\infty}< \sigma_{n,3}^0
  < \cdots \\
  \cdots < \sigma_{n,m_n-1}^{\infty} < \sigma_{n,m_n}^0 < \cdots <
  \sigma_{n,m_n}^{\infty}.
\end{multline*}

As mentionned above one has 
\begin{align*}
  \gamma_n\big(\sigma_{n,j}^\infty\big) \in \glGb \cup \sgGb, 
  j=1, \dots, m_n.
\end{align*}
For $0 \leq s_1 \leq s_2 \leq S_{\max}$ set 
\begin{align*}
  M_n^{[s_1,s_2]} = \big\{ j \in M_n; \ 
  s_1 < \sigma_{n,j}^\infty \ \text{and} \  \sigma_{n,j}^0 < s_2\big\},
\end{align*}
meaning that on the interval $[s_1,s_2]$ one encounters the $j$th
\bbichars for $j \in M_n^{[s_1,s_2]}$. Note that $M_n^{[s_1,s_2]}$ may be empty.

The curve $\gamma_n$ is differentiable away from endpoints of affine
parts and away from points in $B_n$. Write
$\gamma_n(s) = (\pgamma_n(s),\vartheta_n(s))$ using the notation of
Section~\ref{sec: local extension away from boundary}, and using the $(u,z)$ coordinates for
$\pgamma_n(s) \in \pTL$.  

Now set 
\begin{align*}
%\label{eq: def cgammaG}
  \cgamma_n(s) = \begin{cases}
    \cphi  \big(\gamma_n(s)\big) 
    = \gamma_n(s) & \text{if} \  s \notin B_n\\
     \cphi \big( \lim_{s' \to s^-} \gamma_n(s)\big) 
     = \cphi \big( \lim_{s' \to s^+} \gamma_n(s)\big) & \text{if} \  s \in B_n.
  \end{cases}
\end{align*}
The map $\cphi $ is defined in Section~\ref{sec: The compressed cotangent bundle}. 
The curve $\cgamma_n(s)$ is continuous with values in the compressed
cotangent bundle $\cTL$.
Our next goal is to prove the
equicontinuity of $s \mapsto\cgamma_n(s) = \cphi\big(\gamma_n(s)\big)$
for $n \in \N^*$ chosen and $s \in [0,S_{\max}]$.

\medskip
On a piece of $\gamma_n$ given by a maximal \bichar (within a \bbichar), say on $[s_1,s_2]$, one has
\begin{align}
  \label{eq: gen bichar - Ascoli 1}
  \gamma_n(s) = \gamma_n(s_1^+) 
  + \int_{s_1}^s \Hp \big(\gamma_n (\sigma)\big)\, d \sigma.
\end{align}
This yields
\begin{align}
\label{eq: estimate broken bichar}
  \Norm{\pgamma_n(s_2') - \pgamma_n(s_1') }{}
  + \big| \vartheta_n(s_2') - \vartheta_n (s_1') \big| \lesssim s_2' - s_1',
  \qquad s_1 \leq s_1' \leq s_2' \leq s_2.
\end{align}
On a piece given by an affine part of the construction, say on 
$[s_1, s_1 + \delta_n]$, one  has,
from  Proposition~\ref{prop: propagation glancing region}
\begin{align*}  
  \frac{d}{d s}\pgamma_n(s) 
  = \HpG \big(\gamma_n (s_1)\big) + \mathcal{O}( 1 / n)
  = \HpG \big(\gamma_n (s)\big) +  \pe_n(s), 
\end{align*}
where the errors $|\pe_n|$ goes  to zero {\em uniformly}
as $n\rightarrow + \infty$ by the uniform continuity of $\HpG$ in $K^0$,
\begin{align*}  
  |\vartheta_n(s) | \lesssim (s- s_1) \delta_n^{1/2},
\end{align*}
using that $\vartheta_n(s_1)=0$.
 One thus finds
\begin{align}
  \label{eq: gen bichar - Ascoli 2}
  \pgamma_n(s) 
  &=   \pgamma_n(s_1)  +  \int_{s_1}^s \HpG \big(\gamma_n
  (\sigma)\big)\, d \sigma 
  +   \int_{s_1}^s \pe_n(\sigma) \, d \sigma.
\end{align}
Since $s \mapsto \pgamma_n(s)$ is continuous, with \eqref{eq: gen
  bichar - Ascoli 1}--\eqref{eq: gen bichar - Ascoli 2} one obtains 
\begin{align}
  \label{eq: construction gammaG - est y' z}
  \bigNorm{|\pgamma_n (s_2) - \pgamma_n (s_1)}{} \lesssim s_2-s_1,
  \qquad 0\leq s_1 \leq s_2 \leq S_{\max},
\end{align}
uniformly with respect to $n \in \N^*$. 

\medskip
We now proceed with an estimation of the variations of $\vartheta_n (s)$. Note that $s \mapsto \vartheta_n (s)$ is not continuous on the
parts made with a \bbichar as hyperbolic points are encountered. 
For this reason, we are interested by the equicontinuity of
$(\cgamma_n)_{n \in \N^*}$ and 
our goal is to obtain an estimation of
$\cdist\big( \cgamma_n(s_2),
\cgamma_n(s_1) \big)$, for
$0 \leq s_1 \leq s_2 \leq S_{\max}$. 
We refer to Section~\ref{sec: The compressed cotangent bundle} for the definition of the distance $\cdist(.,.)$.

A first case to be considered is $M_n^{[s_1, s_2]} = \emptyset$. Then,
with \eqref{eq: construction gammaG - est y' z} one has 
\begin{align*}
  %\label{eq: estimate distance - Mn empty}
  \cdist\big( \cgamma_n(s_2),\cgamma_n(s_1)\big)
  \lesssim \Norm{\gamma_n (s_2) - \gamma_n (s_1)}{} 
  \lesssim s_2-s_1 + Z_{[s_1,s_2]},
\end{align*} 
where 
\begin{align*}
  %\label{eq: estimate distance - Mn empty2}
   Z_{[s_1,s_2]} =   \norm{\vartheta_n (s_2) - \vartheta_n (s_1)}{} ,
\end{align*}
to be estimated in Lemma~\ref{lemma: construction gammaG - est zeta final} below.

Next, different cases have to be considered if $M_n^{[s_1, s_2]} \neq \emptyset$.
In all cases, a simple yet  key observation is 
\begin{align}
  \label{eq: going from G to G}
    \vartheta_n (\sigma_{n, j}^\infty) =0 
  \ \ \text{for} \ j=1, \dots, m_n,
  \end{align}
since $\gamma_n(\sigma_{n, j}^\infty) \in \glGb \cup \sgGb$.

\medskip
Given $0\leq s_1 \leq s_2 \leq S_{\max}$, one has  $M_n^{[s_1, s_2]} =
\{ j_1, \dots, j_N\}$ for some  $N\leq m_n$ and with $1\leq j_1 < \cdots < j_N$. 
Set 
\begin{align*}
  &\ell_{n,j_1}^{\inf} = \min \{ 0\leq \ell \leq  L_{n,j_1}+1; \ s_1 \leq
  \sigma^{\ell}_{n,j_1} \}, \\
  &\ell_{n,j_N}^{\sup} = \sup \{ 0\leq \ell \leq  L_{n,j_N}+1; \ 
  \sigma^{\ell}_{n,j_N} \leq s_2 \}, 
\end{align*}
and 
\begin{align*}
  \sigma_{n,j_1}^{\inf} =\sigma_{n,j_1}^{\ell_{n,j_1}^{\inf}}
  \ \  \text{and} \ \ \sigma_{n,j_N}^{\sup} = \sup \{ \sigma^{\ell}_{n,j_N}; \ 
  \sigma^{\ell}_{n,j_N} \leq s_2 \}, 
\end{align*}
Recall that $L_{n,j}$ is possibly infinite.

%% Cas 1
If $\sigma_{n,j_N}^{\sup} \leq s_2 < \sigma_{n,j_N}^\infty$ one has
the estimate, using \eqref{eq: distance compressed cotangent bundle}-\eqref{eq: cont quotient map},
\begin{align*}
  &\cdist\big( \cgamma_n(s_2),\cgamma_n(s_1)\big)\\
  &\qquad \quad \lesssim
  \big\|\gamma_n(\sigma_{n,j_1}^{\inf,-} ) - \gamma_n( s_1^+) \big\| 
  + 
  \sum_{\ell_{n,j_1}^{\inf} \leq \ell \leq  L_{n,j_1}}
  \big\|\gamma_n(\sigma_{n,j_1}^{\ell+1, -})
  - \gamma_n(\sigma_{n,j_1}^{\ell,+}) \big\| \\
   &\qquad \quad\quad  + \big\|]
    \gamma_n(\sigma_{n,j_N}^{0, -})- \gamma_n(\sigma_{n,j_1}^{\infty,+}) \big\| 
  \notag\\
  &\qquad \quad\quad +\sum_{0\leq \ell \leq  \ell_{n,j_N}^{\sup}-1} \big\|\gamma_n(\sigma_{n,j_N}^{\ell+1, -})
  - \gamma_n(\sigma_{n,j_N}^{\ell,+}) \big\|   
 + \big\|\gamma_n(s_2^-) - \gamma_n(\sigma_{n,j_N}^{\sup, +}) \big\|.\notag
\end{align*}
With \eqref{eq: construction gammaG - est y' z} and \eqref{eq: going
  from G to G} one finds
\begin{align}
  \label{eq: estimate distance E}
  \cdist\big( \cgamma_n(s_2),\cgamma_n(s_1)\big)
  \lesssim s_2-s_1 + Z_{[s_1,s_2]}, 
\end{align} 
with $Z_{[s_1,s_2]}$ given by 
\begin{align}
  \label{eq: construction gammaG - est zeta 0}
   Z_{[s_1,s_2]} &= \big|\vartheta_n(\sigma_{n,j_1}^{\inf,-} ) - \vartheta_n( s_1^+) \big| 
  + 
  \sum_{\ell_{n,j_1}^{\inf} \leq \ell \leq  L_{n,j_1}}
  \big|\vartheta_n(\sigma_{n,j_1}^{\ell+1, -})
  - \vartheta_n(\sigma_{n,j_1}^{\ell,+}) \big| \\
   &\quad  
  +\big|\vartheta_n(\sigma_{n,j_N}^{0, -}) \big| 
     +\sum_{0\leq \ell \leq  \ell_{n,j_N}^{\sup}-1} \big|\vartheta_n(\sigma_{n,j_N}^{\ell+1, -})
  - \vartheta_n(\sigma_{n,j_N}^{\ell,+}) \big|   \notag\\
 &\quad  + \big|\vartheta_n(s_2^-) - \vartheta_n(\sigma_{n,j_N}^{\sup, +}) \big|.\notag
\end{align}

%% cas 2
If $\sigma_{n,j_N}^{\infty} = s_2 $, one writes 
\begin{align*}
  &\cdist\big( \cgamma_n(s_2),\cgamma_n(s_1)\big)\\
  &\qquad \quad \lesssim
\big\|\gamma_n(\sigma_{n,j_1}^{\inf,-} ) - \gamma_n( s_1^+) \big\| 
  + 
  \sum_{\ell_{n,j_1}^{\inf} \leq \ell \leq  L_{n,j_1}}
  \big\|\gamma_n(\sigma_{n,j_1}^{\ell+1, -})
  - \gamma_n(\sigma_{n,j_1}^{\ell,+}) \big\| \\
   &\qquad \quad \quad  + \big\|\gamma_n(s_2)
    - \gamma_n(\sigma_{n,j_1}^{\infty,+}) \big\| ,\notag
\end{align*}
yielding the same estimate as in \eqref{eq: estimate distance E}
with now $Z_{[s_1,s_2]}$ given by 
\begin{align}
  \label{eq: construction gammaG - est zeta 0-bis}
   Z_{[s_1,s_2]} &= \big|\vartheta_n(\sigma_{n,j_1}^{\inf,-} ) - \vartheta_n( s_1^+) \big| 
  + 
  \sum_{\ell_{n,j_1}^{\inf} \leq \ell \leq  L_{n,j_1}}
  \big|\vartheta_n(\sigma_{n,j_1}^{\ell+1, -})
  - \vartheta_n(\sigma_{n,j_1}^{\ell,+}) \big|.
\end{align}
since $\gamma_n(s_2)\in \glGb \cup \sgGb$ and thus $\vartheta_n(s_2)=0$.

%cas 3
\medskip
If now $\sigma_{n,j_N}^\infty < s_2$, one has 
$\ell_{n,j_N}^{\sup}= L_{n,j_N}+1$ or equivalently
$\sigma_{n,j_N}^{\sup} = \sigma_{n,j_N}^{L_{n,j_N}+1} =
\sigma_{n,j_N}^\infty$. Arguing as above
one has 
\begin{align*}
  &\cdist\big( \cgamma_n(s_2),\cgamma_n(s_1)\big)\\
  &\qquad \quad \lesssim 
\big\|\gamma_n(\sigma_{n,j_1}^{\inf,-} ) - \gamma_n( s_1^+) \big\| 
  + 
  \sum_{\ell_{n,j_1}^{\inf} \leq \ell \leq  L_{n,j_1}}
  \big\|\gamma_n(\sigma_{n,j_1}^{\ell+1, -})
  - \gamma_n(\sigma_{n,j_1}^{\ell,+}) \big\| \\
   &\qquad \quad \quad  + \big\|\gamma_n(\sigma_{n,j_{N}}^{\infty, -})
    - \gamma_n(\sigma_{n,j_1}^{\infty,+}) \big\| 
  + \big\|\gamma_n(s_2^-) - \gamma_n(\sigma_{n,j_N}^{\infty, +}) \big\|,\notag
\end{align*}
yielding the same estimate as in \eqref{eq: estimate distance E}
with now $Z_{[s_1,s_2]}$ given by 
\begin{align}
  \label{eq: construction gammaG - est zeta 0-ter}
   Z_{[s_1,s_2]} &= \big|\vartheta_n(\sigma_{n,j_1}^{\inf,-} ) - \vartheta_n( s_1^+) \big| 
  + 
  \sum_{\ell_{n,j_1}^{\inf} \leq \ell \leq  L_{n,j_1}}
  \big|\vartheta_n(\sigma_{n,j_1}^{\ell+1, -})
  - \vartheta_n(\sigma_{n,j_1}^{\ell,+}) \big|
                   + \big|\vartheta_n(s_2^-) \big|.
\end{align}

%%%%%%%%%%%%%%%%%%%%%%%
% figure
%%%%%%%%%%%%%%%%%%%%%%%
\begin{figure}
  \begin{center}
    \subfigure[Case $s_2 < \sigma_{n,j_N}^\infty$;  \eqref{eq: construction gammaG - est zeta 0} applies. 
    \label{fig: config gamma n-1}]
    {\resizebox{10.0cm}{!}{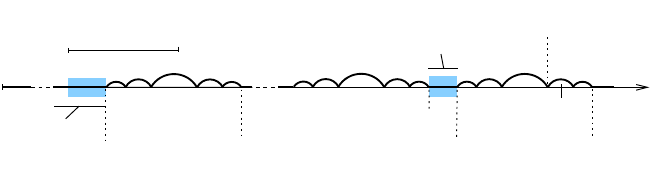}}
    
    \subfigure[Case $\sigma_{n,j_N}^\infty < s_2$; \eqref{eq: construction gammaG - est zeta 0-ter} applies.
  \label{fig: config gamma n-2}]
    {\resizebox{10.0cm}{!}{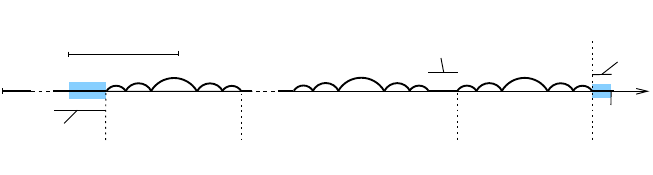}}
    \caption{Two configuratons of $[s_1, s_2]$ with respect to
      $\gamma_n$; `p.a.'~stands for `piecewise affine'. A shaded box
      refers to a location where an estimate making the term
      $\sqrt{s_2 - s_1}$ appear is used in the proof of
      Lemma~\ref{lemma: construction gammaG - est zeta final}. Note
      that discontinuities of \bbichars at hyperbolic points are not
    represented here for convenience.}
  \label{fig: config [s1,s2] gamma n}
  \end{center}
\end{figure}
Some configurations are illustrated in Figure~\ref{fig: config [s1,s2]
  gamma n}. In each case above one has   the following result. 
%%%%%%%%%%%%%%%%%%%%%%%%
% lemma                %
%%%%%%%%%%%%%%%%%%%%%%%%
\begin{lemma}
  \label{lemma: construction gammaG - est zeta final}
  There exists $C>0$, independent of $n \in \N^*$, such that 
  \begin{align*}
    Z_{[s_1,s_2]} \leq C \big( s_2 -s_1 +  \sqrt{s_2 -s_1}\, \big),
  \end{align*}
  for all $0 \leq s_1 \leq s_2 \leq S_{\max}$.  
\end{lemma}
A proof is given below. Then the different estimates above yield
\begin{align*}
  \cdist\big( \cgamma_n(s_2),\cgamma_n(s_1)\big)
  \lesssim s_2 -s_1 +  \sqrt{s_2 -s_1},
\end{align*}
implying the equicontinuity of $(\cgamma_n)_{n \in \N^*}$ 
 on $[0,S_{\max}]$. Since the sequence  is also pointwise bounded one can extract a
subsequence $\big(s \mapsto \cgamma_{n_p}\big)_{p \in \mathbb{N}}$ that
converges uniformly to a curve $\cgamma(s)$, $s \in [0,S_{\max}]$ with
values in $\cTL$ by the Arzel\`a-Ascoli theorem. One has
$\cgamma(0) = \cphi(\y^0)  = \y^0$.  
Set 
\begin{align*}
  B = \{ s \in [0, S_{\max}]; \  \cgamma(s) \in \cphi(\Hb)\}.
\end{align*}
For $s \in [0, S_{\max}] \setminus B$  one defines $\gammaG(s) =
\cphi^{-1} (\cgamma(s))$. One has $\gammaG(0)= \y^0$. 

For any $s \in [0,S_{\max}]$, the convergence of
$\big(\cgamma_{n_p}\big)$ yields a sequence of points of
$\cphi\big(\supp \mu\big)$ that converges to $\cgamma(s)$,
implying that $\cgamma(s)\in \cphi\big(\supp \mu\big)$ since this
set is closed. Thus, for any $s \in [0,S_{\max}] \setminus B$ one has
$\gammaG(s) \in \supp \mu$. 

It now remains to prove that $\gammaG$ is a \gbichar on
$[0,S_{\max}] \setminus B$, including that  $B$ is a discrete set. 

\subsubsection{{\bfseries Proof that the limit curve is a local
    generalized bicharacteristic}}
% \subsubsection{{\bfseries Proof that  $\bld{\gammaG}$ is a local
%     generalized bicharacteristic}}
\label{sec: proof local bichar}
We review (yet again) all possible occurences. Above, we considered a
subsequence $\gamma_{n_p}$ to achieve convergence (by means of the
Arzel\`a-Ascoli theorem). Here, for the sake of simplicity one uses
$\gamma_n$ to denote this subsequence. 

\medskip
\noindent
{\bfseries (1) Points away from $\bld{\dTL}$.}
Assume that $\y^1 = \gammaG(s^0) \in \TL \setminus \dTL$. Since $\y^1$ is away from $\Gb \cup \Hb$, for
$n$ chosen \suff large $\gamma_{n}$ is a piece of \bichar in a \nhd
of $\y^1$. Thus, for some $\delta>0$, one has
\begin{align*}
  \gamma_{n}(s) = \gamma_{n}(s^0) 
  + \int_{s^0}^s \Hp \big(\gamma_{n} (\sigma)\big)\, d \sigma,
  \qquad s \in [s^0 -\delta, s^0 + \delta].
\end{align*}
As convergence is uniform, passing to the limit $n \to \infty$, one finds
\begin{align*}
  \gammaG(s) = \y^1
  + \int_{s^0}^s \Hp \big(\gammaG (\sigma)\big)\, d \sigma,
  \qquad s \in [s^0 -\delta, s^0 + \delta],
\end{align*}
meaning that $\gammaG$ is a \bichar in a \nhd of $\y^1$.

\medskip
\noindent
{\bfseries (2) Hyperbolic points.}
Consider $s^0 \in B$. One has $\cgamma(s^0) \in \cphi(\Hb)$.
Set $\{ \y^{1,+}, \y^{1,-}\} = \cphi^{-1}\big(\{ \cgamma(s^0)\}\big)$ with $\y^{1,\pm} \in \Hb^{\pm}$ and $\y^{1,\pm} = \Sigma(\y^{1,\mp})$.

Note that $\cdist\big(\cgamma(s^0),\glGb\cup
\sgGb)\big)>0$. Hence, there exist $C>0$ and $\delta_0>0$ such that 
\begin{align*}
  \cdist\big(\cgamma(s),\glGb\cup\sgGb)\big)\geq C>0, \qquad 
  s \in [s^0 -\delta_0, s^0 +\delta_0],
\end{align*}
meaning that for some $n_0 \in \N$ chosen \suff large, each term of
the sequence $\gamma_n$ is made of a \bbichar on the interval
$[s^0 -\delta_0, s^0 +\delta_0]$ for $n \geq n_0$.  

For each point $\y^{1,\pm}$, use an open \nhd $V^{0,\pm}$ of
$\y^{1,\pm}$ in $\Hb^\pm$ and $S_0^{\pm}>0$ as given by
Lemma~\ref{lemma: limit of a bichar in H}. Set $0< \delta_1 < \min(S_0^+/2,
S_0^-/2, \delta_0)$.

For $R>0$ set $\mathcal B_R^{+}$ as the open ball of radius $R$ in the
variables $(u, \vartheta)$ centered at the point $\y^{1,\pm} = (u^1, z^1=0, \vartheta^1)$.
Set $C_1 = \vartheta^1/2 >0$.  
There exist $R_0>0$ and $Z>0$ such that 
\begin{align*}
  \mathcal B_{2R_0}^+ \subset V^{0,+} 
  \ \  \text{and} \ \ 
  \Sigma \big(\mathcal B_{2R_0}^+\big) \subset V^{0,-},
\end{align*}
and 
\begin{align*}
  |\vartheta| \geq C_1 \ \text{in} \ W_2^+ \cup W_2^- \ \ \text{with}\ W_2^+ =
  \mathcal B_{2R_0}^+ \times [0,Z[ \ \text{and} \ 
  W_2^- = \Sigma\big(\mathcal B_{2R_0}^+\big) \times [0,Z[.
\end{align*} 
For $N \in \N$, such that $Z_N= C_1 \delta_1  /(N+2) \leq Z$, and $k =1,2$,
further set 
\begin{align*}
  W_{k,N}^+ =
  \mathcal B_{k R_0}^+ \times [0, Z_N [ \ \text{and} \ 
  W_{k,N}^- = \Sigma\big(B_{k R_0}^+\big) \times [0, Z_N[.
\end{align*} 
Note that $W_{k,N}^\pm \cap \{ z=0\} \subset V^{0,\pm} \subset \Hb^\pm$.

In $W_2^+ \cup W_2^-$ one has $|\Hpz| = |\vartheta| \geq C_1$.  
For $N_0\in \N$ chosen \suff large and $N \geq N_0$, observe that any
\bichar initiated in $W_{1,N}^+$ (\resp $W_{1,N}^-$) exits
$W_{2,N}^+$ (\resp $W_{2,N}^-$) at two points: one  located in
$\mathcal B_{2R_0}^+\times \{ 0\}$ and one located in $\mathcal
B_{2R_0}^+\times \{ Z_N\}$ (\resp points located in 
$\Sigma \big(\mathcal B_{2R_0}^+\big)\times \{ 0, Z_N\}$). 
Moreover, this occurs within a $s$-interval of size at most $Z_N /C_1 =
\delta_1 / (N+2)$.
Set
$\delta_2 = \delta_1/(N_0+1)$.

For some  $n_1\geq N_0$, if $n \geq n_1$ one has
$\cgamma_n(s^0)  \in \cphi\big( W_{1,N_0}^+ \cup
W_{1,N_0}^-\big)$. From the discussion above, one finds that  there exists
\begin{align*}
  s_n \in [s^0 -\delta_2, s^0 + \delta_2]
\end{align*}
such that $\cgamma_n(s_n) \in \cphi(\Hb)$.
Since $\gamma_n(s_n^-)\in V^{0,-}$, $\gamma_n(s_n^+)\in V^{0,+}$, and
$2 \delta_2 < \min(S_0^+, S_0^-)$ one concludes with 
Lemma~\ref{lemma: limit of a bichar in H} that $s_n$ is the unique value
of $s \in [s^0 -\delta_2, s^0 + \delta_2]$ such
that $\cgamma_n(s) \in \cphi(\Hb)$.

Similarly, given $N \geq N_0$, for some $n_2 = n_2(N) \geq n_1$, if $n \geq n_2$, one has
$\cgamma_n(s^0)  \in \cphi \big( W_{1,N}^+ \cup
W_{1,N}^-\big)$.
With the same argument as above there exists 
\begin{align*} 
  s_n' \in [s^0 -\delta_1  /(N+1), s^0 + \delta_1  /(N+1)]
  \subset  [s^0 -\delta_2, s^0 + \delta_2],
\end{align*}
such that $\cgamma_n(s_n') \in \cphi(\Hb)$. Naturally,
the uniqueness of $s_n$ obtained above yields $s'_n = s_n$.  One therefore obtain the
convergence of the sequence $(s_n)_{n\geq n_1}$ to $s^0$.

Consider now $0< \eps <\delta_2/2$. For some $n_3\geq n_1$, if $n \geq
n_3$ one has
$s_n \leq s^0 + \eps$.
On the interval $]s_n, s^0 + \delta_2]$,
the curve $\gamma_n$ is a \bichar located in $\TL \setminus \dTL$ and 
\begin{align*}
  \gamma_n(s) = \gamma_n(s^0 + \eps) 
  + \int_{s^0+\eps}^s \Hp \big(\gamma_n(\sigma)\big)\, d \sigma, 
  \quad s \in [s^0 +\eps, s^0+ \delta_2].
\end{align*}
On $]s_n, s^0 + \delta_2]$ one has
$\vartheta_n(s) = \Hpz \big(\gamma_n(s)\big) >0$.  Yet, locally, one has
$\Hpz = \vartheta \geq C_1>0$ and $z_n(s_n)=0$, then
$z_n(s^0+\eps) \geq \eps C_1$. On the interval
$[s^0+\eps, s^0+ \delta_2]$ one thus has $\cgamma_n(s) = \gamma_n(s)$.
With the uniform limit of $(\cgamma_n)_n$ one obtains for
$\gammaG(s) = \big(u(s), z(s), \vartheta(s)\big)$ that
$z(s) \geq z(s^0+ \eps) \geq C_1 \eps$ for
$s \in [s^0+\eps, s^0+ \delta_2]$ 
and 
\begin{align*}
  \gammaG(s) = \gammaG(s^0 + \eps) 
  + \int_{s^0+\eps}^s \Hp \big(\gammaG(\sigma)\big)\, d \sigma, 
  \quad s \in [s^0 +\eps, s^0+ \delta_2].
\end{align*}
 As $\eps>0$ is arbitrary this
implies that $\gammaG(s) \in \TL \setminus \dTL$ for
$s \in ]s^0, s^0+ \delta_2]$.
With the continuity of $\cgamma$ one has
moreover $\lim_{s\to s^{0,+}} \gammaG(s)= \y^{1,+} \in \H^+$ and 
\begin{align*}
  \gammaG(s) = \y^{1,+}
  + \int_{s^0}^s \Hp \big(\gammaG(\sigma)\big)\, d \sigma, 
  \quad s \in ]s^0, s^0+ \delta_2].
\end{align*}
On the interval $ ]s^0, s^0+ \delta_2]$, the curve $\gammaG(s)$ is
thus a piece of bicharacteristic initiated from a point in $\Hb^+$,
here $\y^{1,+}$.  

Similarly, one finds $\gammaG(s) \in \TL
\setminus \dTL$ for $]s^0-\delta_2, s^0[$, $\lim_{s\to s^{0,-}}
\gammaG(s) = \y^{1,-} \in \H^-$, and on the interval $ [s^0-\delta_2, s^0[$, the curve $\gammaG(s)$ is
this a piece of bicharacteristic initiated from a point in $\Hb^-$,
here $\y^{1,-}$.  

Since,
$\y^{1,-} = \Sigma\big(\y^{1,+}\big)$ this proves that 
$\gammaG$ fulfills the required conditions at hyperbolic points for
\bbichars and thus \gbichars; See Definitions~\ref{def: broken bichar} and \ref{def: generalized
  bichar}. In particular $s \in [s^0 - \delta_2, s^0[ \cup ]s^0, s^0+\delta_2[$
one has $\frac{d}{d s} \gammaG (s) = \Hp \big( \gammaG (s) \big)$. 

Moreover, what is above implies that the set $B$ is discrete.

\medskip
\noindent
{\bfseries (3) Points in $\bld{\sdGb}$.}
Suppose $s^0 \in [0,S_{\max}]$ is such that $\y^1 = \gammaG(s^0) \in \sdGb$.
One has $\dist (\gammaG(s^0), \glGb\cup\sgGb\big)>0$. Thus, there exist $C>0$ and $\delta_0>0$ such that 
\begin{align*}
  \dist\big(\gammaG(s), \glGb\cup\sgGb \big)\geq C>0, \qquad 
  s \in [s^0 -\delta_0, s^0 +\delta_0],
\end{align*}
meaning that for some $n_0 \in \N$ chosen \suff large, each term of the sequence
$\gamma_n$ is made of a \bbichar on the interval $[s^0 -\delta_0,
s^0 +\delta_0]$ for $n \geq n_0$. One has $\Hpz (\y^1)=0$ and $\Hp^2
z (\y^1)>0$. Thus, there exists a \nhd $V^0$ of $\y^1$ where
$\Hppz \geq C_2$ for some $C_2>0$. For some  $n_1\geq n_0$ and some
$0 < \delta_1 < \delta_0$, one has $\gamma_n(s) \in V^0$ for $n\geq
n_1$ and $s\in  [s^0 -\delta_1,s^0 +\delta_1]$.

One faces two occurences: either $(1)$ there exists a subsequence
$n_p \to \infty$ such that
$B_{n_p} \cap [s^0 - \delta_1, s^0+\delta_1] = \emptyset$, that is, $\gamma_{n_p}$ does not encounter any
hyperbolic point in this interval  or $(2)$ for
$n_0$ chosen \suff large one has
$B_{n} \cap [s^0 - \delta_1, s^0+\delta_1] \neq \emptyset$ if $n \geq
n_0$. 

\medskip
In case $(1)$, one concudes as for a point away from $\dTL$:
passing to the limit $n_p \to \infty$ one obtains 
\begin{align*}
  \gammaG(s) = \gammaG(s^0) 
  + \int_{s^0}^s \Hp \big(\gammaG (\sigma)\big)\, d \sigma,
  \qquad s \in [s^0 -\delta_1, s^0 + \delta_1],
\end{align*}
meaning that $\gammaG$ is a \bichar in a \nhd of $\y^1$.

\medskip
Consider now case $(2)$ and $n \geq n_0$ and
$s_n  \in B_n \cap [s^0 -\delta_1, s^0 +\delta_1]$, that is, $s_n$
associated with a hyperbolic point of $\gamma_n$ for $s\in [s^0 -\delta_1, s^0 +\delta_1]$. One has
$\Hpz(s_n^+) >0$. As $\gamma_n(s)$ remains in $V^0$ one concludes
that $\Hpz \big( \gamma_n(s)\big) >0$ and $z_n(s)>0$ for
$s\in ]s_n, s^0+\delta_1]$.  Similarly
$\Hpz \big( \gamma_n(s) \big)<0$ and $z_n(s)>0$ for
$s \in [s^0-\delta_1,s_n[$.  Hence, $\gamma_n$ has at most one
isolated hyperbolic point in the interval $[s^0 - \delta_1,
s^0+\delta_1]$. More precisely one has
\begin{align}
  \label{eq: estimation diffractive point z}
  z_n(s) \geq C_2(s-s_n)^2/2, 
  \qquad s\in [s^0 -\delta_1,s^0+\delta_1].
\end{align}
Assume that $\sigma^0$ is an accumulation point of the sequence 
$(s_n)_n$ and $s_{n_p}$ a subsequence that converges to
$\sigma^0$. For $\gammaG(s) = \big(t(s), x'(s),z(s), \tau(s),\xi(s)\big)$, by \eqref{eq: estimation diffractive point z} one obtains
in the limit $n_p \to \infty$
\begin{align}
  \label{eq: estimation diffractive point z2}
  z(s) \geq C_2 (s-\sigma^0)^2/2, 
  \qquad s\in [s^0 -\delta_1,s^0+\delta_1],
\end{align}
implying that $\sigma^0 = s^0$. Thus, one concludes that $s_n \to s^0$ as $n \to \infty$.

Consider $0< \eps < \delta_1/2$. For $n$ \suff large one has $s_n < s^0 + \eps$
and one has 
\begin{align*}
  \gamma_n(s) = \gamma_n(s^0 + 2 \eps) 
  +  \int_{s^0 +  2\eps}^s \Hp \big(\gamma_{n} (\sigma)\big)\, d \sigma,
  \qquad s \in [s^0 + 2 \eps, s^0+\delta_1].
\end{align*}
Since $z_n (s) \geq C_2 (s^0+ 2 \eps - s_n)^2 \geq C_2 \eps^2$,
$\gamma_n(s)$ remains away from $\Hb$ for $s \in [s^0 + 2 \eps,
s^0+\delta_1]$ and 
in the limit $n \to \infty$ one finds
\begin{align*}
  \gammaG(s) = \gammaG(s^0 + 2 \eps) 
  +  \int_{s^0 +  2\eps}^s \Hp \big(\gammaG (\sigma)\big)\, d \sigma,
  \qquad s \in [s^0 +  2\eps, s^0+\delta_1].
\end{align*}
Since $\gammaG(s^0) \in \sdGb$ one has $\gammaG(s^{0,-}) =
\gammaG(s^{0,+}) = \y^1$ as the continuity of $\cgamma$ implies the
continuity of
$\gammaG$ away from points of  $\Hb$. 
Thus, letting $\eps \to 0^+$, one finds
 \begin{align*}
  \gammaG(s) = \y^1
  +  \int_{s^0 }^s \Hp \big(\gammaG (\sigma)\big)\, d \sigma,
  \qquad s \in [s^0, s^0+\delta_1],
\end{align*}
using that $\gammaG(s) \in \TL \setminus \dTL$ by \eqref{eq:
  estimation diffractive point z2} for $s > s^0$.
Similarly, one has 
\begin{align*}
  \gammaG(s) = \y^1
  +  \int_{s^0}^s \Hp \big(\gammaG (\sigma)\big)\, d \sigma,
  \qquad s \in [s^0-\delta_1, s^0].
\end{align*}
As in Case $(1)$ this means that $\gammaG$ is a \bichar in a \nhd of $\y^1$.

\medskip
\noindent
{\bfseries (4) Points in $\bld{\glGb\cup\sgGb}$.}
Suppose $\y^1= \gammaG(s^0) \in \glGb\cup\sgGb$. An example of such a point
is naturally $\y^0 = \gammaG(0)$ where the constuction of $\gammaG$
is initiated.

Applying $\ppi$ to \eqref{eq: gen bichar - Ascoli 1} gives 
\begin{align*}
  \pgamma_n(s) = \pgamma_n(s_1^+) 
  + \int_{s_1}^s \HpG \big(\gamma_n (\sigma)\big)\, d \sigma.
\end{align*}
for $[s_1, s]$ within a maximal bicharateristic and
this extends to any whole maximal \bbichar in the
construction made above. With \eqref{eq: gen bichar - Ascoli 2} one
has  
\begin{align*}
  \pgamma_n(s) 
  &=   \pgamma_n(s_1)  +  \int_{s_1}^s \HpG \big(\gamma_n
  (\sigma)\big)\, d \sigma 
  +   \int_{s_1}^s \pe_n(\sigma) \, d \sigma,
\end{align*}
for any $s, s^1 \in [0,S_{\max}]$, 
with the errors $|\pe_n|$ going  to zero {\em uniformly}
as $n\rightarrow + \infty$. Note that $\HpG \big(\gamma_n
  (\sigma)\big)$ may only be discontinuous at hyperbolic points that
  form a discrete set and thus a countable set for each $n$. Dominated
convergence yields
\begin{align*}
  \ppi\big(\gammaG(s)\big) 
  &=   \ppi\big(\gammaG(s^0)\big)  
    +  \int_{s^0}^s \HpG \big(\gammaG  (\sigma)\big)\, d \sigma,
\end{align*}
for $s \in  [0,S_{\max}]$. From the continuity of the gliding vector
field $\HpG$ and the continuity of
$\gammaG (s)$ at $s=s^0$ one has 
$\HpG \big(\gammaG  (s)\big) = \HpG (\y^1)+ o(1)$
as $s \to s^0$, yielding
\begin{align*}
  \ppi\big(\gammaG(s)\big) 
  &=   \ppi\big(\gammaG(s^0)\big)  
    + (s-s^0) \HpG (\y^1)+ (s-s^0) o(1),
\end{align*}
implying that $ \ppi\big(\gammaG(s)\big) $ is differentiable at
$s=s_0$ and 
\begin{align*}
    \frac{d}{d s}  \ppi\big(\gammaG(s)\big)  (s^0) = \HpG (\y^1).
  \end{align*}
Parts (1), (2), and (3) show that the first
two assumption of Proposition~\ref{prop: limit generalized
  bichar in G-G0} are fulfilled by $s \mapsto \gammaG(s)$. Hence, Proposition~\ref{prop: limit generalized bichar
  in G-G0} applies (one may need to change $s$ into $s^0-s$ depending
if one considers $s > s^0$ or $s< s^0$). Consequently, $\gammaG$
is differentiable at $s =s^0$ and 
\begin{align*}
  \frac{d}{d s} \gammaG (s^0) = \HpG(\y^1).
\end{align*}
This concludes the proof of Lemma~\ref{lemma: construction generalized
    bichar near glancing point}.
\hfill \qedsymbol \endproof

\bigskip
%%%% proof of lemma
\begin{proof}[Proof of Lemma~\ref{lemma: construction gammaG - est
    zeta final}]
First, consider the case  $M_n^{[s_1,s_2]} = \emptyset$, meaning that
all points $\gamma_n(s)$ are in affine parts of the construction for
$s \in [s_1, s_2]$. Thus, there exists $r_0 \in [0,S_{\max}[$ and $N>0$
such that a sequence of $N$ affine parts is initiated at
$\gamma(r_0)$. For some $0 \leq \ell_1 \leq \ell_2 \leq N$, one has 
\begin{align*}
   s_j \in [r_{\ell_j}, r_{\ell_j+1}], \quad j=1,2 
  \ \  \text{for} \ r_\ell = r_0 + \ell
  \delta_n, 
\end{align*}
If $\ell_1 = \ell_2 = \ell$ one as 
\begin{align*}
  \gamma_n (s_j) = \gamma_n(r_\ell) 
  + \frac{s_j - r_\ell}{\delta_n}
  \big( \gamma_n(r_{\ell+1}) -  \gamma_n(r_\ell) \big)
\end{align*}
yielding 
$\vartheta_n (s_j) = 
  \frac{s_j - r_\ell}{\delta_n} \vartheta_n(r_{\ell+1})$,
as $\vartheta_n(r_{\ell})=0$ since $\gamma_n(r_\ell) \in \Gb$. One thus
obtains
\begin{align*}
  \vartheta_n (s_2) - \vartheta_n (s_1) = 
  \frac{s_2 - s_1}{\delta_n} \vartheta_n(r_{\ell+1}).
\end{align*}
As $|\vartheta_n(r_{\ell+1})| \lesssim \delta_n^{1/2}$ by
Proposition~\ref{prop: propagation glancing region} and $0 < s_2 - s_1
\leq \delta_n$ one obtains 
\begin{align}
  \label{eq: estimate Z s1s2 Mn=emptyset 1}
  Z_{[s_1,s_2]} = \norm{\vartheta_n (s_2) - \vartheta_n (s_1)}{} \lesssim \sqrt{s_2 - s_1} .
\end{align}
If now $\ell_1 < \ell_2$ one writes
\begin{align*}
  \vartheta_n(s_2 ) - \vartheta_n (s_1) =  \vartheta_n(s_2 ) - \vartheta_n(r_{\ell_2})
  + \vartheta_n(r_{\ell_1+1}) -  \vartheta_n (s_1) ,
\end{align*}
since $ \vartheta_n(r_{\ell_2}) = \vartheta_n(r_{\ell_1+1})=0$. The argument
that led to the previous
estimate \eqref{eq: estimate Z s1s2 Mn=emptyset 1} gives 
\begin{align*}
 \norm{\vartheta_n(s_2 ) - \vartheta_n(r_{\ell_2})}{}
  + \norm{\vartheta_n(r_{\ell_1+1}) -  \vartheta_n (s_1)}{} 
  \lesssim \sqrt{s_2  - r_{\ell_2}} + \sqrt{r_{\ell_1+1} - s_1},
\end{align*}
yielding also  in this second case
\begin{align*}
  %\label{eq: estimate Z s1s2 Mn=emptyset 2}
  Z_{[s_1,s_2]} = \norm{\vartheta_n (s_2) - \vartheta_n (s_1)}{} \lesssim \sqrt{s_2 - s_1} .
\end{align*}

Second, consider the case $M_n^{[s_1,s_2]} \neq \emptyset$. Start
with the case of $Z_{[s_1,s_2]}$ given by \eqref{eq: construction
  gammaG - est zeta 0}, that is, $\sigma_{n,j_N}^{\sup} \leq s_2 <
\sigma_{n,j_N}^\infty$.

By \eqref{eq: gen bichar - Ascoli 1} and \eqref{eq: estimate broken bichar}, since the following terms only
concern parts of \bichars,  one obtains 
\begin{align}
  \label{eq: construction gammaG - est zeta 1}
  &\sum_{\ell_{n,j_1}^{\inf} \leq \ell \leq  L_{n,j_1}}
  \big|\vartheta_n(\sigma_{n,j_1}^{\ell+1, -})
  - \vartheta_n(\sigma_{n,j_1}^{\ell,+}) \big| 
  + \sum_{0\leq \ell \leq  \ell_{n,j_N}^{\sup}-1} 
  \big|\vartheta_n(\sigma_{n,j_N}^{\ell+1, -})
  - \vartheta_n(\sigma_{n,j_N}^{\ell,+}) \big| \\
  &\quad \lesssim \sum_{\ell_{n,j_1}^{\inf} \leq \ell \leq  L_{n,j_1}}
  \big(\sigma_{n,j_1}^{\ell+1}- \sigma_{n,j_1}^{\ell}\big)
  + \sum_{0\leq \ell \leq  \ell_{n,j_N}^{\sup}-1} 
  \big( \sigma_{n,j_N}^{\ell+1}
  - \sigma_{n,j_N}^{\ell} \big) \notag\\
  &\quad \lesssim 
    \sigma_{n,j_1}^\infty - \sigma_{n,j_1}^{\inf}
  + \sigma_{n,j_N}^{\sup} - \sigma_{n,j_N}^{0}. \notag
\end{align}

Now, consider the interval
$]s_1, \sigma_{n,j_1}^{\inf}[$ and  estimate the term
$\big|\vartheta_n(\sigma_{n,j_1}^{\inf,-} ) - \vartheta_n( s_1^+) \big|$.  In
fact $\gamma_n( s_1^+)$ is either (1) located on the $j_1$th \bbichar, or (2) located on some affine part of
$\gamma_N$ that stands upstream from the initiation of $j_1$th \bbichar. The two situations are illustrated in
Figure~\ref{fig: construction gammaG}.
\begin{figure}
  \begin{center}
    \subfigure[$\gamma_n(s_1)$ in the $j_1$th \bbichar. \label{fig:
      construction gammaG 2}]
    {\resizebox{13.5cm}{!}{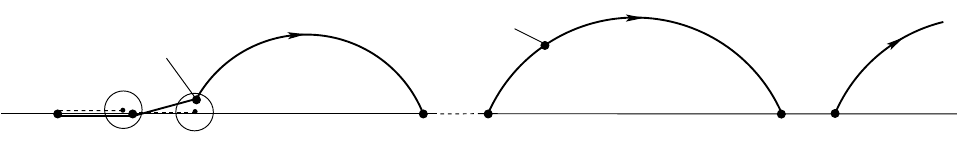}}
    
    \subfigure[$\gamma_n(s_1)$ in some affine part of $\gamma_n$
  upstream from the $j_1$th \bbichar is initiated. Here
  $\sigma_{n,j_1}^{\inf}  = \sigma_{n,j_1}^0$. \label{fig: construction gammaG 1}]
    {\resizebox{13.5cm}{!}{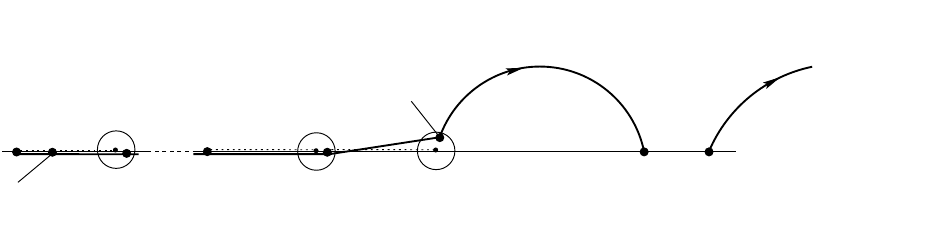}}
    \caption{Two possible locations of $\gamma_n(s_1^+)$, upstream from or
      within 
      the $j_1$th \bbichar.}
  \label{fig: construction gammaG}
  \end{center}
\end{figure}
In the first configuration, arguing as above one finds
\begin{align}
  \label{eq: construction gammaG - est zeta 3}
  \big|\vartheta_n(\sigma_{n,j_1}^{\inf,-} ) - \vartheta_n( s_1^+) \big|
  \lesssim \sigma_{n,j_1}^{\inf} - s_1.
\end{align}
 Let us now perform the
estimation in the second configuration; then 
$\sigma_{n,j_1}^{\inf} = \sigma_{n,j_1}^0$. Suppose the piecewise affine part of $\gamma_n$ upstream from the $j_1$th \bbichar
is composed of at least $m$ pieces.  Assume that 
\begin{align*}
  r_0 \leq s_1 \leq r_1= r_0  + \delta_n < \cdots <r_m = r_0 + m
  \delta_n = \sigma_{n,j_1}^0,
\end{align*}
 with $\gamma_n(r_0), \dots, \gamma_n(r_{m-1}) \in \supp \mu \cap (\glGb \cup \sgGb)$ and
 $\gamma_n(r_m) \in \supp \mu \setminus  (\glGb \cup \sgGb)$, and $\gamma_n$ affine between
 these points. The point $\gamma_n(s_1)$ is here located on the affine
 part joinging $\gamma_n(r_0)$ and $\gamma_n(r_1)$. It is not excluded
 here that there could be additional affine parts upstream from the
 point $\gamma_n(r_0)$. 

The $j_1$th \bbichar is initiated at
 $s=r_m = \sigma_{n,j_1}^0= \sigma_{n,j_1}^{\inf,-}$; See Figure~\ref{fig: construction gammaG 1}.
Estimates obtained above in the case $M_n^{[s_1,s_2]} =\emptyset$ give
\begin{align}
 \label{eq: construction gammaG -
  est zeta 4} 
  \bignorm{\vartheta_n(\sigma_{n,j_1}^{\inf,-} ) - \vartheta_n( s_1^+)}{}
  \lesssim \sqrt{\sigma_{n,j_1}^{\inf} - s_1} \lesssim \sqrt{s_2 - s_1}
\end{align}
To take into account both configurations, the first one given by
\eqref{eq: construction gammaG - est zeta 3} and  the second one given by \eqref{eq: construction gammaG -
  est zeta 4} one writes
\begin{align}
  \label{eq: construction gammaG - est zeta 4-bis}
  \big|\vartheta_n(\sigma_{n,j_1}^{\inf,-} ) - \vartheta_n( s_1^+) \big|
  \lesssim \sigma_{n,j_1}^{\inf} - s_1+ \sqrt{s_2 - s_1}.
\end{align}

Similarly, one finds 
\begin{align}
  \label{eq: construction gammaG - est zeta 5}
  \big|\vartheta_n(\sigma_{n,j_N}^{0, -})\big|  
  = \big|\vartheta_n(\sigma_{n,j_N}^{0, -}) -
  \vartheta_n(\sigma_{n,j_{N-1}}^{\infty,+}) \big| 
  \lesssim \sqrt{\sigma_{n,j_N}^{0} - \sigma_{n,j_{N-1}}^{\infty}}
  \lesssim \sqrt{s_2 - s_1}.
\end{align}
Finally, since $\sigma_{n,j_N}^{\sup} \leq s_2 <
\sigma_{n,j_N}^\infty$ one has as in \eqref{eq: construction gammaG - est zeta 1}
\begin{align}
  \label{eq: construction gammaG - est zeta 6}
  \big|\vartheta_n(s_2^-) - \vartheta_n(\sigma_{n,j_N}^{\sup, +}) \big|
  \lesssim s_2 - \sigma_{n,j_N}^{\sup},
\end{align}
since this term only
concerns parts of \bichars.

Summing estimates \eqref{eq: construction gammaG - est zeta 1},
\eqref{eq: construction gammaG - est zeta 4-bis}, \eqref{eq:
  construction gammaG - est zeta 5}, and \eqref{eq: construction gammaG - est zeta 6} one finds the estimate 
\begin{align*}
  Z_{[s_1,s_2]} \lesssim s_2 -s_1 +  \sqrt{s_2 - s_1}
  + \sigma_{n,j_1}^\infty - \sigma_{n,j_{N}}^0 
  \lesssim s_2 -s_1 +  \sqrt{s_2 - s_1}.
\end{align*}
for $Z_{[s_1,s_2]}$ as given by \eqref{eq: construction gammaG - est
  zeta 0}.

\bigskip Treat now the case of $Z_{[s_1,s_2]}$ given by \eqref{eq: construction
  gammaG - est zeta 0-bis}, that is, if $\sigma_{n,j_N}^{\infty} =
s_2$.
As above one has 
\begin{align}
  \label{eq: construction gammaG - est zeta 1-bis}
  &\sum_{\ell_{n,j_1}^{\inf} \leq \ell \leq  L_{n,j_1}}
  \big|\vartheta_n(\sigma_{n,j_1}^{\ell+1, -})
  - \vartheta_n(\sigma_{n,j_1}^{\ell,+}) \big| 
    \lesssim 
    \sigma_{n,j_1}^\infty - \sigma_{n,j_1}^{\inf}.
\end{align}
 Estimation \eqref{eq: construction gammaG - est zeta 4-bis}
holds and yield
\begin{align*}
  Z_{[s_1,s_2]} \lesssim  \sigma_{n,j_1}^\infty -s_1 +  \sqrt{s_2 - s_1}
  \lesssim s_2 -s_1 + \sqrt{s_2 - s_1},
\end{align*}
for $Z_{[s_1,s_2]}$ as given by \eqref{eq: construction gammaG - est
  zeta 0-bis}.

\bigskip Finally, Treat the case of $Z_{[s_1,s_2]}$ given by \eqref{eq: construction
  gammaG - est zeta 0-ter}, that is, if $\sigma_{n,j_N}^\infty < s_2$.
The $\gamma(s_2)$ lies in some affine part. 
Similarly to \eqref{eq: construction gammaG - est zeta 4} one has 
\begin{align*}
  %\label{eq: construction gammaG - est zeta 3-bis}
\big|\vartheta_n(s_2^-)\big| 
  \lesssim \sqrt{s_2 - s_1}.
\end{align*}
 Estimations \eqref{eq: construction gammaG - est zeta 1-bis} and
 \eqref{eq: construction gammaG - est zeta 4-bis}  hold and yield
\begin{align*}
  Z_{[s_1,s_2]} \lesssim \sigma_{n,j_1}^\infty  -s_1 + \sqrt{s_2 - s_1}
  \lesssim s_2 -s_1 +  \sqrt{s_2 - s_1},
\end{align*}
for $Z_{[s_1,s_2]}$ as given by \eqref{eq: construction gammaG - est
  zeta 0-ter}.
This concludes the proof of Lemma~\ref{lemma: construction gammaG - est zeta final}.
\end{proof}
%%%%%%%%%%%%%%%%%%%%%%%%%%
% section
%%%%%%%%%%%%%%%%%%%%%%%%%%
\section{Mass property of the boundary measure}
\label{sec: proof prop: no mass on Gd G3}
Here, we prove Proposition~\ref{prop: no mass on Gd G3}. Locally, we
first use the quasi-normal geodesic coordinates of
Proposition~\ref{prop: quasi-normal coordinates}.  With
Assumption~\ref{assumption: first properties of the semiclassical
  measure} one has $|\tau|\geq C_0>0$ in $\supp \mu$.  Consider $\y=
(t,x',z=0, \tau, \xi',\zeta=0)\in \pHb$ with $|\tau|< C_0$. Then, 
$\y^{\pm} \notin \supp \mu$.  In a \nhd of $\y^\pm$ Equation~\ref{eq:
  Gerard-Leichtnam equation} reads
\begin{align*}
   \int_{\y \in \pHb} 
    \frac{\delta_{\y^+} - \delta_{\y^-}}
    {\dup{\xi^+- \xi^-}{\nx}_{T_x^*\M, T_x\M} 
    } \ d \nu (\y) = 0,
  \end{align*}
implying that $\y \notin \supp \nu$. Thus
$\supp \nu \cap \pHb \cap \{ |\tau|< C_0\}= \emptyset$.  Consider now
$\y\in \pGb= \Gb$ with $|\tau|< C_0$. With Remark~\ref{remark:
  integrand GL equation on G-intro} one finds $\int_{\pGb} \d_\zeta
a\, d\nu =0$, for any $a\in \Con^1_c(T^* \hL)$ supported near
$\y$. Yet, any $\Con^0$-function on $\pdTL\simeq T^* \d \L$ can take the form $\d_\zeta
a_{|\!{z=0\atop \zeta=0}}$ implying that $\supp \nu \cap \pGb \cap \{ |\tau|< C_0\}=
\emptyset$. This gives the first result.

\medskip
The glancing set $\Gb$ is given by $\{
z=\Hpz=p=0 \} = \{z=\zeta=p=0 \}$. The set $\sdGb \cup \glGb$ is given by $\{
z=\Hpz=p=0 \ \et\ \Hppz\geq 0\}$.
  Define $r(x,\tau,\xi') =  \tau^2 -\sum_{1\leq i,j \leq d-1} g^{i j}(x) \xi_i \xi_j$.
One has 
   \begin{align*}
     p(x, \tau,\xi', \zeta)
     = (1 +  z h_{d d}(x)) \zeta^2
     + z \sum_{1\leq j \leq d-1} h_{j d}(x) \xi_j \zeta
     - r(x, \tau,\xi').
   \end{align*}
   Note that $\Gb$ is also given by $\{ z=r=p=0 \}$.

   As $\d_\tau r = 2 \tau \neq 0$ in $\supp \nu$ by the first part of
   the proposition, one finds that $r(x,\tau,\xi')$ can be used as a
   coordinate on $\{z=\Hpz=0\}$ near $\supp \nu$. Denote by $\sigma \in
   \R^{2d-1}$ coordinates such that $(z,\Hpz, r, \sigma)$ are local
   coordinates of $T^* \hL$.

   Consider $\psi \in \Cinf(\R)$ such that $\psi(s)=1$ for $s\geq 0$
   and $\psi(s)=0$ for $s\leq -1/2$. Consider also $b \in
   \Cinfc(\R^{2d-1})$ such that $b(z,r, \sigma)$ is independent of $z$ and $r$  in a \nhd of $\{
   z=r=0\}$. In $\supp b$, $|\tau|$ is bounded. Hence, $|\zeta|\leq C_1$ in $\supp b \cap \Char p$ for some $C_1>0$. 
   Pick $\varphi \in \Cinfc(\R)$ such that $\varphi(s)=1$
   for $|s|\leq 2C_1+1$ . For $\eps>0$ and $\alpha>0$, set
   \begin{align*}
     a_{\eps,\alpha}(\y) = (\Hpz)\varphi(\Hpz)
     \psi\big(\eps^{-1/4} \phi_\alpha * \Hppz\big)(\y) \, 
     b(\eps^{-1/2}z,\eps^{-1} r, \sigma)
     ,
   \end{align*}
   where $\phi_\alpha(\y) = \alpha^{-2d-2} \phi( \y /\alpha)$ with $\int
   \phi =1$, that is, an approximation to the identity. One has
   \begin{align*}
     \Hp a_{\eps,\alpha} (\y)
     &= (\Hppz)\big(\varphi(\Hpz) \\
     &\quad + (\Hpz)\varphi'(\Hpz) \big)
     \psi\big(\eps^{-1/4} \phi_\alpha * \Hppz  \big)(\y)
     b(\eps^{-1/2}z,\eps^{-1} r, \sigma)
     \\
     &\quad + \eps^{-1/2} (\Hpz)^2 \varphi(\Hpz)
     \psi\big(\eps^{-1/4} \phi_\alpha * \Hppz \big)(\y)
     \d_1 b(\eps^{-1/2}z,\eps^{-1} r, \sigma)
     \\
     &\quad +\eps^{-1} (\Hp r) (\Hpz)\varphi(\Hpz)
     \psi\big(\eps^{-1/4} \phi_\alpha * \Hppz \big)(\y)
      \d_2 b(\eps^{-1/2}z,\eps^{-1} r, \sigma)
     \\
     &\quad + (\Hpz)\varphi(\Hpz)
      \psi\big(\eps^{-1/4} \phi_\alpha * \Hppz \big)(\y) 
     d_3 b(\eps^{-1/2}z,\eps^{-1} r, \sigma)(\Hp \sigma)     
     \\
     &\quad + \eps^{-1/4}
     \big( (\Hp \phi_\alpha) * \Hppz\big)
     (\Hpz) \varphi(\Hpz) \\
     &\qquad \times
     \psi'\big(\eps^{-1/4}  \phi_\alpha * \Hppz \big)(\y)
     b(\eps^{-1/2}z,\eps^{-1} r, \sigma)
     .
   \end{align*}
   One has $|z| \lesssim \eps^{1/2}$, $|r| \lesssim \eps$ in
   $\supp a$. Since $\supp \mu \subset\Char p$, solving $p=0$ for
   $\zeta$ one finds that $|\zeta|\lesssim \eps^{1/2}$. The estimate
   $|\Hp z| = |\d_\zeta p| \lesssim \eps^{1/2}$ follows. As $\Hp r =
   \Hp (p+r)$ one finds $|\Hp r|\lesssim \eps^{1/2}$. These estimates
   and the dominated convergence theorem give
   \begin{align*}
     \dup{\transp \Hp \mu}{a_{\eps,\alpha}}\  \mathop{\longrightarrow}_{\eps \to 0}\ 
     \dup{\mu}{\Hppz(\y) b(0,0,\sigma)
       \unitfunction{\Gb} \unitfunction{A_{\alpha}}},
   \end{align*}
   using that $\supp \mu \subset \Char p$, with
   $A_{\alpha}
   =\{ \phi_\alpha * \Hppz(\y) \geq 0\}$. One also obtains
   \begin{align*}
     \dup{f \mu}{a_{\eps,\alpha}}\
     \mathop{\longrightarrow}_{\eps \to 0}\  0.
   \end{align*}

   \medskip
   If $\y  \in \pHb$ one has  $\Hpz(\y^+) = 2 \zeta^+  = - \Hpz(\y^-) = -2 \zeta^-$
   and $\varphi(2\zeta^+) = \varphi(2\zeta^-) = 1$.  One obtains
   \begin{multline*}
     \Bigdup{\int_{\y \in \pHb } 
       \frac{\delta_{\y^+} - \delta_{\y^-}}
            {\dup{\xi^+- \xi^-}{\nx}_{T_x^*\M, T_x\M}}
            \ d \nu (\y)}{a_{\eps,\alpha}}
     =\int_{\y \in \pHb }
     \frac{a_{\eps,\alpha}(\y^+)- a_{\eps,\alpha}(\y-)}
          {2 \zeta^+}\ d \nu (\y)\\
     =\int_{\y \in \pHb }
          b(0, \eps^{-1} r, \sigma)
          \big( \psi\big(\eps^{-1/2} \phi_\alpha * \Hppz(\y) \big)(\y^+)
          + \psi\big(\eps^{-1/2} \phi_\alpha * \Hppz(\y) \big)(\y^-)
          \big)\ d \nu (\y).
  \end{multline*}
   Since $r>0$ in $\pHb$, with the support property of $b$, by dominated convergence one finds 
   \begin{align*}
     &\Bigdup{\int_{\y \in \pHb } 
       \frac{\delta_{\y^+} - \delta_{\y^-}}
            {\dup{\xi^+- \xi^-}{\nx}_{T_x^*\M, T_x\M}}
            \ d \nu (\y)}{a_{\eps,\alpha}}
     \  \mathop{\longrightarrow}_{\eps \to 0} \  0.
     \end{align*}
   Next, as $\pGb = \Gb$ is also given by $\{ z=r=0\}$ in $\pTL$ one has 
  \begin{align*}
    \Bigdup{\int_{\y \in \pGb} 
    \frac{\delta_{\y^+} - \delta_{\y^-}}
    {\dup{\xi^+- \xi^-}{\nx}_{T_x^*\M, T_x\M} 
    } \ d \nu (\y)}{a_{\eps,\alpha}}
    = \dup{\nu}{2 b(0, 0, \sigma)
          \psi\big(\eps^{-1/2} \phi_\alpha * \Hppz\big)}.
  \end{align*}
  By dominated convergence one finds 
     \begin{align*}
    \Bigdup{\int_{\y \in \pHb \cup \pGb} 
    \frac{\delta_{\y^+} - \delta_{\y^-}}
         {\dup{\xi^+- \xi^-}{\nx}_{T_x^*\M, T_x\M}}
         \ d \nu (\y)}{a_{\eps,\alpha}}
    \  \mathop{\longrightarrow}_{\eps \to 0}\
       2\dup{\nu}{b(0,0,\sigma) \unitfunction{\Gb}
         \unitfunction{A_{\alpha}}}.
     \end{align*}
The measure equation of Assumption~\ref{assumption: Gerard-Leichtnam
  equation}  gives
\begin{align*}
  \dup{\mu}{\Hppz(\y) b(0,0,\sigma)
    \unitfunction{\Gb} \unitfunction{A_{\alpha}}}
     = - 2 \dup{\nu}{b(0,0,\sigma) \unitfunction{\Gb}
         \unitfunction{A_{\alpha}}}.
\end{align*}
Letting $\alpha \to 0$ gives by dominated convergence
\begin{align*}
  \dup{\mu}{\Hppz(\y) b(0,0,\sigma)
    \unitfunction{\Gb} \unitfunction{\Hppz\geq 0}}
     = - 2\dup{\nu}{b(0,0,\sigma) \unitfunction{\Gb}
         \unitfunction{\Hppz\geq 0}}.
\end{align*}
One has $\unitfunction{\Gb} \unitfunction{\Hppz\geq 0}=
\unitfunction{\sdGb \cup \glGb}$.  If $b\geq 0$ one finds that both
sides have opposite signs since $\mu$ and $\nu$ are both nonnegative
measures. Hence, both side vanish if $b \geq 0$. Thus, on the one
hand, $\dup{\nu}{b(0,0,\sigma) \unitfunction{\sdGb \cup \glGb}}=0$ for
any $b \geq 0$, yielding $\dup{\nu}{\unitfunction{\sdGb \cup
    \glGb}}=0$. On the other hand, $\dup{\mu}{\Hppz(\y) b(0,0,\sigma)
  \unitfunction{\sdGb \cup \glGb}}=0$ for any $b\geq 0$.  Since
$\Hppz$ vanishes in $\glGb$ one obtains $\dup{\mu}{\Hppz(\y)
  b(0,0,\sigma) \unitfunction{\sdGb}}=0$ for any $b\geq 0$. One
concludes that $\dup{\mu}{ \unitfunction{\sdGb}}=0$ as $\Hppz>0$ on $\sdGb$.\hfill \qedsymbol \endproof

%%%%%%%%%%%%%%%%%%%%%%%%%%
%%%%%%%%%%%%%%%%%%%%%%%%%%
% Appendices 
%%%%%%%%%%%%%%%%%%%%%%%%%%
\appendix
%%%%%%%%%%%%%%%%%%%%%%%%%%
% section
%%%%%%%%%%%%%%%%%%%%%%%%%%
\section{Mesure associated with a single generalized bicharacteristic}
\label{sec: proof prop mesure single bichar}

Consider a \gbichar $\gammaG: \R \setminus B \to \TL$
according to Definition~\ref{def: generalized bichar-intro} or
equivalently Definition~\ref{def: generalized bichar}.  One wishes to
introduce the measure $\delta_{\GammaG}$ as the linear measure
supported on $\GammaG$ (see Definition~\ref{def: generalized bichar
  2-intro}), that is, 
\begin{align}
  \label{eq: linear measure bichar}
    \dup{\delta_{\GammaG}}{a} = \int_{\R \setminus B}
    a \big(\gammaG(s)\big)\ d s, 
    \qquad a \in \Con^0_c(T^* \hL).
\end{align}
First, observe that the curve $\gammaG$ is noncontinuous at hyperbolic
points and therefore not rectifiable. However, since $B$ is a discrete
set, upon defining the value of $\gammaG(s) = \gammaG(s^-)$ for $s
\in B$, one finds that $\gammaG$ is left continuous and therefore
measurable.
Second, observe that  $I_K = \{ s \in \R \setminus B; \ \gammaG(s) \in K\}$ is
bounded for $K$ a compact set of $\ThL$, using that $\frac{d
  t}{ds}(s) = -2 \tau\neq0$ and $\frac{d \tau}{ds}(s)=0$.
Consequently, given $a\in \Con^0_c(T^* \hL)$ with $\supp a \subset K$
one sees that the integral in \eqref{eq: linear measure
  bichar} is sensible and
\begin{align*}
  |\dup{\delta_{\GammaG}}{a} | \leq \Norm{a}{L^\infty} \, |I_K|,
\end{align*}
meaning that $\delta_{\GammaG}$ is a Radon measure.

\medskip The following theorem states that the measure
$\delta_{\GammaG}$ fulfills a transport equation of the form
considered in the present article.
%%%%%%%%%%%%%%%%%%%%%%%%
% theorem         %
%%%%%%%%%%%%%%%%%%%%%%%%
\begin{theorem}
  \label{theorem: mesure single bichar1}
  Consider a nontrivial \gbichar
  $\gammaG: \R \setminus B \to \TL$. Then \eqref{eq: linear measure bichar}
  yields a measure $\mu = \delta_{\GammaG}$ that fulfills the
  transport equation of \eqref{eq: Gerard-Leichtnam equation} for some
  nonnegative measure $\nu$ on $\pdTL$.
\end{theorem}
\begin{proof}
  To lighten notation write $\gamma$ in place of $\gammaG$.  Having
  $\gamma$ not trivial means that the cotangent component of
  $\gamma(s)$ is not zero, and thus neither $\Hp$ nor $\HpG$ vanish at
  any point of the \bichar.

  Partition $\R \setminus B$  into 
  \begin{align*}
    &L= \{ s \in \R \setminus B; \ \gamma(s) \in \TL \setminus
    (\glGb \cup \sgGb)\}, \\
    &G= \{ s \in \R \setminus B; \ \gamma(s) \in \glGb \cup \sgGb\}.
  \end{align*}

  \medskip
  As equation \eqref{eq: Gerard-Leichtnam equation}  is of geometrical
  nature, it suffices to check that it holds in some local chart. 
  The argument is simple away from the boundary. We thus choose a
  local chart $(\O, \cdiff)$  at the boundary of $\d\L$,
  see Section~\ref{sec: Local coordinates}, and we use the
  quasi-geodesic coordinates of Proposition~\ref{prop: quasi-normal coordinates} to simplify some computations.

  Consider $a \in \Conc^1\big(T^* (\O)\big)$.
  One writes 
  \begin{align*}
    \dup{ \transp{\Hp} \mu}{a} = \dup{\mu}{\Hp a} 
    =  \int_{\R \setminus B} (\Hp a) \big(\gamma(s)\big)\ ds.
  \end{align*}
  If $s \in L$ then 
  \begin{align}
    \label{eq: transport L}
     (\Hp a) \big(\gamma(s)\big) = \frac{d}{d s } a \big(\gamma(s)\big). 
  \end{align}
  If $s \in G$ then 
  \begin{align*}
    \Hp \big(\gamma(s)\big) 
    = \HpG \big(\gamma(s)\big) - \frac{\Hppz}{\Hzzp} \big(\gamma(s)\big)\Hz
    = \HpG \big(\gamma(s)\big) +  \frac12 \Hppz \big(\gamma(s)\big)\d_\zeta,
  \end{align*} 
  by \eqref{sec: def gliding vector field into} (see also Lemma~\ref{lemma: properties
    HpG}) and using that $\Hzzp = 2$ at $z=0$ in the chosen
  quasi-geodesic  coordinates. Thus one finds
  \begin{align}
    \label{eq: transport G}
     (\Hp a) \big(\gamma(s)\big) = \frac{d}{d s } a
    \big(\gamma(s)\big) 
    + \frac12 \big( \Hppz \, \d_\zeta a\big)  \big(\gamma(s)\big) . 
  \end{align}
  With \eqref{eq: transport L}
  and \eqref{eq: transport G} one obtains
  \begin{align}
    \label{eq: transport integral}
    \dup{ \transp{\Hp} \mu}{a} =  \frac12 A_1 + A_2,
  \end{align}
  with
  \begin{align*}
    A_1=\int_{G}
    (\Hppz) \d_\zeta a  \big(\gamma(s)\big)  ds
    \ \ \et \ \
   A_2= \int_{\R\setminus B}  \frac{d}{d s} a  \big(\gamma(s)\big) d s.
  \end{align*}

  \medskip
  First, consider the term $A_1$ in \eqref{eq: transport integral}.
  With \eqref{eq: understanding GL-glancing}, for $\y \in \glGb \cup
  \sgGb$ such as $\gamma(s)$ for $s  \in G$ one can write
  \begin{align*}
   \d_\zeta a (\y)
    =\lim 
    \frac{\dup{\delta_{\y^{(n)+}} - \delta_{\y^{(n)-}}}{a}}
    {\dup{\xi^{(n)+}-\xi^{(n)-}}{\nx}_{T_x^*\M, T_x\M}},
  \end{align*}
   for a sequence  $(\y^{(n)} )_n\subset \pHb$ that converges to $\y$.
   With the notation understanding used for the measure propagation
   equation \eqref{eq: Gerard-Leichtnam equation} for the part of the
   integration performed on $\pGb = \Gb$, one thus finds
  \begin{align}
    \label{eq: app term A1}
    A_1= - \int_{\pHb\cup \pGb}
    \frac{\dup{\delta_{\y^{+}} - \delta_{\y^{-}}}{a}}
    {\dup{\xi^{+}-\xi^{-}}{\nx}_{T_x^*\M, T_x\M}}
    d \nu_G,
  \end{align}
  with $d \nu_G = - \Hppz \, 
     \unitfunction{\glGb \cup \sgGb}\,   \delta_{\GammaG}$, measure on
     $\pdTL$. 
  It is nonnegative as one has 
   $- (\Hppz)
   \unitfunction{\glGb \cup \sgGb}  \geq 0$.

   \medskip Second, before considering the term $A_2$ in \eqref{eq: transport
     integral},   with $\pgamma(s) = \ppi \big(\gamma(s)\big)$ set
  \begin{align*}
    \pHb^\gamma = \{  \pgamma(s);\ s \in B\} \subset \pHb.
  \end{align*}
  Like $B$, this is a discrete set.  Write
  $\gamma(s)$ in the $\Con^1$-variables $(\py,\vartheta)$, that is,
  $\gamma(s) = \big( \pgamma(s), \vartheta(s)\big)$. Recall that
  $\vartheta = \Hpz$ and $\vartheta_{|z=0} = 2 \zeta$ in the
 used quasi-geodesic coordinates.  Set
  $[\vartheta]_s=\vartheta (s^+) - \vartheta (s^-)$ for $s \in
  B$. Note that $[\vartheta]_s = 2 \vartheta (s^+)>0$ since
  $\vartheta (s^+)>0$ and $\vartheta (s^-)= - \vartheta (s^+)$.  One
  needs the following lemma whose proof is given below.
 %%%%%%%%%%%%%%%%%%%%%%%%
% lemma                %
%%%%%%%%%%%%%%%%%%%%%%%%
\begin{lemma}
  \label{lemma: nu H measure}
  \leavevmode  
  \begin{enumerate}
    
  \item Suppose $I$ is a bounded interval.  The series
    $\sum_{s\in B\cap I}\, \vartheta(s^+)$ is absolutely convergent.

  \item The series
  $\nu_H = \sum_{\y \in \pHb^\gamma } \, \dup{\xi^+-
    \xi^-}{\nx}_{T_x^*\M, T_x\M} \ \delta_{\y}$ yields a nonnegative
  Radon measure on $\pdTL$.
  \end{enumerate}
\end{lemma}

   \medskip
   Third, consider the term $A_2$ in \eqref{eq: transport integral}.
   Recall from the main text that 
    $\pgamma(s) $ is defined in
     $\R\setminus B$ and yet can be continuously extended to $\R$ since
     $\pgamma(s^-) = \pgamma(s^+) \in \pHb$ if $s\in B$. Moreover this
     extended function is Lipschitz on $\R$. 

     Express the $\Conc^1$-function $a(\y)$ in terms of the variables
     $(\py, \vartheta)$. The term $A_2$ in \eqref{eq: transport integral}
     reads
     \begin{align*}
       A_2= \int_{\R\setminus B}
       \frac{d}{d s} b (s) d s
       \ \ \avec  \ \
       b(s) =  a \big(\pgamma(s),\vartheta(s)\big)  - a(\pgamma(s),0),
     \end{align*}
     using that
     $\int_{\R\setminus B}
       \frac{d}{d s} a  \big(\pgamma(s), 0\big) d s
       = \int_{\R}
       \frac{d}{d s} a  \big(\pgamma(s), 0\big) d s=0$
     since $a  \big(\pgamma(s), 0\big)$ is Lipschitz thus absolutely
     continuous.

     Note that $s \in \supp b$ implies $\gamma(s)
     \in \supp a$ or $(\pgamma(s),0) \in \supp a$. Hence, there exists an
     open 
     bounded interval $J$ such that $\supp b\subset J$ yielding
     \begin{align*}
       A_2= \int_{J\setminus B}
       \frac{d}{d s} b (s) d s.
     \end{align*}
     One has
     \begin{align*}
       \frac{d}{d s} b \big(\pgamma(s),\vartheta (s)\big)
       &= d_{\py} a \big(\pgamma(s),\vartheta (s)\big)
       \big(\pgamma'(s)\big)
       - d_{\py} a \big(\pgamma(s),0\big)
       \big(\pgamma'(s)\big)\\
       &\quad + \vartheta '(s) \d_\zeta a \big(\pgamma(s),\vartheta (s)\big) \ \text{\pp}
     \end{align*}
     As $(\ovl{B} \setminus B) \subset G$ by Lemma~\ref{lemma: generalized bichar - derivative endpoint}, if
     $s \in \ovl{B}  \setminus B$ one has $\vartheta (s) =
     \Hpz\big(\gamma(s)\big)=0$ and $\vartheta '(s)= \HpG \Hpz\big(\gamma(s)\big)=0$
     yielding
     $\frac{d}{d s} b \big(\pgamma(s),\vartheta (s)\big)=0$.
     One concludes that
     \begin{align*}
       A_2= \int_{J\setminus \ovl{B}}
       \frac{d}{d s} b(s) d s.
     \end{align*}
     If $s \in J\setminus \ovl{B}$, there exists $\eps>0$ such that
     $[s-\eps, s+\eps] \subset J$ and 
     $[s-\eps, s+\eps] \cap \ovl{B}= \emptyset$. Denote by $I_s \subset J$ the largest
     interval such that $s \in I_s$ and $I_s \cap \ovl{B} = \emptyset$. For
     $n \in \N^*$ set
     \begin{align*}
       R_n = \{ s \in J\setminus \ovl{B}; \ |I_s| \geq 1/n\}.
     \end{align*}
     Observe that $R_n$ is a finite union of disjoint intervals
     $I_{n,1}, \dots, I_{n, k(n)}$ all subsets of $J$. Note that $R_n \subset R_{n+1}$,
     and $J\setminus \ovl{B}= \cup_n R_n$, thus writing $J\setminus
     \ovl{B}$ as an at-most-countable union of disjoint intervals.
     Note that
     \begin{align*}
       B \subset E = \{ \inf I_{n,j}, \sup I_{n,j}; \ n \in \N, j = 1, \dots, k(n)\}. 
     \end{align*}
     and
     \begin{align*}
       E \setminus B \subset (\ovl{B} \setminus B)
       \cup \{ \inf J , \sup J\}.
     \end{align*}
     One obtains
     \begin{align}
       \label{eq: app integral all intervals}
       A_2= \lim_{n \to +\infty}\int_{R_n}
       \frac{d}{d s} b(s) d s
       = \lim_{n \to +\infty} \sum_{j=1}^{k(n)} \int_{I_{n,j}}
       \frac{d}{d s} b(s) d s.
     \end{align}
     One computes
     \begin{align}
       \label{eq: app integral single interval}
      \int_{I_{n,j}} \frac{d}{d s} b(s) d s = b(\sup I_{n,j}^-) - b(\inf I_{n,j}^+).
     \end{align}
     Note that if $s\in  (\ovl{B} \setminus B) \cup \{ \inf J, \sup J\}$
     then $b(s) =0$. Hence, non vanishing terms  on the rhs of
     \eqref{eq: app integral single interval} correspond only to the
     cases $s=\sup
     I_{n,j}$ lying in $B \cap J$ and $s=\inf I_{n,j}$ lying in $B\cap
     J$. If $s\in B \cap J$ the term
     $b(s^-)$  appears exactly once on the \rhs of \eqref{eq: app
       integral all intervals} as in \eqref{eq: app integral single interval}. The same holds for the term $- b(s^+)$.
     One has
     \begin{align*}
       |b(s^\pm)|
       = |a \big(\pgamma(s),\vartheta(s^\pm)\big)  - a(\pgamma(s),0)|
       \leq C |\vartheta(s^\pm)|.
     \end{align*}
     By Lemma~\ref{lemma: nu H measure} the series
     $\sum_{s \in B \cap J} |\vartheta(s^+)|$ is absolutely convergent
     using that $J$ is a bounded interval. As
     $|\vartheta(s^-)| = |\vartheta(s^+)|$ the same holds
     for $\sum_{s \in B \cap J} |\vartheta(s^-)|$.  Hence, the series
     $\sum_{s \in B \cap J} b(s^-)$ and $\sum_{s \in B \cap J} b(s^+)$
     are absolutely convergent. Summation order is thus not of
     importance and with \eqref{eq: app integral all intervals} one
     obtains
     \begin{align}
       \label{eq: app term A2}
       A_2&= \sum_{s \in B\cap J} \big( b(s^-)- b(s^+)\big)
       = - \sum_{s \in B\cap J} [a \circ \gamma]_s\\
          &=  -  \int_{\pHb\cup \pGb} \,
            \frac{\dup{\delta_{\y^+} -\delta_{\y^-}}{a}}
            {\dup{\xi^+- \xi^-}{\nx}_{T_x^*\M, T_x\M}  } \ d \nu_H (\y).\notag
     \end{align}
with $\nu_H$ as given in Lemma~\ref{lemma: nu H measure}.
     Combining \eqref{eq: app term A1} and \eqref{eq: app term A2}
     yields the result
     with $\nu = \frac12 \nu_G + \nu_H$.
\end{proof}

\medskip
%%%% proof of lemma
\begin{proof}[Proof of Lemma~\ref{lemma: nu H measure}]

  The first result is trivial if $\#
  I \cap B$ is finite. Assume it is infinite.
  At a point $s\notin B$ one has
  $\vartheta'(s) = \Hppz \big(\gamma(s)\big)$ or
  $\vartheta'(s) = (\HpG \Hpz) \big(\gamma(s)\big)=0$. It implies that $\vartheta'(s)$
  is bounded on $I\setminus B$.

Pick $s\in I \cap B$ with $s< \sup (I \cap \ovl{B})$. There exists
$\tilde{s}\in I \cap \ovl{B}$ with $s< \tilde{s}$ and
$]s, \tilde{s}[ \cap \ovl{B} = \emptyset$. One has $\vartheta(\tilde{s}) \leq
0$ ($<0$ if $\tilde{s}\in B$ and $=0$ if $\tilde{s}\in \ovl{B}$).
Thus one has
\begin{align*}
  \big| \vartheta(s^+) \big|
  \leq \sup_{I \setminus B}|\vartheta'|\,  |\tilde{s} -s|.
\end{align*}
If $s \in I \cap B$ and $s = \sup (I \cap \ovl{B})$ and
$s > \inf (I \cap \ovl{B})$ on picks $\tilde{s} \in I \cap \ovl{B}$
with $\tilde{s}< s$ and argues similarly.  This implies
$\sum_{s\in I \cap B} \big| \vartheta(s^+) \big| \leq  |I| \sup_{I
  \setminus B}|\vartheta'|$.

  \medskip
  If $\y \in \pHb^\gamma$, as
  $\vartheta_{|z=0}=2\zeta$ in the present coordinates, by \eqref{eq: relevement
    H-G-intro} one has 
  \begin{align*}
    \dup{\xi^+-\xi^-}{\nx}_{T_x^*\M, T_x\M}
    =[\Hpz]_{\y},
    \ \ \avec \ \ [\Hpz]_{\y}
    =\Hpz (\y^+) - \Hpz (\y^-).
  \end{align*}
  The series reads
  $\nu_H = \frac12 \sum_{\y \in \pHb^\gamma } [\Hpz]_{\y}\,
  \delta_{\y}$.  As $[\Hpz]_{\y}>0$, this is a positive
  measure. Continuity remains to be proven. Consider a smooth test
  function $a$ with compact support and set $b(s) = a\big(\pgamma(s)\big)$
  and $\vartheta(s) = \Hpz\big(\gamma(s)\big)$ as above.  One has
  $\dup{\nu_H}{a} = \dup{\tnu_H}{b}$ where
  \begin{align*}
    \tnu_H &= \frac12 \sum_{s\in B}  [\vartheta]_s\,\delta_s
             \quad
       \avec \ \  [\vartheta]_s =    \vartheta(s^+) - \vartheta(s^-) >0.
  \end{align*}
  Continuity of $\nu_H$ is equivalent to that of $\tnu_H$. This
  follows from the first part since the series $\sum_{s \in B \cap I}
  \vartheta(s^+)$ and $\sum_{s \in B \cap I}
  \vartheta(s^-)$ converge; recall that $\vartheta(s^-) = - \vartheta(s^+)$.
\end{proof}

\section{Existence and continuity properties of generalized
  bicharacteristics}
\label{sec: Existence and continuity properties of generalized
  bicharacteristics}
\subsection{Existence of generalized
  bicharacteristics}
\label{sec: proof: theorem: existence bichar}
Here, we show how our main result, Theorem~\ref{theorem: measure
  support propagation}, can be used to prove Theorem~\ref{theorem:
  existence bichar}, that is, the existence of \gbichars without
having to carry out a very similar and rather long proof.

As in the main text, the manifold $\TL$ is slightly extended beyond
its boundary and the metric $g$ and the wave symbol $p$ are extended
in a $\Con^1$ manner. With the Liouville measure
$\omega^{d+1}/ (d+1)!$ one identifies a function with a measure on
$\ThL$, for instance $\unitfunction{\TL}$.

Suppose $\y^0= (t^0,x^0,\tau^0, \xi^0) \in \TL \cap \Char p$. If
$\tau^0=0$ then $(\tau^0,\xi^0)=0$ and $\gammaG(s) = \y^0$,
$s \in \R$, is a maximal \gbichar that goes through $\y^0$. Assume now
that $\tau^0\neq 0$.

Note that, away from the zero section, $\Char p$ is a $\Con^1$-submanifold of codimension one of
$\ThL$, and consider the uniform measure
density $\mu_p$ on $\Char p \cap \TL \cap \{ \tau = \tau^0\}$, that
is, $\mu_p = \unitfunction{ \TL}\, \delta_{p=0}\delta_{\tau=\tau^0}$.  Note that this
product of  measures makes sense if one considers the product of
distributions with a wavefront set criterium; see for instance
Section~9.2 in \cite{Hoermander:V1}.

One has $\Hp p =0$, meaning that
the Hamiltonian vector field $\Hp$ is tangent to $\Char p$, and $\Hp\tau=0$, and thus $\Hp \delta_{p=0} \delta_{\tau=\tau^0}=0$ impling
that 
$\Hp \mu_p=0$ away from
$\dTL$.  

With the same arguments used in Proposition~\ref{prop: G submanifolds} one finds that
$\Char p \cap \dTL$ is a submanifold of codimension two,
locally given by $\{ p=0, z=0\}$. One has $\Char p \cap \dTL =
\Hb \cup \Gb$ as expressed below Definition~\ref{def: G H}.

Near the boundary, in a local chart as in \eqref{eq: local chart
  boundary}, use the quasi-normal geodesic coordinates of
Proposition~\ref{prop: quasi-normal coordinates}.
One has
$\mu_p = \unitfunction{ z\geq 0}\,  \delta_{p=0}\delta_{\tau=\tau^0}$.
The Leibnitz rule applies to this distribution product and one finds
\begin{align*}
  \Hp \mu_p = (\Hp \unitfunction{ z\geq 0})\, \delta_{p=0}\delta_{\tau=\tau^0}
  = (\Hpz)  \delta_{z=0} \delta_{p=0}\delta_{\tau=\tau^0}.
\end{align*}
As above the product $\delta_{z=0} \delta_{p=0}\delta_{\tau=\tau^0}$ makes sense under the
wavefront set criterium and $\Hp \mu_p$ is a measure. The measure $\delta_{z=0} \delta_{p=0} $ is
the uniform positive measure $\ell$ for $\Char p \cap \{ z=0\}$
inherited from the Liouville measure.
Note that $\unitfunction{\Gb} \, \Hp \mu_p =0$,
since $\Hpz=0$ on $\Gb$. Hence one can write 
\begin{align*}
  \Hp \mu_p = \unitfunction{\Hb} \Hp \mu_p,
\end{align*}
and 
\begin{align*}
  \dup{\Hp \mu_p}{a}
  = \int_{\Hb} (\Hpz) a_{|{z=0\atop p=0}} (\y) \, d \ell_{\tau^0} (\y),
\end{align*}
with $\ell_{\tau^0} = \delta_{\tau=\tau^0}\, \ell$.
Since $(z,\zeta) \mapsto (-z, -\zeta)$ leaves $\omega^{d+1}$ and thus $\ell$
invariant and moreover exchanges
$\Hb^+$ and $\Hb^-$ (as $\Hpz = 2 \zeta$ at $z=0$ in the
used quasi-geodesic coordinates), 
one finds
\begin{align*}
  \dup{\Hp \mu_p}{a}
  = \int_{\Hb^+}  \Big( (\Hpz)  a_{|{z=0\atop p=0}} (\y)
  + (\Hpz)  a_{|{z=0\atop p=0}} \big(\Sigma(\y)\big) \Big) d \ell_{\tau^0} (\y),
\end{align*}
with $\Sigma$ defined in \eqref{eq: involution-intro}.
Denoting by $\para{\ell}_{\tau^0}$ the pullback of the measure $\ell_{\tau^0}$ 
by the diffeomorphism
\begin{align*}
   \pHb &\to \Hb^+\\
  \y &\mapsto \y^+,
\end{align*}
one obtains
\begin{align*}
  \dup{\Hp \mu_p}{a}
  = \int_{\pHb}  \Big( (\Hpz)  a_{|{z=0\atop p=0}} (\y^+)
  + (\Hpz)  a_{|{z=0\atop p=0}} (\y^-) \Big) d\, \para{\ell}_{\tau^0} (\y),
\end{align*}
where $\y^+ \in \Hb^+$ and $\y^- = \Sigma(\y^+) \in \Hb^-$ with $\ppi(\y^+) =
\ppi(\y^-) = \y$ if $\y \in \pHb$ as in Section~\ref{sec: A partition of the cotangent bundle at the boundary-intro}.
One has  
$0 < \Hpz(\y^+) = - \Hpz(\y^-)$ yielding
\begin{align*}
  \dup{\Hp \mu_p}{a}
  = \int_{\pHb}  (\Hpz) (\y^+)\,
  \bigdup{\delta_{\y^+} - \delta_{\y^-}}{a} \,
  d\, \para{\ell}_{\tau^0} (\y).
\end{align*}
As computed in \eqref{eq: computation denominator GL} one has 
$\dup{\xi^+- \xi^-}{\nx}_{T_x^*\M, T_x\M}  = 2 \alpha(x)\Hpz(\y^+) =
\Hpz(\y^+)$ since $\alpha = 1/2$ at $z=0$ in the present
quasi-geodesic coordinates.
One can thus write 
\begin{align*}
  \Hpz(\y^+) = \frac{\Hpz(\y^+)^2}{\dup{\xi^+-
  \xi^-}{\nx}_{T_x^*\M, T_x\M} },
\end{align*}
implying 
\begin{align*}
  \dup{\Hp \mu_p}{a}= \int_{\pHb} 
    \frac{ \dup{\delta_{\y^+} - \delta_{\y^-}}{a}}
    {\dup{\xi^+- \xi^-}{\nx}_{T_x^*\M, T_x\M} 
    }    d \nu (\y), 
\end{align*}
with $\nu$ a nonnegative measure on $\pHb$ given by
$d \nu (\y) =\Hpz(\y^+)^2  d\, \para{\ell}_{\tau^0} (\y)$.
This is precisely the form of the equation one has in
Assumption~\ref{assumption: Gerard-Leichtnam equation} for the semiclassical
measure. 
The other condition on this measure stated in Assumption~\ref{assumption: first properties of the
  semiclassical measure} is obvious here: 
the measure $\mu_p$ is supported in $\Char p \cap \TL \setminus 0$.
Consequently, all the constructions made in Section~\ref{sec:
  Propagation of the measure  support} can be carried out with $\mu$
replaced by $\mu_p$. Theorem~\ref{theorem: measure support
  propagation} then implies Theorem~\ref{theorem: existence bichar}.
%%%%%%%%%%%%%%%%%%%%%%%%
% remark               %
%%%%%%%%%%%%%%%%%%%%%%%%
\begin{remark}
  \label{remark: measure localization tau}
  Note that one can replace the use of $\delta_{\tau=\tau^0}$ with
  $\unitfunction{\tau^1< |\tau| < \tau^2}$ for some $0< \tau^1<
  |\tau^0| < \tau^2< \infty$. The argument remains the same.
\end{remark}

%%%%%%%%%%%%%%%%%%%%%%%%%%
% sub-section
%%%%%%%%%%%%%%%%%%%%%%%%%%
\subsection{Continuity properties of generalized
  bicharacteristics}
\label{sec: proof: continuity properties bichar}
Here, we prove Proposition~\ref{prop:continuity-flow}.
 One proceeds by contradiction. Then, there exists
  $\y^0 \in \Char p \cap \TL$ and $T>0$ and $\eps^0>0$, and for all
  $n\in \N^*$ there exists $\y_n, \y^1_n$ such that
  $\dist(\y^1_n, \y^0) \leq 1/n$, $\y_n \in \Gamma^T(\y^1_n)$, and
  $\dist \big(\y_n, \Gamma^T(\y^0)\big) \geq \eps^0$.

  Write $\y^1_n = (t^1_n, x^1_n, \tau^1_n, \xi^1_n)$ and $\y_n = (t_n,
  x_n, \tau_n, \xi_n)$. One has $\y_n = \hgammaG_n(s_n)$ for some
  $s_n\in \R$ with
  $\hgammaG_n$ a \gbichar such that $\hgammaG_n(0) = \y^1_n$. Possibly $\y_n = \hgammaG_n(s_n^\pm)$, if $\y_n \in \Hb$.
  For all \gbichars, $\tau$ is constant; thus, $\tau_n = \tau^1_n$.

   If $\y^0=(t^0, x^0, \tau^0, \xi^0)$, the assumption states
   $(\tau^0, \xi^0)\neq 0$. Since $\y^0 \in \Char p$ one has
  $\tau^0 \neq 0$. One has $|\tau^0|/2  \leq |\tau^1_n| \leq 2
  |\tau^0|$ for $n$ \suff large. Set $S= 2 T/ |\tau^0|$.
  Denote $\hgammaG_n(s) = (\hatt_n(s),\hx_n(s),\tau^1_n, \hxi_n(s))$. For
  $|s| >S/2$ one has $|\hatt_n(s) - t^1_n| >  |\tau^1_n| S>  T$.
  This gives $s_n \in  [-S/2,S/2]$ for the value of $s_n$ introduced above.

  Set 
\begin{align*}
  \chgammaG_n(s) = \begin{cases}
    \cphi  \big(\hgammaG_n(s)\big) 
    = \hgammaG_n(s) & \text{if} \  s \notin B_n\\
     \cphi \big( \lim_{s' \to s^-} \hgammaG_n(s)\big) 
     = \cphi \big( \lim_{s' \to s^+} \hgammaG_n(s)\big) & \text{if} \  s \in B_n.
  \end{cases}
\end{align*}
The map $\cphi$ is defined in Section~\ref{sec: The compressed cotangent bundle}. 
The curve $\chgammaG_n(s)$ is continuous with values in the compressed
cotangent bundle $\cTL$.

  The sequence of functions ${\chgammaG_n}_{|[-S,S]}$ is equicontinuous. By
  the Arzel\`a-Ascoli theorem one can extract a subsequence
  $(s \mapsto \chgammaG_{n_p})_{p \in \mathbb{N}}$ that converges
  uniformly to a curve $\chgammaG(s)$ in $\cTL$ for $s \in [-S,S]$ and $s_{n_p}$
  converges to some $S' \in [-S/2,S/2]$. Write $n$ in place of
  $n_p$ for simplicity.
  In particular, one has $\chgammaG(0) = \cphi(\y^0)$.
  Set 
\begin{align*}
  B = \{ s \in [-S, S]; \  \chgammaG(s) \in \cphi(\Hb)\}.
\end{align*}
For $s \in [-S, S] \setminus B$  one defines $\hgammaG(s) =
\cphi^{-1} (\chgammaG(s))$. One has $\hgammaG(0)= \y^0$, with the understanding that  $\hgammaG(0^\pm) = \y^0$ if $\y^0 \in \Hb^\pm$.
With the same sequence
  of arguments as in Section~\ref{sec: proof local bichar} one obtains
  that $\hgammaG(s)$ is a \gbichar for $s\in[-S,S]$. This \gbichar can
  be extended to a maximal \gbichar by Theorem~\ref{theorem: existence bichar}.

  Assume first $S' \notin B$.  One has $\chgammaG_n(s_n) \to
  \chgammaG(S')$ with the equicontinuity property and the uniform
  convergence. Then, $\y_n = \hgammaG_n(s_n) \to \hgammaG(S')$. One
  has $t_n - t^1_n = - 2 \tau^1_n s_n$. With $|t_n - t^1_n|\leq T$
  this gives $|\tau^1_n s_n| \leq T/2$.  In the limit, this gives
  $|\tau^0 S'| \leq T/2$, implying $|\hatt(S')- t^0| \leq T$. Hence,
  $\y_n$ converges to a point of $\Gamma^T(\y^0)$. A contradiction.

  Assume second that $S'\in B$. One has $\chgammaG_n(s_n) \to
  \chgammaG(S')$. This means that one can find a subsequence of $\y^n$ that converges either to $\hgammaG(S^{\prime +})$ or to $\hgammaG(S^{\prime -})$. Then, one reaches the same contradiction.
\hfill \qedsymbol \endproof
%%%%%%%%%%%%%%%%%%%%%%%%%%
% section
%%%%%%%%%%%%%%%%%%%%%%%%%%
\section{Quasi-normal geodesic coordinates}
\label{sec proof prop: quasi-normal coordinates}

This section is devoted to the proof of Proposition~\ref{prop:
  quasi-normal coordinates}. In fact, we prove the following more
general result that is adatpted to various levels of regularity.
%%%%%%%%%%%%%%%%%%%%%%%%
% proposition          %
%%%%%%%%%%%%%%%%%%%%%%%%
\begin{proposition}[quasi-normal coordinates]
  \label{prop: quasi-normal coordinates-bis}
  Suppose  $k \geq 0$ and $\M$ is a $d$-dimensional manifold of
  class $\Con^{1+k}$ (\resp $W^{1+k,\infty}$) equipped
  with a metric of class $\Con^k$ (\resp $W^{k,\infty}$). Suppose
  $\mathcal N$ is a submanifold of class $\Con^{1+k}$ (\resp
  $W^{1+k,\infty}$) and codimension one and 
$m^0 \in \mathcal N$. There exists a
  local chart $(\hO, \chdiff)$ with regularity consistent
  with that of $\M$, that is, $\Con^{1+k}$ (\resp
  $W^{1+k,\infty}$), such that $m^0 \in \hO$,  $\chdiff(m) = (x',z)$ with $x' \in \R^{d-1}$ and $z \in \R$ and 
  \begin{enumerate}
    \item the vector field $\d_z$ is transverse to $\mathcal N$; 
    \item $\chdiff(\hO\cap \mathcal N) =  \{ z= 0\} \cap \chdiff(\hO)$;
    \item on $\mathcal N$ the representative of the metric reads 
      \begin{align*}
        g(x',0) = \sum_{1\leq i, j \leq d-1} g_{i j}(0,x')  dx^i \otimes dx^j
        + |dz|^2.
      \end{align*}
  \end{enumerate}
\end{proposition}
 If $k \geq 1$ one  has
  \begin{align*}
    g_{j d}(x',z)  = z h_{j d}(x',z)  \ \ \text{and} \ \ 
    g_{d d}(x',z)  = 1 +  z h_{d d}(x',z),
  \end{align*}
  for some $h_{j d}$, $j=1, \dots, d$, of class $\Con^{k-1}$ (\resp
  $W^{k-1, \infty}$).
\begin{proof}
  Consider a local chart $(\hO^0, \chdiff^0)$ such that
  $m^0 \in \hO^0$, with coordinates $(x',x_d)$ with
  $\chdiff^0(\hO^0\cap \mathcal N)$ given by $\{ x_d = 0\}$. Without
  loss of generality we may assume $x^0 = \chdiff^0(m^0) =0$.  Denote
  by $g^0$ the representative of the metric in this local chart.  In
  these coordinates, at $(x',0)$ the one form
  $\nu_{x'} = \big(g^{0,dd}(x',0)\big)^{-1/2}dx^d$ is unitary with respect
  to $g^*_x = \big(g^{ij}(x)\big) = \big(g_{ij}(x)\big)^{-1}$ and
  $n_{x'} = \nu_{x'}^\sharp$ is a unitary vector field with respect to
  $g_x=\big(g_{ij}(x)\big)$ and orthogonal to $T_{(x',0)} \mathcal N$
  in the sense of the metric $g_x$, thus transverse to $\mathcal
  N$. Its coodinates are given by
 \begin{align*}
   n_{x'}^j  =  (g^{0,dd})^{-1/2} g^{0, jd}(x',0).
 \end{align*}
Note that  the vector field $n$ defined along $\mathcal N$ has the
regularity of the metric, that is, $\Con^k$ (resp. $W^{k,\infty}$).

Suppose $R>0$ and $E>0$ are  such that $B'(0, 2R) \times [-E,E] \subset
\chdiff^0(\hO^0)$, where $B'(0, a)$ denotes the Euclidean ball of radius $a$ centered at
$x^0=0$ in the $x'$ variables. Suppose $\chi \in \Cinfc(\R^{d-1})$ with
$\chi=1$ in a \nhd of $B'(0, R)$ and $\supp \chi \subset B'(0, 2R)$. 
Introduce the vector field
\begin{align*}
  m(x',z) = e^{1 - \jp{z D_{x'}}} (\chi n)_{x'}, 
\end{align*}
with $z$ acting as a parameter at this stage and 
where $e^{1 - \jp{z D_{x'}}}$ stands as the operator  associated with the Fourier
multiplier $e^{1 - \jp{z \xi'}}$ with the usual notation $\jp{u} =
\sqrt{1+|u|^2}$ for $u \in \R^{d-1}$. 

%%%%%%%%%%%%%%%%%%%%%%%%
% lemma                %
%%%%%%%%%%%%%%%%%%%%%%%%
\begin{lemma}
  \label{lemma: regularity m}
  The function $(x',z) \mapsto m(x',z)$ is of class $\Con^{k}$
  (\resp $W^{k, \infty}$).
  Moreover  the function $(x',z) \mapsto z m(x',z)$ is of class $\Con^{1+k}$
  (\resp $W^{1+k, \infty}$).
\end{lemma}
A proof of this lemma is given below.

\medskip
Consider now the $\Con^{1+k}$ (\resp $W^{1+k, \infty}$)  function 
\begin{align*}
  \phi (x',z) = (x',0) + z m(x',z).
\end{align*}
Observe that $\d_{x_j} \phi (x',0) = e_j$, $j=1, \dots, d-1$, the
Euclidean unit
vector in the $j$-direction  and $\d_z \phi
(x',0) = m(x',0) = \chi n_{x'} = n_{x'} $ if $x' \in B'(0, R)$.  The Jacobian matrix of $\phi$ is thus full
rank in a \nhd of $x^0$. Thus, $\phi$ is a local  $\Con^{1+k}$ (\resp $W^{1+k, \infty}$) diffeomorphism,   and $(x',z)$ provides local
coordinates for a \nhd of $x^0$ in $\R^{d}$ of the form $(x',z)  \in B'(0,r)\times ]-e,e[$ for
some $0< r < R$ and 
$0< e< E$ chosen \suff small, and thus coordinates for an open \nhd of
$m^0$ in $\M$.  These coordinates have the announced regularity.   

We claim that $\d_z$ is the representative of $n_{x'}$ at $x=(x',z=0)$ in
the $(x',z)$ coordinates. In fact, consider a function on
$\M$ with $f$ as its representative in the original coordinates
$(x',x_d)$ and $\tilde{f}$ its representative in the $(x',z)$
coordinates. Then $\tilde{f} = f \circ \phi$ and one has 
\begin{align*}
  \d_z \tilde{f} (x',z) 
  = d f \big( \phi(x',z) \big) (\d_z \phi(x',z)).
\end{align*}
It was seen above that $\d_z \phi (x',0) = n_{x'}$ for $x' \in
B'(0,r)$ implying 
\begin{align*}
  \d_z \tilde{f} (x',0) 
  = d f(x',0) (n_{x'})
  = n_{x'}(f ), 
\end{align*}
which is precisely our claim. 

Note that $\phi(x',0) = (x',0)$ meaning that $\mathcal N$ is locally
given by $\{z=0\}$ and one sees that $\d_z$ is transverse to $\mathcal
N$ and $\{z>0\}$ coincides locally with $\{ x_d>0\}$.

We prove that the metric has the announced structure at a point of
$\mathcal N$ in the coordinates
$(x',z)$. By abuse of notation, still denote the representative of
the metric by $g$. One has
\begin{align*}
  &g_{ij}(x',0) = g_{(x',0)} (\d_{x_i}, \d_{x_j}), \ \ 1\leq i, j\leq
  d-1, \\
  &g_{dj}(x',0) = g_{(x',0)} (\d_z, \d_{x_j}), \ \ 1\leq j\leq
  d-1, \ \ \text{and} \\
   &g_{dd}(x',0) = g_{(x',0)} (\d_z, \d_z).
\end{align*}
Since $\d_z$ is the representative of $n_{x'}$ at $x=(x',z=0)$ in
the $(x',z)$ coordinates it follows that 
\begin{align*}
  &g_{dj}(x',0) = 0, \ \ 1\leq j\leq
  d-1, \ \ \text{and}  \ \
   g_{dd}(x',0) = 1,
\end{align*}
since $n_{x'}$ is orthogonal to $T_{(x',0)} \mathcal N$ and unitary in the
sense of the metric $g$. This concludes the proof of Proposition~\ref{prop: quasi-normal coordinates}.
\end{proof}

%%%% proof of lemma
\begin{proof}[Proof of Lemma~\ref{lemma: regularity m}]

  Preliminay observation.
 Suppose $p: \R^{d-1} \to \R$ and $q: \R \to \R$ are both polynomial.
 Observe that for $z>0$ the operator
 \begin{align*}
   p(z D_{x'}) q(\jp{z D_{x'}}) \exp(1-\jp{z
     D_{x'}})
 \end{align*}
 acts as a convolution with the function $k_z(x') = z^{1-d} \ell(x'/z)$ where
  \begin{align*}
    \ell(x') =\frac{1}{(2\pi)^{d-1}} \int e^{\scp{ix'}{\xi'}}
    p(\xi') q(\jp{\xi'}) e^{1-\jp{ \xi'}} d\xi'.
  \end{align*}
  Note that  $\ell \in \mathscr S(\R^{d-1}) \subset L^1(\R^{d-1})$ as the
  inverse Fourier transform of a Schwartz function, yielding $k_z
  \in L^1(\R^{d-1})$ uniformly in $z>0$ and moreover  $k_z \to C \delta$
  in the sense of measures with $C  = \int \ell$. Here, the Dirac
  measure acts in the $x'$ variable. If one considers a $\Con^0$-function $h(x')$ one thus obtains that $p(z D_{x'}) q(\jp{z D_{x'}}) \exp(1-\jp{z
    D_{x'}}) h$ is a $\Con^0$-function in both variables $x'$ and $z$.

  Recall the Fa\`a di Bruno formula for the repeated
  differentiation of the composition of two functions of one variable:
  \begin{align}
    \label{eq: Faa di Bruno}
   \frac{d^n}{d z^n} (f \circ h) (z) = \sum \frac{n!}{r_1! \cdots r_n!} 
  f^{(r_1+ \cdots +r_n)} \big( h(z) \big)
  \prod_{j=1}^n \Big( \frac{h^{(j)}(z)}{j!}\Big)^{r_j},
\end{align}
where the sum is carried out for $r_1 + 2 r_2 + \cdots +
n r_n=n$. Applying \eqref{eq: Faa di Bruno} with $f (z) = \exp(z)$ and $h(z) = 1-\jp{z \xi'}$ one
finds 
\begin{align}
    \label{eq: Faa di Bruno2}
   \frac{d^n}{d z^n} e^{1 - \jp{z \xi'}}= \sum \frac{n!}{r_1! \cdots r_n!} 
  e^{1 - \jp{z \xi'}}
  \prod_{j=1}^n \Big( - \frac{1}{j!}\frac{d^j}{d z^j} \jp{z \xi'}\Big)^{r_j},
\end{align}
and applying \eqref{eq: Faa di Bruno} with $f(z) = \sqrt{z}$ and $h(z) =1+ |\xi'|^2 z^2$,
since $h'(z) = 2|\xi'|^2 z$, $h''(z) = 2|\xi'|^2$ and $h^{(3)}=0$ one
finds for some $\alpha_{m_1, m_2}>0$
\begin{align}
   \label{eq: Faa di Bruno3}
  \frac{d^j}{d z^j} \jp{z \xi'} 
  &= \sum_{m_1 +2 m_2=j} \alpha_{m_1, m_2}\,
    z^{m_1} |\xi'|^{2(m_1+m_2)} \jp{z \xi'}^{1 - 2(m_1+m_2)}.
\end{align}
Expanding $|\xi'|^{2(m_1+m_2)} = (\xi_1^2 + \cdots +
\xi_{d-1}^2)^{m_1+m_2}$ one writes $z^{m_1} |\xi'|^{2(m_1+m_2)}$
as a linear combination  of terms of the form
\begin{align*}
  (z \xi')^{\beta_1^j} (\xi')^{\beta_2^j},
  \ \ \text{with} \ |\beta_1^j| =m_1 \ \text{and} \ 
  |\beta_2^j| =m_1 + 2 m_2 = j.
\end{align*}
Combining \eqref{eq: Faa di Bruno2} and \eqref{eq: Faa di Bruno3} one
obtains $\frac{d^n}{d z^n} \exp(1 - \jp{z \xi'})$ as a linear  combination  of terms of the form
\begin{align*}
 (\xi')^{\beta} p(z \xi') q(\jp{z \xi'}) \exp(1-\jp{z \xi'}),
\end{align*}
with
$p: \R^{d-1}\to \R$ and $q: \R \to \R$ both polynomial, and with
$|\beta|=n$ using that $r_1 + 2 r_2 + \cdots + n r_n=n$. 
One thus obtains $\d_{x'}^\gamma \d_z^n m$ as a linear  combination  of terms of the form
\begin{align*}
 p(z D_{x'}) q(\jp{z D_{x'}}) e^{1-\jp{z D_{x'}}}
  (\d_{x'})^{\beta + \gamma} (\chi n_{x'}),
\end{align*}
with $|\beta| =n$. 
Since $\chi n_{x'}$ is $\Con^k$  (\resp $W^{k,\infty}$)  with respect to $x'$, 
with the preliminary observation made above,  one
concludes that  $m(x',z)$ is $\Con^k$  (\resp $W^{k,\infty}$) with
respect to $(x',z)$.

Observe now that for $n\geq 1$
\begin{align*}
 \d_{x'}^\gamma \d_z^n \big(z m\big) 
  = n \d_{x'}^\gamma \d_z^{n-1}m  + z  \d_{x'}^\gamma \d_z^n  m
\end{align*}
One obtains $\d_{x'}^\gamma \d_z^n \big(z m\big) $ as a linear  combination  of terms of the form
\begin{align*}
 p(z D_{x'}) q(\jp{z D_{x'}}) e^{1-\jp{z D_{x'}}}
  (\d_{x'})^{\tilde{\beta} + \gamma}  (\chi n_{x'}),
\end{align*}
with
$\tilde{p}: \R^{d-1}\to \R$ and $q: \R \to \R$ both polynomial and with
$|\tilde{\beta}|=n-1$. Arguing as above one concludes that $z m(x',z)$ is 
$\Con^{k+1}$  (\resp $W^{k+1,\infty}$) with
respect to $(x',z)$. 
\end{proof}

% references
\bibliographystyle{plain}

\end{document}